\documentclass[reqno,10pt]{amsart}
\usepackage{amsmath}
\usepackage{amssymb}
\usepackage{amsthm}
\usepackage{amsfonts}
\usepackage{graphicx}
\usepackage{ esint }
\usepackage[abs]{overpic}
\usepackage{caption}
\usepackage{wrapfig}
\usepackage{float}
\usepackage{graphicx}
\usepackage{enumitem}

\newtheorem{theorem}{Theorem}[section]
\newtheorem{lemma}[theorem]{Lemma}

\newtheorem{corollary}[theorem]{Corollary}
\newtheorem{definition}[theorem]{Definition}
\newtheorem{example}[theorem]{Example}

\newcommand{\Proof}{\par\noindent{\em Proof. }}
\newcommand{\eop}{\nopagebreak\hspace*{\fill}$\Box$\smallskip}

\numberwithin{equation}{section}

\newtheorem{rem}[theorem]{Remark}

\newcommand{\N}{\Bbb N}
\newcommand{\Z}{\Bbb Z}
\newcommand{\R}{\Bbb R}

\def\sat{{\rm sat}}

\def\eps{\varepsilon}

\def\e{\mathbf{e}}

\def\dist{\operatorname{dist}}

\def\XXint#1#2#3{{\setbox0=\hbox{$#1{#2#3}{\int}$ }
\vcenter{\hbox{$#2#3$ }}\kern-.6\wd0}}

\usepackage{color}

\newcommand{\BBB}{\color{black}} 
 
\newcommand{\EEE}{\color{black}} 
\newcommand{\UUU}{\color{black}} 
\newcommand{\GGG}{\color{black}}
\begin{document}

\title[A piecewise Korn inequality in SBD]{A piecewise Korn inequality in $SBD$ and applications to embedding and density results}

\author{Manuel Friedrich}
\address[Manuel Friedrich]{Applied Mathematics M\"unster, University of M\"unster\\
Einsteinstrasse 62, 48149 M\"unster, Germany.}
\email{manuel.friedrich@uni-muenster.de}
\urladdr{https://www.uni-muenster.de/AMM/Friedrich/index.shtml}


\begin{abstract} 
We present a piecewise Korn inequality for generalized special functions of bounded deformation ($GSBD^2$) in a planar setting generalizing the classical result in elasticity theory to the setting of functions with jump discontinuities.  We show  that for every configuration there is a partition of the domain such that on each component of the cracked body the distance of the function from an infinitesimal rigid motion can be controlled solely in terms of the linear elastic strain. In particular, the result implies  that $GSBD^2$ functions have bounded variation after subtraction of a piecewise infinitesimal rigid motion. As an application we prove a density result in $GSBD^2$ \BBB in dimension two. \EEE Moreover, for all $d \ge 2$ we show $GSBD^2(\Omega) \subset (GBV(\Omega;\R))^d$ and the embedding $SBD^2(\Omega) \cap L^\infty(\Omega;\R^d) \hookrightarrow SBV(\Omega;\R^d)$ into the space of special functions of bounded variation ($SBV$). Finally, we present a Korn-Poincar\'e inequality for functions with small jump sets in arbitrary space dimension.  
\end{abstract}

\maketitle
                  
              \centerline{\today}                       %
      
      
\pagestyle{myheadings} 

 \section{Introduction}\label{sec: intro}

A natural framework for the investigation of damage and fracture models in a geometrically linear setting is given by the space of \emph{special functions of bounded deformation} investigated in \cite{Ambrosio-Coscia-Dal Maso:1997, Bellettini-Coscia-DalMaso:98}. The space $SBD^p(\Omega)$ consists of all functions $u \in L^1(\Omega;\R^d)$ whose  symmetrized distributional derivative $Eu := \frac{1}{2}((Du)^T + Du)$ is a finite $\R^{d \times d}_{\rm sym}$-valued Radon measure which can be written as the sum of an $L^p(\Omega;\R^{d \times d}_{\rm sym})$ function $e(u):= \frac{1}{2}((\nabla u)^T + \nabla u)$ and a part concentrated on a rectifiable set $J_u$ with finite $\mathcal{H}^{d-1}$ measure. Starting with the seminal paper   \cite{Francfort-Marigo:1998},  various variational problems for fracture mechanics in the realm of linearized elasticity have recently been investigated in the literature  (see e.g. \cite{Bourdin-Francfort-Marigo:2008, Chambolle:2003, Chambolle:2004, Focardi-Iurlano:13, Friedrich:15-2, SchmidtFraternaliOrtiz:2009}), where the common ground of all these models is that the main energy term is essentially of the form 
\begin{align}\label{rig-eq: general energy}
\int_\Omega |e(u)|^2 \,dx + \mathcal{H}^{d-1}(J_u)
\end{align}
for $u \in SBD^2(\Omega)$.  These so-called Griffith functionals comprise elastic bulk contributions for the unfractured regions of the body represented by the linear elastic strain $e(u)$ and surface terms that assign energy contributions on the crack paths comparable to the size of the `jump set' $J_u$. 

 For technical reasons models are often formulated in $SBD^2(\Omega) \cap L^\infty(\Omega; \R^d)$, e.g. in \cite{Bellettini-Coscia-DalMaso:98,Chambolle:2004,SchmidtFraternaliOrtiz:2009}, since for this setting a compactness result in $SBD$ was proved in \cite[Theorem 1.1]{Bellettini-Coscia-DalMaso:98} and hereby the existence of solutions to minimization problems related to \eqref{rig-eq: general energy} is guaranteed. To overcome this restriction, the space of \emph{generalized special functions of bounded deformation}, denoted by $GSBD(\Omega)$, was introduced in \cite{DalMaso:13}, admitting a compactness result under weaker assumptions and leading to the investigation of fracture models \cite{Friedrich:15-2, FriedrichSolombrino, Iurlano:13} in a more general setting. (For an exact definition of $GSBD$, which is based on properties of one-dimensional slices, we refer to Section \ref{rig-sec: sub, bd}.)

On the one hand,  models of the form \eqref{rig-eq: general energy} with linearized elastic energies  are in general easier to treat than their nonlinear counterparts  since in  the regime of finite elasticity the energy density of the elastic contributions is genuinely geometrically nonlinear due to frame indifference rendering the problem highly non-convex. On the other hand, a major difficulty of these models in contrast to nonlinear problems  formulated in the space  $SBV$ of \emph{special functions of bounded variation} (see \cite{Ambrosio-Fusco-Pallara:2000}) is given by the fact that  one controls only the linear elastic strain. 

Indeed, as discussed in \cite{Conti-Iurlano:15}, it appears that various properties being well established in $SBV$ are only poorly understood in $SBD$ due to the lack of control on the skew symmetric part of the distributional derivative $(Du)^T - Du$. Especially, to the best of our knowledge, for the \emph{coarea formula} in $BV$ (see \cite[Theorem 3.40]{Ambrosio-Fusco-Pallara:2000}), being useful in various applications as \cite{Ambrosio:90, BraidesPiat:96, Braides-Defranceschi, DalMaso-Francfort-Toader:2005, DeGiorgiCarrieroLeaci:1989, Francfort-Larsen:2003}, no $BD$ analog has been obtained in the literature.  

Therefore, it is a natural and highly desirable issue to gain a deeper understanding of the relation between $SBD$ and $SBV$ or even to show that in certain circumstances $SBD$ functions have bounded variation. In fact, one may expect that hereby various problems in $SBD$ could be solvable by a reduction to corresponding results in $SBV$. 

Clearly, a generic function of bounded deformation does not necessarily lie in $SBV$ as one can already construct a function with $J_u  = \emptyset$ such that  $e(u) \in L^1(\Omega)$, but $\nabla u \notin L^1(\Omega)$ (cf. \cite{ContiFaracoMaggi:2005,  Ornstein}). Moreover, following the examples in \cite{Ambrosio-Coscia-Dal Maso:1997, DalMaso:13}, one can define functions in $GSBD^2(\Omega)$ not having bounded variation (see Example \ref{eq: ex} and Lemma \ref{lemma: lq} below), i.e., also the higher  integrability for the elastic strain and the finiteness of the energy \eqref{rig-eq: general energy} is in general not sufficient. However, the examples involve unbounded configurations and it has therefore been conjectured that $SBD^2(\Omega) \cap L^\infty(\Omega;\R^d)$ is a subspace of $SBV(\Omega;\R^d)$.

Profound understanding of the connection between $SBV$ and $SBD$ functions is directly related to the validity of a Korn-type inequality in the setting of functions exhibiting  jump discontinuities. In elasticity theory, Korn's inequality is the key estimate providing a relation between the symmetric and the full part of the gradient (see e.g. \cite{Nitsche}). It states that the distance of $\nabla u$ from a skew symmetric matrix can be controlled solely in terms of the linear elastic strain $e(u)$. In the recently appeared contributions \cite{Conti-Iurlano:15.2, Friedrich:15-3}, this well-know estimate has been generalized in a planar setting to configurations in $SBD$ with small jump set controlling $\nabla u$   away from a small exceptional set. 

However, for general functions of  bounded deformation Korn's inequality in its basic form is doomed to fail. In fact, the domain may be disconnected into various parts by the jump set with completely different behavior on each component. In the special case that a brittle material does not store elastic energy,  i.e., $e(u)=0$, Chambolle, Giacomini, and Ponsiglione \cite{Chambolle-Giacomini-Ponsiglione:2007}  showed that  the body behaves piecewise rigidly, i.e.,  the only possibility that $u$ is not an affine mapping is that the body is divided into at most countably many  parts each of which subject to a different infinitesimal rigid motion.

Consequently, this observation already shows  that an analogous statement of Korn's inequality in $(G)SBD$ has to be formulated in a considerably more complex way involving a partition of the domain. The problem is related to the result in \cite{Friedrich-Schmidt:15}, where a quantitative piecewise  rigidity result in a geometrically nonlinear setting is established stating that in the planar case the distance of the deformation gradient from a piecewise  rigid motion can be controlled. 

The main goal of the present work is the derivation of an analogous result in the geometrically linear setting which we call a \emph{piecewise Korn inequality}. We show that in \BBB \emph{two dimensions} \EEE for each $u\in GSBD^2(\Omega)$ there is an associated partition whose boundary length is controlled by $\mathcal{H}^1(J_u)$ and a corresponding piecewise infinitesimal rigid motion $a$, being constant on each connected component of the cracked body, such that the distance of $u$ and $\nabla u$  from $a$ can be estimated in terms of $\Vert e(u) \Vert_{L^2(\Omega)}.$  This result for configurations storing elastic energy and exhibiting cracks may be seen as a suitable combination of Korn's inequality for elastic materials and the qualitative result in \cite{Chambolle-Giacomini-Ponsiglione:2007}.

 The estimate  proves to be useful to gain a deeper understanding of the relation between $SBV$ and $SBD$ functions. Although, as discussed above, $GSBD^2$ functions do not have bounded variation in general, the piecewise Korn inequality yields \BBB (in dimension two) \EEE that after subtraction of a piecewise infinitesimal rigid motion the function lies in the space $SBV^p$ for all $p<2$.  Hereby we particularly derive that each function in $GSBD^2$ has bounded variation away from an  at most countable union of sets of finite perimeter with arbitrarily small Lebesgue measure.
 
It turns out that for the subspace of bounded functions this observation can be substantially improved. In this case it is possible to control the piecewise infinitesimal rigid motion in terms of $\Vert u \Vert_\infty$ and $\mathcal{H}^{d-1}(J_u)$ and we indeed obtain the embedding 
\begin{align}\label{eq:em}
SBD^2(\Omega) \cap L^\infty(\Omega; \R^d) \hookrightarrow SBV(\Omega;\R^d),
\end{align}
which, in contrast to the piecewise Korn inequality, is proved in \BBB \emph{arbitrary space dimension} \EEE using a slicing technique. Moreover, similar arguments yield that without the $L^\infty$-bound one may  derive the inclusion $GSBD^2(\Omega) \subset (GBV(\Omega; \R))^d$ for $d \ge 2$.   (See Section \ref{rig-sec: sub, bd} below for the  definition of \emph{generalized functions of bounded variation}.) In particular, from \eqref{eq:em} we deduce that in the space of bounded $SBD^2$ functions, being a natural and widely adopted space in the investigation of linear fracture models, the \emph{coarea formula} is applicable.

Apart from the derivation of embeddings a major motivation for the derivation of the piecewise Korn inequality are applications to Griffith models. We show that the embedding result and the coarea formula allow to derive a  Korn-Poincar\'e inequality in $SBD$ stating that for a $GSBD^2$ function with small jump set the distance from a single infinitesimal rigid motion can be controlled outside a small exceptional set. This estimate, established in arbitrary space dimension, enhances a result recently obtained by Chambolle, Conti, and Francfort \cite{Chambolle-Conti-Francfort:2014} in the sense that also the length of the boundary of the exceptional set can be bounded in terms of $\mathcal{H}^{d-1}(J_u)$ and  therefore compactness results in $GSBD$ (see \cite{DalMaso:13}) are applicable.

We also present an approximation result for $GSBD^2$ functions in a \emph{planar setting} improving \cite{Iurlano:13}  in the sense that no $L^2$-bound on the function is needed. Hereby we can complete the $\Gamma$-convergence result for the linearization of Griffith energies \cite{Friedrich:15-2} by providing recovery sequences for every $GSBD^2$ function. \BBB (We also refer to the recent contribution \cite{Crismale} where a  generalization to arbitrary space dimension has been obtained.) \EEE Moreover, the main statement of this paper will be a fundamental ingredient to prove a general compactness theorem and to analyze quasistatic crack growth  for energies of the form \eqref{rig-eq: general energy} (see \cite{FriedrichSolombrino}).

The major difficulties in the derivation of the result come from the fact we treat   a full free discontinuity problem in the language of Ambrosio and De Giorgi \cite{DeGiorgi-Ambrosio:1988} deriving  an estimate without  any a priori assumptions on the crack geometry. Already simplified situations, in which the jump set decomposes the body into a finite number of sets with Lipschitz boundary, are subtle since  there are no uniform bounds on the constants in  Poincar\'e's and Korn's inequality. In fact, in the basic case of simply connected sets, in \cite{Friedrich:15-5} a lot of effort was needed to derive a decomposition result into domains for which uniform bounds on the constants can be derived. Even more challenging difficulties occur for highly irregular jump sets forming, e.g., infinite
crack patterns on various mesoscopic scales.

At first sight the derivation of the piecewise Korn inequality appears to be easier than the related problem \cite{Friedrich-Schmidt:15} due to the geometrical linearity. However, whereas in \cite{Friedrich-Schmidt:15} the statement is established only for a suitable,  arbitrarily small modification of the configuration, in the present context the estimates hold for the original function, where this stronger result is indispensable to prove the announced embedding and approximation results. In particular, the modification scheme for the deformation and jump set, which was iteratively applied in \cite{Friedrich-Schmidt:15} on mesoscopic scales becoming gradually larger, is not adequate in the present context and we need to use comparably different proof techniques.

The general strategy will be to first derive an auxiliary piecewise Korn inequality, which similarly as in \cite{Conti-Iurlano:15.2, Friedrich:15-3}  holds up to a small exceptional set. Then the main result follows by an iterative application of the estimate on various mesocopic scales. For the derivation of the auxiliary statement we first identify the regions where the jump set is too large. Hereby we can (1) construct a partition of the domain into simply connected sets and (2) use the Korn inequality \cite{Conti-Iurlano:15.2, Friedrich:15-3} for functions with small jump sets to find \BBB an approximation \EEE of the configuration, which is smooth on  each component  of the domain.
To control the shape of the components we apply the main result of the paper \cite{Friedrich:15-5} which allows to pass to a refined partition consisting of \emph{John domains} with uniformly controlled John constant. The auxiliary statement then follows by using a Korn  inequality for John domains (see e.g. \cite{Acosta}). 

Let us remark that the piecewise Korn inequality is only proved in a planar setting since also the major proof ingredients, the Korn inequality for functions with small jump sets  \cite{Conti-Iurlano:15.2, Friedrich:15-3} and the decomposition result \cite{Friedrich:15-5}, have only been shown in dimension two. Although the derivation of a higher dimensional analog  seems currently out of reach, the majority of the applications in this contribution, namely the embedding \eqref{eq:em} and the Korn-Poincar\'e inequality, hold in \emph{arbitrary dimensions} due to slicing arguments.

The paper is organized as follows. In Section \ref{sec: main} we present the main results including the embedding and approximation for $GSBD^2$ functions. In Example \ref{eq: ex} we also recall a standard construction for an $SBD$ function not having bounded variation, which is based on the idea to cut out small balls from the bulk part. Hereby we see that in view of the piecewise Korn inequality, examples of this type  essentially represent the only possibility to define $SBD^2$ functions not lying in $SBV$. Moreover, the construction shows that the embedding \eqref{eq:em} is sharp in the sense that the result is false if (1) $SBV$ is replaced by $SBV^p$, $p>1$, and if (2) $L^\infty$ is replaced by $L^q$, $q < \infty$.

Section \ref{sec: prep} is devoted to some preliminaries. We first recall the definition of John domains and formulate the decomposition result established in  \cite{Friedrich:15-5}. Then we state some properties of affine mappings being  especially important for the derivation of \eqref{eq:em}, where we need fine estimates how affine mappings and its derivatives can by controlled in terms of the volume and diameter of the underlying set. Afterwards, we recall the definition of $(G)SBV$ and $(G)SBD$ functions and discuss basic properties.

In Section \ref{sec: small set} we prove a version of the piecewise Korn inequality outside a small exceptional set. Together with the construction of a partition into simply connected sets we also provide a covering of Whitney-type such that the jump set $J_u$ in each element of the covering is small. This then allows us to apply the inequality \cite{Friedrich:15-3}  in order to define \BBB an approximation \EEE of the deformation being smooth in each part of the domain.  We discuss properties  of the covering and the exceptional set being crucial for the iterative application of the auxiliary result in Section 
 \ref{sec: main*}.

Section \ref{sec: main-proof} then contains the proof of the piecewise  Korn inequality. In Section 
 \ref{sec: main*} we first prove the result  for configurations in $SBD$ with a regular jump set consisting of a finite number of segments. Then in Section \ref{sec: main***} we derive the general case by an approximation argument similar to the one in \cite{Friedrich:15-3}. Here we also prove our main approximation result.
 
 Finally, Section \ref{sec: main**} contains the proof of the embedding results and the Korn-Poincar\'e inequality. As an ingredient for the derivation of \eqref{eq:em} we also provide  a \emph{piecewise Poincar\'e inequality}, which will allow to control the distance of configurations from piecewise infinitesimal rigid motions in terms of the $L^\infty$-norm.  For the latter as well as for the Korn-Poincar\'e inequality, the  coarea formula in $BV$ will turn out to be a key ingredient.

 \section{Main results}\label{sec: main}

 \subsection{Piecewise Korn inequality}\label{sec: main1}

Before we formulate the main results of this article, let us introduce some notions. We say a partition  $(P_j)_{j=1}^\infty$ of $\Omega \subset \R^d$ is a \emph{Caccioppoli partition} if each $P_j$ is a set of finite perimeter such that $\sum_{j=1}^\infty\mathcal{H}^{d-1}(\partial^* P_j)< + \infty$ (see Section \ref{rig-sec: sub, bd} for details). Given corresponding  affine mappings $a_j= a_{A_j,b_j}: \R^d \to \R^d$ with $a_j(x) = A_j\,x  + b_j$ for $A_j \in \R^{d \times d}_{\rm skew}$, $b_j \in \R^d$, we say $a_j$ is an \emph{infinitesimal rigid motion} and the function $\sum_{j=1}^\infty a_j \chi_{P_j}$ is a  \emph{piecewise infinitesimal rigid motion}. 

For the definition and properties of special functions of bounded variation ($SBV$) and deformation ($SBD$) we   refer to Section \ref{rig-sec: sub, bd}. In particular, for $u \in SBD^2(\Omega)$ we denote by $e(u)$ the part of the strain $Eu = \frac{1}{2}((Du)^T + Du)$ which is absolutely continuous with respect to $\mathcal{L}^d$ and let $J_u$ be the jump set. Our main goal is to show the following result.

\begin{theorem}\label{th: main korn}
Let $\Omega \subset \R^2$ open, bounded with Lipschitz boundary and $p \in [1,2)$.  Then there is a  constant $c=c(p)>0$ and a constant $C>0$ only depending on ${\rm diam}(\Omega)$ and $p$ such that the following holds: for  each $u \in GSBD^2(\Omega)$ there is a Caccioppoli partition $(P_j)^\infty_{j=1}$ of  $\Omega$ with 
 \begin{align}\label{eq: main1}
 \sum\nolimits^\infty_{j=1} \mathcal{H}^1( \partial^* P_j) \le c(\mathcal{H}^1(J_u) + \mathcal{H}^1(\partial \Omega))
 \end{align}
and corresponding infinitesimal rigid motions $(a_j)_j = (a_{A_j,b_j})_j$ such that 
\begin{align*}
v:= u - \sum\nolimits_{j=1}^\infty a_j \chi_{P_j} \in SBV^p(\Omega;\R^2) \cap L^\infty(\Omega;\R^2)
\end{align*}
and 
\begin{align}\label{eq: main korn2} 
\begin{split}
(i)& \ \ \Vert v \Vert_{L^\infty(\Omega)}\le C(\mathcal{H}^1(J_u) + \mathcal{H}^1(\partial \Omega))^{-1} \Vert e(u) \Vert_{L^2(\Omega)} ,\\ (ii)& \ \  \Vert \nabla v \Vert_{L^p(\Omega)} \le C \Vert e(u) \Vert_{L^2(\Omega)}.
\end{split}
\end{align}
\end{theorem}

In its most general form the result holds in the generalized space $GSBD^2(\Omega)$, which loosely speaking arises from $SBD^2(\Omega)$ by requiring that one-dimensional slices  have bounded variation. (We again refer to Section \ref{rig-sec: sub, bd} for the  exact definition.) The reader not interested in the generalized space may readily replace $GSBD^2(\Omega)$ by $SBD^2(\Omega)$ here and in the following (except for Lemma \ref{lemma: lq}). The result will first be proved for configurations with a regular jump set consisting of a finite number of segments. Only at the very end we pass to the general case by using a density result (see Theorem \ref{th: iurlano}, Lemma \ref{lemma: iurlano*}) and Ambrosio's compactness theorem in $SBV$. Consequently, the article is in large part accessible for readers without familiarity with the properties of the above mentioned function spaces. 
 
\begin{rem}\label{rem: remark}
{\normalfont

Let us comment the results of Theorem \ref{th: main korn}. 
\begin{itemize}
\item[(i)] An important special case is given by functions not storing linearized elastic energy, i.e., $e(u) = 0$ a.e. Then \eqref{eq: main korn2}  gives that $u$ is a piecewise infinitesimal rigid motion and thus Theorem \ref{th: main korn} can be interpreted as a \emph{piecewise rigidity result}. It shows that the only way that global rigidity may fail is that the body is divided into at most countably many parts each of which subject to a different infinitesimal rigid motion. This result has been discussed for $SBD$ functions in \cite{Chambolle-Giacomini-Ponsiglione:2007} and Theorem \ref{th: main korn} may be seen as a quantititive extension in the $GSBD$ setting.
\item[(ii)] It is a well known fact that in general $\nabla u$ and $u$ are not summable for $u \in  GSBD^2(\Omega)$.   The result shows that one may gain control over the function and its derivative after subtraction of a piecewise infinitesimal rigid motion. We hope that such an estimate will be useful to reduce various problems from the $(G)SBD$ to the $SBV$ setting.
\item[(iii)] In contrast to  \cite{Friedrich-Schmidt:15}, where a related result in a geometrically nonlinear setting is established, in  Theorem \ref{th: main korn} no modification of the configuration $u$ is necessary. This will be crucial for the embedding results announced in Section \ref{sec: main2}. 
\item[(iv)] In general one has $\sum\nolimits^\infty_{j=1} \mathcal{H}^1( \partial^* P_j \setminus J_u) >0$ even if $J_u$ already induces a partition $(P'_j)_j$ of $\Omega$ with $J_u = \bigcup_j \partial^* P'_j \cap \Omega$ up to a set of negligible $\mathcal{H}^1$-measure. Indeed, it is well known that in irregular domains, e.g. in domains with external cusps, Korn's inequality may fail (cf. \cite{Geymonat-Gilardi}). Also in the case that each $P'_j$ is a Lipschitz set, $\sum\nolimits^\infty_{j=1} \mathcal{H}^1( \partial^* P_j \setminus J_u) >0$ might be necessary since possibly the constants $C_{\rm korn}(P'_j)$ involved in the $W^{1,p}$ Korn inequality tend to infinity as $j \to \infty$. 
 
\item[(v)] Note that the constant in \eqref{eq: main korn2} depends on $\Omega$ only through its diameter and $p \in [1,2)$. Consequently, the above result may be also interesting in the case of $u \in H^1(\Omega)$ and varying domains $\Omega$. 
\item[(vi)] At some points of the proof we have to pass to a slightly smaller exponent by H\"older's inequality and therefore we derive the result only for $p<2$. In particular, we do not know if a similar result holds in the critical case $p=2$. Note, however, that for the applications we have in mind, it is sufficient that \eqref{eq: main korn2}  holds for $p=1$.
\item[(vii)] \UUU  In revising this paper, we remark that a generalized version of the result in the space $GSBD^q(\Omega)$ for $1 < q < \infty$ can be derived along similar lines, see Remark \ref{rem: general version} for some details. \EEE
\end{itemize}

}

\end{rem} 

 The result is first shown up to a small exceptional set (see Theorem \ref{th: korn-small set}) and then Theorem \ref{th: main korn} for functions with regular jump set follows from an iterative application of the estimate in Theorem \ref{th: korn-small set} (see Section \ref{sec: main*}). Finally, the general version is proved by approximation arguments given in Section \ref{sec: main***}. For a more detailed outline of the proof we refer to the beginning of Section \ref{sec: small set} and Section \ref{sec: main-proof}.

In Remark \ref{rem: remark}(ii) we have discussed that in general $GSBD$ functions are  not summable. A similar problem occurs in  $GSBV$  (see Section  \ref{rig-sec: sub, bd} for the exact definition). However, we get that  $GSBV$ functions are integrable after subtraction of a \emph{piecewise constant function}.

 \begin{theorem}\label{th: main poinc}
Let $\Omega \subset \R^d$ open, bounded with Lipschitz boundary. Let $\rho>0$ and $m \in \N$. For  each $u \in (GSBV(\Omega; \R))^m$ with $\Vert  \nabla u\Vert_{L^1(\Omega; \R^{d\times m})} + \mathcal{H}^{d-1}(J_u) < + \infty$  there is a Caccioppoli partition $(P_j)_{j=1}^\infty$ of $\Omega$ and corresponding translations $(b_j)_{j=1}^\infty \subset \R^m$ such that  $v:= u - \sum_{j=1}^\infty b_j \chi_{P_j} \in SBV(\Omega;\R^m) \cap L^\infty(\Omega;\R^m)$ and 
\begin{align}\label{eq: kornpoinsharp0}
\begin{split}
(i) & \ \   \sum\nolimits_{j=1}^\infty \mathcal{H}^{d-1}( (\partial^* P_j \cap \Omega) \setminus J_u ) \le  c\rho,\\
(ii) &\ \ \Vert v \Vert_{L^{\infty}(\Omega; \R^m)} \le c \rho^{-1}  \Vert  \nabla u\Vert_{L^1(\Omega; \R^{d\times m})}
\end{split}
\end{align}
for a dimensional constant $c=c(m)>0$.
\end{theorem}

Note that for the special choice $\rho = \theta (\mathcal{H}^{d-1}(J_u) + \mathcal{H}^{d-1}(\partial \Omega))$,   $\theta>0$, \eqref{eq: kornpoinsharp0}(ii) is similar to \eqref{eq: main korn2}(i). In particular, in the proofs we will use Theorem \ref{th: main poinc} to derive \eqref{eq: main korn2}(i) once \eqref{eq: main korn2}(ii) is established. The estimate is essentially based on the \emph{coarea formula} in $BV$ (see \cite[Theorem 3.40]{Ambrosio-Fusco-Pallara:2000}) and  the argumentation in the proof goes back to \cite[Theorem 3.3]{Ambrosio:90} and \cite[Proposition 6.2]{Braides-Defranceschi}. We remark that, in contrast to Theorem \ref{th: main korn},  the result is derived in arbitrary space dimension. The proof will be given in Section \ref{subsec: main**1}.

\begin{rem}\label{rem: better}
{\normalfont
Before we proceed with some applications of Theorem \ref{th: main korn}, let us comment on the estimates \eqref{eq: main1} and \eqref{eq: kornpoinsharp0}(i). In \eqref{eq: kornpoinsharp0}(i) for the choice $\rho = \theta (\mathcal{H}^{d-1}(J_u) + \mathcal{H}^{d-1}(\partial \Omega))$ with $\theta$ small we obtain a fine estimate in the sense that the additional surface energy associated to the partition may be taken arbitrarily small with respect to the original jump set, where the constant in \eqref{eq: kornpoinsharp0}(ii) blows up for $\theta \to 0$. 

In contrast to the embedding and approximation results addressed in this article, for certain applications of the piecewise Korn inequality it is crucial to have also a refined version of \eqref{th: main korn}, e.g. for a jump transfer theorem in $SBD$ or a general compactness and existence result for Griffith energies in the realm of linearized elasticity (see \cite{FriedrichSolombrino}). The derivation of such an estimate needs additional techniques and we refer to \cite{FriedrichSolombrino} for a deeper analysis.  
}
\end{rem}

 \subsection{Approximation of GSBD functions}\label{sec: main3}

For each $u \in GSBD^2(\Omega)$ we can consider the sequence
\begin{align}\label{eq: sequence}
v_n = u - \sum\nolimits_{j\ge n} a_j \chi_{P_j} \in SBV^p(\Omega; \R^2) \cap L^\infty(\Omega; \R^2)
\end{align}
with $(P_j)_j$  and $(a_j)_j$ as in Theorem \ref{th: main korn}. Hereby  we can approximate $GSBD^2$ functions by bounded functions of bounded variation which coincide with the original function up to a set of arbitrarily small measure. In particular, combining this observation with the density result proved in \cite{Iurlano:13} (see Theorem \ref{th: iurlano}) we obtain the following improved approximation result for $GSBD^2$ functions in a planar setting (see Section \ref{sec: main***} for the proof).

\begin{theorem}\label{th: fried-iurlano}
Let $\Omega \subset \R^2$ open, bounded with Lipschitz boundary. Let $u \in GSBD^2(\Omega)$. Then there exists a sequence $(u_k)_k \subset  SBV^2(\Omega; \R^2)$ such that each $J_{u_k}$ is the union of  a finite number of closed connected pieces of $C^1$-curves, each $u_k$ belongs to
$W^{1,\infty}(\Omega \setminus J_{u_k};\R^2)$ and the following properties hold:
\begin{align}\label{eq: dense}
\begin{split}
(i) & \ \ u_k \to  u  \text{ in measure on } \Omega,\\
(ii) & \ \ \Vert e(u_k) - e(u) \Vert_{L^2(\Omega)} \to 0,\\
(iii) &  \ \ \mathcal{H}^{1}(J_{u_k} \triangle J_u) \to  0.
\end{split}
\end{align}
\end{theorem}

Here the symbol $\triangle$ denotes the symmetric difference of two sets. Observe that in \cite{Chambolle:2004, Iurlano:13} the additional assumption $u \in L^2(\Omega)$ was needed which \BBB we circumvent here by \EEE  (1) first approximating $u$ by $(v_n)_n \subset GSBD^2(\Omega) \cap L^2(\Omega)$ as given in \eqref{eq: sequence} and then (2) applying the original density result \cite{Chambolle:2004, Iurlano:13}. In fact, in various applications in fracture mechanics  approximation results are essential, e.g.  for the construction of recovery sequences. Consequently, to establish complete $\Gamma$-limits it is indispensable to have a density result for all admissible functions. In particular, Theorem \ref{th: fried-iurlano} allows us to complete the picture in the linearization result for Griffith energies presented in \cite{Friedrich:15-2}. \BBB We also refer to the recent contribution \cite{Crismale} where with other techniques a generalization of the density result to $GSBD^p$ functions, $1 < p < \infty$, in arbitrary space dimension has been obtained.  \EEE 

Recall that we establish Theorem \ref{th: main korn} first for  functions with a regular jump set and the general case then follows by approximation of any $u \in GSBD^2(\Omega)$. Although Theorem \ref{th: fried-iurlano} provides such a density result, in the  proof we have to argue differently since we already use the general version of Theorem \ref{th: main korn} to prove Theorem \ref{th: fried-iurlano}. The strategy in the proof of Theorem \ref{th: main korn} is to provide an adaption of the arguments in \cite{Chambolle:2004, Iurlano:13} and to approximate each $u\in GSBD^2(\Omega)$  (without $L^2(\Omega)$ assumption) by functions $(u_k)_k$ with regular jump set such that \eqref{eq: dense}(i) holds and $\Vert e(u_k) \Vert_{L^2(\Omega)} \le c\Vert e(u) \Vert_{L^2(\Omega)}$, $\mathcal{H}^1(J_{u_k})\le c\mathcal{H}^1(J_{u})$ (see Lemma \ref{lemma: iurlano*}), which is sufficient for \eqref{eq: main1} and \eqref{eq: main korn2}.

 \subsection{Relation between SBV and SBD: Embedding results}\label{sec: main2}

The standard examples for functions of bounded deformation not having bounded variation are given by configurations where small balls are cut out from the bulk part with an appropriate choice of the functions on these specific sets (see e.g.  \cite{Ambrosio-Coscia-Dal Maso:1997, Conti-Iurlano:15, DalMaso:13}). We include the construction here for convenience of the reader. 

\begin{example}\label{eq: ex}
{\normalfont
Let $B_k = B_{r_k}(x_k) \subset B_1(0)\subset \R^d$ be pairwise disjoint balls for $k \in \N$ with $x_k$ converging to some limiting point $x_\infty$. For $k\in \N$ let $A_k =(a^k_{ij})\in \R^{d \times d}_{\rm skew}$ with $a^k_{12} = -a^k_{21} = d_k \in \R$ and $a^k_{ij} = 0$ otherwise. Define the piecewise infinitesimal rigid motion
\begin{align}\label{eq: cif}
u(x) = \sum\nolimits_{k=1}^\infty (A_k \,(x-x_k)) \chi_{B_k}(x).
\end{align}
With the choice $r_k = (\frac{1}{k (\log(k))^2})^{\frac{1}{d-1}}$ and $d_k = r_k^{-1}$  we get  $u \in  L^\infty(B_1(0);\R^d) $  and $u \in SBV(B_1(0);\R^d) \cap SBD(B_1(0))$ with $e(u) = 0$ a.e. and $\mathcal{H}^{d-1}(J_u)< + \infty$. However,  $\nabla u \notin L^p(B_1(0);\R^{d \times d})$ for $p>1$.  
}
\end{example}

In view  of Theorem \ref{th: main korn}, we see that the above construction is essentially the only way to  define   $(G)SBD$ functions with $e(u) \in L^2(\Omega)$ and $\mathcal{H}^1(J_u) < +\infty$ not having bounded variation. In particular, by considering $u \chi_{\bigcup_{j=1}^n P_j}$, $n \in \N$, with $(P_j)_j$ as in \eqref{eq: main1}, we observe that each  $u \in GSBD^2(\Omega)$  has bounded variation away from an  at most countable union of sets of finite perimeter with arbitrarily small Lebesgue measure.

Observe that $\nabla u \in L^1(B_1(0);\R^{d \times d})$ for the function given in Example \ref{eq: ex}. An   $L^\infty$ bound indeed  implies that functions are of bounded variation as the following embedding result shows.

\begin{theorem}\label{th: embed}
Let $\Omega \subset \R^d$ open, bounded  with Lipschitz boundary. Then   $SBD^2(\Omega) \cap L^\infty(\Omega; \R^d) \hookrightarrow SBV(\Omega; \R^d)$. More precisely, there is a constant $C>0$ only depending on $\Omega$ such that
\begin{align}\label{eq: embed}
\Vert \nabla u\Vert_{L^1(\Omega)} \le  C\Vert e(u) \Vert_{L^2(\Omega)} + C\Vert u \Vert_{L^\infty(\Omega)}(\mathcal{H}^{d-1}(J_u)+1).
\end{align}
\end{theorem}

Note that in contrast to Theorem \ref{th: main korn}, our main embedding result as well as the following theorems are valid in arbitrary space dimension, where we apply a slicing technique to extend the estimate from the planar to the $d$-dimensional setting. As a direct consequence, we get that one may apply the  \emph{coarea formula} in $BV$  for functions in $SBD^2(\Omega) \cap L^\infty(\Omega;\R^d)$. This estimate can be potentially applied in many situations and will hopefully reveal useful.

Note that the property $u \in SBV(\Omega;\R^d)$ does not follow directly from \eqref{eq: embed}. However, by first deriving Theorem \ref{th: embed} for functions with regular jump set and then using a density argument (see Theorem \ref{th: iurlano}), we can show $SBD^2(\Omega) \cap L^\infty(\Omega; \R^d) \subset BV(\Omega; \R^d)$. The claim then follows from Alberti's rank one property in $BV$ (see \cite{Alberti:93}) implying $|D^c u|(\Omega) \le \sqrt{2} |E^c u|(\Omega)$.

The first term on the right hand side in \eqref{eq: embed}  controls the distance of $u$ from a piecewise infinitesimal rigid motion $(a_j)_j = (a_{A_j,b_j})_j$ (cf. \eqref{eq: main korn2}(ii)) and the last term allows to estimate the skew symmetric matrices $(A_j)_j$. Let us note that the above inequality is optimal in following sense: by Example \ref{eq: ex} the $L^1$-norm on the left cannot be replaced by any other $L^p$-norm.  Moreover, the bound $\mathcal{H}^{d-1}(J_u) < + \infty$ is essential by \cite[Theorem 3.1]{Conti-Iurlano:15}, where a function in $u \in  SBD(\Omega) \cap L^\infty(\Omega;\R^d)$ is constructed  with $e(u) = 0$ a.e., $\mathcal{H}^{d-1}(J_u) = \infty$ and $\nabla u \notin L^1(\Omega;\R^{d \times d})$. Finally, the following lemma based on the construction given in Example \ref{eq: ex} shows that the $L^\infty$-norm cannot be replaced by any other $L^q$-norm. 

\begin{lemma}\label{lemma: lq}
Let $\Omega \subset \R^d$ open, bounded. Then for all $1 \le q < \infty$ we find some $u \in GSBD^2(\Omega) \cap L^q(\Omega;\R^d)$ such that the approximate differential $\nabla u$ exists a.e. and $\nabla u \notin L^1(\Omega; \R^{d \times d})$.  
\end{lemma}

 Consequently, $GSBD^2(\Omega) \cap L^q(\Omega;\R^d) \not\subset BV(\Omega;\R^d)$ for all $q < \infty$. In Theorem \ref{th: main korn} we have seen that $v = u - a \in SBV^p(\Omega;\R^2) \subset GBV(\Omega;\R^2)$ for a piecewise infinitesimal rigid motion $a:=\sum_{j=1}^\infty a_j\chi_{P_j} $. Moreover, one has that $a \in GBV(\Omega;\R^2)$ (see \cite[Theorem 2.2]{Conti-Iurlano:15}). Of course, herefrom we cannot directly deduce that $u = v + a \in GBV(\Omega;\R^2)$ since $GBV(\Omega; \R^2)$ is not a vector space. However, the property holds as the following inclusion shows.

\begin{theorem}\label{th: embed2}
Let $\Omega \subset \R^d$ open, bounded  with Lipschitz boundary. Then  $GSBD^2(\Omega) \subset (GBV(\Omega;\R))^d\subset GBV(\Omega;\R^d)$. More precisely, there is a constant $C>0$ only depending on $\Omega$ such that   for all $i=1,\ldots,d$ and   all $M>0$ 
  \begin{align}\label{eq: embed0}
|Du^M_i|(\Omega) &\le   CM(\mathcal{H}^{d-1}(J_u)+1) + C\Vert e(u) \Vert_{L^2(\Omega)}.
 \end{align}

\end{theorem}
Here for  $u \in GSBD^2(\Omega)$  the functions $u_1,\ldots,u_d$ denote the components of $u$ and $u_i^M := \min \lbrace\max\lbrace u_i,-M\rbrace,M\rbrace$ for $M>0$. Theorem  \ref{th: embed2} shows that the coarea formula in $GBV$ (see \cite[Theorem 4.34]{Ambrosio-Fusco-Pallara:2000}) is applicable in $GSBD^2(\Omega)$. The results announced in this section will be proved in Section \ref{subsec: main**2}.

 \subsection{A Korn-Poincar\'e inequality for functions with small jump set}\label{sec: main4}

As an application of our embedding results and the coarea formula, we present a Korn-Poincar\'e inequality for functions with small jump set.  

\begin{theorem}\label{th: kornpoin-d}
Let $Q\subset \R^d$ be a cube. Then there is a constant  $C>0$ only depending on ${\rm diam}(Q)$ and $d \in \N$  such that for all $u \in GSBD^2(Q)$ the following holds: there is a set of finite perimeter $E \subset Q$ with 
\begin{align}\label{eq: R2mainx^2}
\mathcal{H}^{d-1}(\partial^* E \cap Q) \le C(\mathcal{H}^{d-1}(J_u))^{\frac{1}{2}}, \ \ \ \ |E| \le C(\mathcal{H}^{d-1}(J_u))^{\frac{d}{d-1}}
\end{align}
and an infinitesimal rigid motion $a = a_{A,b}$ such that
\begin{align}\label{eq: main estmainx^2}
\begin{split}
(i) & \ \ \Vert u - a \Vert_{L^{\frac{2d}{d-1}}(Q \setminus E)}\le C \Vert e(u) \Vert_{L^2(Q)},\\ 
(ii) & \ \ \Vert u - a \Vert_{L^\infty(Q \setminus E)} \le C(\mathcal{H}^{d-1}(J_u))^{-\frac{1}{2}} \Vert e(u) \Vert_{L^2(Q)}.
\end{split}
\end{align}
Moreover, if $\mathcal{H}^{d-1}(J_u) >0$, $\bar{u} := (u-a) \chi_{Q \setminus E} \in SBV(Q;\R^d) \cap SBD^2(Q) \cap L^\infty(Q;\R^d)$ and
\begin{align}\label{eq: main estmain2X}
\begin{split}
|D\bar{u}|(Q) \le C\Vert e(u) \Vert_{L^2(Q)}.
\end{split}
\end{align}
\end{theorem}

Note that the exceptional set  $E$ is associated to the parts of $Q$ being detached from the bulk part of $Q$ by $J_u$. A variant of the proof shows that $(\mathcal{H}^{d-1}(J_u))^{\frac{1}{2}}$ in \eqref{eq: R2mainx^2}, \eqref{eq: main estmainx^2}(ii) can be replaced by $(\mathcal{H}^{d-1}(J_u))^{\frac{p}{2}}$ if in \eqref{eq: main estmainx^2}(i) we replace $\frac{2d}{d-1}$ by $\frac{2d}{p(d-1)}$ for $1 \le p \le 2$.

This result improves an estimate recently obtained by Chambolle, Conti, and Francfort \cite{Chambolle-Conti-Francfort:2014} in the sense that $\mathcal{H}^{d-1}(\partial^* E)$ can be controlled and  therefore compactness results in $GSBD$ (see \cite{DalMaso:13}) are applicable. In addition, we  provide a Korn-type estimate in \eqref{eq: main estmain2X}. Similar results in a planar setting have been investigated in \cite{Conti-Iurlano:15.2, Friedrich:15-3}.

 The statement essentially follows by combining (1) the result in  \cite {Chambolle-Conti-Francfort:2014} and (2) applying the coarea formula in $BV$ together with \eqref{eq: embed0}. In particular, the argument uses a truncation at a specific level set (cf. \eqref{eq: main estmainx^2}(ii)) which is reminiscent of the Poincar\'e inequality in $SBV$  due to De Giorgi, Carriero, and Leaci (see \cite{DeGiorgiCarrieroLeaci:1989}). The proof will be given in Section \ref{subsec: main**3}.

\section{Preparations}\label{sec: prep}

\subsection{John domains}

A key step in our analysis will be the construction of a partition and a corresponding \BBB approximation of the configuration \EEE being in $W^{1,p}$ on each component of the partition. Then in the  application of  Poincar\'e's and Korn's inequality on each component, it is essential to provide uniform bounds for the constants involved in the inequalities. To this end, we introduce the notion of \emph{John domains}.

\begin{definition}\label{def: chain} 
{\normalfont

Let $\Omega \subset \R^d$ be a bounded domain and let $x_0 \in \Omega$. We say $\Omega$ is a $\varrho$-\emph{John domain} with respect to the \emph{John center} $x_0$ and with the constant
$\varrho \ge 1$ if for all $x \in \Omega$ there exists a rectifiable curve $\gamma: [0,l_\gamma] \to \Omega$, parametrized by arc length, such that $\gamma(0) = x$, $\gamma(l_\gamma) = x_0$ and 
$t \le \varrho\dist(\gamma(t),\partial \Omega)$ for all $t \in [0,l_\gamma]$. 
}
 \end{definition}
 
 Domains of this form were introduced by John  \cite{John:1961} to study problems in elasticity theory and the term was first used by Martio and Sarvas  \cite{Martio-Sarvas:1978}. Roughly speaking, a domain is a John domain if it is possible to connect two arbitrary points without getting too close to the boundary of the set. This class is much larger than the class of Lipschitz domains and contains sets which may possess fractal boundaries or internal cusps (external cusps are excluded), e.g. the interior of Koch's \emph{snow flake} is a John domain. Although in the following we will only consider domains with Lipschitz boundary, it is convenient to introduce the much  more general notion of John domains as the constants in  Poincar\'e's and Korn's inequality only depend on the John constant. More precisely, we have the following statement (see e.g. \cite{Acosta, BuckleyKoskela, Diening}).

\begin{theorem}\label{th: kornsobo}
Let ­ $\Omega \subset \R^d$ be a $c$-John domain. Let $p \in (1,\infty)$ and $q \in (1,d)$. Then there is a constant $C=C(c,p,q)>0$ such that for all $u \in W^{1,p}(\Omega)$ there is some $A \in \R^{d \times d}_{\rm skew}$ such that
$$\Vert \nabla u -A\Vert_{L^p(\Omega)} \le C\Vert e(u) \Vert_{L^p(\Omega)}.$$
Moreover,  for all $u \in W^{1,q}(\Omega)$ there is some $b \in \R^d$ such that
$$\Vert u -b\Vert_{L^{q^*}(\Omega)} \le C\Vert \nabla u \Vert_{L^q(\Omega)},$$
where $q^* = \frac{dq}{d-q}$. The constant is invariant under rescaling of the domain. 
\end{theorem}

We recall a result about the decomposition of sets into John domains.

 \begin{theorem}\label{th: main part2}
There is a universal constant $c>0$ such that for all simply connected, bounded domains $\Omega \subset \R^2$ with Lipschitz boundary and all $\eps >0$ there is a  partition  $\Omega = \Omega_0 \cup \ldots \cup \Omega_N$ (up to a set of negligible measure) with $|\Omega_0|  \le \eps$ and the sets $\Omega_1,\ldots,\Omega_{N}$ are $c$-John domains with Lipschitz boundary satisfying 
\begin{align}\label{eq:parti2}
\sum\nolimits^N_{j=0}\mathcal{H}^1(\partial \Omega_j) \le c\mathcal{H}^1(\partial \Omega).
\end{align}
\end{theorem}

Observe that in general it is necessary to introduce an (arbitrarily) small exceptional set $\Omega_0$ as can be seen, e.g., by  considering polygons with very acute interior angles. The statement is given explicitly in \cite[Theorem 6.4] {Friedrich:15-5}, but can also easily be derived from the main result of that paper (\cite[Theorem 1.1] {Friedrich:15-5}) since Lipschitz sets can be approximated from within by smooth sets. 

 The essential step in the proof of Theorem \ref{th: main part2} is to consider polygonal domains. One observes that convex polygons can be (iteratively) separated into convex polygons such that \eqref{eq:parti2} holds and each set contains a ball whose size is comparable to the diameter of the set, whereby Definition \ref{def: chain} (for fixed $\varrho>0$)  can be confirmed. For general, nonconvex polygons the property in Definition \ref{def: chain}  may be violated if a concave vertex is `too close to the opposite part of the boundary'. It is shown that in this case the set may be (iteratively) separated into smaller polygons with less concave vertices, for which Definition \ref{def: chain} holds. We refer to \cite{Friedrich:15-5} for more details.  

\UUU 
We will also use the following simple property (see e.g. \cite[Lemma 2.3]{Vaisala:2000}).

\begin{lemma}\label{lemma: plump}
Let $\Omega$ be a $\varrho$-John domain. Then for each $x \in \Omega$ and $r >0$ with $\Omega \setminus B_r(x) \neq \emptyset$, there is $z \in \overline{B_r(x)}$ with $B_{r/(2\varrho)}(z) \subset \Omega$. Particularly, we have $|\Omega| \ge  |B_1(0)| (2\varrho)^{-d}  ({\rm diam}(\Omega))^d$.
\end{lemma}
 
\BBB

\subsection{An elementary lemma about the distance of sets}

In this short section we formulate an elementary lemma which will be needed in the construction of a partition in the proof of Theorem \ref{th: bad part}.

\begin{lemma}\label{lemma: topo}
Let $(F_j)_{j=1}^n \subset \R^2$, $(R_j)_{j=1}^m \subset \R^2$ be families of pairwise disjoint, closed and connected  sets. Let  $F= \bigcup_{j=1}^n F_j$ and $R= \bigcup_{j=1}^m R_j$. Let $G \subset \R^2$ be a closed, path connected set satisfying $G \supset F$ and  
\begin{align}\label{eq: violatXXXXX}
\dist(x, F) \le d \ \ \ \text{for all $x \in G \setminus R$}
\end{align}
for some $d>0$. Then for each $F_j$ we find another component $F_k$ such that (a) $\dist(F_j,F_k) \le 4d$ or (b) there is a corresponding set $R_l$ such that 
\begin{align}\label{eq: desprop}
\dist(F_j,R_l) \le 4d \ \ \ \ \text{and} \ \ \ \ \dist(F_k,R_l) \le 4d.
\end{align}

\end{lemma}

\UUU  In the proof of Theorem \ref{th: bad part}, we will use this lemma to show that two connected components of $F$ can be connected by a line segment of length at most $4d$ (case (a)) or two connected components of $F$ can be connected by a curve lying inside a component $R_l$ and two line segments of length at most $4d$  (case (b)). \EEE

\Proof
Fix $F_{j}$ and assume that  (a) does not hold, i.e., $\dist(F_j,F_k)>4d$ for all $k \neq j$. As $G$ is path connected, we can find $F_k$, $j \neq k$, and a (closed)   curve $\gamma$ in $\overline{G\setminus F}$ such that $F_j \cup F_k \cup \gamma$ is connected. We let $\gamma'$ be a connected component of
$$\gamma \cap \lbrace x \in \R^2: \dist(x,F_j) \ge  2d \rbrace $$
having nonempty intersection with $F_k$.  By \eqref{eq: violatXXXXX} and the fact that (a) does not hold, we then find some $x \in R \cap \gamma'$ with $2d\le \dist(x,F_j) \le 4d$.

Let  $R_l$ be the component  containing $x$ and define $\gamma'' = \gamma' \cap R_l$. If $\gamma''$ intersects $F_k$, we see that $R_l$ has the desired property \eqref{eq: desprop}. Otherwise, we find $y \in \gamma''$ with  $y \in \partial R_l$ and therefore by \eqref{eq: violatXXXXX} and the definition of $\gamma'$ we get some $F_{k'}$, $k' \neq j$ (and possibly $k'\neq k$), so that $\dist(y,F_{k'}) \le d$. Herefrom we again deduce \eqref{eq: desprop} with $F_{k'}$ in place of $F_k$.\eop

\EEE

\subsection{Properties of infinitesimal rigid motions}

In this section we collect some properties of infinitesimal rigid motions. As before $a = a_{A,b}$ stands for the mapping $a(x) = A\,x+b$ with $A \in \R^{2 \times 2}_{\rm skew}$ and $b \in \R^2$. The following lemma is shown in \cite[Lemma 2.3]{FriedrichSolombrino}.

\begin{lemma}\label{lemma: rigid motionXXX}
Let $M>0$, $\delta>0$,  $F \subset \R^2$ bounded, measurable  and  $\psi \colon \R^+\to \R^+$ a continuous nondecreasing function satisfying $\lim_{s \to \infty} \psi(s) = + \infty$. Then there is a constant $C=C(M, \delta, \psi, F)$ such that for every Borel set $E \subset F$ with $|E| \ge \delta$ and every infinitesimal rigid motion $a_{A,b}$ one has 
$$
\int_E \psi(|Ax+b|)\,dx \le M  \ \ \ \Rightarrow \ \ \  |A|+|b|\le C. 
$$

\end{lemma}

In the following we denote by $d(E)$   the diameter of a set $E \subset \R^2$.

\begin{lemma}\label{lemma: rigid motion}
Let $q \in [1,\infty)$. There is a constant $c = c(q)>0$ such that for every Borel set $E \subset \R^2$  and every infinitesimal rigid motion $a=a_{A,b}$ one has for all $x \in E$ 
\begin{align}\label{eq: extra thing}
\begin{split}
(i)& \ \ |A|  \le c  |E|^{-\frac{1}{2}- \frac{1}{q}}\Vert a\Vert_{L^q(E)}, \ \ \ |A|  \le c |E|^{-\frac{1}{2}}\Vert a\Vert_{L^\infty(E)},\\
(ii)& \ \  |A\,x+b| \le c|E|^{-\frac{1}{2}- \frac{1}{q}} d(E)\Vert a\Vert_{L^q(E)}.
\end{split}
\end{align}
\end{lemma} 
 
\Proof Without restriction assume that $A\neq 0$ as otherwise the statement is clear. The assumption $A \in \R^{2 \times 2}_{\rm skew}$ implies that $A$ is invertible and that $|Ay|=\frac{\sqrt2}{2}|A| |y|$ for all $y \in \R^{2}$. Setting $z:= - A^{-1}b$ we find for all $x \in E$
\begin{align}\label{eq: franc}
\tfrac{1}{\sqrt2}|A| |x-z| = |a(x)| \le \Vert a\Vert_{L^\infty(E)}.
\end{align}
\UUU This together with  $|E \setminus B_\rho(z)|\ge \frac{1}{2}|E|$ for $\rho = (\frac{|E|}{2\pi})^{\frac{1}{2}}$ yields $\frac{1}{2}|E| \rho^q|A|^q \le (\sqrt{2})^q\Vert a \Vert^q_{L^q(E \setminus B_\rho(z))}   \le (\sqrt{2})^q\Vert a\Vert^q_{L^q(E)}$ for $q < \infty$. \EEE This implies the first part of \eqref{eq: extra thing}(i). In the case  $q = \infty$ we immediately get $|A|  \le c|E|^{-\frac{1}{2}}\Vert a\Vert_{L^\infty(E)}$ from \eqref{eq: franc}. As there is some $x_0 \in E$ with $|A\,x_0 + b| \le |E|^{-\frac{1}{q}} \Vert a\Vert_{L^q(E)}$, we conclude for all $x \in E$
$$|A\,x + b| \le |A\,x_0+b| + \tfrac{1}{\sqrt{2}}|A||x-x_0| \le |E|^{-\frac{1}{q}} \Vert a\Vert_{L^q(E)} + cd(E)|E|^{-\frac{1}{2}- \frac{1}{q}}\Vert a\Vert_{L^q(E)}.$$
This concludes the proof since $|E|^{\frac{1}{2}} \le d(E)$ by the isodiametric inequality.    \eop

\begin{rem}\label{rem: slice}
{\normalfont

Analogous estimates hold on lines. If, e.g., $\Gamma \subset \R \times \lbrace 0\rbrace$ with $l:=\mathcal{H}^1(\Gamma) < \infty$, then for a constant $c=c(q)>0$ we have for all infinitesimal rigid motions $a=a_{A,b}$
$$l^q|A|^q   \le cl^{-1}\int_{\Gamma} |a(x)|^q \, d\mathcal{H}^1(x).$$

}
\end{rem}

With Lemma \ref{lemma: rigid motion} at hand, we now see that the $L^q$-norm on larger sets can be controlled.

\begin{lemma}\label{lemma: rigid}
Let $q\in [1,\infty)$. Then there is a constant $c=c(q)>0$ such that for all $y \in \R^2$, $R >0$, Borel sets $E \subset Q^y_R :=y+(-R,R)^2$, and $a=a_{A,b}$ one has
$$\Vert a \Vert_{L^q(Q^y_R)}\le  c (R^2|E|^{-1})^{\frac{1}{2} + \frac{1}{q}} \Vert a \Vert_{L^q(E)}.$$
\end{lemma}

\Proof Define $F = Q^y_R \supset E$. We repeat the last estimate of the previous proof for each $x \in F$ and obtain 
$$|A\,x + b| \le   |E|^{-\frac{1}{q}} \Vert a\Vert_{L^q(E)} + cd(F)  |E|^{-\frac{1}{2} - \frac{1}{q}}\Vert a\Vert_{L^q(E)}.$$
This implies with $|F| = 4R^2$ and $d(F) = 2\sqrt{2} R$  
$$\Vert a \Vert^q_{L^q(F)} \le cR^2|E|^{-1}  \Vert a\Vert^q_{L^q(E)} + cR^{2+q} |E|^{-\frac{q}{2} - 1}\Vert a\Vert^q_{L^q(E)}.$$
In view of $R^q|E|^{-\frac{q}{2}} \ge 4^{-\frac{q}{2}}$, this concludes the proof.  \eop

\subsection{(G)SBV and (G)SBD functions}\label{rig-sec: sub, bd}

In this section we collect the definitions and fundamental properties of the function spaces needed in this article. In the following let $\Omega \subset \R^d$ be an open, bounded  set. 

\textbf{SBV- and GSBV-functions.} The space $BV (\Omega; \mathbb R^d)$ consists of the functions $u\in L^1(\Omega; \mathbb R^d)$ such that the distributional gradient $Du$ is a $\R^{d\times d}$-valued finite Radon measure on $\Omega$. $BV$-functions have an approximate differential $\nabla u(x)$ at $\mathcal L^{d}$-a.e.\ $x\in \Omega$ (\cite[Theorem 3.83]{Ambrosio-Fusco-Pallara:2000}) and their jump set $J_u$ is $\mathcal H^{d-1}$-rectifiable in the sense of \cite[Definition 2.57]{Ambrosio-Fusco-Pallara:2000}. The space $SBV (\Omega; \mathbb R^d)$, often abbreviated  hereafter as $SBV(\Omega)$,  of \emph{special functions of bounded variation}  consists of those $u\in BV(\Omega; \mathbb R^d)$ such that
$$
Du= \nabla u\mathcal L^{d}+[u]\otimes \nu_u \mathcal H^{d-1}\lfloor J_u,
$$ 
where $\nu_u$ is a normal of $J_u$ and $[u] = u^+ - u^-$ (the `crack opening') with $u^{\pm}$ being the one-sided limits of $u$ at $J_u$. If in addition $\nabla u \in L^p(\Omega)$ for \BBB some \EEE $1 < p < \infty$ and $\mathcal{H}^{d-1}(J_u) < + \infty$, we write $u \in SBV^p(\Omega)$. Moreover, $(S)BV_{\rm loc}(\Omega)$ denotes the space of functions which belong to $(S)BV(\Omega')$ for every open set $\Omega' \subset \subset \Omega$.

We define the space $GBV(\Omega;\R^d)$  of \emph{generalized functions of bounded variation} consisting of all $\mathcal{L}^d$-measurable functions $u: \Omega \to \R^d$ such that for every $\phi \in C^1(\R^d;\R^d)$ with the support of $\nabla \phi$ compact, the composition $\phi \circ u $ belongs to $BV_{\rm loc}(\Omega)$ (see \cite{DeGiorgi-Ambrosio:1988}). Likewise, we say $u \in GSBV(\Omega;\R^d)$ if $\phi \circ u $ belongs to $SBV_{\rm loc}(\Omega)$  and  $u \in GSBV^p(\Omega;\R^d)$ for $u \in GSBV(\Omega;\R^d)$ if $\nabla u \in L^p(\Omega)$ and $\mathcal{H}^{d-1}(J_u) < + \infty$. As usual we write for shorthand $GSBV(\Omega) := GSBV(\Omega;\R)$. See \cite{Ambrosio-Fusco-Pallara:2000} for the basic properties of these function spaces. 

We now state a version of Ambrosio's compactness theorem in $GSBV$ adapted for our purposes (see e.g. \cite{Ambrosio-Fusco-Pallara:2000, DalMaso-Francfort-Toader:2005}):

\begin{theorem}\label{clea-th: compact}
Let $1 < p < \infty$.  Let $(u_k)_k$ be a sequence in $GSBV^p(\Omega;\R^d)$ with
$$\Vert \nabla u_k \Vert_{L^p(\Omega)}+ \mathcal{H}^{d-1}(J_{u_{k}}) + \Vert u_k \Vert_{L^1(\Omega)} \le C  $$
for some  $C>0$ not depending on $k$. Then there is a subsequence (not relabeled) and a function $u \in GSBV^p(\Omega;\R^d)$ such that $u_k \to u$ a.e., $\nabla u_k \rightharpoonup \nabla u$ weakly in $L^p(\Omega)$. If in addition $\Vert u_k\Vert_{\infty} \le C$ for all $k \in \N$, we find $u \in SBV^p(\Omega) \cap L^\infty(\Omega)$.
\end{theorem}

\textbf{Caccioppoli-partitions.} We say a partition $\mathcal{P} = (P_j)_{j=1}^\infty$ of $\Omega$ is a \textit{Caccioppoli partition} of $\Omega$ if $\sum_{j=1}^\infty \mathcal{H}^{d-1}(\partial^* P_j) < + \infty$, where $\partial^* P_j$ denotes the \textit{essential boundary} of $P_j$ (see \cite[Definition 3.60]{Ambrosio-Fusco-Pallara:2000}).  We say a partition is \textit{ordered} if $|P_i| \ge |P_j|$ for $i \le j$.   We now state a compactness result for ordered Caccioppoli partitions (see \cite[Theorem 4.19, Remark 4.20]{Ambrosio-Fusco-Pallara:2000}).

\begin{theorem}\label{th: comp cacciop}
Let $\Omega \subset \R^d$ open, bounded with Lipschitz boundary.  Let $\mathcal{P}_i = (P_{j,i})_{j=1}^\infty$, $i \in \N$, be a sequence of ordered Caccioppoli partitions of $\Omega$ such that $\sup_i \sum_{j=1}^\infty \mathcal{H}^{d-1}(\partial^* P_{j,i} \cap \Omega) < + \infty$. Then there exists a Caccioppoli partition $\mathcal{P} = (P_j)_{j=1}^\infty$ and a not relabeled subsequence such that  $\chi_{P_{j,i}} \to \chi_{P_j}$ in measure for all $j \in \N$ as $i \to \infty$.
\end{theorem}

\textbf{SBD- and GSBD-functions.} We say that a function $u \in L^1(\Omega; \R^d)$ is in $BD(\Omega)$ if the  symmetrized distributional derivative $Eu := \frac{1}{2}((Du)^T + Du)$ 
is a finite $\R^{d \times d}_{\rm sym}$-valued Radon measure. Likewise, we say $u$ is  a \emph{special  function of bounded deformation} if $Eu$ has vanishing Cantor part $E^c u$. Then $Eu$ can be decomposed as 
\begin{align*}
 Eu  = e(u) \mathcal{L}^d + [u] \odot \nu_u \mathcal{H}^{d-1}|_{J_u},
 \end{align*}
where $e(u)$ is the absolutely continuous part of $Eu$ with respect to the Lebesgue measure $\mathcal{L}^d$, $[u]$, $\nu_u$, $J_u$ as before and $\odot$ denotes the symmetrized tensor product. If in addition $e(u) \in L^q(\Omega)$ \UUU for some $1<q<\infty$ \EEE and $\mathcal{H}^{d-1}(J_u) < +\infty$, we write \UUU $u \in SBD^q(\Omega)$. \EEE For basic properties of this function space we refer to \cite{Ambrosio-Coscia-Dal Maso:1997,  Bellettini-Coscia-DalMaso:98}.

We now introduce the space of \emph{generalized functions of bounded variation}. Observe that it is not possible to follow the approach in the definition of $GSBV$ since for $u \in SBD(\Omega)$ the composition $\phi \circ u$ typically does not lie in $SBD(\Omega)$. In \cite{DalMaso:13}, another approach is suggested which is based on certain properties of one-dimensional slices. 
For fixed $\xi \in S^{d-1}$ we set
\begin{align*}
\Pi^\xi:=\{y \in \R^{d}: y\cdot \xi=0\}\,,&\quad \Omega^\xi_y:=\{t \in \R: y+t\xi\in \Omega\}\hbox{ for }y\in \Pi^\xi\,,\\
&\quad \Omega^\xi:=\{y\in \Pi^\xi: \Omega^\xi_y\neq \emptyset\}\,.
\end{align*}

\begin{definition}
An $\mathcal L^{d}$-measurable function $u:\Omega\to \R^{d}$ belongs to $GBD(\Omega)$ if there exists a positive bounded Radon measure $\lambda_u$ such that, for all $\tau \in C^{1}(\R^{d})$ with $-\frac12 \le \tau \le \frac12$ and $0\le \tau'\le 1$, and all $\xi \in S^{d-1}$, the distributional derivative $D_\xi (\tau(u\cdot \xi))$ is a bounded Radon measure on $\Omega$ whose total variation satisfies
$$
\left|D_\xi (\tau(u\cdot \xi))\right|(B)\le \lambda_u(B)
$$
for every Borel subset $B$ of $\Omega$. A function $u \in GBD(\Omega)$ belongs to the subset $GSBD(\Omega)$   if in addition for every $\xi \in S^{d-1}$ and $\mathcal H^{d-1}$-a.e.\ $y \in \Pi^\xi$, the function $u^\xi_y(t):=u(y+t\xi)$ belongs to $SBV_{\mathrm{loc}}(\Omega^\xi_y)$.
\end{definition}

As before we say   \UUU $u \in GSBD^q(\Omega)$ \EEE if in addition \UUU $e(u) \in L^q(\Omega)$ \EEE and $\mathcal{H}^{d-1}(J_u) < +\infty$. For later purpose we note that there is a compactness result in $GSBD^q(\Omega)$ similar to Theorem \ref{clea-th: compact} (see \cite[Theorem 11.3]{DalMaso:13}).

\textbf{Density results.} We recall a density result in $GSBD^2$   (see \cite{Iurlano:13} and also \cite[Theorem 3, Remark 5.3]{Chambolle:2004}).

\begin{theorem}\label{th: iurlano}
Let $\Omega \subset \R^d$ open, bounded with Lipschitz boundary. Let $u \in GSBD^2(\Omega) \cap L^2(\Omega;\R^d)$.    Then there exists a sequence $(u_k)_k \subset  SBV^2(\Omega; \R^d)$ such that each $J_{u_k}$ is the union of  a finite number of closed connected pieces of $C^{1}$-hypersurfaces, each $u_k$ belongs to
$W^{1,\infty}(\Omega \setminus J_{u_k};\R^d)$ and
\begin{align}\label{eq: iuiu} 
\begin{split}
(i) & \ \ \Vert u_k - u \Vert_{L^2(\Omega)} \to 0,\\
(ii) & \ \ \Vert e(u_k) - e(u) \Vert_{L^2(\Omega)} \to 0,\\
(iii) &  \ \ \mathcal{H}^{d-1}(J_{u_k} \triangle J_u) \to  0.
\end{split}
\end{align}
If in addition $u \in L^\infty(\Omega)$, one can ensure $\Vert u_k \Vert_\infty \le \Vert u \Vert_\infty$ for all $k \in \N$.
\end{theorem}

\begin{rem}\label{rem: iur}
{\normalfont The result together with \cite{Cortesani-Toader:1999} shows that the approximating sequence $(u_k)_k$ can also be chosen such that $J_{u_k}$ is a finite union of closed $(d-1)$-simplices intersected with $\Omega$. In this case, \eqref{eq: iuiu}(iii) has to be replaced by $\mathcal{H}^{d-1}(J_{u_k}) \to \mathcal{H}^{d-1}(J_u).$
}
\end{rem}

\textbf{Korn inequality in SBD.} As a final preparation, we recall a Korn inequality in $SBD^2$ in a planar setting for functions with small jump set.

\begin{theorem}\label{th: kornSBDsmall}
Let $p \in [1,2]$. Then there is a universal constant $c>0$ such that for all squares $Q_\mu = (-\mu,\mu)^2$, $\mu>0$, and all  $u \in SBD^2(Q_{\mu})$ there is a set of finite perimeter $E \subset Q_\mu$ with 
\begin{align}\label{eq: R2main}
\mathcal{H}^1(\partial E) \le c\mathcal{H}^1(J_u), \ \ \ \ |E| \le c(\mathcal{H}^1(J_u))^2
\end{align}
and $A \in \R^{2 \times 2}_{\rm skew}$, $b \in \R^2$ such that
\begin{align}\label{eq: main estmain}
\begin{split} 
(i) & \ \ \Vert u - (A\,\cdot + b) \Vert_{L^p(Q_\mu \setminus E)}\le c\mu^{\frac{2}{p}} \Vert e(u) \Vert_{L^2(Q_\mu)},\\
(ii) & \ \ \Vert \nabla u - A \Vert_{L^p(Q_\mu \setminus E)}\le c\mu^{\frac{2}{p}-1} \Vert e(u) \Vert_{L^2(Q_\mu)}.
\end{split}
\end{align}
Moreover, there is a Borel set $\Gamma \subset \partial Q_\mu$ such that
\begin{align}\label{eq: main estmainXX}
\mathcal{H}^1(\Gamma) \le c\mathcal{H}^1(J_u), \ \ \ \ \int_{\partial Q_\mu \setminus \Gamma} |Tu - (A\,x + b)|^2\, d\mathcal{H}^1(x) \le c\mu\Vert e(u) \Vert^2_{L^2(Q_\mu)},
\end{align}
where $Tu$ denotes the trace of $u$ on $\partial Q_\mu$. 
\end{theorem}

The fact that the constant is independent of $\mu$ and $p$ follows from a standard scaling argument and H\"older's inequality. More generally, the result holds on connected, bounded Lipschitz sets $\Omega$  for a constant $C$ also depending on $\Omega$. \UUU Note that the result is indeed only relevant for functions with sufficiently small jump set, as otherwise one can choose $E = Q_\mu$, $\Gamma= \partial Q_\mu$ and \eqref{eq: main estmain}-\eqref{eq: main estmainXX} trivially hold. \EEE 

\Proof The result stated in \eqref{eq: R2main}-\eqref{eq: main estmain} has first been derived in \cite{Friedrich:15-3} for $p \in [1,2)$ and was then improved in \cite{Conti-Iurlano:15.2} by showing the estimate for $p=2$.  In \cite{Chambolle-Conti-Francfort:2014}, the trace estimate has been established. We need to confirm that in both estimates \eqref{eq: main estmain}, \eqref{eq: main estmainXX} one can indeed take the same infinitesimal rigid motion. 

Let $E$ and $a_{A,b}$ be given such that  \eqref{eq: R2main}-\eqref{eq: main estmain} hold for $p=2$. By \cite{Chambolle-Conti-Francfort:2014} we find Borel sets $F \subset Q_\mu$, $\Gamma \subset \partial Q_\mu$ and $a' = a_{A',b'}$ such that \eqref{eq: main estmainXX} holds with $a'$ in place of $a$ and $|F|\le c(\mathcal{H}^1(J_u))^2$ as well as $\Vert u - a' \Vert_{L^2(Q_\mu \setminus F)} \le c\mu\Vert e(u) \Vert_{L^2(Q_\mu)}$. We can assume that $\mathcal{H}^1(J_u)$ is so small that   $|E \cup F| \le \BBB 2 \EEE c(\mathcal{H}^1(J_u))^2 \le  \frac{1}{2} |Q_\mu| $ since otherwise passing to a larger $c>0$ in \eqref{eq: R2main} we could choose $E = Q_\mu$ and the statement was trivially satisfied. Therefore, we have $\Vert a - a' \Vert_{L^2(Q_\mu \setminus (E \cup F))} \le c\mu\Vert e(u) \Vert_{L^2(Q_\mu)}$. By Lemma \ref{lemma: rigid} we get $\Vert a - a' \Vert_{L^2(Q_\mu)} \le c\mu\Vert e(u) \Vert_{L^2(Q_\mu)}$   and thus by Lemma  \ref{lemma: rigid motion} $\Vert a-a' \Vert_{L^\infty(Q_\mu)} \le c\Vert e(u) \Vert_{L^2(Q_\mu)}$. This yields \eqref{eq: main estmainXX} for  the function $a$. \eop

 \section{A piecewise Korn inequality up to a small exceptional set}\label{sec: small set}

As a key step for the proof of Theorem \ref{th: main korn}, we first derive a piecewise Korn inequality up to a small exceptional set. By iterative application of this result in Section \ref{sec: main-proof}   we then derive the main inequality. In this section we will consider configurations on a square with regular jump set consisting of a finite number of segments and defer the general case also to Section \ref{sec: main-proof}, where we  will make use of a density result (see Theorem \ref{th: iurlano} and Lemma \ref{lemma: iurlano*}). 

For $\Omega \subset \R^d$ open, bounded we let $\mathcal{W}(\Omega) \subset SBV^2(\Omega)$ be the functions $u$ such that $J_{u} = \bigcup^n_{j=1} \Gamma^u_j$  is a finite union of closed $(d-1)$-simplices intersected with $\Omega$ and  $u|_{\Omega\setminus J_u} \in W^{1,\infty}(\Omega \setminus J_u)$. (In dimension $d=2$, each $\Gamma^u_j$ is a closed segment.) In the following we say $(P_j)_{j=1}^n$ is a partition of $\Omega$ if the sets are open, pairwise disjoint and satisfy $|\Omega \setminus \bigcup^n_{j=1} P_j| = 0$. For convenience we often write $\Omega = \bigcup^n_{j=1} P_j$ although the identity only holds up to a set of negligible $\mathcal{L}^2$-measure.

\begin{theorem}\label{th: korn-small set}
Let  $p \in [1,2)$ and  $\theta>0$. Then there is a universal constant $c>0$ and some $C=C(p,\theta)>0$ such that the following holds: 

(1) For  each square $Q_\mu = (-\mu,\mu)^2$ for $\mu >0$ and each $u \in \mathcal{W}(Q_\mu)$   one finds an exceptional set $E \subset Q_\mu$ with
\begin{align}\label{eq: except}
|E| \le c\mu   \theta^2 \, \mathcal{H}^1(J_u), \ \ \ \ \ \mathcal{H}^1(\partial E) \le C\mathcal{H}^1(J_u),
\end{align} 
a partition $Q_\mu = \bigcup^n_{j=1} P_j$ and corresponding infinitesimal rigid motions $(a_j)_j = (a_{A_j,b_j})_j$ such that  the function $v:= u - \sum_{j=1}^n a_j \chi_{P_j}$ satisfies
\begin{align}\label{eq: small set main}
\begin{split}
(i) & \ \ \mathcal{H}^1\big(\bigcup\nolimits^n_{j=1}\partial P_j \cap Q_\mu \big) \le  C\mathcal{H}^1\big(\bigcup\nolimits^n_{j=1}\partial P_j \cap J_u\big),\\
(ii) &\ \  \Vert \nabla v \Vert_{L^p(Q_\mu \setminus E)} \le C\mu^{\frac{2}{p}-1} \Vert  e(u) \Vert_{L^2(Q_\mu)}.\\
\end{split}
\end{align}

(2) One can choose $E$, $(P_j)_{j=1}^n$  and $v$ such that we also have
\begin{align}\label{eq: small set mainXXX}
\Vert v \Vert_{L^\infty(Q_\mu)} \le C\Vert u \Vert_{L^\infty(Q_\mu)}.
\end{align}
\end{theorem}

Note that \eqref{eq: small set main}(ii) is similar to \eqref{eq: main korn2}(ii) with the correct scaling, \UUU where the estimate only holds outside of an exceptional set  which, roughly speaking, covers the jump set $J_u$ in a specific sense. \EEE Since \eqref{eq: main korn2}(i) will eventually follow from \eqref{eq: main korn2}(ii) by Theorem \ref{th: main poinc}, it is not necessary at this stage to derive an analog of \eqref{eq: main korn2}(i). However, we establish \eqref{eq: small set mainXXX} which may be of independent interest.

Estimate \eqref{eq: small set main}(i) differs from \eqref{eq: main1} in the sense that the length of the boundary of the partition  can be controlled solely by the part of $J_u$ contained in $\bigcup_{j=1}^n\partial P_j$. This property will be crucial for the iterative application of the arguments in Section \ref{sec: main-proof} since hereby a blow up of the length of the boundary of the partition  can be avoided. Also the fact the volume of the exceptional set in \eqref{eq: except} is controlled in terms of a (small) parameter $\theta>0$  will be convenient for the subsequent analysis. Note that therefore it will not be  restrictive to concentrate on the case  $\theta \le 1$ and
\begin{align}\label{eq: concentration}
\mu\theta^{2} \le \mathcal{H}^1(J_u) \le 2\sqrt{2}\mu\theta^{-2}.
\end{align} 
In fact, for functions with smaller jump sets, Theorem \ref{th: korn-small set} follows directly  from Theorem \ref{th: kornSBDsmall} for a partition only consisting of one element, and if the jump set is too large, one can choose $E = Q_\mu$ and \eqref{eq: small set main} trivially holds.

\BBB Let us mention that for the iterative scheme in  Section \ref{sec: main-proof} we will make use of a \emph{refined version} of Theorem \ref{th: korn-small set} which particularly takes a  specific structure of the exceptional set $E$ into account. As its formulation needs some of the notions and concepts developed during the proof of Theorem \ref{th: korn-small set}, we defer the exact statement to Lemma \ref{rem: large components} in Section \ref{sec: refined}. For a detailed outline of the iterative application of  Lemma \ref{rem: large components}, which is necessary to derive the main result without an exceptional set, we refer to the beginning of Section \ref{sec: main*}. 

We now concentrate on the proof  of Theorem \ref{th: korn-small set}. We start by describing the basic strategy. \EEE We first identify so-called \emph{bad squares} $Q$ of various mesoscopic sizes where $\mathcal{H}^1(J_u \cap Q)$ is `too large', i.e., comparable to the diameter of the square.  Suitably combining regions formed by such bad squares with segments, we will be able to construct an \BBB auxiliary \EEE partition of $Q_\mu$  consisting of simply connected sets. Moreover, we obtain a corresponding Whitney covering such that $\mathcal{H}^1(J_u \cap Q) \ll d(Q)$ for the squares of the covering (see Section \ref{sec: 3sub1}). Then applying the Korn inequality for functions with small jump set (Theorem \ref{th: kornSBDsmall}) on each of these squares, we can replace  the configuration $u$ by a \BBB \emph{piecewise smooth approximation} \EEE $\bar{u}$ which is smooth on each component of  the auxiliary partition  and whose distance from $u$ can be controlled outside a small exceptional set (see Section \ref{sec: 3sub2}). Afterwards, we use Theorem \ref{th: main part2} to find another refined partition consisting of John domains whose John constant can be uniformly controlled. \BBB We then apply the Korn inequality in John domains (Theorem \ref{th: kornsobo}) for the auxiliary function $\bar{u}$, which is smooth in each component of the partition. Hereby, as an intermediate step,   we  can establish an estimate on $\nabla \bar{u}$. Finally,  Theorem \ref{th: korn-small set} follows from  the fact that $\bar{u}$ is a suitable approximation of $u$, in particular  $\nabla u - \nabla \bar{u}$ can be controlled outside a small exceptional set  (see Section \ref{sec: 3sub5}). \EEE

Before we start with the construction of the auxiliary partition, we introduce some general notation which will be used throughout the paper. For $s>0$ we partition $\R^2$ up to a set of measure zero into squares $Q^s(p) = p + s(-1,1)^2$ for $p \in s(1, 1) + 2s\Z^2$ and write  
\begin{align}\label{eq: enlarged squaresXXX}
 \tilde{\mathcal{Q}}^s = \lbrace Q =Q^s(p): \ p \in s(1, 1) + 2s\Z^2 \rbrace.
 \end{align}
Let $\theta>0$ be given with $\theta \in 2^{-\N}$. For all $i \in \N_0$ we define  $s_i = \mu {\theta}^{i}$ and we consider the coverings of dyadic squares \BBB $\mathcal{Q}^i := \tilde{\mathcal{Q}}^{s_i}$. All following statements  will be formulated under the assumption that $\theta>0$ is   \emph{small}. More precisely, this means that there exists a universal constant $\theta_0>0$ such that all statements hold for $\theta \in 2^{-\N} \cap (0,\theta_0]$. \EEE

 For each square $Q$ of the above form we also introduce the corresponding enlarged squares $Q  \subset Q' \subset Q'' \subset Q'''$ defined by 
\begin{align}\label{eq: enlarged squares}
 Q' = \tfrac{3}{2}Q, \ \ \ Q'' = 3Q, \ \ \   Q''' = 5Q,
 \end{align}
where $\lambda Q$ denotes the square with the same center and orientation and $\lambda$-times the sidelength of $Q$. By $d(A)$ we again indicate the diameter of a set  $A \subset \R^2$ and  by $\dist(A,B)$ we denote the Euclidian distance between $A,B \subset \R^2$.  Finally, we introduce the \emph{saturation} of a bounded set $A \subset \R^2$ by  $\sat(A) = \R^2 \setminus A'$, where $A'$ denotes the (unique)  unbounded connected component of $\R^2 \setminus A$. Loosely speaking, $\sat(A)$ arises from $A$ by `filling the holes of $A$'.

\subsection{Construction of an auxiliary  partition}\label{sec: 3sub1}

\BBB In this section we start the proof of Theorem \ref{th: korn-small set} by  constructing an auxiliary partition consisting of simply connected sets and a covering of Whitney-type, both associated to the jump set \UUU of the function. At this stage the constructions only involve the structure and geometry of the jump set. We therefore formulate the result in terms of a general set $J = \bigcup_{j=1}^n \Gamma_j \subset Q_\mu$ which consists of a finite number of closed segments and satisfies (cf. \eqref{eq: concentration})
\begin{align}\label{eq: concentration-new}
\mu\theta^{2} \le \mathcal{H}^1(J) \le 2\sqrt{2}\mu\theta^{-2}.
\end{align} 
(For convenience, we will still call $J$ the \emph{jump set} in the following.) \EEE Estimates for $\nabla u$ will only play a role starting from Section \ref{sec: 3sub2}. We begin with some definitions. \EEE

\smallskip

\emph{Auxiliary jump set:} To ensure that we provide an estimate which is valid up to a small exceptional set (cf. \eqref{eq: except}), it will be convenient to decompose $Q_\mu$ into smaller squares. To this end, we introduce the \emph{auxiliary jump set} 
\begin{align}\label{eq: Ju*}
\UUU J^* \EEE := \UUU J \EEE \cup J_0 \ \ \  \text{with  } J_0 = \bigcup\nolimits_{Q \in \mathcal{Q}^{7}, Q \subset Q_\mu} \partial Q
\end{align}
and observe that together with \UUU \eqref{eq: concentration-new} \EEE we find 
\begin{align}\label{eq: Ju*XX}
\mathcal{H}^1(\UUU J^* \EEE ) \le \mathcal{H}^1( \UUU J \EEE ) + c\theta^{-7}\mu \le c\theta^{-9}\mathcal{H}^1( \UUU J \EEE )
\end{align}
for a universal constant $c>0$. 

\smallskip

\begin{figure}[H]
\centering
\begin{overpic}[width=0.88\linewidth,clip]{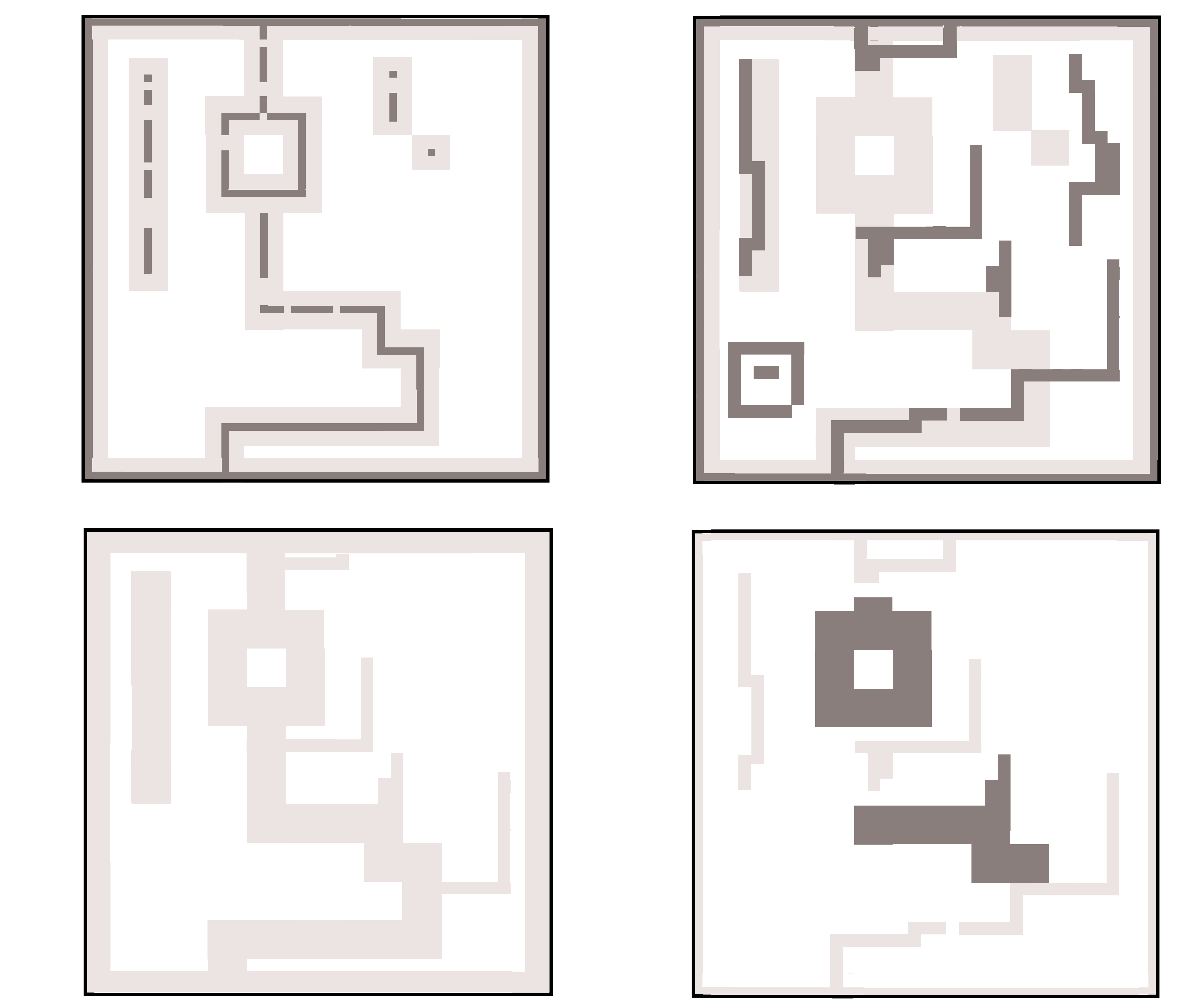}
\put(8,297){{$(a)$}}
\put(195,297){{$(b)$}}
\put(8,137){{$(c)$}}
\put(195,137){{$(d)$}}

\put(132,280){{$B^{I-1}_1$}}
\put(125,286) {\line(1,0){6}}

\put(93,230){{$B^{I-1}_2$}}
\put(86,236) {\line(1,0){6}}

\put(45,205){{$B^{I-1}_3$}}
\put(48,215) {\line(0,1){6}}

\put(36,185){{$J_0$}}
\put(26,191) {\line(1,0){10}}

\put(251,191){{$B^{I}_2$}}
\put(237,195) {\line(1,0){13}}

\put(234,170){{$B^{I}_1$}}
\put(234,177) {\line(0,1){6}}

\put(310,245){{$B^{I}_4$}}
\put(323,251) {\line(1,0){6}}

\put(314,220){{$B^{I}_3$}}
\put(308,226) {\line(1,0){6}}

\put(38,45){{$D^{I-1}_1$}}
\put(46,57) {\line(0,1){7}}

\put(121,100){{$D^{I-1}_2$}}
\put(114,106) {\line(1,0){6}}

\put(291,112){{$R^{I-1}_2$}}
\put(284,118) {\line(1,0){6}}

\put(315,67){{$R^{I-1}_1$}}
\put(308,73) {\line(1,0){6}}

\end{overpic}
\caption{(a) \UUU A square of $\mathcal{Q}^7$ contained in \EEE $Q_\mu$ is depicted with $B^{I-1} = \bigcup_{k=1}^3 B^{I-1}_k$ in light gray and the squares $\mathcal{Q}_{\rm bad}^{I-1}$ in dark gray. Note that $B^{I-1}_1 \in \mathcal{B}^{I-1}_{\rm iso}$ and $J_0 \subset B^{I-1}_2$.  \ (b) In light gray we see $B^{I-1}$ and in dark gray $B^I$. The sets $B^I_1, B^I_2, B^I_3$ are contained in $\mathcal{B}^I_{\rm iso}$ with $B_2^I\subset \sat(B_1^I)$. We suppose $B_4^I \notin \mathcal{B}^I_{\rm iso}$ and observe $B_4^I \cap B^{I-1} = \emptyset$. \ (c) The set $E^{I-1} = D^{I-1}_1 \cup D^{I-1}_2$ is sketched. Note that $B_4^I \cap E^{I-1} = \emptyset$ although $B_4^I \notin \mathcal{B}^I_{\rm iso}$.  \ (d) In light gray \BBB $E^I$ \EEE is depicted and in dark gray $R^{I-1} = R^{I-1}_1 \cup R^{I-1}_2$. \BBB Observe that $E^I \subset E^{I-1}$. \EEE } \label{korn1}
\end{figure}
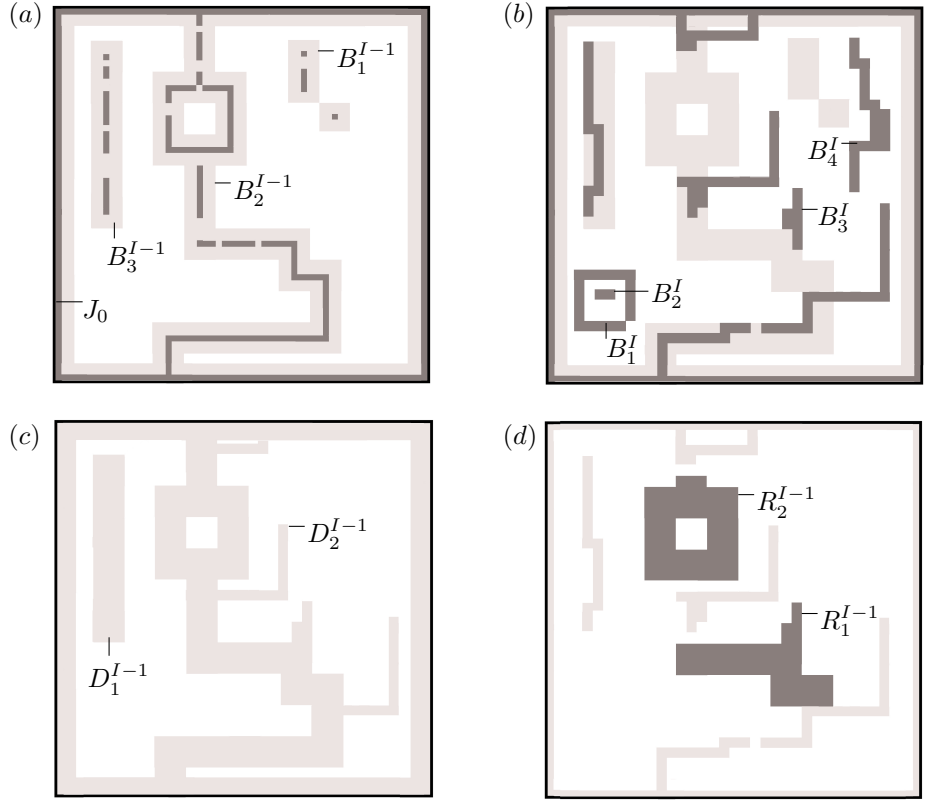

\emph{Bad squares:} One of the main strategies will be the application of Theorem \ref{th: kornSBDsmall}. As the result only provides an estimate for functions with small jump set, we introduce for  $i \ge 1$  the set of \emph{bad squares} 
\begin{align}\label{eq: bad1}
\begin{split}
\mathcal{Q}^{i}_{\rm bad} &= \lbrace Q \in \mathcal{Q}^i: \mathcal{H}^1( \UUU J^* \EEE \cap Q') \ge \theta^{3} s_i \rbrace.
\end{split}
\end{align}
\BBB The squares $\mathcal{Q}^i \setminus \mathcal{Q}^{i}_{\rm bad}$ will sometimes be called \emph{good squares}. \EEE  In particular, note that each square with $Q' \cap J_0 \neq \emptyset$  satisfies $\mathcal{H}^1(Q' \cap \UUU J^* \EEE) \ge 2 s_i$ and thus lies in $\mathcal{Q}^{i}_{\rm bad}$. \BBB For $i \ge 8$ \EEE we let 
\begin{align}\label{NNNNN} 
B^i = \bigcup\nolimits_{Q \in \mathcal{Q}^i_{\rm bad}} \overline{Q'''}
\end{align}
and observe that $B^i$ is \BBB a \EEE union of squares in $\mathcal{Q}^i$ up to a set of negligible measure (see Figure \ref{korn1}(a) below). \BBB Since $B^i$ is defined in terms of  enlarged squares, the squares $Q \in \mathcal{Q}^i$ lying inside $B^i$ and touching the boundary  $\partial B^i$ are good squares. This property will be crucial for the construction of the Whitney covering in Theorem \ref{th: bad part}. \EEE Moreover, we will frequently use the following lemma about unions of `bad sets'.

\begin{lemma}\label{lemma: cup}
There is a universal constant $c>0$ such that for all $8 \le i_1 \le i_2 < \infty$ and for all $\mathcal{R}^j \subset \mathcal{Q}^j_{\rm bad}$ we have for $R := \bigcup_{j=i_1}^{i_2} \bigcup_{Q \in \mathcal{R}^j} \overline{Q'''}$
$$\mathcal{H}^1(\partial R ) \le c\theta^{-3}\mathcal{H}^1( \UUU J^* \EEE \cap R).$$
Moreover, there is a set $\Gamma_R \supset \partial R$ being a finite union of closed segments  with  $\mathcal{H}^1(\Gamma_R) \le c\theta^{-3}\mathcal{H}^1( \UUU J^* \EEE \cap R)$ such that for each connected component $R_k$ of $R$ the set $\Gamma_R \cap R_k$ is connected.  
\end{lemma}

\Proof For $i_1 \le j \le i_2$ let $\hat{\mathcal{R}}^j = \lbrace Q \in \mathcal{R}^j: \overline{Q'''} \not\subset \bigcup^{j-1}_{i = i_1} \bigcup_{Q \in \mathcal{R}^i} \overline{Q'''}\rbrace$ and $\hat{R}^j = \bigcup_{Q \in \hat{\mathcal{R}}^j} Q'$. Note that in view of \eqref{eq: enlarged squares}, for $\theta$ small the sets $(\hat{R}^j)_{j=i_1}^{i_2}$ are pairwise disjoint and by \eqref{eq: bad1} we have for all $i_1 \le j  \le i_2$
\begin{align*}
\# \hat{\mathcal{R}}^j \le \theta^{-3}s_j^{-1}\sum\nolimits_{Q \in \hat{\mathcal{R}}^j} \mathcal{H}^1(\UUU J^* \EEE \cap Q') \le  4\theta^{-3}s_{j}^{-1} \mathcal{H}^1(\UUU J^* \EEE \cap \hat{R}^j),
\end{align*}
where we used that each $x \in \R^2$ is contained in at most four $Q'$, $Q \in \mathcal{Q}^j$. Since $\bigcup\nolimits_{j=i_1}^{i_2} \hat{R}^j  \subset R$,   $(\hat{R}^j)_j$ are pairwise disjoint and $\mathcal{H}^1(\partial Q''') = 40s_j$ for $Q \in \hat{\mathcal{R}}^j$, we then derive
\begin{align*}
\mathcal{H}^1(\partial R ) &\le \sum\nolimits^{i_2}_{j=i_1} \sum\nolimits_{Q \in \hat{\mathcal{R}}^j} \mathcal{H}^1(\partial Q''') \le \sum\nolimits^{i_2}_{j=i_1}c\theta^{-3}  \mathcal{H}^1(\UUU J^* \EEE \cap \hat{R}^j)  \le c\theta^{-3}\mathcal{H}^1( \UUU J^* \EEE \cap R).
\end{align*}
To see the second part of the statement, we define $\Gamma_R = \bigcup^{i_2}_{j=i_1} \bigcup_{Q \in \hat{\mathcal{R}}^j} \partial Q'''$ and note that $\Gamma_R \cap R_k$ is connected for each connected component $R_k$ of $R$. \eop  

\smallskip

\emph{Components and isolated components:} Let $\mathcal{B}^i := (B^i_k)_k$ be the connected components of $B^i$ for \BBB $i \ge 8$. \EEE By the remark below \eqref{eq: bad1}, for each \BBB $i \ge 8$ \EEE there is a component $B^i_k$ with $J_0 \subset B^i_k$ (see also Figure \ref{korn1}(a) below). \UUU A basic strategy in the construction of the auxiliary partition  will lie in inductively combining such connected components with line segments. Fix $r \in (0,\frac{1}{8})$. \EEE  We say a set $B_k^i \in \mathcal{B}^i$ is an \emph{isolated component} if   
\begin{align}\label{eq: crack bound}
d(B_k^i) \le   \theta^{-ir} s_{i},
\end{align}
\BBB where, as before, $d(\cdot)$ denotes the diameter of a set. \EEE  Let $\mathcal{B}^i_{\rm iso} \subset \mathcal{B}^i$ be the subset consisting of the isolated components. \BBB Due to their small diameter, these components will play a special role in the construction of the auxiliary partition. For more details, in particular concerning the necessity to introduce these objects, we refer to Lemma \ref{lemma: bad sets} below and the comments thereafter. At this point, let us just introduce the associated sets    \EEE 
\begin{align}\label{eq: V}
U_i :=  \sat\Big(\bigcup\nolimits_{B^i_k \in \mathcal{B}^i_{\rm iso} } B^i_k\Big), \ \ \BBB i \ge 8, \EEE
\end{align}
where $\sat(\cdot)$ was defined below  \eqref{eq: enlarged squares}. \UUU Using the saturation ensures that the connected components of $U_i$ are simply connected sets. \EEE

\smallskip

\emph{Definition of a final iteration step:} \BBB For the iterative scheme applied below, we have to introduce a final step (or equivalently a smallest length scale) at which we stop the construction. \EEE \UUU Recall that $J^{*}$ is  a finite union of closed segments \EEE and denote the connected components of $\UUU J^* \EEE$ by $(\Gamma^*_j)_{j=1}^N$, each of which being a union of finitely many closed segments. \EEE Choosing $I \in \N$ sufficiently large, we see that for $i\ge I$ each set $(B^i_k)_k$ contains exactly one connected component of $\UUU J^* \EEE$. In particular, we can assume that  $I \BBB =I(u,\theta, \UUU r \EEE)  \in \N$ is so large that 
\begin{align}\label{eq: distance}
\begin{split}
(i)& \ \ \dist(\Gamma^*_{j_1}, \Gamma^*_{j_2}) \ge \theta^{-1} s_I \ \ \ \text{for} \ \ \ 1 \le j_1 < j_2 \le N,\\
(ii)& \ \ \mathcal{B}^i_{\rm iso} = \emptyset \ \ \ \text{for} \ \ i \ge I,\\
(iii) & \  \ \BBB N \le s_I^{-1} \mathcal{H}^1(J^{*}_u). \EEE
\end{split}
\end{align} 
Indeed, for (ii) we observe $d(B^i_k) \ge d(\UUU J^* \EEE \cap B_k^i) \ge  \min\nolimits_{j=1,\ldots,N} \, d(\Gamma_j^*)>0$ and $  \theta^{-ir} s_{i} \to 0$ for $i \to \infty$.

We now formulate the main result of this section. We show that $Q_\mu$ can be partitioned into simply connected components $(P'_j)_j$ with a corresponding covering of Whitney-type in terms of squares in $\bigcup_{i \ge 8} \mathcal{Q}^i$ such that the jump set is controllable in the squares not associated to the isolated components $U_i$.

\begin{theorem}\label{th: bad part}
Let $\mu>0$, $\theta>0$ \BBB small \EEE and \UUU $r \in (0,\frac{1}{8})$. \EEE Then there is a universal constant $c>0$ and a constant  $C=C(\theta,r)>0$  such that for each \UUU $J \subset Q_\mu$ which is a finite union of closed segments with \UUU \eqref{eq: concentration-new} \EEE there is a partition $(P'_j)_{j=1}^m$ of $Q_\mu$ consisting of open, simply connected sets satisfying
\begin{align}\label{eq: bad length}
\mathcal{H}^1\big(\bigcup\nolimits_{j=1}^m \partial P'_j \big) \le  C \mathcal{H}^1( \UUU J \EEE )
\end{align}
with the following properties: the set $\bigcup^m_{j=1} \partial P'_j$ is a finite union of closed segments. Moreover, there is a covering $\mathcal{C}   \subset \bigcup^{\infty}_{i=8} \mathcal{Q}^i$ of  $Q_\mu$ \UUU (up to a set of negligible measure) \EEE  with pairwise disjoint dyadic squares,  a closed set $Z \subset Q_\mu$ and an index  $I=I(u,\theta,\UUU r \EEE ) \in \N$   such that
\begin{align}\label{eq: U}
Z\subset \bigcup\nolimits_{i = 8}^I U_i,
\end{align}
for $U_i$ as in \eqref{eq: V}, and with $\mathcal{C}^i := \mathcal{C} \cap \mathcal{Q}^i$  we have
\begin{align}\label{eq: main pro}
(i) & \ \  Q_\mu \setminus \big( \UUU J \EEE \cap \bigcup\nolimits^m_{j=1} \partial P'_j \big) \subset \bigcup\nolimits_{Q \in \mathcal{C}} Q' \subset Q_\mu, \notag \\
(ii) & \ \ Q_1' \cap Q_2' \neq \emptyset \ \text{ for } \ Q_1,Q_2 \in \mathcal{C} \  \Rightarrow \ \theta d(Q_1) \le d(Q_2) \le \theta^{-1}d(Q_1), \notag \\
(iii) & \ \ \# \lbrace Q \in \mathcal{C}: x \in Q'\rbrace \le 12 \ \ \text{ for all } \ x \in Q_\mu,\\
(iv) & \ \ \mathcal{H}^1\big( \UUU J \EEE \cap Q'\big) \le \theta^{2}  s_i \ \ \text{ for all } \ Q \in \mathcal{C}^i \ \text{with} \ Q'' \not\subset Z,  \notag  \\
(v) & \ \ \mathcal{H}^1\big( \UUU J \EEE \cap Q'\big) \le \theta^{2}  s_i \ \ \text{ for all } \ Q \in \mathcal{C}^i    \text{ with } Q'' \cap \hat{Q}'' \neq \emptyset  \text{ for some } \hat{Q} \in \mathcal{C}^{i-1} \cup \mathcal{C}^{i+1}, \notag \\
(vi) & \ \ \UUU J \EEE   \   \cap  \    \bigcup\nolimits^\infty_{i \ge I+1}\bigcup\nolimits_{Q \in \mathcal{C}^i} Q' \ = \emptyset,\notag
\end{align}
where $Q',Q''$ denote the enlarged squares corresponding to $Q \in \mathcal{C}$ (cf. \eqref{eq: enlarged squares}).

\end{theorem}

The fact that the auxiliary partition consists of simply connected sets is essential since  then  Theorem \ref{th: main part2} is applicable.  \BBB Hereby, we will be able  to find a refined partition (the one appearing in Theorem \ref{th: korn-small set})    consisting of John domains with uniformly controlled John constant, see Section \ref{sec: 3sub5}.  In this context, also note that \eqref{eq: bad length} is crucial for the proof of \eqref{eq: small set main}(i). \EEE

Observe that \eqref{eq: main pro}(ii),(iii) are the typical properties of a Whitney covering. \UUU The main relation between the covering and the partition is given in \eqref{eq: main pro}(i). Specifically, the (enlarged) squares  \EEE of $\mathcal{C}$  do possibly not cover the part of $\UUU J \EEE$ contained in the boundary of the partition. In particular, \UUU note that the squares $\mathcal{C}$ cover $Q_\mu$ only up to set of negligible measure. For simplicity, we will nevertheless call $\mathcal{C}$ \emph{a covering of $Q_\mu$}, referring to \eqref{eq: main pro}(i) for the exact property. \EEE

The crucial conditions \eqref{eq: main pro}(iv)-(vi) state that outside of $Z$ the jump set in each square is small compared to its diameter. \BBB (In particular, \eqref{eq: main pro}(vi) shows that squares of sufficiently small size do not intersect the jump at all, which is related to the definition of $I$ in \eqref{eq: distance}.) \UUU Note that, to achieve this, it is in general indispensable that the enlarged squares do not cover certain parts of the jump set $J$, cf. \eqref{eq: main pro}(i).  \EEE The conditions will allow us to apply the Korn inequality  for functions with small jump set (Theorem \ref{th: kornSBDsmall}) and to approximate  the configuration $u$  on each square outside of $Z$ by an infinitesimal rigid motion. By using a partition of unity associated to $\mathcal{C}$ , this will be the starting point to construct a piecewise smooth approximation  of $u$, \UUU which is smooth on each component $P'_j$ \EEE  (see Section \ref{sec: 3sub2}).  

In the set $Z$  we will have to apply other techniques. \UUU Note that its definition is related to the isolated components $U_i$, but in general the inclusion \eqref{eq: U} is strict. \EEE  Before we proceed with the proof of Theorem \ref{th: bad part}, we give a more precise meaning to $Z$ and \BBB discuss its \UUU role in the construction of the partition. \EEE

\begin{lemma}\label{lemma: bad sets}
Let be given the situation of Theorem \ref{th: bad part} with $(P'_j)_j$ and $Z$ as in \eqref{eq: bad length}-\eqref{eq: U}. Then the covering $\mathcal{C}$ can be chosen such that \eqref{eq: main pro} holds and that there is a partition $Z = \bigcup_{i=8}^I Z^i$ into pairwise disjoint, closed sets with the property that each $Z^i$ is \BBB a \EEE union of squares in $\mathcal{C}^i = \mathcal{C} \cap \mathcal{Q}^i$ up to a set of measure zero. We have 
\begin{align}\label{eq: main proXXX1}
\begin{split}
(i)&\ \ \big|Z\big| \le c\mu\theta^5\mathcal{H}^1( \UUU J \EEE ), \ \ \ \
 |Z^i| \le C\theta^{-ir}s_i \mathcal{H}^1( \UUU J \EEE ),\\
 (ii)& \ \ \sum\nolimits_{i=8 }^I \mathcal{H}^1(\partial Z^i) \le C \mathcal{H}^1( \UUU J \EEE )
\end{split}
\end{align}
for a universal $c>0$ and  $C=C(\theta)>0$ only depending on $\theta$. Denote by   $(X^i_k)_k$ the connected components of $Z^i$ and let $\mathcal{X}^i_k  = \lbrace Q \in \mathcal{C}^i: Q \subset X^i_k, \ \partial Q \cap \partial X^i_k \neq \emptyset \rbrace$. \UUU Each $X^i_k$ is a simply connected set and we have \EEE
\begin{align}\label{eq: main proXXX}
\begin{split}
(i)& \ \ d(X^i_k) \le \theta^{-ir}s_i, \ \ \ \# \mathcal{X}^i_k \le c  \theta^{-2ir}, \\
(ii) & \ \ Q' \cap Z \subset X^i_k \ \ \ \text{for all} \ \ \ Q  \in \mathcal{X}^i_k,\\
(iii) & \ \ \sum\nolimits_k\# \mathcal{X}^i_k \le C \mathcal{H}^1(\UUU J \EEE ) s_{i}^{-1}.
\end{split}
\end{align}

\end{lemma}

\BBB 
Since  the jump set in squares of $\mathcal{C}$ lying in $Z$  is possibly large, we cannot use Theorem \ref{th: kornSBDsmall} to approximate $u$ on $Z$ by infinitesimal rigid motions. Consequently, $Z$ will be a part of the exceptional set of Theorem \ref{th: korn-small set}. The control on the area and the boundary of $Z$ stated in \eqref{eq: main proXXX1} will be needed to show \eqref{eq: except}. The fact that each connected component of $Z$ is covered by squares in $\mathcal{C}$ of a specific size will be important for the  formulation of a refined version of Theorem \ref{th: korn-small set} (see Lemma \ref{rem: large components}) and its application in Section \ref{sec: main-proof}.

Let us comment on the necessity of $Z$. The construction of the auxiliary partition in Theorem \ref{th: bad part} is based on a suitable combination of different connected components of $B^i$, $i \ge 8$,  with small segments (recall \eqref{NNNNN}). In this context, it is essential that the total length of the added segments is controlled in terms of $\mathcal{H}^1( \UUU J \EEE )$, cf. \eqref{eq: bad length}. As the total length of the added segments crucially depends on the number of connected components  which have to be combined, we are only allowed to combine components with sufficiently large diameter since their number can be  controlled suitably.  Components with diameter smaller than $\theta^{-ir}s_i$ (see \eqref{eq: crack bound} and \eqref{eq: main proXXX}(i)), called \emph{isolated components}, are neglected in the construction of the partition.

On the one hand, this neglection allows us to obtain \eqref{eq: bad length}. On the other hand, inside these isolated components we cannot confirm the crucial property that the jump set in squares of the Whitney covering  is small (cf. \eqref{eq: main pro}(iv)). Consequently, Theorem \ref{th: kornSBDsmall} cannot be applied inside $Z$. Therefore, in Section \ref{sec: 3sub2} we will obtain a piecewise smooth approximation $\bar{u}$ of $u$ first only outside $Z$. Then other techniques will have to be applied to extend $\bar{u}$ inside $Z$. In this context, it will be fundamental that the diameter of each connected component of $Z$ is controlled (see \eqref{eq: main proXXX}(i)), that the components are \UUU simply connected and \EEE well separated in the sense of \eqref{eq: main proXXX}(ii). Finally, we will also exploit that the set of `boundary squares' $\mathcal{X}^i_k$ consists of good squares (see \eqref{eq: main pro}(iv)).

\EEE
 
 \smallskip

We now proceed with the proofs of  Theorem \ref{th: bad part} and Lemma \ref{lemma: bad sets}.

\smallskip

 \noindent {\em Proof of Theorem \ref{th: bad part}.} \BBB For the proof recall the notions introduced in \eqref{eq: Ju*}-\eqref{eq: distance}. We start by describing the proof strategy.  
 
 \smallskip
\textbf{Outline of the proof:} \emph{Step I:} We first construct a sequence of `bad sets' $(E^i)_{i = 8}^I$ being the central object in the proof. In iteration step $i$ the starting point of the construction is the set $B^i$ defined in \eqref{NNNNN}. We consider a suitable union of $B^i$ with  components of the set $\bigcup_{l=i+1}^I B_l$, which contains  bad squares of smaller length scales. Considering only connected components with sufficiently large diameter (i.e. not isolated components), we then define iteratively sets $E^i$ such that the sequence $(E^i)_{i=8}^I$ is \emph{decreasing} in $i$. In this context, it is essential to consider bad squares of various length scales in order to obtain a decreasing sequence of sets. The fundamental property of each $E^i$ is the fact that a small layer around $E^i$ consists only of good squares $\mathcal{Q}^i \setminus \mathcal{Q}^i_{\rm bad}$ (we refer to \eqref{eq: nocrack} for details). This together with the monotonicity of $(E^i)_{i=8}^I$ will be essential in the definition of the Whitney covering in Step IV where regions around $E^i$ are covered with squares of $\mathcal{Q}^i$.

\smallskip

\emph{Step II and III:} The sets $(E^i)_{i=8}^I$ are also the starting point for the definition of the partition   $(P'_j)_j$. In fact, we construct a family of sets consisting of closed segments, denoted by $\mathcal{S}_I$, such that $E^I \cup \bigcup_{S \in \mathcal{S}_I} S$ is connected and $\sum_{S \in \mathcal{S}_I} \mathcal{H}^1(S)$ is controlled in terms of $\mathcal{H}^1( \UUU J^* \EEE )$. As $E^I$ contains $\partial Q_\mu$, we find that the connected components of  $Q_\mu \setminus (E^I \cup \bigcup_{S \in \mathcal{S}_I}S)$ are simply connected. It turns out that each connected component of $E^I$ contains exactly one component of $\UUU J^* \EEE$ (see \eqref{eq: elarge3}). Using these components of \UUU $J^*$ \EEE and the sets $\mathcal{S}_I$, we then construct the partition  $(P'_j)_j$ consisting of simply connected sets (see Step III).

The iterative construction of the sets $\mathcal{S}_I$ is performed in Step II. Here,  in each iteration step $i$, the essential idea is to combine the various connected components of $E^i$ with small segments of length $\sim s_i$. Due to a suitable bound on the number of connected components, which is related to the definition of isolated components (see \eqref{eq: length m}), we are able to show that the length of added segments in iteration step $i$ is of the order $\theta^{ir}   \mathcal{H}^1(\UUU J^* \EEE)$ (see \eqref{eq: length addi}). This guarantees a control on the total length of the segments which are added in all iterations steps (see left contribution in \eqref{eq: one propi}).

In the construction an additional difficulty has to be faced: if $E^{i-1} \setminus E^{i}$ contains `large' connected sets, so-called \emph{removed sets} as defined in \eqref{eq: E large}, possibly longer segments have to be added. The length of these segments, however, can be controlled in terms of the original jump set contained in the (pairwise disjoint) removed sets, see the right contribution in \eqref{eq: one propi}.

\smallskip

\emph{Step IV:} In the final step of the proof we define a covering of  Whitney-type by tessellating regions around $E^i$ with squares of $\mathcal{Q}^i$.  We observe that the area of the sets $(E^i)_{i=8}^I$ becomes gradually smaller and eventually each component of $E^I$ is a small neighborhood of a  component of  $\UUU J^*\EEE$. We then finally tessellate $E^I \setminus \UUU J \EEE$ with squares of $\bigcup_{l \ge I+1}^\infty\mathcal{Q}^l$. This allows us to show \eqref{eq: main pro}(i),(vi). For properties \eqref{eq: main pro}(ii),(iii) we follow   the usual construction of Whitney coverings. Finally, for \eqref{eq: main pro}(iv),(v) we fundamentally use that small layers around each $E^i$ consist only of good squares (see \eqref{eq: nocrack}). To establish the latter property, it is crucial that the sets $B^i$   are defined in terms of  enlarged squares (cf. \eqref{eq: enlarged squares} and \eqref{NNNNN}). 

In this step of the proof we also introduce the set $Z$ and comment on its definition.

 \smallskip
 
We now start with the proof. In the following  $c>0$ stands for a universal constant. We may  suppose that $\theta$ is small, i.e., $\theta \le \theta_0$ for some small universal constant $\theta_0 \in (0,\frac{1}{16})$ chosen in dependence of the appearing universal constants in the proof (cf. \eqref{eq: NNN} and \eqref{eq: new constant}). Having fixed the value of $\theta$, we let  $I = I(u,\theta, \UUU r \EEE ) \in \N$   be the index introduced in \eqref{eq: distance}. \EEE

\smallskip
\textbf{Step I (Definition  of bad sets):}   We define sets forming the starting point for the construction of the partition and the Whitney covering.   Recall the definition of $B^i$ in \eqref{NNNNN}.  For $\BBB 8 \EEE \le i \le I$ let $\mathcal{D}^i =(D^i_k)_k$ be the connected components of  $\bigcup_{l=i}^I B^l$ having nonempty intersection with $B^i$ and satisfying 
\begin{align}\label{eq: bound2}
d(D_k^i) >   \theta^{-ir} s_{i}.
\end{align} 
Observe that by definition of $\mathcal{B}^i_{\rm iso}$, in particular \eqref{eq: crack bound}, each $B^i_k \in \mathcal{B}^i \setminus \mathcal{B}^i_{\rm iso}$ is contained in a component of  $\mathcal{D}^i$. As each connected component is a finite union of squares, we find $d(D^i_k) \le \mathcal{H}^1(\partial D^i_k)$ and thus by Lemma \ref{lemma: cup} (for $\mathcal{R}^j = \mathcal{Q}^j_{\rm bad}$, $i \le j \le I$) we obtain 
\begin{align}\label{eq: cup}
\sum\nolimits_k d(D^i_k) \le \mathcal{H}^1\Big(\partial \, \bigcup\nolimits_{l=i}^I B^l \Big)\le c\theta^{-3} \mathcal{H}^1\big(\UUU J^* \EEE \cap \bigcup\nolimits_{l=i}^I B^l\big)\le c\theta^{-3} \mathcal{H}^1(\UUU J^* \EEE)
\end{align} 
independently of $ \BBB 8 \EEE \le i \le I$.

\BBB For notational reasons, it is convenient to start with the set $E^{7} := \bigcup\nolimits_{Q \in \mathcal{Q}^{7}, Q \subset Q_\mu} \overline{Q'''} $ containing $Q_\mu$. \EEE Assuming that $E^k$, $ 7 \le k \le i-1$, have been constructed, we let 
\begin{align}\label{eq: E large2}
\mathcal{E}^i = \lbrace D^i_k: D^i_k \cap E^{i-1} \neq \emptyset \rbrace \ \ \  \ \text{and} \ \ \ \  E^i= \bigcup\nolimits_{D^i_k \in \mathcal{E}^i} D^i_k.
\end{align}
(See also Figure \ref{korn1}(b)-(d).) \BBB As mentioned above, it will turn out that for each $8 \le i \le I$ a small layer around $E^i$ consists of good squares $\mathcal{Q}^i \setminus \mathcal{Q}^i_{\rm bad}$. (This property is essential for the construction of the Whitney covering in Step IV and we refer to \eqref{eq: nocrack} for details.) Consequently, in the definition of $E^i$ it suffices to consider the components $D^i_k$ having nonempty intersection with $E^{i-1}$. This is indeed crucial to obtain a sequence of sets which is decreasing in $i$. In fact, we have   \EEE
\begin{align}\label{eq: nested}
E^i \subset E^j \ \ \ \ \ \text{for $7\le j <i \le I$.}
\end{align}
\BBB First, $E^8 \subset E^7$ is clear. Let $i \ge 9$. Note that each $D^i_k \in \mathcal{E}^i$ intersects $E^{i-1}$ and   thus we find some $D^{i-1}_{k'} \in \mathcal{E}^{i-1}$ with $D^i_k \cap D^{i-1}_{k'} \neq \emptyset$. As  $D^{i-1}_{k'}$ is a connected component of $\bigcup_{l=i-1}^I B^l$ and $D^i_k$ is contained in a connected component  of $\bigcup_{l=i-1}^I B^l$, we get $D^i_k \subset D^{i-1}_{k'}$. \EEE  Then the definition of $E^{i-1}$ implies $D^i_k \subset D^{i-1}_{k'} \subset E^{i-1}$.  This yields $E^i \subset E^{i-1}$, as desired. 

Moreover, note that each $E^i$ contains $J_0$ (recall \eqref{eq: Ju*}) and therefore particularly 
\begin{align}\label{eq: nested2}
\partial Q_\mu \subset J_0 \subset E^i \ \ \ \ \text{ for } \ \ \  \BBB 8 \le i \le I \EEE.
\end{align}
 \BBB In fact, by the remark below \eqref{eq: bad1}, for each $i \ge 8$  there is a component $D^i_k$ satisfying  $J_0 \subset D^i_k$.\EEE

\UUU  For the construction of the partition we will also need \EEE the   associated \emph{removed sets}    
\begin{align}\label{eq: E large}
\begin{split}
&\mathcal{P}^i =  \Big\{ Q \in \bigcup\nolimits_{l=i}^I\mathcal{Q}^l_{\rm bad}: \overline{Q'''}\subset E^i, \  \overline{Q'''} \cap E^{i+1} = \emptyset \Big\}, \ \ \  R^i = \bigcup\nolimits_{Q \in \mathcal{P}^i} \overline{Q'''}
\end{split} 
\end{align}
for $ 7 \le i \le I-1$, \BBB which will play a pivotal role for  the construction of the partition in Step II. We \EEE observe that 
\begin{align}\label{eq: E not shrink cond}
R^j \cap R^{i} = \emptyset \  \ \ \text{ for } 7 \le j < i \le I-1.
\end{align} 
To see this, consider $j < i$ and note that  $\overline{Q'''} \cap E^i = \emptyset$ for $Q \in \mathcal{P}^j$ since $E^i \subset E^{j+1}$. Thus, we find  $R^j \cap E^i = \emptyset$  and therefore \eqref{eq: E not shrink cond} holds as $R^i \subset E^i$.

 Essentially, $E^{i}$ arises from $E^{i-1}$ by removing the squares $\mathcal{P}^{i-1}$. Moreover, $E^i$ still intersects the squares $\lbrace Q \in \mathcal{Q}^{i-1} \setminus \mathcal{P}^{i-1}: \overline{Q'''} \subset E^{i-1} \rbrace$, but in general covers a smaller portion compared to $E^{i-1}$ (cf. Figure \ref{korn1}(d)). More precisely, we get
\begin{align}\label{eq: violat2}
\dist(x,E^i)\le  10\sqrt{2}s_{i-1} \ \ \ \text{for all} \ \ x \in E^{i-1} \setminus R^{i-1}.
\end{align}
In fact, since $x \in E^{i-1}$, by \eqref{eq: E large2} we find $Q \in \bigcup^I_{l=i-1} \mathcal{Q}^l_{\rm bad}$ such that $x \in \overline{Q'''} \subset E^{i-1}$. As $x \notin R^{i-1}$, we obtain $Q \notin \mathcal{P}^{i-1}$ and therefore $\overline{Q'''} \cap E^i \neq  \emptyset$, which implies \eqref{eq: violat2}.

\smallskip

\textbf{Step II (Construction of the partition, induction step):}  \BBB Note that the sets $(E^i)_{i=8}^I$ defined in Step I   are in general not connected. The goal of this step is to construct    a family of closed, connected sets $\mathcal{S}_I$ such that each $S \in \mathcal{S}_I$ is a finite union of closed segments, the set $E^I \cup \bigcup_{S \in \mathcal{S}_{I}} S$ is connected, and we have
\begin{align}\label{eq: newnew}
\# \mathcal{S}_I \le \theta^{-4}s_I^{-1} \mathcal{H}^1(\UUU J^* \EEE), \ \ \ \ \  \mathcal{H}^1\big(\bigcup\nolimits_{S \in \mathcal{S}_{I}}S \big) \le C\mathcal{H}^1(\UUU J^* \EEE),
\end{align}
where   $C=C(\theta, \UUU r \EEE )>0$.   In particular, since $\partial Q_\mu \subset E^I$, this implies that the connected components of $Q_\mu \setminus (E^I \cup \bigcup_{S \in \mathcal{S}_{I}} S)$ are simply connected. In Step III we will use this fact to construct the partition $(P_j')_{j=1}^m$ consisting of simply connected sets. \BBB
The family $\mathcal{S}_I$ is defined inductively, where in each iteration step $i$ the strategy is to  connect the components of $E^i$ by suitable segments lying in $E^{i-1}$. The argument is mainly based on a suitable control on the number of components induced by \eqref{eq: bound2}, the definition of the removed sets \eqref{eq: E large}, and an elementary lemma about the distance of sets (Lemma \ref{lemma: topo}). As the details are not needed for the subsequent parts of the proof, the reader may prefer to  proceed directly with Step III.  \EEE

We will show \UUU by induction \EEE that for \BBB $i \ge 8$ \EEE there is a family of closed, connected sets  $\mathcal{S}_{i}$ with 
\begin{align}\label{eq: one propi2}
\# \mathcal{S}_{i} \le \theta^{-4}\theta^{ir} s_{i}^{-1}  \mathcal{H}^1(\UUU J^* \EEE),
\end{align} 
such that each $S \in \mathcal{S}_i$ is a finite union of closed segments and $E^{i}_{\mathcal{S}} := E^i \cup \bigcup_{S \in \mathcal{S}_{i}} S$ is connected.   Moreover, we show  (recall \eqref{eq: E large})
\begin{align}\label{eq: one propi}
\begin{split}
 \mathcal{H}^1\big(\bigcup\nolimits_{S \in \mathcal{S}_{i}}S \big) \le c\theta^{-5}\mathcal{H}^1(\UUU J^* \EEE) \sum\nolimits_{j=8}^{i} \theta^{jr} + c \theta^{-3}\sum\nolimits^{i-1}_{j=7}\mathcal{H}^1(\UUU J^* \EEE \cap R^j).
\end{split}
\end{align}
\BBB Herefrom we then obtain \eqref{eq: newnew}. Indeed, it suffices to apply \eqref{eq: one propi2}-\eqref{eq: one propi} for step $I$ and to recall that the sets $(R^{j})_j$ are pairwise disjoint (see  \eqref{eq: E not shrink cond}). \EEE

Let us assume that the above assertions, in particular \eqref{eq: one propi2}-\eqref{eq: one propi}, hold for $i-1$, $i \ge 8$. \BBB (Note that the assertions hold for $i = 7$ with $\mathcal{S}_7 := \emptyset$ since $E^7$ is connected, see before \eqref{eq: E large2}.)  We now construct $\mathcal{S}_i$ and confirm \eqref{eq: one propi2}-\eqref{eq: one propi}.  \EEE

\smallskip

\emph{Components:} We denote the connected components  of $E^i \cup \bigcup_{S \in \mathcal{S}_{i-1}} S$ by $\mathcal{F}^i:= (F^i_k)^{M_i}_{k=1}$. The goal is to connect these sets  by suitable segments.   Observe that each $F^i_k$ is a union of closed sets with Hausdorff-dimension one or two. Recall that each connected component $D^i_k \in \mathcal{E}^i$ satisfies $d(D^i_k) \ge   \theta^{-ir} s_{i}$ by \eqref{eq: bound2}. Consequently,  by \eqref{eq: cup} we get for $\theta$ small   
\begin{align}\label{eq: NNN}
\begin{split}
\# \mathcal{E}^i \le  \theta^{ir}s^{-1}_{i}\sum\nolimits_k d(D^i_k)\le c\theta^{-3} \theta^{ir} s^{-1}_{i}\mathcal{H}^1(\UUU J^* \EEE)
 \le \tfrac{1}{4}\theta^{-4}\theta^{ir} s^{-1}_{i}\mathcal{H}^1(\UUU J^* \EEE). 
 \end{split}
\end{align}
  By \eqref{eq: one propi2} (for $i-1$) and \eqref{eq: NNN} we have 
\begin{align}\label{eq: length m}
\begin{split}
M_i= \# \mathcal{F}^i &\le \tfrac{1}{4}\theta^{-4}\theta^{ir} s^{-1}_{i}\mathcal{H}^1(\UUU J^* \EEE) + \theta^{-4}\theta^{(i-1)r} s_{i-1}^{-1} \mathcal{H}^1(\UUU J^* \EEE) \le \tfrac{1}{2}\theta^{-4}\theta^{ir} s^{-1}_{i}\mathcal{H}^1(\UUU J^* \EEE)
\end{split}
\end{align} 
recalling that \UUU $r < \frac{1}{8}$, \EEE \BBB $\theta \le \frac{1}{16}$ \EEE  and $s_i = \theta s_{i-1}$.  By the definition of $\mathcal{F}^i$, $E^{i-1}_{\mathcal{S}} = E^{i-1} \cup \bigcup_{S \in \mathcal{S}_{i-1}} S$ and the fact that $E^i \subset E^{i-1}$ (see   \eqref{eq: nested}) we get
\begin{align}\label{eq: connected}
E^{i-1}_{\mathcal{S}} = \bigcup\nolimits^{M_i}_{k=1} F^i_k \cup  (E^{i-1} \setminus E^i).  
\end{align}
We denote the connected components of $R^{i-1}$ (see \eqref{eq: E large}) by $\mathcal{R}^{i-1} =  (R^{i-1}_k)_k$ and let $\bar{c} = 10\sqrt{2}$ for shorthand. Observe that by \eqref{eq: violat2}, \BBB \eqref{eq: connected} \EEE and the fact that $E^i \subset \bigcup\nolimits_{k=1}^{M_i} F_k^i$ we obtain
\begin{align}\label{eq: violat}
\dist\big(x, \bigcup\nolimits_{k = 1}^{M_i} F_k^i \big)\le   \bar{c} s_{i-1} \ \ \ \text{for all} \ \ x \in E^{i-1}_{\BBB \mathcal{S}\EEE} \setminus R^{i-1}.
\end{align}
\BBB Recall that $E^{i-1}_{\mathcal{S}}$ is (path) connected by hypothesis.  Using \eqref{eq: connected}-\eqref{eq: violat} we now apply Lemma \ref{lemma: topo} for the families $\mathcal{F}^i$ and $\mathcal{R}^{i-1}$ with $G = E^{i-1}_{\mathcal{S}}$ and $d =  \bar{c} s_{i-1} $. We get that for each \EEE $F^i_{j} \in \mathcal{F}^i$ there is another component $F^i_{k} \in \mathcal{F}^i$ such that (a) $\dist(F^i_{j},F^i_{k}) \le 4\bar{c} s_{i-1}$ or (b) there is a corresponding $R^{i-1}_l \in \mathcal{R}^{i-1}$ such that 
$$
\dist(F^i_{j},R_l^{i-1}) \le 4\bar{c} s_{i-1}, \ \ \ \ \dist(F^i_{k},R_l^{i-1}) \le 4\bar{c} s_{i-1}.
$$
%
%
Due to \eqref{eq: violat} and the fact that  $E^{i-1}_{\mathcal{S}}$ is connected, we further note that for each pair $F^i_{j_1}, F^i_{j_2} \in \mathcal{F}^i$ there is a chain $j_1 = k_1,k_2,\ldots,k_n = j_2$ such that each pair  $F^i_{k_l}, F^i_{k_{l+1}}$, $l=1,\ldots,n-1$, satisfies (a) or (b).


\smallskip
\emph{Construction of segments:} Consider $F^i_{j}, F^i_{k}\in \mathcal{F}^i$  such that (a) or (b) holds. We observe that in case (a) one can choose a (closed) segment $S$ with $\mathcal{H}^1(S) \le 4\bar{c}s_{i-1}$ such that $F^i_{j} \cup S \cup F^i_{k}$ is connected and in case (b) we can find two segments $S^1$, $S^2$ with $\mathcal{H}^1(S^1 \cup S^2) \le 8\bar{c}s_{i-1}$ such that $S^1 \cup S^2 \cup R^{i-1}_l \cup F^i_{j} \cup F_{k}^i$ is connected, where the component $R^{i-1}_l$ is chosen as above. Thus, for a universal $c>0$ large enough we can find $M_i-1$ sets $(\hat{R}^{i-1}_l)^{M_i-1}_{l=1} \subset \mathcal{R}^{i-1} \cup \lbrace \emptyset \rbrace$ and sets $(S_l)_{l=1}^{M_i-1}$ with 
\begin{align}\label{eq: sk}
\mathcal{H}^1(S_l) \le cs_{i-1},
\end{align}
where each $S_l$ consists of at most two segments, so that 
$\bigcup^{M_i}_{k=1} F_k^i \cup \bigcup^{M_i-1}_{l=1} (\hat{R}^{i-1}_l \cup S_l)$ is connected \BBB and each $\hat{R}^{i-1}_l \cup S_l$ is connected. \EEE Indeed, the construction of sets connecting two different components $(F^i_k)_k$ has been addressed above and the fact that it suffices to consider $M_i-1 = \# \mathcal{F}^i-1$ sets may be seen by induction over the number of components $M_i$. 

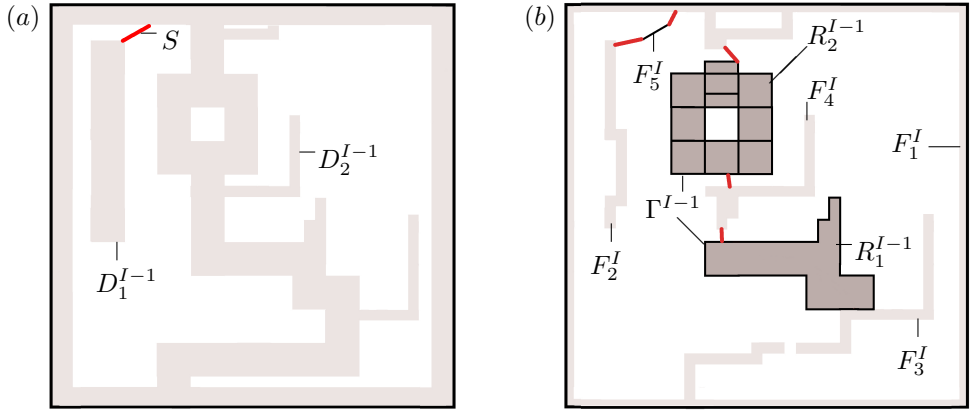
\begin{figure}[H]
\centering
\hspace{0.5cm}
\begin{overpic}[width=0.88\linewidth,clip]{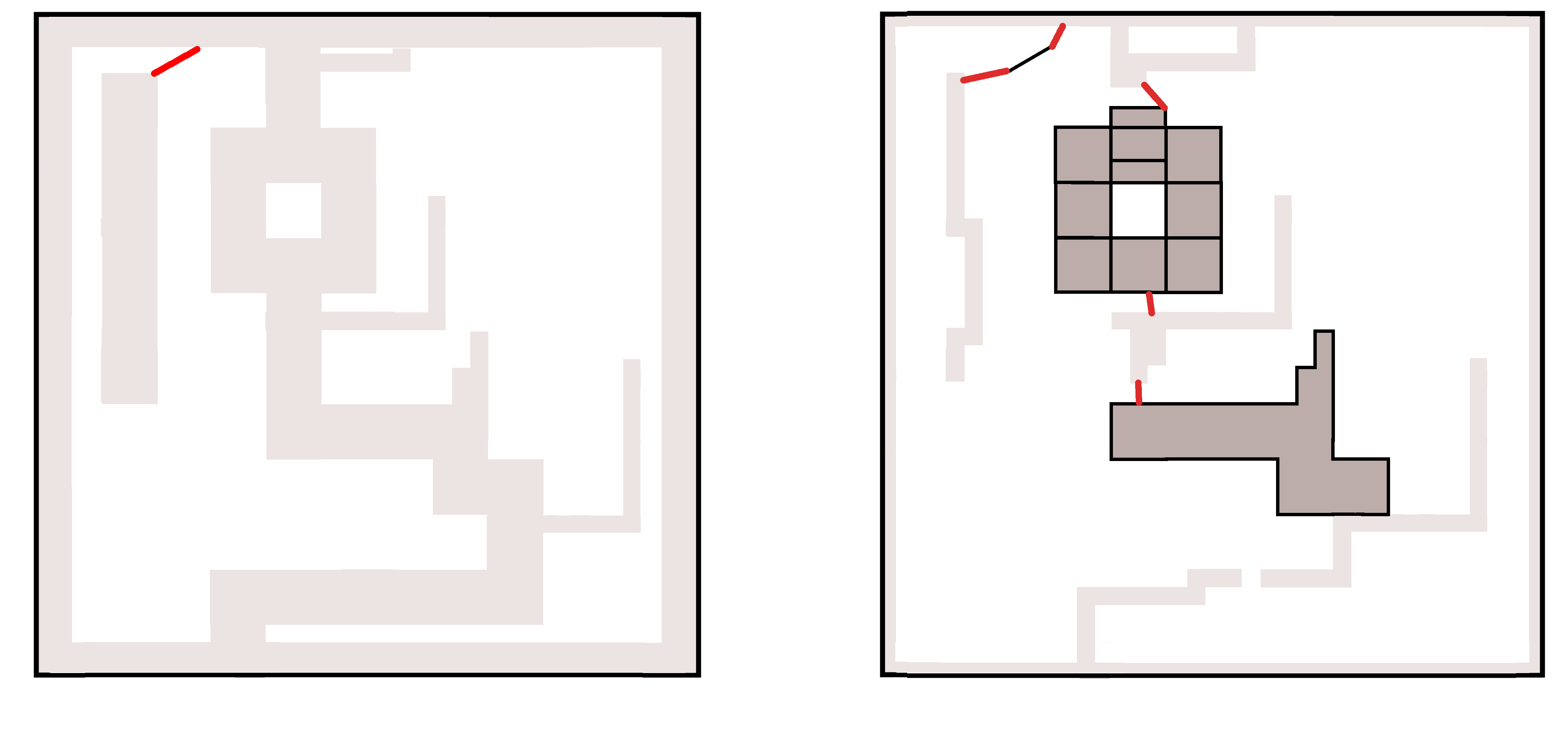}

\put(-9,159){{$(a)$}}
\put(187,159){{$(b)$}}

\put(25,58){{$D^{I-1}_1$}}
\put(33,70) {\line(0,1){7}}

\put(109,105){{$D^{I-1}_2$}}
\put(102,111) {\line(1,0){6}}

\put(50,150){{$S$}}
\put(42,156) {\line(1,0){6}}

\put(212,65){{$F^{I}_2$}}
\put(220,77) {\line(0,1){7}}

\put(228,137){{$F^{I}_5$}}
\put(236,149) {\line(0,1){7}}

\put(294,131){{$F^{I}_4$}}
\put(294,123) {\line(0,1){7}}

\put(326,111){{$F^{I}_1$}}
\put(340,113) {\line(1,0){12}}

\put(328,29){{$F^{I}_3$}}
\put(336,41) {\line(0,1){7}}

\put(233,86){{$\Gamma^{I-1}$}}
\put(247,96) {\line(0,1){6}}
\put(245,86) {\line(1,-1){10}}

\put(294,152){{$R^{I-1}_2$}}
\put(294,151) {\line(-1,-1){14}}

\put(312,70){{$R^{I-1}_1$}}
\put(305,76) {\line(1,0){6}}

\end{overpic}
\caption{(a) \BBB The set $E^{I-1}$ consists of two components  $D_1^{I-1}$ and $D_2^{I-1}$ which have been connected by the red segment $S$. Thus, $E^{I-1}_{\mathcal{S}} = E^{I-1} \cup S$ is connected. \ (b) The set $E^I \cup S$ consists of five  components $\mathcal{F}^{I} = \lbrace F^I_1,\ldots, F^I_5 \rbrace$ where $F^I_5 = S$.  \EEE The components of $\mathcal{S}_*$ are depicted, where the (union of) segments $S_l$, $l=1,\ldots,4$, connecting the sets $\mathcal{F}^{I}$ are highlighed in red and (a possible choice of) the set $\Gamma^{I-1}$ is illustrated in black. While in the example we can choose $\Gamma^{I-1} \cap R_1^{I-1} = \partial R_1^{I-1}$, in $R_2^{I-1}$ a more elaborated definition of $\Gamma^{I-1}$ is necessary (cf. Lemma \ref{lemma: cup}) since $R_2^{I-1}$ is not simply connected.  }  \label{korn2-1}
\end{figure}

We recall \eqref{eq: E large} and note $R^{i-1} = \bigcup_{j=i-1}^I\bigcup_{Q \in \mathcal{P}^{i-1} \cap \mathcal{Q}^j} \overline{Q'''}$. We apply Lemma \ref{lemma: cup} for the squares $\mathcal{P}^{i-1} \cap \mathcal{Q}^j \subset \mathcal{Q}^j_{\rm bad}$, $i-1 \le j \le I$, and obtain a set $\Gamma^{i-1}\supset \partial R^{i-1}$, being a finite union of closed segments, such that $\Gamma^{i-1} \cap \hat{R}^{i-1}_l$ is connected for $l=1,\ldots, M_i-1$ and 
\begin{align}\label{eq:prepull}
\mathcal{H}^1\big(\bigcup\nolimits^{M_i-1}_{l = 1}  \hat{R}^{i-1}_l \cap \Gamma^{i-1}\big) \le  \mathcal{H}^1\big(\Gamma^{i-1}\big)\le  c\theta^{-3}\mathcal{H}^1( R^{i-1} \cap \UUU J^* \EEE). 
\end{align}
Then we introduce the new jump components  
$$\mathcal{S}_* = \big(S_l \cup ( \hat{R}_l^{i-1} \cap \Gamma^{i-1})\big)^{M_i-1}_{l=1}, \ \ \ \ \ \ \mathcal{S}_i = \mathcal{S}_{i-1} \cup  \mathcal{S}_*. $$
We are now in the position to confirm \eqref{eq: one propi2}-\eqref{eq: one propi}. First, each component in $\mathcal{S}_i$ is  connected and a finite union of closed segments.  Moreover,  we observe that $E^{i}_{\mathcal{S}} := E^i \cup \bigcup_{S \in \mathcal{S}_{i}} S$ is connected by the definition of $(F^i_k)_k$ and  $\mathcal{S}_*$ (cf. Figure \ref{korn2-1}(b)).  By  \eqref{eq: one propi2} (for $i-1$) and  \eqref{eq: length m} we find
\begin{align*}
\# \mathcal{S}_i &\le \# \mathcal{S}_{i-1} +  \#\mathcal{F}^i-1 \le  \theta^{-4}\theta^{r(i-1)} s_{i-1}^{-1}  \mathcal{H}^1(\UUU J^* \EEE) + \tfrac{1}{2} \theta^{-4}\theta^{ir} s^{-1}_{i} \  \mathcal{H}^1(\UUU J^* \EEE)  \le  \theta^{-4}\theta^{ir} s_{i}^{-1} \mathcal{H}^1(\UUU J^* \EEE)
\end{align*}
\BBB as \UUU $r < \frac{1}{8}$ \EEE and  $\theta \le \frac{1}{16}$. This gives \eqref{eq: one propi2}. Moreover,  we use \eqref{eq: length m} as well as \eqref{eq: sk} to find
\begin{align}\label{eq: length addi}
\sum\nolimits_{l=1}^{M_i-1} \mathcal{H}^1(S_l) &\le c(M_i-1) s_{i-1} \le   c\theta^{-5} \theta^{ir}   \mathcal{H}^1(\UUU J^* \EEE).
\end{align}
This together with \eqref{eq:prepull} and    \eqref{eq: one propi} for step $i-1$  yield \eqref{eq: one propi} for step $i$.

\smallskip

\textbf{Step III (Construction of the partition, final step):} \BBB In Step II we have constructed    a family of closed, connected sets $\mathcal{S}_I$ such that each $S \in \mathcal{S}_I$ is a finite union of closed segments, the set $E^I \cup \bigcup_{S \in \mathcal{S}_{I}} S$ is connected, and \eqref{eq: newnew} holds. We now  define the partition $(P'_j)_{j=1}^m$ with the desired properties. The argument will be based on the observation that, due to the choice of $I$ in \eqref{eq: distance}, each connected component of $E^I$ contains exactly one component of $\UUU J^* \EEE$. These components will be connected suitably with $\bigcup_{S \in \mathcal{S}_{I}} S$.\EEE

First, note that by the definition of $\mathcal{D}^I$ in \eqref{eq: bound2} and  \eqref{eq: distance}(ii) we have $\mathcal{D}^I = \mathcal{B}^I$. Consequently,  recalling \eqref{eq: E large2} we see that $E^I$ consists of the connected components of ${B}^I$ having nonempty intersection with $E^{I-1}$. As before, by $\mathcal{E}^I \subset \mathcal{D}^I = \mathcal{B}^I$ we denote the connected components $D^I_k$ with $D^I_k \cap E^{I-1} \neq \emptyset$.

\begin{figure}[H]
\centering
\hspace{0.5cm}
\begin{overpic}[width=0.88\linewidth,clip]{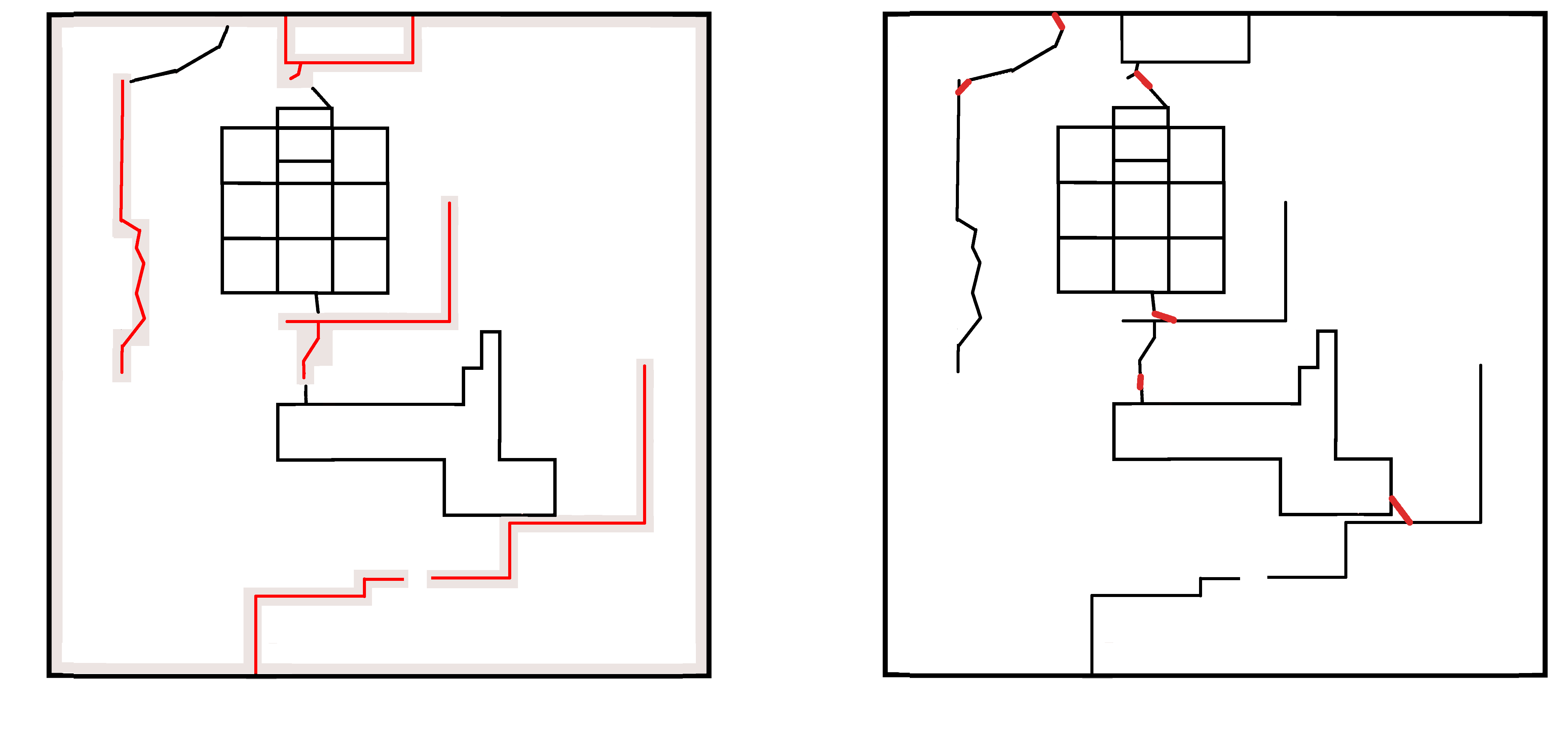}

\end{overpic}
\caption{(a) \BBB The segments $\mathcal{S}_I$ are depicted in black and the set $E^{I}$ in gray, where $E^I \cup \bigcup_{S \in \mathcal{S}_{I}} S$ is connected. Moreover, the components $(\Gamma^*_{j_k})_k$ contained in $E^I$ are illustrated in red. (b) The  components of $\mathcal{S}_I$ and $(\Gamma^*_{j_k})_k$ are combined to a connected set by the red segments.  \EEE   }  \label{korn2-2}
\end{figure}

In view of the remark before \eqref{eq: distance},   we observe that each $D^I_k$  contains exactly one connected component $\Gamma_{j_k}^*$ of $\UUU J^* \EEE$  (cf. Figure \ref{korn2-2}(a)). More precisely, by \eqref{eq: bad1}  and \eqref{NNNNN} we have 
\begin{align}\label{eq: elarge3}
\dist(x, \Gamma_{j_k}^*) \le \BBB \bar{c} \EEE s_I  \ \ \ \text{ for all } x \in D^I_k,
\end{align}
\BBB where as before we set $\bar{c} = 10\sqrt{2}$ for brevity. \EEE Moreover,  since $\# \mathcal{E}^I$ is bounded from above by the number of components of $\UUU J^* \EEE$, \BBB we find  $\# \mathcal{E}^I \le s_I^{-1} \mathcal{H}^1(\UUU J^* \EEE)$ by \eqref{eq: distance}(iii). \EEE

\enlargethispage{\baselineskip}

We now proceed similarly as in Step II to connect the components $\mathcal{S}_I$ and $(\Gamma_{j_k}^*)_k$ with segments (cf. Figure \ref{korn2-2}(b)). \BBB We apply Lemma \ref{lemma: topo} on the family $\mathcal{F} := \mathcal{S}_I \cup \lbrace \Gamma_{j_k}^*: D^I_k \in \mathcal{E}^I \rbrace $ with $R= 0$, $G= E^I \cup \bigcup_{S \in \mathcal{S}_{I}} S$ and $d = \bar{c} s_I$, where we note that \eqref{eq: elarge3} implies \eqref{eq: violatXXXXX}. We see that for each set in $\mathcal{F}$ there is another set in $\mathcal{F}$ with distance at most $4\bar{c}s_I$. Thus, \EEE we can find a family of closed segments $\hat{\mathcal{S}}$ with $\# \hat{\mathcal{S}} = \BBB \# \mathcal{F} - 1 \EEE$ and $\mathcal{H}^1(\hat{S}) \le 4\bar{c}s_I$ for all $\hat{S} \in \hat{\mathcal{S}}$ such that 
\begin{align*}
\bigcup\nolimits_{D^I_k \in \mathcal{E}^I} \Gamma_{j_k}^*  \cup \bigcup\nolimits_{\hat{S} \in \hat{\mathcal{S}}} \hat{S} \cup \bigcup\nolimits_{S \in  \mathcal{S}_{I}} S
\end{align*}
is connected, where as before $\Gamma_{j_k}^*$ denotes the component contained in $D^I_k$. Note that $\partial Q_\mu \subset \bigcup\nolimits_{D^I_k \in \mathcal{E}^I} \Gamma_{j_k}^* $ since $\partial Q_\mu \subset J_0 \subset  E^I$ \BBB by \eqref{eq: nested2}. \EEE Thus, letting    
\begin{align}\label{eq: S const}
\mathcal{S} = \mathcal{S}_{I} \cup \hat{\mathcal{S}} \cup \lbrace \Gamma_{j_k}^*: D^I_k \in \mathcal{E}^I \rbrace 
\end{align}
we finally find that the connected components of $Q_\mu \setminus \bigcup_{S \in \mathcal{S}} S$, denoted by $(P'_j)^m_{j=1}$, are simply connected and form a partition of $Q_\mu$ (up to a set of negligible measure). By construction $\bigcup^m_{j=1} \partial P'_j$ is a finite union of closed segments. \BBB To conclude this step of the proof, it remains to show \eqref{eq: bad length}. Recall $\# \mathcal{E}^I \le s_I^{-1} \mathcal{H}^1(\UUU J^* \EEE)$ by \eqref{eq: distance}(iii). This together with \eqref{eq: newnew} gives
$$ \sum\nolimits_{\hat{S} \in \hat{\mathcal{S}}} \mathcal{H}^1(\hat{S}) \le 4\bar{c}s_I (\# \mathcal{F}-1) \le 4\bar{c}s_I(\# \mathcal{S}_I + \# \mathcal{E}^I) \le C \mathcal{H}^1(\UUU J^* \EEE)$$
for  $C=C(\theta,\UUU r \EEE )>0$.   Then recalling \eqref{eq: S const} and using again \eqref{eq: newnew} we conclude \EEE
\begin{align}\label{eq: NNNNN}
 \mathcal{H}^1\big(\bigcup\nolimits_{j=1}^m \partial P'_j \big) \le \sum\nolimits_{j=1}^m \mathcal{H}^1(\partial P'_j) \le C\mathcal{H}^1(\UUU J^* \EEE)  \le C\mathcal{H}^1(\UUU J \EEE),
\end{align}
where in the last step we employed \eqref{eq: Ju*XX}.

\smallskip

\textbf{Step IV (Covering):} \BBB After having completed the construction of the partition $(P'_j)^m_{j=1}$, we now finally define an associated Whitney  covering of $Q_\mu$, introduce the set $Z$, and confirm \eqref{eq: main pro}. For this part of the proof, recall  the definition of the sets $E^i$ in \eqref{eq: E large2} and the properties \eqref{eq: nested}-\eqref{eq: nested2}.  Moreover, recall the relation between $E^I$, the components of $\UUU J^* \EEE$, and the partition $(P'_j)^m_{j=1}$ in \eqref{eq: elarge3}-\eqref{eq: S const}.

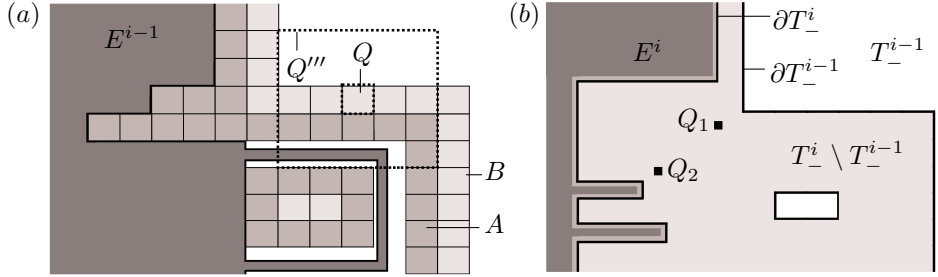
\begin{figure}[H]
\centering
\begin{overpic}[width=0.855\linewidth,clip]{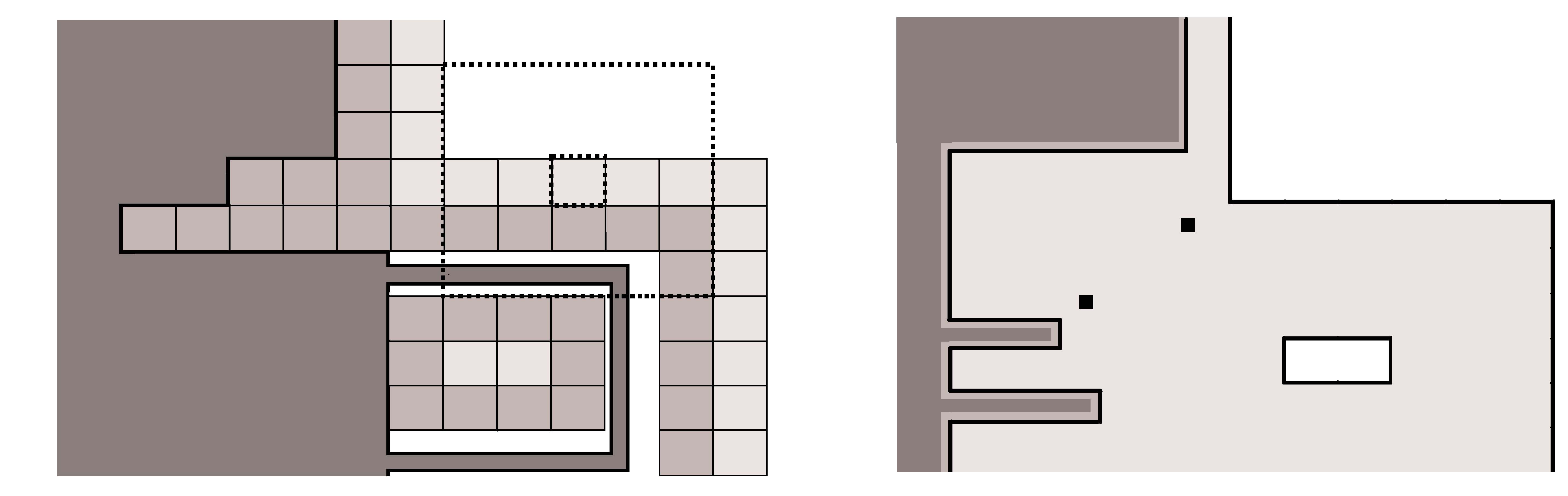}

\put(-4,100){{$(a)$}}
\put(186,100){{$(b)$}}

\put(33,90){{$E^{i-1}$}}
\put(102,79){{$Q'''$}}
 \put(127,85){{$Q$}}
 \put(106,89) {\line(0,1){7}}
  \put(130,73) {\line(0,1){9}}

\put(177,40){{$B$}}
\put(170,42) {\line(1,0){6}}

\put(177,20){{$A$}}
\put(155,22) {\line(1,0){20}}

  \put(250,60){{$Q_1$}}
  \put(246,41){{$Q_2$}}

\put(233,85){{$E^{i}$}}

\put(286,97){{$\partial T^{i}_-$}}
\put(285,77){{$\partial T^{i-1}_-$}}
\put(265,102) {\line(1,0){20}}
\put(275,82) {\line(1,0){8}}

\put(323,85){{$T^{i-1}_-$}}
\put(293,46){{$T^{i}_- \setminus T^{i-1}_-$}}

\end{overpic}
\caption{(a) The set $E^{i-1}$ and the `layers' $A = T^{i-1} \setminus  T^{i-1}_-$, $B = T^{i-1}_- \setminus  T^{i-1}_{--}$ are sketched, where in general $E^{i-1}$ is not a union of squares in $\mathcal{Q}^{i-1}$. \BBB For the proof of \eqref{eq: nocrack}, \EEE note that a square $Q$ in $A \cup B$ is not in $\mathcal{Q}^{i-1}_{\rm bad}$ since $\overline{Q'''}$ would then belong to $E^{i-1}$.   \ (b) In dark gray we have depicted  $E^{i}$, which is a subset of $E^{i-1}$. The sets $T^{i-1}_-$ and $T^{i}_- \setminus T^{i-1}_-$ are illustrated in white and light gray, respectively. Note that $\partial T_-^{i-1}$, $\partial T_-^i$ satisfy \eqref{eq: gooddisti}. \BBB For the proof of \eqref{eq: main pro}(iv), case (a),  note that  the black square $Q_1 \in \mathcal{Q}^i$ lies in $A$ and thus satisfies  \eqref{eq: badspec}. For case (c), observe that $Q_2 \subset E^{i-1}$ and thus $Q_2 \notin\mathcal{Q}^{i}_{\rm bad}$ or $Q_2'' \subset U_i$ since otherwise it would be contained in $E^{i}$. \EEE  } \label{korn3}
\end{figure}

As motivated above, we will cover certain regions around $E^i$  with squares of $\mathcal{Q}^i$. Let us first introduce some relevant notation. \EEE For each $8 \le i \le I$ we define the sets
\begin{align}\label{eq: Ts} 
\begin{split}
&T^i = \bigcup\nolimits_{Q \in \mathcal{T}^i}   \overline{Q}    \ \ \ \, \ \ \ \text{ with } \ \ \ \mathcal{T}^i = \lbrace Q\in    \mathcal{Q}^i: \ Q \subset Q_\mu \setminus E^i \rbrace,\\
&T^i_- = \bigcup\nolimits_{Q \in \mathcal{T}^i_-} \overline{Q} \ \ \ \ \ \  \text{ with } \ \ \  \mathcal{T}^i_- = \lbrace Q\in    \mathcal{T}^i: Q'' \subset T^i \rbrace,\\
&T^i_{--} = \bigcup\nolimits_{Q \in \mathcal{T}^i_{--}} \overline{Q} \ \ \ \text{ with } \ \ \  \mathcal{T}^i_{--} = \lbrace Q\in    \mathcal{T}^i_-: Q'' \subset T^i_- \rbrace,
\end{split}
\end{align}
where, loosely speaking,  $T^i_-, T^i_{--}$ arise from $T^i$ by  removing `layers'  of squares. (We refer to Figure \ref{korn3} for an illustration.) As a preparation, we  observe that $T^j \supset T^i$ for $8 \le i < j \le I$ since $E^j \subset E^i$ (see \eqref {eq: nested}) and $s_j < s_i$. This particularly implies
\begin{align}\label{eq: gooddisti}
T^{i-1}_- \subset T^{i}_-, \ \ \ \ \dist(\partial T^{i-1}_-, \, \partial T^{i}_-)\ge s_{i-1}
\end{align}
for $9 \le i \le I$. \BBB As mentioned above, for the proof of \eqref{eq: main pro} it will be crucial  that a small layer around $E^i$ does not contain bad squares. Indeed, we have \EEE 
\begin{align}\label{eq: nocrack}
\big( T^i \setminus T^i_{--} \big) \cap \bigcup\nolimits_{Q \in \mathcal{Q}^i_{\rm bad}} Q = \emptyset.
\end{align}
\BBB To see this, we note as a preparation that each $Q \in \mathcal{T}^i \setminus \mathcal{T}^i_{--}$ satisfies $\overline{Q'''} \cap E^i \neq \emptyset$. Indeed, since $Q \notin \mathcal{T}^i_{--}$,  we find some $Q_* \in \mathcal{T}^i \setminus \mathcal{T}^i_{-}$ such that $Q_* \subset Q''$. As $Q_* \notin \mathcal{T}^i_{-}$, we find some $Q_{**} \in \mathcal{Q}^i \setminus \mathcal{T}^i$ with $Q_{**} \subset Q_*''$. Thus, $\overline{Q_*''} \cap E^i \neq \emptyset$ and it then suffices to note that $\overline{Q'''} \supset \overline{Q_*''}$ (we also refer to Figure \ref{korn3}(a)).  

Now suppose by contradiction that $Q \in \mathcal{Q}^i_{\rm bad} \cap (\mathcal{T}^i \setminus \mathcal{T}^i_{--})$. Since $\overline{Q'''} \cap E^i \neq \emptyset$, by \eqref{eq: E large2} we find some $D^i_k \in \mathcal{E}^i$ such that $ \overline{Q'''} \cap D^i_k \neq \emptyset$. As $D^i_k$ is a connected component  of $\bigcup_{l=i}^I B^l$  (see before \eqref{eq: bound2}), we get $\overline{Q'''} \subset D^i_k$ and thus  $\overline{Q'''}  \subset E^i$ by \eqref{eq: E large2}. This gives a contradiction as clearly $\overline{Q'''} \not\subset E^i$ due to $Q \in \mathcal{T}^i$.   \EEE

\smallskip

\emph{Definition of the covering and the set $Z$:} We now introduce the covering as follows. Recall that each set $T^i_- \setminus  T^{i-1}_-$ for $8 \le i\le I$ is \BBB a \EEE union of squares in $\mathcal{Q}^i$ (up to a negligible set), where we set $T^{7}_ - = \emptyset$. We  define for $8 \le i \le I$  
\begin{align}\label{eq:covi}
\mathcal{C}^i = \lbrace Q \in \mathcal{Q}^i: Q \subset T^i_- \setminus T^{i-1}_-\rbrace.
\end{align}
Hereby, we obtain a covering of $T^I_-$. More precisely, we have
\begin{align}\label{eq:covi*}
T^I_- \subset \bigcup\nolimits_{Q \in \bigcup\nolimits_{8 \le k \le I}\mathcal{C}^k} Q' \subset Q_\mu,
\end{align}
where the second inclusion follows from \eqref{eq: Ts}  and the fact that $\partial Q_\mu \subset E^i$ for all $i \in \N$ by \eqref{eq: nested2}.  To  complete the covering, we introduce for all $k \ge I+1$
\begin{align}\label{eq: fkk}
\mathcal{G}_k = \lbrace Q \in \mathcal{Q}^k: \  Q'' \subset Q_\mu, \ \  Q'' \cap (E^I \cap \UUU J \EEE) = \emptyset \rbrace.
\end{align}
 Assuming that we have already constructed $\mathcal{C}^I, \ldots, \mathcal{C}^n$ for $n \ge I$ we define
\begin{align}\label{eq: cov def}
\mathcal{C}^{n+1} = \big\{Q \in \mathcal{G}_{n+1}: Q \cap \bigcup\nolimits^{n}_{k=8} \bigcup\nolimits_{\hat{Q} \in \mathcal{C}^k} \hat{Q} = \emptyset  \big\}.
\end{align} 
Clearly, $\mathcal{C} := \bigcup_{k=8}^{\infty} \mathcal{C}^k$ is a covering of $Q_\mu$ (up to a set of negligible measure) consisting of pairwise disjoint, dyadic squares. Finally, for $8 \le i \le I$ we define with $U_i$ as given in \eqref{eq: V} 
\begin{align}\label{eq: bad sets-new}
\mathcal{Y}^i= \lbrace Q \in \mathcal{C}^i: Q'' \subset U_i, \  \mathcal{H}^1(\UUU J \EEE \cap Q') > \theta^2 s_i \rbrace, \ \ \ \ \  Y^i=\sat\Big(\bigcup\nolimits_{Q \in \mathcal{Y}^i} \overline{Q''}\Big),
\end{align}
\BBB where $\sat(\cdot)$ was defined below  \eqref{eq: enlarged squares}. \EEE Let $Z = \bigcup\nolimits_{8 \le i \le I} Y^i$ and observe that \eqref{eq: U} holds since $Y^i \subset U_i$. \UUU For some explanation and motivation about the definition of the set $Z$ we refer to Remark \ref{rem: Z}(ii) below. \EEE

\smallskip

\emph{Proof of \eqref{eq: main pro}(i):}  First, we have already noticed \eqref{eq:covi*}. Then taking \eqref{eq: fkk}  into account, using $E^I \cap \UUU J \EEE \BBB = \EEE \bigcup\nolimits_{D^I_k \in \mathcal{E}^I} \Gamma_{j_k}^* \cap \UUU J \EEE \subset \bigcup\nolimits_{j=1}^m \partial P'_j\cap \UUU J \EEE$ (see \eqref{eq: S const}) and arguing as in the construction of a Whitney covering for open sets, we obtain \eqref{eq: main pro}(i).

\smallskip

\emph{Proof of \eqref{eq: main pro}(ii):}  We fix $Q_1 \in \mathcal{C}^i$ and $Q_2 \in \mathcal{C}^j$ with $j \ge i$ and $Q_1' \cap Q_2' \neq \emptyset$. (a) If $i \ge I+1$,
 we deduce $j \le i+1$ in view of \eqref{eq: fkk}-\eqref{eq: cov def}, where one argues as in the construction of a Whitney covering. (b) If $i \le I-1$, we suppose that we had $j \ge i+2$. Then by \eqref{eq:covi} we get $Q_1 \subset T^i_-$ and $Q_2 \cap T^{i+1}_- = \emptyset$. Consequently, using \eqref{eq: gooddisti} we obtain $\dist(Q_1,Q_2) \ge  s_i$ and thus the contradiction $Q_1' \cap Q_2' = \emptyset$. 
 
(c) Finally, for the case $i= I$ it suffices to show that all $Q \in \mathcal{Q}^{I+1}$ with $Q  \subset Q_\mu \setminus T^I_-$ and $Q_1' \cap Q' \neq \emptyset$ fulfill $Q \in \mathcal{C}^{I+1}$. First, since $Q_1' \cap Q' \neq \emptyset$, $Q$ satisfies $Q'' \subset T^I$ by \eqref{eq: Ts} and \eqref{eq:covi}.  Then $Q'' \cap E^I = \emptyset$ by \eqref{eq: Ts}, which also implies $Q'' \subset Q_\mu$ since $\partial Q_\mu \subset E^I$  \BBB (see \eqref{eq: nested2}). \EEE Consequently, $Q \in \mathcal{G}_{I+1}$ and thus  $Q \in \mathcal{C}^{I+1}$ as $Q  \subset Q_\mu \setminus T^I_-$.

 \smallskip

\emph{Proof of \eqref{eq: main pro}(iii):} Property  \eqref{eq: main pro}(iii) follows directly  from \eqref{eq: main pro}(ii) recalling the fact that each $x \in Q_\mu$ is contained in at most four sets $Q'$, $Q \in \mathcal{Q}^i$.

\smallskip

\emph{Proof of \eqref{eq: main pro}(vi):}  We note that each $Q \in \mathcal{C}^i$ for $i \ge I+1$ satisfies $Q \cap (Q_\mu \setminus T^I_-) \neq \emptyset$ and thus $\dist(E^I, Q') \le cs_I$ by \eqref{eq: Ts}. Recalling \eqref{eq: elarge3}, $E^I = \bigcup_{D^I_k \in \mathcal{E}^I} D^I_k$, and \eqref{eq: fkk} we derive  
\begin{align}\label{eq: new constant}
\dist(Q',\bigcup\nolimits_{D^I_k \in \mathcal{E}^I} \Gamma_{j_k}^* ) \le cs_I, \ \ \ \ \ \ \ \ \dist(Q',\bigcup\nolimits_{D^I_k \in \mathcal{E}^I} \Gamma_{j_k}^* \cap \UUU J \EEE) >0
\end{align}
for some universal $c>0$. Consequently, \eqref{eq: distance}(i) for $\theta$ small enough   yields (vi). It now remains to prove \eqref{eq: main pro}(iv),(v), where by \eqref{eq: main pro}(vi) it suffices to consider $8  \le i \le I$.

\smallskip

\emph{Proof of \eqref{eq: main pro}(iv):}  Fix $Q \in \mathcal{C}^i$, $8 \le i \le I$, with $Q'' \not\subset Z$. By \eqref{eq: bad sets-new} we can even assume that $Q'' \not\subset U_i$ since $Q'' \not\subset Z$ together with $Q'' \subset U_i$ already implies \eqref{eq: main pro}(iv). \BBB Moreover, we suppose that $Q \in  \mathcal{Q}^i_{\rm bad}$ since otherwise the assertion follows directly from \eqref{eq: bad1}. \EEE Let $Q_* \in \mathcal{Q}^{i-1}$ with $Q \subset Q_*$.  \BBB We distinguish three cases and refer  also to Figure \ref{korn3}(b) for an illustration: \EEE

\BBB (a) Suppose that  $Q \subset T^{i-1} \setminus T^{i-1}_-$. As $T^{i-1}$ consists of squares in $\mathcal{Q}^{i-1}$, \EEE we have $Q_* \subset T^{i-1} \setminus T^{i-1}_-$ and  then \EEE $Q_* \notin \mathcal{Q}^{i-1}_{\rm bad}$ by \eqref{eq: nocrack}. Thus, by   \eqref{eq: bad1} we get 
\begin{align}\label{eq: badspec}
\mathcal{H}^1(\UUU J \EEE \cap Q') \le \mathcal{H}^1(\UUU J^* \EEE \cap Q_*') \le \theta^{3} s_{i-1} = \theta^2  s_i.
\end{align}

\BBB (b) Suppose that  $Q \subset T^i_- \setminus T^{i-1}$ and   $Q_*  \notin \mathcal{Q}^{i-1}_{\rm bad}$. In this case, the assertion follows by repeating the argument in \eqref{eq: badspec}. \EEE  

\BBB (c) Finally, suppose that $Q \subset T^i_- \setminus T^{i-1}$ and $Q_* \in\mathcal{Q}^{i-1}_{\rm bad}$. We will show that this assumption leads to a contradiction. First, note that  $\overline{Q_*'''} \subset  B^{i-1}$ by \eqref{NNNNN}. \EEE As $T^{i-1}$ consists of squares in $\mathcal{Q}^{i-1}$, we get $Q_* \cap  T^{i-1} = \emptyset$ and then $Q_*  \cap E^{i-1} \neq \emptyset$ by \eqref{eq: Ts}. Consequently, recalling \eqref{eq: E large2} (for $i-1$) we get that $\overline{Q_*'''}$ is contained in a component of $\mathcal{E}^{i-1}$. This implies
\begin{align}\label{eq: impli}
Q \subset \overline{Q_*'''} \subset E^{i-1}.
\end{align}
\BBB Recall  that $Q \in \mathcal{Q}^i_{\rm bad}$ and $Q'' \not\subset U_i$. \EEE Consequently, $\overline{Q'''}$ is contained in a component $\mathcal{B}^i \setminus \mathcal{B}^i_{\rm iso}$ and therefore by the  remark below \eqref{eq: bound2},   $\overline{Q'''}$ is contained in a component of $\mathcal{D}^i$. Since $Q \subset E^{i-1}$ \BBB by \eqref{eq: impli}, \EEE  by \eqref{eq: E large2} (for $i$) we derive  $\overline{Q'''} \subset E^i$.  This, however, contradicts the fact that  $Q \in \mathcal{C}^i$. \BBB Indeed, by \eqref{eq:covi} we have  $Q \subset T^i_-$, and since   $T^i_- \cap E^i = \emptyset$, we get $Q \cap E^i = \emptyset$. \EEE This concludes the proof of (iv). 
 
 \smallskip

\emph{Proof of \eqref{eq: main pro}(v):} We suppose that for $Q \in \mathcal{C}^i$, $8 \le i \le I$, there is a neighbor $\hat{Q} \in \mathcal{C}^{i-1}$ such that $Q'' \cap \hat{Q}'' \neq \emptyset$. Then by construction we find $\hat{Q} \subset T^{i-1}_- \setminus  T^{i-1}_{--}$ and a further square $Q_* \in \mathcal{Q}^{i-1}$ with $Q \subset Q_* \subset T^{i-1} \setminus T^{i-1}_-$. In view of \eqref{eq: nocrack}, the desired property follows as in \eqref{eq: badspec}.   Likewise, if there is a neighbor $\hat{Q} \in \mathcal{C}^{i+1}$ with $Q'' \cap \hat{Q}'' \neq \emptyset$, we see $Q \subset T^{i}_- \setminus  T^{i}_{--}$ and then again by \eqref{eq: nocrack} we get the claim from \eqref{eq: bad1}. \eop

 \begin{rem}\label{rem: Z}
{\normalfont \UUU 

(i) For later reference, we remark that assumption \eqref{eq: concentration-new} was necessary to obtain the estimate  \eqref{eq: bad length} on the length of the boundary of the partition (see \eqref{eq: NNNNN} and \eqref{eq: Ju*XX}).

(ii) Let us comment  briefly on the definition of $Z$ in \eqref{eq: bad sets-new}. \EEE In the definition of $\mathcal{Y}^i$, the condition on the amount of jump contained in a square, which  differs from \eqref{eq: bad1}, is chosen in such a way that we will be able to use   \eqref{eq: main pro}(v) in the proof of Lemma  \ref{lemma: bad sets} below. In contrast to the definition of $U_i$, we use squares $\overline{Q''}$ instead of $\overline{Q'''}$ in the definition of $Y^i$. This will be essential to prove property \eqref{eq: main proXXX}(ii). \UUU We again use the saturation of sets to ensure that the connected components of $Y^i$ are simply connected. \EEE  Note that $Y^i$ has empty intersection with the set $E^i$, which is the main object for the construction of the partition. Therefore, we call the components of $Y^i$ \emph{isolated}.  
}
\end{rem}

 \noindent {\em Proof of Lemma \ref{lemma: bad sets}.} Let $\mathcal{C}$ be given satisfying \eqref{eq: main pro} and $I$ as in \eqref{eq: distance}. \BBB Recall the definition of $U_i$ in \eqref{eq: V} and the definitions of $Y^i$, $\mathcal{Y}^i$, and $Z = \bigcup\nolimits_{8 \le i \le I} Y^i$ in  \eqref{eq: bad sets-new}. In the following proof we will particularly need property  \eqref{eq: main pro}(v). Moreover, at the end we will redefine the Whitney covering on the components of $Z$ and in this context, we will exploit the construction of $\mathcal{C}$ in \eqref{eq:covi}-\eqref{eq: cov def}. \EEE We denote the (closed) connected components of $Y^i$ by $(Y^i_k)_k$.

 \smallskip
\emph{Properties of $(Y^i_k)_k$:}  We first show that
 \begin{align}\label{eq: la2}
{Q} \in \mathcal{C}^i,   {Q} \subset Y^i_k  \ \  \text{ for all } \ \  Q \in  \hat{\mathcal{Y}}^i_k :=  \lbrace {Q} \in \mathcal{C}:  {Q} \cap Y^i_k \neq \emptyset, \  \partial  {Q} \cap \partial Y^i_k \neq \emptyset\rbrace.
\end{align}
 In fact, for each ${Q} \in \hat{\mathcal{Y}}^i_k$  there is $\tilde{Q} \in \mathcal{Y}^i$, $\tilde{Q}'' \subset Y^i_k$ such that $|\tilde{Q}'' \cap  {Q}|>0$. \BBB Since $\tilde{Q}$ satisfies $\mathcal{H}^1(\UUU J \EEE \cap \tilde{Q}')>\theta^2s_i$ by \eqref{eq: bad sets-new}, \EEE property \eqref{eq: main pro}(v) yields $\tilde{Q}'' \subset \bigcup\nolimits_{{Q}_* \in \mathcal{C}^i} \overline{{Q}_*}$. As $\mathcal{C}$ consists of pairwise disjoint dyadic squares, this gives ${Q} \in \mathcal{C}^i$ and ${Q} \subset \tilde{Q}'' \subset Y^i_k$, as desired.  
 
\begin{figure}[H]
\centering
\begin{overpic}[width=0.89\linewidth,clip]{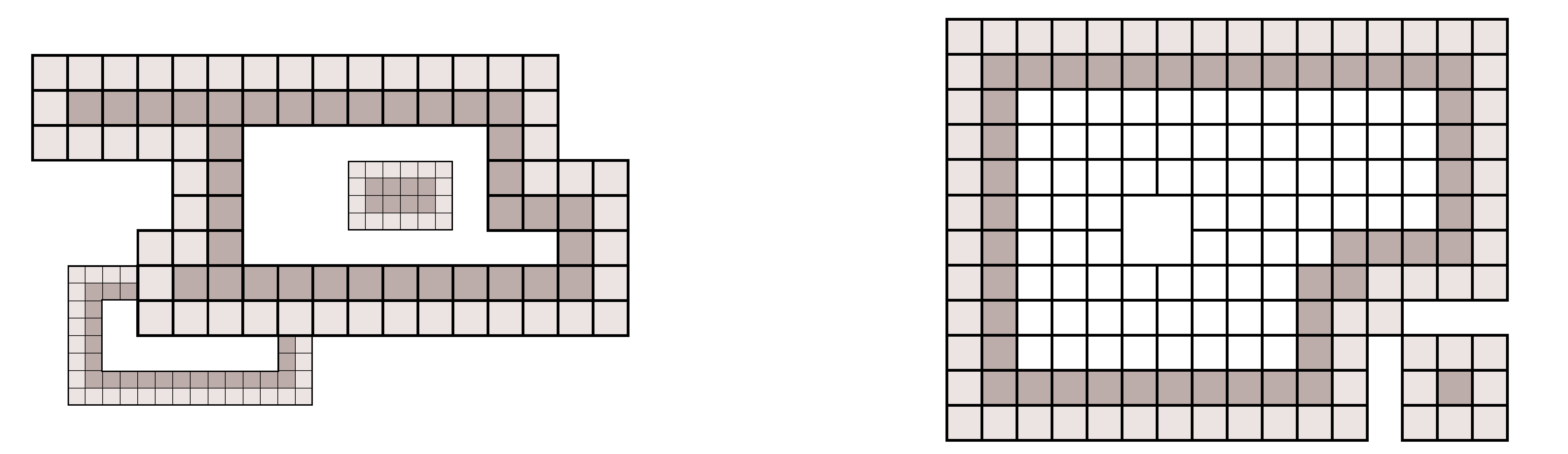}

\put(-8,86){{$(a)$}}
\put(205,86){{$(b)$}}

\put(136,80){{$\partial X^i_k$}}
\put(130,86) {\line(1,0){6}}

\put(150,60){{$\mathcal{X}^i_k = \hat{\mathcal{X}^i_k}$}}
\put(143,64) {\line(1,0){6}}

\put(155,35){{$\mathcal{Y}^i$}}
\put(133,41) {\line(1,0){21}}

\put(66,57){{$Y^j_1$}}
\put(75,12){{$Y^j_2$}}

\end{overpic}
\caption{(a) The squares $\mathcal{X}^i_k$ are depicted in light gray and the dark gray squares are contained in $\mathcal{Y}^i$. It is possible that $Y^j_1 \subset X^i_k$, whereas $Y^j_2$, which overlaps with $X^i_k$, cannot exist since this contradics \eqref{eq: main pro}(v). \ (b) We have  illustrated the part of $\mathcal{C}$ contained in $X^i_k$, which may (in contrast to $\mathcal{C}_*$) contain squares in $\mathcal{Q}^j$, $j < i$.} \label{korn4}
\end{figure}
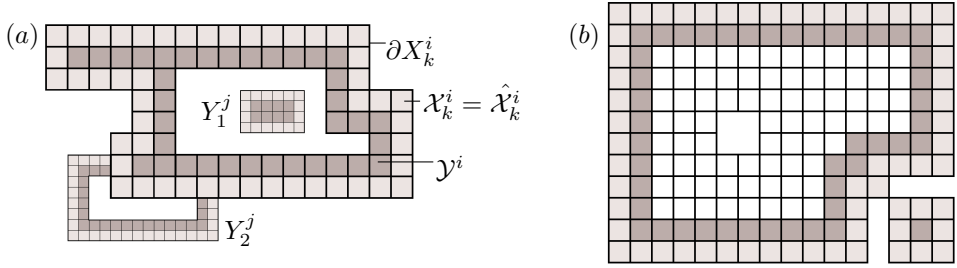

 We now see that   each pair $Y^i_k$, $Y^j_l$, $i < j$, either satisfies
\begin{align}\label{eq: pairwisedisjoint}
\dist_\infty(Y^i_k,Y^j_l) \ge \tfrac{1}{2}s_i
\end{align}
  or one set is contained in the other. (Here $\dist_\infty$ denotes the distance with respect to the maximum norm.) If not, we would find $Q_i \in \hat{\mathcal{Y}}^i_k$ and $Q_j \in \hat{\mathcal{Y}}^j_l$  such that $Q'_j \cap Q_i' \neq \emptyset$. Then $j=i+1$  by  \eqref{eq: main pro}(ii) and \eqref{eq: la2}.  Moreover, \BBB $Q_j \cap Q_i'' \neq \emptyset$ and \EEE there would be $Q_* \in \mathcal{Y}^j$, $Q_*'' \subset Y^j_l$, such that $Q''_* \cap Q_i'' \neq \emptyset$. This, however, contradicts \eqref{eq: main pro}(v) and the fact that $\mathcal{H}^1(\UUU J \EEE \cap {Q}_*') > \theta^2 s_j$ by \eqref{eq: bad sets-new}. (We also refer to Figure \ref{korn4}(a).)
    
    \smallskip
\emph{Components $(X^i_k)_k$:} Now for each $8 \le i \le I$, let $(X^i_k)_k \subset (Y^i_k)_k$ be the components such that for all $j \neq i$ we have
$X^i_k \not\subset Y^j_l$  for all components  $Y^j_l.$ Accordingly, we denote the sets defined in \eqref{eq: la2} by $\hat{\mathcal{X}}^i_k$. Define $Z^i = \bigcup_k X^i_k$ and note that by \eqref{eq: pairwisedisjoint} the sets $(Z^i)_i$ are pairwise disjoint with $Z  = \bigcup_{i =8}^I Z^i$.  Recalling the definition of $\mathcal{X}^i_k$ before \eqref{eq: main proXXX}, we then see that $\hat{\mathcal{X}}^i_k = {\mathcal{X}}^i_k$. \UUU As the sets $Y^i$ are defined using the saturation (see \eqref{eq: bad sets-new}), we also get that the components $(X^i_k)_k$ are simply connected. \EEE We now show \eqref{eq: main proXXX1} and \eqref{eq: main proXXX}.

 \smallskip
 
 \emph{Proof of \eqref{eq: main proXXX}:}  As \BBB by \eqref{eq: bad sets-new} \EEE each component  $(X^i_k)_k$ of $Z^i$ is contained in a component of $U_i$, \eqref{eq: main proXXX}(i) follows directly from \eqref{eq: crack bound}. By \eqref{eq: pairwisedisjoint} and the fact that $\dist_\infty(X^i_k,X^i_l) \ge 2s_i$ for all $i$ and $k\neq l$, we obtain \eqref{eq: main proXXX}(ii).  The proof of \eqref{eq: la2} showed that for each square in $\mathcal{X}^i_k$ there is an adjacent square in $\mathcal{Y}^i$ (cf. Figure \ref{korn4}(a)). Consequently,  we obtain by \eqref{eq: bad sets-new} 
\begin{align}\label{eq: repeat estimate}
\sum\nolimits_k \# \mathcal{X}^i_k &\le 8 \sum\nolimits_k\# \lbrace Q \in \mathcal{Y}^i: Q'' \subset X^i_k \rbrace\le c\theta^{-2}s_i^{-1}\sum\nolimits_{Q \in \mathcal{Y}^i} \mathcal{H}^1(Q' \cap Z^i\cap \UUU J \EEE )\notag \\& \le c\theta^{-2}s_i^{-1}\mathcal{H}^1(\UUU J \EEE \cap Z^i),
\end{align}
for $c>0$ universal, where we  used that each $x \in Q_\mu$ is contained in at most four enlarged squares. This gives \eqref{eq: main proXXX}(iii).  

\smallskip
\emph{Proof of \eqref{eq: main proXXX1}:} Since   $(Z^i)_i$ are pairwise disjoint,  \eqref{eq: repeat estimate} also yields
\begin{align}\label{eq: EEEE} 
\sum\nolimits_{i = 8}^I \mathcal{H}^1(\partial Z^i) & \le \sum\nolimits_{i = 8}^I \sum\nolimits_k \mathcal{H}^1(\partial X^i_k) \le \sum\nolimits_{i = 8}^I  \big(  \sum\nolimits_{Q \in \bigcup_k\mathcal{X}^i_k}  \mathcal{H}^1(\partial Q)\big) \notag\\&\le \sum\nolimits_{i = 8}^I  c\theta^{-2}\mathcal{H}^1(\UUU J \EEE \cap Z^i)\le  c\theta^{-2}\mathcal{H}^1(\UUU J \EEE),
\end{align}
i.e., \eqref{eq: main proXXX1}(ii) holds. Moreover, by \eqref{eq: main proXXX}(i) and the fact that $d(X^i_k) \le \mathcal{H}^1(\partial X^i_k)$  
 $$|Z^i| =  \sum\nolimits_k |X^i_k| \le  \sum\nolimits_k d(X^i_k)^2 \le \theta^{-ir} s_i  \sum\nolimits_k \mathcal{H}^1(\partial X^i_k) \le c\theta^{-ir}s_i \theta^{-2}\mathcal{H}^1(J).$$
Likewise, \BBB note that \eqref{eq: main proXXX}(i) and the fact that  \UUU $r < \frac{1}{8}$ \EEE yield  $d(X^i_k) \le \mu \theta^7$. We then conclude the proof of \eqref{eq: main proXXX1} by \BBB using   \eqref{eq: EEEE} \EEE and calculating
 $$|Z| = \big|\bigcup\nolimits_{i = 8}^I Z^i\big| \le  \mu \theta^7 \sum\nolimits_{i = 8}^I\sum\nolimits_k \mathcal{H}^1(\partial X^i_k) \le c\mu\theta^{5}\mathcal{H}^1(\UUU J \EEE).$$
 
\smallskip

\emph{Redefinition of the covering:} Finally, we define a modification $\mathcal{C}_*$ of $\mathcal{C}$ such that each $Z^i$ consists of squares $\mathcal{C}^i_* := \mathcal{C}_* \cap \mathcal{Q}^i$. To this end, consider a component $X^i_k$. By \eqref{eq:covi} and \eqref{eq: la2} we have $\bigcup_{Q \in \hat{\mathcal{X}}^i_k} Q \subset T^i_- \setminus T^{i-1}_-$, which by \eqref{eq: Ts} has empty intersection with $E^i$. Since the diameter of each connected component of $E^i$ exceeds $s_i\theta^{-ir}$ (see \eqref{eq: bound2}), \eqref{eq: main proXXX}(i) then yields $X^i_k \cap E^i = \emptyset$.  Consequently, $X^i_k \subset T^i_-$ and again by \eqref 
{eq:covi} we then obtain $X^i_k \subset \bigcup_{j=8}^i \bigcup_{Q \in \mathcal{C}^j}\overline{Q}$ (cf. Figure \ref{korn4}(b)). This together with the fact that $\hat{\mathcal{X}}^i_k \subset \mathcal{Q}^i$ by \eqref{eq: la2}, shows that we can replace the covering $\mathcal{C}$ in each component $X^i_k$ by a covering $\mathcal{C}_*$ consisting exclusively of squares in $\mathcal{Q}^i$ such that all conditions in \eqref{eq: main pro} remain true.   \eop

\vspace{-0.1cm}
 
 \vspace{-0.35cm} 
\BBB

\subsection{Construction of a piecewise smooth approximation}\label{sec: 3sub2}

In Theorem \ref{th: bad part}, \UUU given a finite union of closed segments $J$, \EEE  we have constructed an auxiliary partition $(P'_j)_{j=1}^m$ consisting of simply connected sets and an associated  covering $\mathcal{C}$ of Whitney-type such that in the squares of the covering the \UUU $\mathcal{H}^1$-norm of $J$ \EEE is small compared to their diameter. Recall that this construction was exclusively based on the geometry of \UUU $J$. \EEE In this section, \UUU applying Theorem \ref{th: bad part} for the jump set $J=J_u$ of $u \in \mathcal{W}(Q_\mu)$, \EEE we will  see that, with the help of the Korn inequality  for functions with small jump set (Theorem \ref{th: kornSBDsmall}), we can approximate $u$  on the squares of $\mathcal{C}$ by  infinitesimal rigid motions. Proceeding similarly as in  \cite{Friedrich:15-3}, by using a partition of unity associated to $\mathcal{C}$, this will allow us to define a piecewise smooth approximation $\bar{u}$ of $u$ which is smooth on every component $P_j'$. 

Note that an additional difficulty has to be faced: due to technical reasons in the construction of Theorem \ref{th: bad part}, there exists a set $Z$ of isolated components given in Lemma \ref{lemma: bad sets} where Theorem \ref{th: kornSBDsmall} is not applicable and $\bar{u}$ has to be defined differently. Consequently, besides the main properties of the covering  \eqref{eq: main pro}, we will also use the structure of $Z$ described in \eqref{eq: main proXXX1}-\eqref{eq: main proXXX}. Note, however, that no other details of the construction from the previous section will be needed.

\begin{theorem}\label{th: modifica}
Let $\mu>0, \theta > 0$ small and   $p \in [1,2)$. Then there are a universal $c>0$ and a constant    $C=C(\theta,p)>0$ such that for all $u \in \mathcal{W}(Q_\mu)$ with \eqref{eq: concentration} the following holds: let  $(P'_j)_{j=1}^m$ be the partition and $\mathcal{C}$ be the covering  given by Theorem \ref{th: bad part} and Lemma \ref{lemma: bad sets} \UUU applied for $J = J_u$ and $r= \frac{1}{24}(2-p)$. \EEE Then there is \BBB an approximation \EEE $\bar{u}: Q_\mu \to \R^2$, being smooth in  $\bigcup\nolimits_{Q \in \mathcal{C}} Q' \supset Q_\mu \setminus (J_u \cap \bigcup_{j=1}^m \partial P'_j)$,  and an exceptional set $F \subset Q_\mu$ with
\begin{align}\label{eq: except***}
|F| \le c\mu   \theta^5  \, \mathcal{H}^1(J_u), \ \ \ \ \  \mathcal{H}^1(\partial F) \le C\mathcal{H}^1(J_u) 
\end{align} 
such that $\bar{u} = u$ on $\partial Q_\mu$ (in the sense of traces) and 
\begin{align}\label{eq: modifica prop**}
\begin{split}
(i) & \ \ \Vert e(\bar{u}) \Vert_{L^p(Q_\mu)} \le C\mu^{\frac{2}{p}-1} \Vert e(u) \Vert_{L^2(Q_\mu)},\\
(ii) & \ \ \Vert \nabla \bar{u}  - \nabla u \Vert_{L^p(Q_\mu \setminus F)} \le C\mu^{\frac{2}{p}-1}  \Vert e(u) \Vert_{L^2(Q_\mu)},\\
(iii) & \ \ \Vert  \bar{u} - u  \Vert_{L^p(Q_\mu\setminus F)}\le C\mu^{\frac{2}{p}} \Vert e(u) \Vert_{L^2(Q_\mu)}.
\end{split}
\end{align}

\end{theorem}

\BBB
\UUU The approximation $\bar{u}$ is constructed by smoothing a piecewise affine approximation of $u$ given on the covering $\mathcal{C}$ and it is therefore smooth on $\bigcup\nolimits_{Q \in \mathcal{C}} Q'$. Recall $\bigcup\nolimits_{Q \in \mathcal{C}} Q' \supset Q_\mu \setminus (J_u \cap \bigcup_{j=1}^m \partial P'_j)$ by \eqref{eq: main pro}(i). This implies  that a connected component of $\bigcup\nolimits_{Q \in \mathcal{C}} Q'$ (i.e., a connected component of the set where $\bar{u}$ is smooth) contains    sets $(P'_j)^m_{j=1}$. (It may contain indeed several of these sets if they are not completely disconnected by the jump set $J_u$.) Particularly, $\bar{u}$ is smooth on each of the sets $P_j'$, which will be the fundamental property used in the following. To highlight this fact, we will call $\bar{u}$ a \emph{piecewise smooth approximation} of $u$, where the word \emph{piecewise} refers to the partition $(P'_j)^m_{j=1}$. 

\EEE

Note that   $\bar{u}$ plays the role of an auxiliary function and does not appear in the statement of Theorem \ref{th: korn-small set}. Later in Section \ref{sec: 3sub5}, after passing to a refined partition consisting of John domains with uniformly controlled John constant,   we will   use the  Korn inequality in John domains (Theorem \ref{th: kornsobo}) to obtain an estimate on $\nabla \bar{u}$. Finally, \eqref{eq: modifica prop**}(ii)  will then yield control on $\nabla u$ outside of $F$. 

In particular, $F$ will be a part of the exceptional set of Theorem \ref{th: korn-small set} and \eqref{eq: except***} will be needed to show \eqref{eq: except}. (We refer to Remark \ref{rem:modi}(i) below for more details on the structure of $F$.) Note that the assumption $p < 2$ is essential since at some point of the proof we need to reduce  the exponent $2$ in a H\"older-type estimate.   \EEE

 \smallskip

\Proof  We apply Theorem \ref{th: bad part} and Lemma \ref{lemma: bad sets} \UUU for $J = J_u$ and $r= \frac{1}{24}(2-p)$. \EEE Let the partition $(P'_j)_{j=1}^m$, the isolated components $Z = \bigcup^I_{l=8} Z^l = \bigcup^I_{l=8} \bigcup_k X^l_k$ and the Whitney-type covering $\mathcal{C}$ of $Q_\mu \setminus (J_u \cap \bigcup_{j=1}^m \partial P'_j)$ be given such that \eqref{eq: main pro} holds and each $X^l_k$ is a union of squares in $\mathcal{C}^l = \mathcal{C} \cap \mathcal{Q}^l$. By $\mathcal{X}^l_k$ we denote the squares at the boundary of $X^l_k$, see before \eqref{eq: main proXXX}.

\smallskip

\BBB
\textbf{Outline of the proof:}  We first use \eqref{eq: main pro}(iv) and  Theorem \ref{th: kornSBDsmall} to approximate $u$ by infinitesimal rigid motions $a_Q$ (also called affine mappings in the following) on the family of squares $Q \in \mathcal{C}$ with $Q'' \not\subset Z$, which particularly includes the squares $\mathcal{X}^l_k$ at the boundary of the isolated components. Inside each isolated component  $X^l_k \subset Z$ we choose a single affine mapping $a^l_k$ with $a^l_k = a_Q$ for some (arbitrary) $Q \in \mathcal{X}^l_k$.

We then estimate the difference of the infinitesimal rigid motions on adjacent squares. For the isolated components $X^l_k$, we have to face the additional difficulty that  the affine mappings given on \emph{all} $\mathcal{X}^l_k$ have to be compared against each other. To do so, we draw some ideas from  \cite[Section 4]{FrieseckeJamesMueller:02}. Hereby we see that, due to the control $d(X^l_k) \le \theta^{-lr} s_l$ (see  \eqref{eq: main proXXX}(i)), the involved constant may  blow up in $l$, but is controlled in terms of $\theta^{-2(p+1)lr}$ (see \eqref{eq: Q1,Q2}).

  Finally, using a partition of unity associated to $\mathcal{C}$, we will construct a function being smooth in $\bigcup\nolimits_{Q \in \mathcal{C}} Q'$. To confirm \eqref{eq: modifica prop**}, we exploit the fact that the difference of the affine mappings on adjacent squares is controlled. In this context, we explicitly use $p<2$ and H\"older's inequality, which together with the  `smallness' of  $|Z^l|$ (see \eqref{eq: main proXXX1}(i))  allows us to compensate the fact that for isolated components the constants blow up in terms of $\theta^{-2(p+1)lr}$  (cf. \eqref{eq: main holder}-\eqref{eq: main holderXXXX}).  \EEE
  
  \smallskip
 
We now start with the proof. In the following $C=C(\theta,p)>0$ denotes a generic constant and $c>0$ is universal. We can suppose that $\theta$ is small with respect to $c$. 

\smallskip

\textbf{Step I (Definition of infinitesimal rigid motions):} \BBB In this step we define  infinitesimal rigid motions for each square of the Whitney covering. \EEE

\enlargethispage{\baselineskip}
\smallskip

\emph{Good squares:} Let us first consider the subset of \emph{good squares}
\begin{align}\label{eq: hatcal def}
{\mathcal{C}}_{\rm g} :=    \big\{ Q \in \mathcal{C}: Q'' \not\subset Z\rbrace \subset  \bigcup\nolimits_{i\ge \BBB 8 \EEE }\big\{ Q \in \mathcal{C}^i: \mathcal{H}^1(J_u \cap Q') \le \theta^2 s_i \big\},
\end{align} 
where the inclusion follows from {\eqref{eq: main pro}(iv). We apply Theorem \ref{th: kornSBDsmall}  on $Q'$, $Q \in {\mathcal{C}}_{\rm g}$,  to find  infinitesimal rigid motions $a_Q = a_{A_Q,b_Q}$ and exceptional sets $E_Q \subset Q'$ such that   
\begin{align}\label{eq: poincare estim}
\begin{split}
\hspace{-0.2cm}  d(Q)^{-\frac{2}{p}}\Vert u - a_Q \Vert_{L^{p}(Q' \setminus E_Q)}   + d(Q)^{1-\frac{2}{p}} \Vert \nabla u - A_Q \Vert_{L^{p}(Q' \setminus E_Q)} \le c  \Vert e(u)\Vert_{L^2(Q')}.
\end{split}
\end{align}
Moreover, taking  \eqref{eq: hatcal def} into account we find for all $Q \in {\mathcal{C}}_{\rm g}$
\begin{align}\label{eq: excpt}
|E_Q| \le c d(Q)\theta^2 \mathcal{H}^1(J_u \cap Q') \le c\theta^{4} d(Q)^2, \  \  \ \ \mathcal{H}^1(\partial E_Q) \le c \mathcal{H}^1(J_u \cap Q').
\end{align} 
As a preparation for Step II, we now estimate the difference of the infinitesimal rigid motions for neighboring squares in ${\mathcal{C}}_{\rm g}$.  For $Q \in {\mathcal{C}}_{\rm g}$ we let 
\begin{align}\label{eq: hat neigh}
{\mathcal{N}}_{\rm g}(Q) =  \lbrace \hat{Q} \in {\mathcal{C}}_{\rm g} \setminus \lbrace Q \rbrace: Q' \cap \hat{Q}' \neq \emptyset \rbrace
\end{align}
and observe that by \eqref{eq: main pro}(ii) we have $\theta d(\hat{Q}) \le d(Q) \le \theta^{-1}d(\hat{Q})$ for all $\hat{Q} \in {\mathcal{N}}_{\rm g}(Q)$. This also implies $\# {\mathcal{N}}_{\rm g}(Q) \le c\theta^{-2}$. Moreover, as the covering consists of dyadic squares, $Q' \cap \hat{Q}'$ contains a ball $B$ with radius larger than $c\min \lbrace d(Q), d(\hat{Q})\rbrace \ge c\theta d(Q)$ for some small universal $c>0$. Consequently, by \eqref{eq: excpt} we find $|\hat{E} \cap B|\le \frac{1}{2}|B|$ for $\theta$ sufficiently small, where $\hat{E} = E_{\hat{Q}} \cup E_Q$.  Thus, $|B \setminus \hat{E}| \ge \frac{1}{2}|B|\ge c\theta^2d(Q')^2$ and \BBB then we can apply  Lemma \ref{lemma: rigid} to find
$$\Vert a_{Q} - a_{\hat{Q}}\Vert^{p}_{L^{p}(Q')} \le c' (d(Q')^2|B \setminus \hat{E}|^{-1})^{\frac{p}{2}+1} \Vert a_{Q} - a_{\hat{Q}}\Vert^{p}_{L^{p}(B \setminus \hat{E})}  \le  c' \theta^{-(p+2)}\Vert a_{Q} - a_{\hat{Q}}\Vert^{p}_{L^{p}(B \setminus \hat{E})}  $$
for all $\hat{Q} \in {\mathcal{N}}_{\rm g}(Q)$, where $c' = c'(p)>0$.  Then by \eqref{eq: poincare estim}, $\theta d(\hat{Q}) \le d(Q)$, and the triangle inequality we find
\begin{align}\label{eq: diff1*} 
\Vert a_{Q} - a_{\hat{Q}}\Vert^{p}_{L^{p}(Q')}   \le c c' \theta^{-(p+2)} \theta^{-2} d(Q)^{2}  \Vert e(u) \Vert^{p}_{L^2(Q' \cup \hat{Q}')} \le  C d(Q)^{2}  \Vert e(u) \Vert^{p}_{L^2(Q' \cup \hat{Q}')}
\end{align} 
for all $\hat{Q} \in {\mathcal{N}}_{\rm g}(Q)$, where $C=C(\theta,p)>0$. \EEE Therefore, by $\# {\mathcal{N}}_{\rm g}(Q) \le c\theta^{-2}$ we get with $N_Q := \bigcup_{\hat{Q} \in {\mathcal{N}}_{\rm g}(Q) \cup \lbrace Q \rbrace} \hat{Q}'$
\begin{align}\label{eq: diff1}
\sum\nolimits_{\hat{Q} \in {\mathcal{N}}_{\rm g}(Q)}\Vert a_{Q} - a_{\hat{Q}}\Vert^{p}_{L^{p}(Q')} \le C d(Q)^{2}  \Vert e(u) \Vert^p_{L^{2}(N_Q)}.
\end{align} 

\smallskip

\emph{Squares in isolated components:} To conclude Step I, it remains to define infinitesimal rigid motions for the squares in $\mathcal{C} \setminus {\mathcal{C}}_{\rm g}$.  Consequently, in view of  \eqref{eq: hatcal def} and  Lemma \ref{lemma: bad sets}, we have to concern ourselves with the isolated components $(X^l_k)_k$ for  $8 \le l \le I$. Let $\mathcal{X}^l_k$ as introduced before \eqref{eq: main proXXX}. For each $Q \in \mathcal{X}^l_k$ we have $Q' \not\subset X^l_k$ and thus $Q' \not\subset Z$ by \eqref{eq: main proXXX}(ii). Therefore, $\mathcal{X}^l_k \subset \mathcal{C}_{\rm g}$ by the definition in \eqref{eq: hatcal def}. Thus, \eqref{eq: diff1} holds for all $Q \in \mathcal{X}^l_k$. 

Define $N^l_k = \BBB X^l_k \cup \EEE \bigcup\nolimits_{Q \in \mathcal{X}^l_k}  Q' $. Since $d(X^l_k) \le \theta^{-lr} s_l$ by \eqref{eq: main proXXX}(i), we have by Lemma \ref{lemma: rigid} for all $\tilde{Q} \in \mathcal{X}^l_k$ and $\hat{Q} \in {\mathcal{N}}_{\rm g}(\tilde{Q}) \cap \mathcal{X}^l_k$
$$\Vert a_{\tilde{Q}} - a_{\hat{Q}} \Vert_{L^p(N^l_k)} \le C\theta^{-(1 + \frac{2}{p})lr}\Vert a_{\tilde{Q}} - a_{\hat{Q}} \Vert_{L^p(\tilde{Q}')}.
$$
As $X^l_k$ is \UUU simply \EEE connected, for all $Q_1,Q_2 \in \mathcal{X}^l_k $ we find a chain of squares in $\mathcal{X}^l_k$ connecting $Q_1$ with $Q_2$  \BBB (cf. Figure \ref{korn4}). \EEE The previous estimate and the triangle inequality  then yield
\begin{align*}
\max_{Q_1,Q_2 \in \mathcal{X}^l_k} \Vert a_{Q_1} - a_{Q_2} \Vert_{L^p(N^l_k)} & \le  \BBB \sum\nolimits_{\tilde{Q} \in \mathcal{X}^l_k} \sum\nolimits_{\hat{Q} \in {\mathcal{N}}_{\rm g}(\tilde{Q}) \cap \mathcal{X}^l_k}\Vert a_{\tilde{Q}} - a_{\hat{Q}}\Vert_{L^p(N^l_k)} \EEE \\
& \le C\theta^{-(1 + \frac{2}{p})lr}\sum\nolimits_{\tilde{Q} \in \mathcal{X}^l_k} \sum\nolimits_{\hat{Q} \in {\mathcal{N}}_{\rm g}(\tilde{Q}) \cap \mathcal{X}^l_k}\Vert a_{\tilde{Q}} - a_{\hat{Q}}\Vert_{L^p(\tilde{Q}')}.
\end{align*}
Therefore,   the discrete H\"older inequality together with  $\# \mathcal{X}^l_k \le c  \theta^{-2lr}$ (see \eqref{eq: main proXXX}(i)), $\# {\mathcal{N}}_{\rm g}(\tilde{Q}) \cap \mathcal{X}^l_k \le 8$ for all $\tilde{Q} \in \mathcal{X}^l_k$  and    \eqref{eq: diff1*}  implies
\begin{align}\label{eq: Q1,Q2}
\max_{Q_1,Q_2 \in \mathcal{X}^l_k} \Vert a_{Q_1} - a_{Q_2} \Vert^p_{L^p(N^l_k)} & \le C\theta^{\BBB -2(p+1)lr \EEE }\Big(\sum\nolimits_{\tilde{Q} \in \mathcal{X}^l_k} \sum\nolimits_{\hat{Q} \in {\mathcal{N}}_{\rm g}(\tilde{Q}) \cap \mathcal{X}^l_k}\Vert a_{\tilde{Q}} - a_{\hat{Q}}\Vert^2_{L^p(\tilde{Q}')}\Big)^{\frac{p}{2}} \notag  \\
& \BBB \le Cs_l^2\theta^{-2(p+1)lr}\Big(\sum\nolimits_{\tilde{Q} \in \mathcal{X}^l_k} \sum\nolimits_{\hat{Q} \in {\mathcal{N}}_{\rm g}(\tilde{Q}) \cap \mathcal{X}^l_k}\Vert e(u)\Vert^2_{L^2(\tilde{Q}' \cup \hat{Q}')}\Big)^{\frac{p}{2}} \EEE \notag  \\
& \le C s_l^2\theta^{-2(p+1)lr}  \Vert e(u) \Vert^p_{L^2(N_k^l)}.
\end{align}
(See \cite[Section 4]{FrieseckeJamesMueller:02} for similar arguments.) For each $X^l_k$ we define an infinitesimal rigid motion by setting $a^l_k = a_Q$ for an arbitrary $Q \in \mathcal{X}^l_k$.

 We can now define affine mappings associated to each $Q \in \mathcal{C}\setminus \mathcal{C}_g$. Given $Q \in \mathcal{C}^l$ with $Q'' \subset Z^l$ we choose the  component $X^l_k$ with $Q'' \subset X^l_k$ and define $a_Q = a^l_k$. 
 
 \smallskip

\textbf{Step II (Difference of affine mappings on neighboring squares):} \BBB In Step I we have already defined  infinitesimal rigid motions for each square of the Whitney covering. In particular,  $a_Q = a^l_k$ is constant for all $Q \in \mathcal{C}$ with $Q'' \subset X^l_k$, i.e., for squares lying well inside an isolated component.  For good squares $\mathcal{C}_g$ defined in  \eqref{eq: hatcal def} we refer to \eqref{eq: poincare estim}.  \EEE Similarly to \eqref{eq: hat neigh}, we  introduce the neighbors of each $Q \in \mathcal{C}$ defined by  
 $$\mathcal{N}(Q) = \lbrace \hat{Q} \in \mathcal{C} \setminus \lbrace Q\rbrace: Q' \cap \hat{Q}' \neq \emptyset \rbrace.$$
\BBB Note that $\# {\mathcal{N}}(Q) \le c\theta^{-2}$ by \eqref{eq: main pro}(ii). We now compare the  infinitesimal rigid motions on neighboring squares by  using \eqref{eq: diff1} and  \eqref{eq: Q1,Q2}. To this end, we distinguish three cases:\EEE
 
\BBB (a) \EEE First, consider $Q \in \mathcal{C}$ with   $Q \cap Z = \emptyset$. Then   $\mathcal{N}(Q) \subset {\mathcal{C}}_{\rm g}$ since each $\hat{Q} \in \mathcal{N}(Q) $ satisfies  $\hat{Q}'' \not\subset Z$. Consequently, the results of Step I are applicable and we find by \eqref{eq: diff1}  
\begin{align}\label{eq:11}
\sum\nolimits_{\hat{Q} \in \mathcal{N}(Q)}\Vert a_{Q} - a_{\hat{Q}}\Vert^{p}_{L^{p}(Q')} \le C d(Q)^2  \Vert e(u) \Vert^p_{L^{2}(N_Q)}.
\end{align}

\BBB (b) Consider $Q \in \mathcal{C}$ with $Q \subset X^l_k$ and $Q''' \not\subset X^l_k$ (recall \eqref{eq: enlarged squares}). Then either $Q \in \mathcal{X}^l_k$, with $\mathcal{X}^l_k$ as introduced before \eqref{eq: main proXXX}, or  $Q$ is a neighbor of a square in $\mathcal{X}^l_k$ (see  light and dark gray squares in  Figure \ref{korn4}(b)).  \EEE  Recall that $a_Q$ as in \eqref{eq: poincare estim} if $Q \in \mathcal{X}^l_k$ and  $a_Q = a^l_k$  otherwise. Moreover, by \eqref{eq: main proXXX}(ii) we get $a_{\hat{Q}} = a^l_k$ for all $\hat{Q} \in \mathcal{N}(Q) \setminus \mathcal{N}_{\rm g}(Q)$.   We apply  \eqref{eq: diff1}-\eqref{eq: Q1,Q2}  and find using $\#\mathcal{N}(Q) \le c\theta^{-2}$ and the definition of $a^l_k$
\begin{align}\label{eq: nonumber?} 
\sum\nolimits_{\hat{Q} \in \mathcal{N}(Q)}\Vert a_{Q} - a_{\hat{Q}}\Vert^{p}_{L^{p}(Q')}  \le Cs_l^{2} \Vert e(u) \Vert^p_{L^2(N_Q)}+ Cs_l^{2} \theta^{-2(p+1)lr}\Vert e(u) \Vert^p_{L^2(N^l_k)},
\end{align}
\BBB where we recall $N_Q = \bigcup_{\hat{Q} \in {\mathcal{N}}_{\rm g}(Q) \cup \lbrace Q \rbrace} \hat{Q}'$ and $N^l_k = X^l_k \cup \bigcup\nolimits_{Q \in \mathcal{X}^l_k} Q'$. \EEE

\BBB (c) Finally, if $Q \subset X^l_k$, $Q''' \subset X^l_k$ \EEE (see the white squares in Figure \ref{korn4}(b)), we observe that all neighbors lie in the same isolated component and thus by definition 
\begin{align}\label{eq:33}
\sum\nolimits_{\hat{Q} \in \mathcal{N}(Q)}\Vert a_{Q} - a_{\hat{Q}}\Vert^{p}_{L^{p}(Q')} =0.
\end{align}
We now sum  over all squares and obtain collecting \eqref{eq:11}-\eqref{eq:33}
\begin{align*}
H&:= \sum\nolimits_{Q \in \mathcal{C} } \sum\nolimits_{\hat{Q} \in \mathcal{N}(Q)} d(Q)^{-p} \Vert a_Q - a_{\hat{Q}} \Vert^p_{L^p(Q')} \\
&\le C\sum\nolimits_{l \ge 8} \sum\nolimits_k \# {\mathcal{X}}^l_k \ s_l^{2-p} \theta^{-2(p+1)lr} \Vert e(u)\Vert^p_{L^2(N^l_k)} +  C\sum\nolimits_{Q \in {\mathcal{C}}} d(Q)^{2-p}  \Vert e(u) \Vert^p_{L^2(N_Q)},
\end{align*}
\BBB where we used that the number of squares in $X^l_k$ of type (b) is bounded by $9\#{\mathcal{X}}^l_k$. \EEE
Since $\#\mathcal{X}^l_k \le c\theta^{-2lr}$ by \eqref{eq: main proXXX}(i), we get
$$H \le C\sum\nolimits_{l \ge 8} \sum\nolimits_k (\# {\mathcal{X}}^l_k)^{1-\frac{p}{2}} \ s_l^{2-p} \theta^{-(3p+2)lr} \Vert e(u)\Vert^p_{L^2(N^l_k)} +  C\sum\nolimits_{Q \in {\mathcal{C}}} d(Q)^{2-p}  \Vert e(u) \Vert^p_{L^2(N_Q)}.  $$ 
Taking \eqref{eq: main pro}\BBB (ii),\EEE (iii) into account we see that each $x \in Q_\mu$ is contained in a bounded number of different neighborhoods $N_Q=\bigcup_{\hat{Q} \in {\mathcal{N}}_{\rm g}(Q) \cup \lbrace Q \rbrace} \hat{Q}'$. Moreover, we note that for each $8 \le l \le I$ the sets $(N^l_k)_{k}$ are pairwise disjoint.  We observe that $M_l \BBB := \EEE \sum_k
\# \mathcal{X}^l_k  \le Cs^{-1}_l \mu$  by \eqref{eq: main proXXX}(iii) and $\mathcal{H}^1(J_u) \le 2\sqrt{2}\theta^{-2}\mu$ (see \eqref{eq: concentration}). Thus, using the discrete H\"older inequality  we derive 
\begin{align}\label{eq: main holder}
H& \le  C\sum_{ l \ge 8} M_l^{1-\frac{p}{2}} s_l^{2-p} \theta^{-(3p+2)lr} \big(\sum_{k}\Vert e(u)\Vert^2_{L^2(N^l_k)}\big)^{\frac{p}{2}} +  C\big(\sum_{Q \in {\mathcal{C}}} d(Q)^{2}\big)^{1-\frac{p}{2}}\BBB \big( \sum_{Q \in {\mathcal{C}}} \Vert e(u) \Vert^2_{L^2(N_Q)}\big)^{\frac{p}{2}}\EEE \notag \\
&\le C\sum\nolimits_{l \ge 8} (\mu s_l)^{1-\frac{p}{2}}  \theta^{- (3p+2)lr} \Vert e(u)\Vert^p_{L^2(Q_\mu)} +  C\mu^{2-p}  \Vert e(u) \Vert^p_{L^2(Q_\mu)},
\end{align}
where   we used that $\sum_{Q \in \mathcal{C}} d(Q)^2\le c|Q_\mu| \le  c\mu^2$. Recalling $s_l = \mu\theta^l$,  using $1-\frac{p}{2}  \BBB  = 12r \ge (4p+2)r \EEE$ as $r = \frac{1}{24}(2-p)$, and $\sum_{l \ge 8} \theta^{plr} \le C=C(\theta, p)$, we conclude 
\begin{align}\label{eq: main holderXXXX}
 H \le    C\mu^{2-p} \big(\sum\nolimits_{l \ge 8} \theta^{plr} + 1 \big) \Vert e(u) \Vert^p_{L^2(Q_\mu)} \le C\mu^{2-p}  \Vert e(u) \Vert^p_{L^2(Q_\mu)}.
 \end{align}
 
\smallskip

\textbf{Step III (Definition of the \BBB approximation): \EEE } \BBB We are now in the position to define the approximation $\bar{u}$ and the exceptional set $F$. Then we will use \eqref{eq: poincare estim} and \eqref{eq: main holderXXXX} to show \eqref{eq: modifica prop**}. At the end   we prove that the traces of $u$ and the approximation coincide on $\partial Q_\mu$. \EEE 

\smallskip

\emph{Definition of $\bar{u}$ and $F$:}  First, we choose a  partition of unity $(\varphi_Q)_{Q \in \mathcal{C}} \subset C^\infty(\R^2)$ with $\sum_{Q \in \mathcal{C}} \varphi_Q(x) = 1$ for $x \in \bigcup\nolimits_{Q \in \mathcal{C}} Q'$ and  
\begin{align}\label{eq: dingens}
\begin{split}
(i) & \ \ {\rm supp}(\varphi_Q) \subset Q' \text{ for all } Q \in \mathcal{C}, \\
(ii) & \ \ \Vert \nabla \varphi_Q \Vert_\infty \le c d(Q)^{-1} \text{ for all } Q \in \mathcal{C}.
\end{split}
\end{align}
As the proof of the existence of such a partition is very similar to the construction of a partition of unity for Whitney coverings (see  \cite{Federer:1969, Stein:1970}), we omit it here. Let us just briefly mention that the idea it to take a cut-off function $\bar{\varphi} \in C^\infty_c((0,1)^2)$ and to define the functions $\varphi_Q$ as suitably rescaled versions of $\bar{\varphi}$ taking the fact into account that each point is contained in only a bounded number of different squares $Q'$, $Q \in \mathcal{C}$. We now define the \BBB approximation \EEE  $\bar{u}: Q_\mu \to \R^2$ by 
\begin{align}\label{eq: function} 
\bar{u} = \sum\nolimits_{Q \in \mathcal{C}} \varphi_Q a_Q   = \sum\nolimits_{Q \in \mathcal{C}} \varphi_Q (A_Q \, \cdot +b_Q)
\end{align}
in $Q_\mu$. First, observe that $\bar{u}$ is smooth in $\bigcup\nolimits_{Q \in \mathcal{C}} Q'$ and that  $\bigcup\nolimits_{Q \in \mathcal{C}} Q' \supset Q_\mu \setminus (J_u \cap \bigcup_{j=1}^m \partial P'_j)$ by \eqref{eq: main pro}(i).  \BBB Recalling \eqref{eq: excpt} we \EEE let 
\begin{align}\label{eq: FFdef}
F :=   \bigcup\nolimits_{Q \in {\mathcal{C}}_{\rm g}} E_Q \cup Z.
\end{align} 
By  \eqref{eq: main pro}(iii), \eqref{eq: main proXXX1}, \eqref{eq: excpt} and $d(Q) \le  2\sqrt{2}\mu\theta^8$ for all $Q \in \mathcal{C} \subset \bigcup_{j\ge 8}\mathcal{Q}^j$ we derive $|F| \le c \mu\theta^5 \mathcal{H}^1(J_u)$ and  $\mathcal{H}^1(\partial F) \le C\mathcal{H}^1(J_u)$. This gives \eqref{eq: except***}. 

\smallskip

\emph{Proof of \eqref{eq: modifica prop**}:} We obtain $\nabla \bar{u} =\sum\nolimits_{\hat{Q} \in \mathcal{C}}\big( \varphi_{\hat{Q}} A_{\hat{Q}}     +  a_{\hat{Q}} \otimes\nabla \varphi_{\hat{Q}}\big)$ and using that $\nabla (\sum_{\hat{Q} \in \mathcal{C}} \varphi_{\hat{Q}})=0$ \BBB on $\bigcup\nolimits_{Q \in \mathcal{C}} Q'$ \EEE  we find for $x \in Q$, $Q \in\mathcal{C}$ 
\begin{align*}
\nabla \bar{u}(x) =  \sum\nolimits_{\hat{Q} \in \mathcal{C}}\big( \varphi_{\hat{Q}}(x) A_{\hat{Q}}   + ( a_{\hat{Q}}(x) - a_Q(x)) \otimes\nabla \varphi_{\hat{Q}}(x)\big).
\end{align*}
Then we  get by \eqref{eq: main pro}(ii),(iii), \eqref{eq: main holderXXXX}, \eqref{eq: dingens}, and the discrete H\"older inequality
$$
\Vert e(\bar{u}) \Vert^p_{L^p(Q_\mu)} \le C\sum\nolimits_{Q \in \mathcal{C} } \sum\nolimits_{\hat{Q} \in \mathcal{N}(Q)} d(Q)^{-p} \Vert a_Q - a_{\hat{Q}} \Vert^p_{L^p(Q')} \le C\mu^{2-p}   \Vert e(u) \Vert^p_{L^2(Q_\mu)},
$$
which implies  \eqref{eq: modifica prop**}(i). Similarly, we find using \eqref{eq: main pro}(ii),(iii), \eqref{eq: poincare estim}, \eqref{eq: main holderXXXX}, \eqref{eq: dingens}, \BBB \eqref{eq: FFdef}, \EEE and the fact that $\bigcup_{Q \in {\mathcal{C}}_{\rm g}} Q' \supset Q_\mu \setminus F$ 
\begin{align*}
\Vert \nabla \bar{u}  - \nabla u \Vert^p_{L^p(Q_\mu\setminus F)} &\le C\sum\nolimits_{Q \in {\mathcal{C}}_{\rm g}}\Vert \nabla u - A_Q \Vert^p_{L^p(Q'\setminus F)}  + CH\\
&\le  C\sum\nolimits_{Q \in {\mathcal{C}}_{\rm g}}d(Q)^{2-p}\Vert e(u) \Vert^p_{L^2(Q')}  +  CH   \le  C\mu^{2-p} \Vert e(u) \Vert^p_{L^2(Q_\mu)},
\end{align*}
where we repeated the H\"older-type estimate in \eqref{eq: main holder}. Let $S^k = \lbrace x \in Q_\mu: \dist(x,\partial Q_\mu) \le s_k\rbrace$ and observe $d(Q) \le 2\sqrt{2}s_k$ for all $Q \in \mathcal{C}$ with $Q \cap S^k  \neq  \emptyset$ by \eqref{eq: main pro}(i). Then we recall \eqref{eq: function} and get by \eqref{eq: main pro}(iii), \eqref{eq: poincare estim}, and again by  a H\"older-type estimate (cf. \eqref{eq: main holder})   
\begin{align}\label{eq: iurarg1}
\begin{split}
\Vert \bar{u}  - u \Vert^p_{L^p(S^k\setminus F)} &\le C\sum\nolimits_{l \ge k}\sum\nolimits_{Q \in {\mathcal{C}}_{\rm g} \cap \mathcal{C}^l}\Vert  u - a_Q \Vert^p_{L^p(Q'\setminus F)} 
\\&
\le  C s_k^p\sum\nolimits_{Q \in \mathcal{C}_{\rm g}} d(Q)^{2-p}  \Vert e(u) \Vert^p_{L^2(Q')} \le  Cs_k^p \mu^{2-p} \Vert e(u) \Vert^p_{L^2(Q_\mu)}.  
 \end{split}
\end{align}
In particular, for $k=0$ this implies \eqref{eq: modifica prop**}(iii).  

\smallskip

\emph{Traces of $u$ and $\bar{u}$:}  Finally, to show that $\bar{u} = u$ on $\partial Q_\mu$ one may argue, e.g., as in the proof of \cite[Theorem 2.1]{Conti-Iurlano:15.2}. We briefly sketch the argument for the reader's convenience. Choose $\psi_k \in C^\infty(Q_\mu)$ such that $\psi_k = 1$ in a neighborhood of $\partial Q_\mu$ and $\psi_k = 0$ on $Q_\mu \setminus S^k$ with $\Vert \nabla \psi_k \Vert_\infty \le c s_k^{-1}$. Note that \eqref{eq: main pro}(vi), \BBB \eqref{eq: excpt} \EEE and \eqref{eq: FFdef} imply $F \cap S_k = \emptyset$ for $k$ large enough. We define $v_k = \psi_k(u - \bar{u}) \in SBD(Q_\mu)$ and show $v_k \to 0$ strongly in $BD$ which implies $v_k|_{\partial Q_\mu} \to 0$ in $L^1$ and implies the assertion.  In fact, by \eqref{eq: iurarg1} and H\"older's inequality we get for $k$ sufficiently large
$$\Vert \bar{u}  - u \Vert_{L^1(S^k)} \le C (\mu s_k)^{1-\frac{1}{p}}\Vert \overline{u}  - u \Vert_{L^p(S^k)} \le C \mu^{\frac{1}{p}} s_k^{2-\frac{1}{p}} \Vert e(u) \Vert_{L^2(Q_\mu)}.$$
Now we derive that
$$|E v_k|(Q_\mu) \le |E u|(S^k) +  |E\bar{u}|(S^k) + cs_k^{-1} \Vert \bar{u}  - u \Vert_{L^1(S^k)} $$
vanishes for $k \to \infty$ and likewise $\Vert v_k \Vert_{L^1(Q_\mu)} \le \Vert u - \bar{u} \Vert_{L^1( S^k)} \to 0$. This implies $v_k \to 0$ in $BD$, as desired. \eop

 \begin{rem}\label{rem:modi}
 {\normalfont
 (i) For later reference, we recall that the exceptional set $F$ consists of the isolated components $Z$ and the exceptional sets $E_Q$ for $Q \in {\mathcal{C}}_{\rm g}$ (cf. \eqref{eq: FFdef}). In particular, by \eqref{eq: main pro}(ii),(iii) and \eqref{eq: excpt} we find (for $\theta$ small) that $|F \cap Q'|\le \theta|Q|$ for all $Q \in \mathcal{C}$ with $Q \cap Z = \emptyset$. 
 
 (ii) Moreover, in view of \eqref{eq: hatcal def}, each $Q \in \mathcal{C}$ with $Q \cap Z = \emptyset$ \BBB satisfies (see \eqref{eq: poincare estim})
 $$ \Vert \nabla u - A_Q \Vert^p_{L^{p}(Q' \setminus F)} \le c d(Q)^{2-p} \Vert e(u)\Vert^p_{L^2(Q')}. $$
 for some $A_Q \in \R^{2 \times 2}_{\rm skew}$. \EEE

 }
 \end{rem}
 
 We close this section with the observation that also $\Vert \bar{u} \Vert_\infty$ can be controlled, which is necessary for the proof of \eqref{eq: small set mainXXX}. The reader not interested in the derivation of \eqref{eq: small set mainXXX} (which will indeed not be needed in the sequel for the proof of Theorem \ref{th: main korn}) may readily skip the following lemma.

\begin{lemma}\label{lemma: modifica}
Let be given the situation of Theorem \ref{th: modifica}. Then the \BBB approximation \EEE  $\bar{u}: Q_\mu \to \R^2$  can be chosen such that \eqref{eq: except***}-\eqref{eq: modifica prop**} hold and  $\Vert \bar{u} \Vert_{L^\infty(Q_\mu \setminus F)} \le C\Vert u \Vert_\infty$ for $C=C(\theta,p)>0$.
\end{lemma}

\Proof The goal is to show that for each $Q \in {\mathcal{C}}_{\rm g}$ there is an infinitesimal rigid motion $\hat{a}_Q$ such that  
\begin{align}\label{eq: pp1}
\Vert  \hat{a}_Q \Vert_{L^\infty(Q')} \le C\Vert u \Vert_{\infty} 
\end{align}
and \eqref{eq: poincare estim} still holds for a possibly larger constant with $\hat{a}_Q$ in place of $a_Q$.  Then we can repeat the previous proof with $\hat{a}_Q$ in place of $a_Q$. Since for each $x \in Q_\mu \setminus F$ we have $Q \in {\mathcal{C}}_{\rm g}$ for all $Q \in \mathcal{C}$ with $x \in Q'$ (cf. \eqref{eq: hatcal def}), we obtain $|\bar{u}(x)| \le C\Vert u \Vert_\infty$ in view of  \eqref{eq: function}. Let us now show \eqref{eq: pp1}. Fix $Q \in {\mathcal{C}}_{\rm g}$. Using \eqref{eq: poincare estim} for the affine mapping $a_Q$ we find by \eqref{eq: excpt} and Lemma \ref{lemma: rigid}  
\begin{align*}
\Vert  a_Q \Vert^p_{L^{p}(Q')}  &\le C\Vert a_Q\Vert^p_{L^{p}(Q'\setminus E_Q)}  \le    Cd(Q)^{2} \Vert e(u)\Vert^p_{L^2(Q')} + C|Q'| \Vert u \Vert^p_{\infty}.
\end{align*}
Applying Lemma \ref{lemma: rigid motion} we find for all $x \in Q'$
\begin{align}\label{eq: poincare estim+}
 d(Q)|A_Q| + |A_Q \,x + b_Q|  \le C(\Vert e(u)\Vert_{L^2(Q')}+ \Vert u \Vert_{\infty}).
\end{align} 
If now $\Vert e(u)\Vert_{L^2(Q')} \le \Vert u \Vert_{\infty}$, then \eqref{eq: pp1}  follows directly with $\hat{a}_Q = a_Q$. Otherwise,  \eqref{eq: pp1} is clearly satisfied with $\hat{a}_Q := 0$ and in view of \eqref{eq: poincare estim+}, also \eqref{eq: poincare estim} holds for a larger constant since  
$$\Vert a_Q \Vert_{L^p(Q')} \le Cd(Q)^{\frac{2}{p}}\Vert e(u)\Vert_{L^2(Q')}, \ \ \ \Vert A_Q \Vert_{L^p(Q')} \le Cd(Q)^{\frac{2}{p} - 1}\Vert e(u)\Vert_{L^2(Q')}.$$

 \eop

\subsection{Partitions into John domains and proof of Theorem \ref{th: korn-small set}}\label{sec: 3sub5}

\BBB

In Section \ref{sec: 3sub1}  we have constructed an auxiliary partition $(P'_j)_{j=1}^m$ consisting of simply connected sets. Subsequently, in Section \ref{sec: 3sub2} we have seen that the configuration $u$ can be approximated, outside a small exceptional set, by an auxiliary function $\bar{u}$ which is smooth on each component of the partition $(P'_j)_{j=1}^m$.   We will  now apply Theorem \ref{th: main part2} to obtain a refined partition (the one appearing in Theorem \ref{th: korn-small set}) consisting of  John domains such that we can apply  Korn's inequality (see Theorem \ref{th: kornsobo}) with uniform constants. This allows to control the distance of the approximation $\bar{u}$ from an associated  infinitesimal rigid motion on each component (see Theorem \ref{th: partiparti}). Afterwards, using the properties of the approximation stated in \eqref{eq: modifica prop**}, we can give the proof of Theorem \ref{th: korn-small set}.  

In the following we essentially need the properties of  $(P'_j)_{j=1}^m$ and $\bar{u}$ given in Theorem \ref{th: bad part}  and Theorem \ref{th: modifica}, respectively. Details of other previously introduced  objects, such as the covering of Whitney-type or the isolated components, are not relevant in this section.  

\EEE

\begin{theorem}\label{th: partiparti}
Let $\mu,  \eta>0$, $\theta >0$ small and \UUU $p \in (1,2)$. \EEE There are a universal $c>0$ and a constant $C=C(\theta,p)>0$ such that for all $u \in \mathcal{W}(Q_\mu)$  with \eqref{eq: concentration} and for the corresponding partition $(P'_j)_{j=1}^m$ and \BBB approximation \EEE $\bar{u}$ as constructed in  Theorem \ref{th: bad part} and Theorem \ref{th: modifica},  respectively, the following holds: \\ There is a partition $(P_j)_{j=1}^n$ of $Q_\mu$ and a Borel set $R \subset Q_\mu$ such that each $P_j$ with $P_j \not\subset  R$ is a $c$-John domain with Lipschitz boundary. We have
\begin{align}\label{eq: partial W}
\begin{split}
(i)& \ \ \mathcal{H}^1\big(\bigcup\nolimits_{j=1}^n \partial P_j\big) \le  C\mathcal{H}^1(J_u),  \ \ \ \mathcal{H}^1\big(\bigcup\nolimits^m_{j=1} \partial P'_j \setminus \bigcup\nolimits^n_{j=1} \partial P_j \big) = 0, \\
(ii)&  \ \  |R| \le \eta,  \  \ \ \ \mathcal{H}^1(\partial R)  \le C\mathcal{H}^1(J_u),  
\end{split}
\end{align}
and there are infinitesimal rigid motions $a_j$, $j=1,\ldots,n$, such that  
\begin{align}\label{eq: toomany**}
\begin{split}
&\Vert \nabla \bar{v} \Vert_{L^{p}(Q_\mu \setminus R)}  \le  C\mu^{\frac{2}{p}-1}\Vert e(u) \Vert_{L^2(Q_\mu)},  \  \ \   \ \ \Vert  \bar{v}\Vert_{L^{p}(P_j \setminus R)} \le Cd(P_j)\Vert  \nabla \bar{v}\Vert_{L^{p}(P_j)},
\end{split}
\end{align}
where $\bar{v}:= \sum^n_{j=1} \chi_{P_j} (\bar{u}- a_j)$.  
\end{theorem}

\BBB Observe that $\theta$ appears in the statement as the partition $(P'_j)_{j=1}^m$ and the approximation $\bar{u}$ have been constructed in dependence of $\theta$. \EEE The occurrence of a \emph{rest set} $R$ is unavoidable as discussed below Theorem \ref{th: main part2}. \BBB Note, however, that its area can be made arbitrarily small.  Since in \eqref{eq: toomany**} we do not obtain any control on $R$,  this set will be a part of the exceptional set of Theorem \ref{th: korn-small set} and \eqref{eq: partial W}(ii) will be needed to show \eqref{eq: except}. \EEE

\smallskip

\Proof We apply Theorem \ref{th: bad part} \UUU (for $J = J_u$ and $r= \frac{1}{24}(2-p)$) \EEE and Theorem \ref{th: modifica}  to obtain a partition  $(P'_j)_{j=1}^m$ and an \BBB approximation \EEE $\bar{u}$ associated to $u$. Let us first note that Theorem \ref{th: main part2} cannot be applied directly since the sets are possibly not Lipschitz. However, we can introduce a refined partition  as follows. (Note that this is just a technical point and not the core of the argument.)

For fixed $i \in \N$ let $\mathcal{R}^i$ be the squares $ \mathcal{Q}^i$ whose closure have nonempty intersection with $\bigcup_{j=1}^m \partial P_j'$. Since $\bigcup_{j=1}^m \partial P'_j$ consists of finitely many closed segments, we clearly get  for $I \in \N$ large enough  that  $R_1 := \bigcup_{Q \in \mathcal{R}^I} \overline{Q''}$ satisfies $|R_1| \le \frac{\eta}{2}$,  $\mathcal{H}^1(\partial R_1) \le c\mathcal{H}^1(\bigcup_{j=1}^m \partial P'_j)$ and there is a partition $(P''_j)^M_{j=1}$ of $Q_\mu$ such that 
 \begin{align}\label{eq: newpart}
 \begin{split}
   (i) & \ \   \bigcup\nolimits^M_{j=1} \partial P''_j = \bigcup\nolimits^m_{j=1} \partial P'_j \cup \partial R_1, \\
  (ii) & \ \  \mathcal{H}^1\big(\bigcup\nolimits^M_{j=1} \partial P''_j \big) \le c \mathcal{H}^1\big(\bigcup\nolimits^m_{j=1} \partial P'_j \big)
 \end{split}
 \end{align} 
 for a universal $c>0$. We order the partition such that  $P''_j = P_j' \setminus R_1$ for $j = 1,\ldots,m$, and note that (again for $I \in \N$ large enough) the sets $(P''_j)_{j=1}^m$ are simply connected  domains with Lipschitz boundary, in particular   unions of squares. Moreover, we have that $R_1 \cup \bigcup_{j=1}^m P''_j$ forms a partition of $Q_\mu$.
 
Now we apply Theorem \ref{th: main part2} for $\eps =  \frac{\eta}{2m}$ on each set in $(P''_j)_{j=1}^m$ and find a partition $Q_\mu = R_1 \cup R_2 \cup P'''_1 \cup \ldots \cup P'''_N$ (up to a set of negligible measure) with $R_2 = \bigcup_{j=1}^m\Omega_0^j$ such that $|\Omega_0^j| \le   \frac{\eta}{2m}$, $\mathcal{H}^1(\partial \Omega_0^j) \le c\mathcal{H}^1(\partial P_j'')$ for $j=1,\ldots,m$ and  $P'''_1,\ldots,P'''_{N}$ are   $c$-John domains with Lipschitz boundary. Then $|R_2| \le \frac{\eta}{2}$. \BBB Applying  \eqref{eq:parti2} for each $P_j''$, $j=1,\ldots,m$, and summing over $j=1,\ldots,m$, we obtain 
\begin{align*}
\mathcal{H}^1\big(\bigcup\nolimits_{j=1}^m \partial \Omega_0^j \big) + \mathcal{H}^1\big(\bigcup\nolimits_{k=1}^N\partial P'''_k\big) &\le c\sum\nolimits_{j=1}^m\mathcal{H}^1\big( \partial P''_j \big).
\end{align*}
Since $\sum_{j=1}^m \mathcal{H}^1(\partial P_j'') \le 2 \mathcal{H}^1(\bigcup_{j=1}^m \partial P_j'')$, we get by \eqref{eq: bad length} and  \eqref{eq: newpart}(ii)  
\begin{align*}
\mathcal{H}^1\big(\bigcup\nolimits_{j=1}^m \partial \Omega_0^j \big) + \mathcal{H}^1\big(\bigcup\nolimits_{k=1}^N\partial P'''_k\big) &\le c\mathcal{H}^1\big(\bigcup\nolimits_{j=1}^M \partial P''_j \big)  \le c\mathcal{H}^1\big(\bigcup\nolimits_{j=1}^m \partial P'_j \big)  \le C\mathcal{H}^1(J_u)
\end{align*}
\UUU for some $C=C(\theta,r) = C(\theta,p)>0$. \EEE By $(P_j)_{j=1}^n$ we denote the partition consisting of $(P'''_k)_{k=1}^N$, $(\Omega_0^j)_{j=1}^m$, and $(P''_j)_{j=m+1}^M$. Then in view of \eqref{eq: newpart} and the previous estimate, \eqref{eq: partial W}(i) follows. We let  $R  = R_1 \cup R_2$ and observe $|R| \le \eta$ as well as $\mathcal{H}^1(\partial R) \le C\mathcal{H}^1(J_u)$, i.e., \eqref{eq: partial W}(ii) holds.  Recall that $\bar{u}$ is smooth on each $P'''_k$, $k=1,\ldots,N$, by Theorem \ref{th: modifica}. Applying Theorem \ref{th: kornsobo}   on each $P'''_k$ we obtain  infinitesimal rigid motions $(a_k)_{k=1}^N$  such that by  H\"older's inequality and the fact that $Q_\mu \setminus R = \bigcup_{k=1}^N P_k'''$
\begin{align}\label{eq: before}
  \Vert \nabla \bar{v}\Vert_{L^{p}(Q_\mu \setminus R)} \le C \Vert e(\bar{u}) \Vert_{L^{p}(Q_\mu)} \le C \mu^{\frac{2}{p}-1}\Vert e(u) \Vert_{L^2(Q_\mu)},
\end{align}
where $\bar{v}:= \bar{u} - \sum^N_{k=1} \chi_{P'''_k} a_k$.  By Poincar\'e's inequality  and a scaling argument we also have  $\Vert  \bar{v}\Vert_{L^{p}(P'''_k)} \le Cd(P'''_k)\Vert  \nabla \bar{v}\Vert_{L^{p}(P'''_k)}$ for $k=1,\ldots, N$. With $a_j =0$ for all other components we get that $\bar{v}$ can be written as $\bar{v} = \sum^n_{j=1} \chi_{P_j} (\bar{u}- a_j)$. \eop

We are now in a position to give the proof of Theorem  \ref{th: korn-small set}. The reader not interested in the derivation of \eqref{eq: small set mainXXX}  may readily skip the second part of the proof. 

\smallskip

\noindent {\em Proof of Theorem \ref{th: korn-small set}.} (1) \UUU We can suppose that $p \in (1,2)$. Once the result is proved in this case, the case $p=1$  follows readily by H\"older's inequality. \EEE Moreover, we may without restriction assume that  $\mathcal{H}^1(J_u) \le 2\sqrt{2}\theta^{-2}\mu$ as otherwise the claim is trivially satisfied (cf. \eqref{eq: concentration}).  \BBB If also  $\mathcal{H}^1(J_{u}) \ge \mu \theta^2$, \EEE we let  $(P'_j)_{j=1}^m$   be the partition  constructed in Theorem \ref{th: bad part} \UUU (for $J = J_u$ and $r= \frac{1}{24}(2-p)$) \EEE and let $\bar{u}$ be the \BBB approximation \EEE given by    Theorem \ref{th: modifica}. We distinguish \BBB three \EEE cases:

(a) \BBB Suppose first that $\mathcal{H}^1(J_{u}) \ge \mu \theta^2$ and  $\mathcal{H}^1(J_{u} \cap \bigcup_{j=1}^m\partial P'_j)  \ge \mu \theta^2$. \EEE  We apply Theorem \ref{th: partiparti} to obtain the partition $Q_\mu =   \bigcup_{j=1}^n P_j$ such that by $\mathcal{H}^1(J_u) \le 2\sqrt{2}\theta^{-2}\mu$ and  \eqref{eq: partial W}(i)  
\begin{align}\label{eq:again-new}
\begin{split}
\mathcal{H}^1\big(\bigcup\nolimits_{j=1}^n\partial P_j\big) & \le C\mathcal{H}^1(J_u) \le C\mu \le C \mathcal{H}^1(J_{u} \cap \bigcup\nolimits_{j=1}^m\partial P'_j)  \le C\mathcal{H}^1(J_{u} \cap \bigcup\nolimits_{j=1}^n\partial P_j), 
\end{split}
\end{align}
where in the last step we used $\bigcup_{j=1}^m \partial P'_j \subset \bigcup_{j=1}^n\partial P_j$ up to a set of negligible $\mathcal{H}^1$-measure. This gives \eqref{eq: small set main}(i). Let $F$ be the exceptional set derived in Theorem \ref{th: modifica} (see \eqref{eq: except***}), let $R$ be the set in \eqref{eq: partial W}(ii) and define 
\begin{align}\label{eq:Efed}
E=F \cup R.
\end{align}
Then   \eqref{eq: except} follows directly from \eqref{eq: except***} and \eqref{eq: partial W}(ii) since $\eta$ can be chosen arbitrarily small. To see \eqref{eq: small set main}(ii), we apply \eqref{eq: modifica prop**}(ii) and \eqref{eq: toomany**} to find
\begin{align}\label{eq:doch}
\begin{split}
\Vert \nabla v \Vert_{L^p(Q_\mu \setminus E)}& \le \Vert \nabla \bar{v} \Vert_{L^p(Q_\mu \setminus R)} + \Vert \nabla u - \nabla \bar{u} \Vert_{L^p(Q_\mu \setminus F)}  \le C\mu^{\frac{2}{p}-1}\Vert e(u) \Vert_{L^2(Q_\mu)}, 
\end{split}
\end{align}
where $v = u - \sum_{j=1}^n \chi_{P_j} a_j$ and $\bar{v}$ as in \eqref{eq: toomany**}  for infinitesimal rigid motions $(a_j)_{j=1}^n$.

\UUU (b)(i) \EEE  We now suppose that  \BBB $\mathcal{H}^1(J_{u}) \ge \mu \theta^2$ and  $ \BBB \mathcal{H}^1(J_{u} \cap \bigcup_{j=1}^m\partial P'_j)  < \mu \theta^2$.   As   $\bar{u}$ is smooth on the complement of $J_{u} \cap \bigcup_{j=1}^m\partial P'_j$, \EEE this implies $\mathcal{H}^1(J_{\bar{u}}) < \mu\theta^2$. Consequently, we can  apply   Theorem \ref{th: kornSBDsmall} on $\bar{u}$ and define  $E = \tilde{E} \cup F$ with $\tilde{E}$ as in \eqref{eq: R2main} and  $F$ as in \eqref{eq: except***}. Observe that also in this case \eqref{eq: except} holds \BBB since particularly $|\tilde{E}| \le c (\mathcal{H}^1(J_{\bar{u}}))^2 \le c\mu\theta^2 \mathcal{H}^1(J_{{u}})$. Moreover, \EEE \eqref{eq: small set main}(i) is trivially satisfied for the partition consisting only of $Q_\mu$. Property \eqref{eq: small set main}(ii) follows from \eqref{eq: main estmain} and \eqref{eq: modifica prop**}(ii).

\UUU (b)(ii) \EEE If $\mathcal{H}^1(J_{u}) < \mu \theta^2$, we proceed as in (b)(i), but apply Theorem \ref{th: kornSBDsmall} directly  on $u$. \EEE

\smallskip
(2) We now show that we can find an exceptional set $E'$ and a configuration $v' = u - \sum_{j=1}^n a_j' \chi_{P_j}$ such that \eqref{eq: except}-\eqref{eq: small set mainXXX} hold. We only treat case (a) where $\mathcal{H}^1(J_{u} \cap \bigcup_{j=1}^m\partial P'_j)  \ge \mu \theta^2$ since \UUU (b)(i),(ii) \EEE are similar. With $E$ as in part (1) we define  $\mathcal{P} = \lbrace P_j: |P_j \setminus E| \ge \frac{1}{2} |P_j| \rbrace$ and let $E' = E \cup \bigcup_{P \notin \mathcal{P}} P$. Observe that $\mathcal{H}^1(\partial E') \le C\mathcal{H}^1(J_u)$ by \eqref{eq: except}, \eqref{eq: partial W}(i) and $|E'| \le 2|E|$. This yields \eqref{eq: except} for $E'$. 

\BBB In view of \eqref{eq:Efed}, we find \EEE $\mathcal{P} \subset   (P'''_k)_{k=1}^N$ for the components considered before \eqref{eq: before}. Since $P'''_k$ is a $c$-John domain, we get $|P'''_k| \ge  c'd(P'''_k)^2$ for $c'=c'(c)$, \UUU see Lemma \ref{lemma: plump}. \EEE Therefore,   $|P_j \setminus E| \ge \frac{1}{2}c'(d(P_j))^2$ for all $P_j \in \mathcal{P}$ and by Lemma \ref{lemma: rigid} we get $\Vert  a_j\Vert_{L^{p}(Q)} \le C\Vert  a_j\Vert_{L^{p}(P_j \setminus E)}$, where $Q$ is a square containing $P_j$ with $|P_j \setminus E| \ge c'|Q|$ for a possibly smaller $c'$. By  Lemma \ref{lemma: rigid motion}, \eqref{eq: toomany**}, and the fact that $\Vert \bar{u} \Vert_{L^\infty(Q_\mu \setminus F)} \le C\Vert u \Vert_\infty$ (see Lemma \ref{lemma: modifica}) we then derive
\begin{align*}
\begin{split}
d(P_j)|A_j| + \Vert a_j \Vert_{L^\infty(P_j)}   &\le Cd(P_j)|P_j|^{-\frac{1}{2}-\frac{1}{p}}\Vert  a_j\Vert_{L^{p}(P_j)} \le C|P_j|^{-\frac{1}{p}}\Vert  a_j\Vert_{L^{p}(P_j)} \\  
& \le C|P_j|^{-\frac{1}{p}}\Vert  a_j\Vert_{L^{p}(P_j \setminus E)} \le C|P_j|^{-\frac{1}{p}}\Vert \bar{u}- a_j\Vert_{L^{p}(P_j \setminus E)}   + C\Vert \bar{u} \Vert_{L^\infty(Q_\mu \setminus F)}  \\&\le C(d(P_j))^{ 1 -\frac{2}{p}}\Vert  \nabla \bar{v}\Vert_{L^{p}(P_j)} + C\Vert u \Vert_\infty. 
\end{split}
\end{align*}
 If  $(d(P_j))^{ 1 -\frac{2}{p}}\Vert  \nabla \bar{v}\Vert_{L^{p}(P_j)} \le \Vert u \Vert_\infty$, we indeed obtain $\Vert v \Vert_{L^\infty(P_j)} \le C\Vert u \Vert_\infty$ and set $v' = v$ on $P_j$, i.e., $a_j' = a_j$. Otherwise, we set $v' = u$ on $P_j$, i.e., $a_j' = 0$, and derive by a short calculation using the previous estimate
$$  \Vert \nabla v'\Vert^p_{L^{p}(P_j \setminus E)} \le C\Vert \nabla v\Vert^p_{L^{p}(P_j \setminus E)} + C|P_j||A_j|^p  \le C\Vert \nabla v\Vert^p_{L^{p}(P_j \setminus E)} + C\Vert \nabla \bar{v}\Vert^p_{L^{p}(P_j)}.$$
\BBB For $P_j \notin \mathcal{P}$, let $a_j' = 0$. \EEE Thus, summing over all components  $P_j \in \mathcal{P}$ and using \eqref{eq: toomany**}, \eqref{eq:doch} we see that also $v'$ satisfies \eqref{eq: small set main}(ii) and additionally $\Vert v'\Vert_\infty \le C\Vert u \Vert_\infty$, which gives \eqref{eq: small set mainXXX}.  \eop

 \UUU
 
 \begin{rem}\label{rem: different}
{\normalfont

Let us comment on the distinction of the  cases in the proof: (b) addresses the case where the jump set of the approximation $\bar{u}$ (or already of $u$ itself, see (b)(ii)) is small in a certain sense and (a) deals with the case of larger jump sets.

Whereas in case (a) we establish a piecewise inequality by means of Theorem \ref{th: partiparti}, the particularity of case (b)  is that we use Theorem \ref{th: kornSBDsmall} to provide a Korn-type estimate for a \emph{single} infinitesimal rigid motion. This is the case even if small pieces of the domain might be detached from the bulk part by $J_u$ and the partition provided by Theorem \ref{th: partiparti} (in case (b)(i)) might be nontrivial. Intuitively, by the application of the Korn inequality for functions with small jump set (Theorem \ref{th: kornSBDsmall}), these small detached pieces of the domain form part of the exceptional set given in \eqref {eq: except} and \eqref{eq: small set main}(ii) does not provide any information in theses pieces. (We refer to \cite{Chambolle-Conti-Francfort:2014, Conti-Iurlano:15.2, Friedrich:15-3} for details.)

It may appear more natural to apply Theorem \ref{th: partiparti} also in case (b)(i) as by means of a piecewise inequality involving different infinitesimal rigid motions also the behavior in such small detached pieces could be controlled suitably. The distinction performed above, however, accounts not only for the Korn-type estimate \eqref{eq: small set main}(ii) but also for the fact that we need  a bound on the partition \eqref{eq: small set main}(i). In view of \eqref{eq:again-new}, it appears to be impossible to derive such a bound in case (b)(i), where $\mathcal{H}^1(J_{u} \cap \bigcup_{j=1}^m\partial P'_j)$ (and thus $\mathcal{H}^1(J_{\bar{u}})$) is small.  Therefore, in case (b)(i) it is more convenient to derive an estimate for a trivial partition consisting only of one component as then  \eqref{eq: small set main}(i) follows immediately. 

The problem is even more apparent in case (b)(ii): if already the jump set of the original function is small, Theorem \ref{th: bad part}, Theorem \ref{th: modifica} and Theorem \ref{th: partiparti} are not applicable. Indeed, in this situation it would be difficult to construct a partition whose boundary  is controlled suitably (cf. Remark \ref{rem: Z}(i)).

}
\end{rem}

\EEE

 \BBB
 
 \begin{rem}\label{rem: excep}
 {\normalfont
 For later purpose, we recall the structure of the exceptional set $E$ in case (a): by Remark \ref{rem:modi} and \eqref{eq:Efed} we have that $E = F \cup R$ consists of the isolated components $Z$ introduced in Lemma \ref{lemma: bad sets}, the exceptional sets related to the application of Theorem \ref{th: kornSBDsmall} on good squares, and the rest set $R$ from Theorem \ref{th: partiparti}. Note that the latter is less delicate since the area of $R$ can be chosen arbitrarily small.
 }
 \end{rem}
 \EEE

\BBB

\subsection{A refined version of Theorem \ref{th: korn-small set}}\label{sec: refined}

In the previous section we have concluded the proof of Theorem \ref{th: korn-small set}. We have seen that for $u \in \mathcal{W}(Q_\mu)$ we can find a partition $(P_j)_j$ of $Q_\mu$  such that on each component an estimate of Korn-type \eqref{eq: small set main}(ii) holds outside a small exceptional set $E$. The derivation of the main result without exceptional set  will be based on an iterative application of  Theorem \ref{th: korn-small set}   on various mesoscopic scales. More precisely, one step of the algorithmic procedure will consist in considering the covering of Whitney-type $\mathcal{C}$ constructed in Section \ref{sec: 3sub1} (see Theorem \ref{th: bad part}) and  using Theorem \ref{th: korn-small set} on every square of the covering. (For details on the algorithm we refer  to the beginning of Section \ref{sec: main*}.) 

In order to derive the main estimate \eqref{eq: main korn2}(ii), it is fundamental to investigate the relation of infinitesimal rigid motions obtained by the application of Theorem \ref{th: korn-small set} in different iteration steps. Therefore, as mentioned at the beginning of Section \ref{sec: small set}, we will need a refined version of Theorem \ref{th: korn-small set} taking the following aspects into account: 

\smallskip

(A) It is essential that the exceptional set $E$ of Theorem \ref{th: korn-small set} covers only a small portion of each square of $\mathcal{C}$. To this end, we have to investigate the \emph{structure} of the exceptional set. 

\smallskip

(B) It is essential that the estimate of Korn-type  \eqref{eq: small set main}(ii) does not only hold on each component $P_j$, but also on each square of the covering $\mathcal{C}$ intersecting $P_j$, i.e., on a set being in general larger than $P_j$. To this end, we possibly need to redefine the partition. 

\smallskip

Before we come to the details, let us recall the relevant notions for this section.  Besides the statement of Theorem \ref{th: korn-small set}, we need some properties of the Whitney covering $\mathcal{C}$ stated in Theorem \ref{th: bad part}. Moreover, recall the structure of the exceptional set $E = F \cup R$ described in Remark \ref{rem: excep}. We also recall that if the partition $(P_j)_j$ in Theorem \ref{th: korn-small set} is not trivial (case (a) in the proof), it coincides with the one from  Theorem \ref{th: partiparti} and we will use that the components  with $P_j \not\subset R$ are $c$-John domains. Finally, at some point we need the approximation constructed in Theorem \ref{th: modifica}. We now address (A) and (B):

\smallskip

(A) By Remark \ref{rem:modi}(i) we get that $F$ covers only a small portion of each square outside the isolated components $Z$. Although $|R|$ is small, its portion inside a square of the covering is not necessarily controllable. To this end, we introduce  \EEE
\begin{align}\label{eq.prep}
R'= R \cup \bigcup\nolimits_{Q \in \mathcal{C}: |Q' \cap R| \ge \theta|Q|} Q
\end{align}
and note that $|R'| \le \eta+ 12\theta^{-1}\eta$ by  \eqref{eq: partial W} and \eqref{eq: main pro}(iii). \BBB Although (A) does not hold on the squares contained in $R'$, this will not affect our analysis as $|R'|$ can be chosen arbitrarily small. (Indeed, in Section \ref{sec: main-proof}  we will see that Theorem \ref{th: main korn} already holds if it holds up to an exceptional set of arbitrarily small size.) Finally, we observe that isolated components always belong to the exceptional set $E$ and (A) is not achievable for squares contained in $Z$. Here, other ideas will have to be applied using the properties in Lemma \ref{lemma: bad sets}  and we refer to the description of the iteration scheme in Section \ref{sec: main*} for more details. \EEE

\smallskip

\BBB
(B) For the generalization of the Korn-type estimate it is essential that all squares of $\mathcal{C}$ intersecting $P_j$ are not `too large' with respect to the diameter of $P_j$. To this end, we have to redefine the partition consisting of \emph{good components} (having exactly the desired property), \emph{small components}, and \emph{rest components}, which are contained in the set $R'$ defined in \eqref{eq.prep}.  \EEE

 \begin{lemma}\label{cor: partiparti}
 Let $(P_j)_{j=1}^n$ be the partition from Theorem \ref{th: partiparti} and \BBB $\mathcal{C}$ be the covering from Theorem \ref{th: bad part}, Lemma \ref{lemma: bad sets}. Let  $E$ as in Remark \ref{rem: excep} and \EEE $R'$ as in \eqref{eq.prep}. Then there is another partition $(P^*_j)^N_{j=1} = \mathcal{P}_R \cup \mathcal{P}_{\rm good} \cup \mathcal{P}_{\rm small}$ of $Q_{\mu}$ and $(A_j)_{j=1}^N \subset \R^{2 \times 2}_{\rm skew}$ such that $\mathcal{P}_{\rm good}$ consists of $c$-John domains, $\mathcal{H}^1\big(\bigcup\nolimits_{j=1}^N \partial P^*_j\big) \le  C\mathcal{H}^1(J_u)$,  and  with $P^*_{j, {\mathcal{C}}} := \bigcup_{Q \in \mathcal{C}, Q \cap P_j^* \neq \emptyset} Q$ we get 
\begin{align}\label{eq:lastlast}
\begin{split}
(i) & \ \ \BBB \sum\nolimits^N_{j=1} \Vert \nabla u - A_j \Vert^p_{L^p(P^*_j \setminus E)} \EEE \le C\mu^{2-p} \Vert e(u) \Vert_{L^2(Q_\mu)}^p,\\
(ii)& \ \  P^*_j \subset R' \text{ for all } P_j^* \in \mathcal{P}_R, \\
(iii) & \ \  d(Q) \le 4\sqrt{2}d(P^*_j)  \text{ for all } P_j^* \in \mathcal{P}_{\rm good}  \text{ and for all }  Q \in \mathcal{C} \text{ with } Q \cap P_j^* \neq \emptyset , \\
(iv)& \ \  \text{for all } P_k^* \in \mathcal{P}_{\rm small}  \text{ there is } P^*_j \in  \mathcal{P}_{\rm good} \text{ with }  P_{k,\mathcal{C}}^* \subset P^*_{j,\mathcal{C}}. 
\end{split}
\end{align}
 \end{lemma}
 
 \BBB
Below in the proof of Lemma \ref{rem: large components} we will exploit \eqref{eq:lastlast}(iii) to obtain refinement (B) and each small component will be combined with a good component by using \eqref{eq:lastlast}(iv). \EEE

 \noindent {\em Proof of Lemma \ref{cor: partiparti}.} With $(P_j)_{j=1}^n$ from Theorem \ref{th: partiparti} and $\mathcal{C}$ from   Theorem \ref{th: bad part}   we let  
\begin{align}\label{eq:rec}
\mathcal{P}'_{\rm good}  = \lbrace P_j: d(Q) \le 4\sqrt{2}d(P_j) \ \text{ for all } \  Q \in \mathcal{C}: Q \cap P_j \neq \emptyset \rbrace.
\end{align}
By  $\mathcal{P}'_{\rm small}$ we denote the remaining components. Recalling \eqref{eq: enlarged squares} we find for each $P_j \in \mathcal{P}'_{\rm small}$ some $Q \in \mathcal{C}$ such that $P_j \subset Q'$.  Then \eqref{eq: main pro}(ii) shows that  $P_j \in \mathcal{P}'_{\rm small}$ intersects at most $c\theta^{-2}$ different squares of $\mathcal{C}$. Consequently, indicating by $\mathcal{P}''_{\rm small}$ the nonempty sets in $\lbrace Q \cap P_j: P_j \in \mathcal{P}'_{\rm small}, \ Q \in \mathcal{C} \rbrace$ we derive for $c>0$ large enough  
\begin{align}\label{eq:lll2}
\sum\nolimits_{P \in \mathcal{P}''_{\rm small}} \mathcal{H}^1(\partial P) \le c\theta^{-2}\sum\nolimits_{P_j \in \mathcal{P}'_{\rm small}} \mathcal{H}^1(\partial P_j).
\end{align}
Note that $\mathcal{P}'_{\rm good} \cup \mathcal{P}''_{\rm small}$ forms a partition of $Q_\mu$. We let $\mathcal{P}_{R} = \lbrace P \in \mathcal{P}'_{\rm good} \cup \mathcal{P}''_{\rm small}: P \subset R'\rbrace$ and let $\mathcal{P}''_{\rm good} = \mathcal{P}'_{\rm good} \setminus \mathcal{P}_{R}$, $\mathcal{P}'''_{\rm small}  = \mathcal{P}''_{\rm small}  \setminus \mathcal{P}_{R}$.

 By $\mathcal{C}_{\rm small}$ we denote the squares $Q \in    \mathcal{C}$ with $Q \not\subset R'$ and $Q \cap \bigcup_{P_j \in \mathcal{P}''_{\rm good}} P_{j} =  \emptyset$. As then $|Q \cap \bigcup_{P \in \mathcal{P}'''_{\rm small}}P| \ge |Q \setminus R'| \ge (1-\theta)|Q| = \frac{1-\theta}{16}(\mathcal{H}^1(\partial Q))^2$ for $Q \in \mathcal{C}_{\rm small}$  by \eqref{eq.prep}, the isoperimetric inequality and \eqref{eq:lll2} yield for $C=C(\theta)>0$
\begin{align}\label{eq:lll}
\sum_{Q \in \mathcal{C}_{\rm small}}\mathcal{H}^1(\partial Q) \le C\sum_{Q \in \mathcal{C}_{\rm small}}\sum_{P \in \mathcal{P}'''_{\rm small}} \mathcal{H}^1(\partial P \cap \overline{Q}) \le C\sum_{P_j \in \mathcal{P}'_{\rm small}} \mathcal{H}^1(\partial P_j).
\end{align}
We define $\mathcal{P}_{\rm good} = \mathcal{P}''_{\rm good} \cup (Q)_{Q \in \mathcal{C}_{\rm small}}$ and let $\mathcal{P}_{\rm small} \subset \mathcal{P}'''_{\rm small}$ be the components not contained in some $Q, Q \in \mathcal{C}_{\rm small}$. We denote the sets of the partition $\mathcal{P}_R \cup \mathcal{P}_{\rm good} \cup \mathcal{P}_{\rm small}$ by $(P^*_j)_{j=1}^N$. Note that by Theorem \ref{th: partiparti} and the definition of $\mathcal{P}_R$, the family $\mathcal{P}_{\rm good}$ consists of $c$-John domains. The assertion $\mathcal{H}^1\big(\bigcup\nolimits_{j=1}^N \partial P^*_j\big) \le  C\mathcal{H}^1(J_u)$ follows  from  \eqref{eq: partial W}(i) and \eqref{eq:lll2}-\eqref{eq:lll}. 

\BBB We show \eqref{eq:lastlast}. First, \EEE \eqref{eq:lastlast}(ii) follows directly from the definition of $\mathcal{P}_{R}$, \eqref{eq:lastlast}(iii) is obvious for  $(Q)_{Q \in \mathcal{C}_{\rm small}}$ and otherwise we recall \eqref{eq:rec}. Likewise, \eqref{eq:lastlast}(iv) follows from the construction of $\mathcal{P}_{\rm small}$, in particular the definition of $\mathcal{C}_{\rm small}$. Indeed,  we recall that each set in $\mathcal{P}_{\rm small}$ is contained in some $Q \in \mathcal{C} \setminus \mathcal{C}_{\rm small}$ with $Q \not\subset R'$ and \BBB thus $Q \cap \bigcup_{P_j \in \mathcal{P}''_{\rm good}} P_{j} \neq \emptyset$. \EEE

\BBB Finally, we show \eqref{eq:lastlast}(i). Recall the approximation $\bar{u}$ in Theorem \ref{th: modifica}. \EEE By \eqref{eq: modifica prop**}(i), the fact that $\bar{u}$ is smooth in $Q, Q\in \mathcal{C}_{\rm small}$, and Korn's inequality there are matrices $A_Q \in \R^{2 \times 2}_{\rm skew}$ such that
$$\sum\nolimits_{Q \in \mathcal{C}_{\rm small}} \Vert \nabla \bar{u} - A_Q \Vert^p_{L^p(Q)} \le C\mu^{2-p} \Vert e(u) \Vert^p_{L^2(Q_\mu)}. $$
\BBB Then by \eqref{eq: modifica prop**}(ii) we also get  (see  \eqref{eq:doch} for a similar argument)
$$\sum\nolimits_{Q \in \mathcal{C}_{\rm small}} \Vert \nabla u - A_Q \Vert^p_{L^p(Q \setminus E)} \le C\mu^{2-p} \Vert e(u) \Vert^p_{L^2(Q_\mu)}$$
with the exceptional set $E$ from Remark \ref{rem: excep}.  For the  components  $\mathcal{P}_R \cup \mathcal{P}''_{\rm good} \cup \mathcal{P}_{\rm small}$, each of which being contained in a component of the original partition $(P_j)_{j=1}^n$, we apply  \eqref{eq: small set main}(ii) \EEE and then get that there are matrices $(A_j)_{j=1}^N$ such that \eqref{eq:lastlast}(i) holds.  \eop

Given a covering $\mathcal{C}$ as in Theorem \ref{th: bad part}, a partition $(P_j)_j$ as in Theorem \ref{th: partiparti} (or Lemma \ref{cor: partiparti}) and    an exceptional set $R'$ as in  \eqref{eq.prep} we define
\begin{align}\label{eq: specialdef}
\begin{split}
&\mathcal{Q}(P_j; \mathcal{C}, R') = \lbrace Q \in \mathcal{C}: Q\cap P_j \neq \emptyset, \ Q \cap P_j\not\subset R'  \rbrace
\end{split}
\end{align}
and the \BBB \emph{covering-adapted sets} \EEE $P_{j,{\rm cov}} =  \bigcup\nolimits_{Q   \in \mathcal{Q}(P_j;  \mathcal{C}, R')} Q$, which will play a crucial role  in Section \ref{sec: main-proof} \BBB (see refinement (B) above). \EEE Note that the definition of  $P_{j,{\rm cov}}$ always depends on $\mathcal{C}, R'$ although not made explicit in the notation.

\BBB
We now   formulate the refined version of Theorem \ref{th: korn-small set}. Note that in Section  \ref{sec: main-proof} we will exclusively use Lemma \ref{rem: large components} and no other details from Section \ref{sec: small set} will be needed.  
\EEE

\begin{lemma}\label{rem: large components}
Let $\mu, \eta>0$, $\theta>0$ small,    $p  \in [1,2)$ and $r= \frac{1}{24}(2-p)$. Then there are a universal constant $c>0$ and $C=C(\theta,p)>0$   such that  the following holds: 

(1) For each $u \in \mathcal{W}(Q_\mu)$ with $\mathcal{H}^1(J_u) \le 2\sqrt{2}\theta^{-2}\mu$   there is a partition $Q_\mu =  \bigcup_{j=0}^n P_j$  
with 
\begin{align}\label{eq:either} 
P_0 =  Q_\mu \ \ \ \text{or} \ \ \ P_0 = \emptyset,
\end{align} 
a covering $\mathcal{C}_u \subset \bigcup_{i\ge 1} \mathcal{Q}^i$ of $Q_\mu$ \BBB consisting of pairwise disjoint dyadic squares, \EEE and \BBB pairwise disjoint, closed \EEE sets $Z^l_u$, $l \ge 8$,  such that   with $S_u := (J_{u} \cap Q_\mu)  \setminus \bigcup_{Q \in \mathcal{C}_u} Q'$ we have
\begin{align}\label{eq: sharpi2}
\begin{split}
(i)& \ \ \mathcal{H}^1\big(\bigcup\nolimits_{j=1}^n\partial P_j \big)  \le C \mathcal{H}^1(S_u),\\
(ii) & \ \ Q \subset Q_\mu \setminus S_u    \ \text{for all}   \ Q \in \mathcal{C}_u, \ \  \overline{Q} \subset Q_\mu \setminus S_u    \ \text{for all}   \ Q \subset \bigcup\nolimits_{l\ge 8} Z^l_u,\\
(iii)& \ \ \partial X^l_k \BBB = X^l_k \cap \bigcup\nolimits_{j=1}^n \partial P_j  \EEE \subset \bigcup\nolimits_{j=1}^n \partial P_j \ \ \text{ for all $X^{l}_k$,}
\end{split}
\end{align}
where  $(X^{l}_k)_k$ denote the connected components of $Z^l_u$. Each $Z^l_u$ is \BBB a \EEE union of squares in $\mathcal{Q}^l \cap \mathcal{C}_u$ up to a set of measure zero.  Moreover, there are exceptional sets $E_u, R_u$ such that for  $l \ge 8$ one has
\begin{align}\label{eq: sharpi3}
\begin{split}
(i) & \ \ |E_u| \le c\mu \theta^2 \mathcal{H}^1(J_u) \le cd(Q_{\mu}) \theta^2 \mathcal{H}^1(J_u \cap Q_\mu), \ \ \ \  \  \ \ |R_u| \le \eta, \\
(ii) & \ \ |Q' \cap E_u | \le c\theta |Q|  \text{ for all } Q \in \mathcal{C}_u \ \text{with} \  \ Q\not\subset R_u,  \ \ Q \cap   \bigcup\nolimits_{l \ge 8} Z^l_u =\emptyset,\\
(iii)& \ \  d(X^{l}_k) \le \theta^{-rl} s_l \ \ \text{ for all $X^{l}_k$,}\\
(iv)& \ \ |Z^l_u | \le C \theta^{-lr} s_l  \mathcal{H}^1(J_u) \le C \theta^{-rl} s_l \mu,\\
\end{split}
\end{align}
and there are infinitesimal rigid motions $a_j = a_{A_j,b_j}$, $j=0,\ldots,n$, such that
\begin{align}\label{eq: sharpi2.b}
\begin{split}
(i)& \ \ \sum\nolimits_{j=0}^n\Vert \nabla u - A_j\Vert^p_{L^{p}(P_{j,{\rm cov}} \setminus E_u)}  \le C\mu^{ 2-p} \Vert e(u) \Vert^p_{L^2(Q_\mu)}, \\
(ii) & \ \ \#\lbrace P_{j,{\rm cov}}: x \in P_{j,{\rm cov}}  \rbrace \le c \ \ \text{ for all } \ \ x \in Q_\mu,
\end{split}
\end{align}
where $P_{j,{\rm cov}}$ as defined below \eqref{eq: specialdef} with respect to $\mathcal{C}_u,R_u$. 

(2) If $P_0 = Q_\mu$, then $P_{0,{\rm cov}} = Q_\mu$, \UUU $\mathcal{C}_u = \lbrace Q \in \mathcal{Q}^1: Q \subset Q_\mu \rbrace$, \EEE $S_u= \emptyset$, $\bigcup_{l \ge 8} Z_u^l = \emptyset$, $R_u = \emptyset$ and  there is a set $\Gamma_u \subset \partial Q_\mu$ with $\mathcal{H}^1(\Gamma_u) \le c\theta\mu$ such that
\begin{align}\label{eq: traceestimate}
\int_{\partial Q_\mu \setminus \Gamma_u} |Tu - a_0|^2\, d\mathcal{H}^1 \le c\mu\Vert e(u) \Vert^2_{L^2(Q_\mu)},
\end{align}
where $Tu$ denotes the trace of $u$ on $\partial Q_\mu$. 
\end{lemma}

\BBB 
 
 \begin{rem}\label{rem: case}
 {\normalfont
 
Two cases are distinguished depending on whether the partition is trivial or not. This is reflected in \eqref{eq:either}. In Section \ref{sec: main-proof} we will say that the square is of \emph{type 1} if $P_0 = \emptyset$ and of \emph{type 2} otherwise. \UUU This distinction corresponds to the two different cases in the proof of Theorem \ref{th: korn-small set} and we refer to Remark \ref{rem: different} for more explanation in that direction. \EEE If $P_0 = Q_\mu$, the situation is actually  much easier and only \eqref{eq: sharpi3}(i), \eqref{eq: sharpi3}(ii) \UUU for a quite simple covering \EEE as well as \eqref{eq: sharpi2.b}-\eqref{eq: traceestimate} are relevant, where \eqref{eq: sharpi2.b}(i) then simplifies to $\Vert \nabla u - A_0\Vert^p_{L^{p}(Q_\mu\setminus E_u)}  \le C\mu^{ 2-p} \Vert e(u) \Vert^p_{L^2(Q_\mu)}$.  However, we prefer to present both cases in one single statement since this is convenient for the application of the lemma in Section \ref{sec: main-proof}.
 
 } 
 \end{rem}

Property \eqref{eq: sharpi2}(i) is a refinement of \eqref{eq: small set main}(i) and shows that only a specific subset $S_u$ of the jump set $J_u$ is needed to control the  length of the boundary of the partition. In particular, the squares of the covering do not intersect $S_u$ (see \eqref{eq: sharpi2}(ii)) which will be crucial in the algorithmic procedure in Section \ref{sec: main-proof}  to show that each part of $J_u$ is only `used' in \emph{one single iteration step} to control the partition. \UUU We emphasize that estimate \eqref{eq: sharpi2}(i) is the essential reason why a distinction into two different cases (1) and (2) is necessary (we refer to Remark \ref{rem: different} for an explanation)\EEE . 

Property \eqref{eq: sharpi3}(ii) reflects  refinement (A) above and \eqref{eq: sharpi3}(i) is related to  \eqref{eq: except}. Properties \eqref{eq: sharpi3}(iii),(iv) on the isolated components have already been established in Lemma \ref{lemma: bad sets}. Note that   \eqref{eq: sharpi2}(iii) states that  each isolated component is a component of the partition. Finally, the fact that the Korn-type estimate  \eqref{eq: sharpi2.b} holds on the covering-adapted sets is a refinement of \eqref{eq: small set main}(ii) and has has been addressed in (B) above. 

Point (2) treats the special case of a trivial partition. In Section \ref{sec: main-proof}  we will  use the trace estimate \eqref{eq: traceestimate} to relate the different behavior of $u$ on adjacent squares of type 2.

\EEE

\smallskip

\Proof Similarly as in the proof of Theorem \ref{th: korn-small set}, we distinguish different cases depending on whether the jump set is large or not.

(1) First assume \BBB $\mathcal{H}^1(J_{u}) \ge \mu \theta^2$ and \EEE   $\mathcal{H}^1(J_{u} \setminus \bigcup\nolimits_{Q \in \mathcal{C}} Q') \ge \mu \theta^2$, where $\mathcal{C}$ is the covering from Theorem \ref{th: bad part},  Lemma \ref{lemma: bad sets} \UUU (applied for $J = J_u$) \EEE  satisfying \eqref{eq: main pro}.  We define $\mathcal{C}_u = \mathcal{C}$ and let $(Z^l_u)_{l \ge 8} =(Z^l)_{l \ge 8}$ be the isolated components from Lemma \ref{lemma: bad sets} \BBB satisfying \eqref{eq: main proXXX1}-\eqref{eq: main proXXX}. \EEE   We let $E_u$ as in Theorem \ref{th: korn-small set} (see also \BBB Remark \ref{rem: excep} \EEE for its definition)  and $R_u = R'$ as in   \eqref{eq.prep}. Denote the partition given in Lemma \ref{cor: partiparti} by $(P_j^*)_{j=1}^N = \mathcal{P}_{R} \cup \mathcal{P}_{\rm good} \cup \mathcal{P}_{\rm small}$ and note that $\mathcal{H}^1\big(\bigcup\nolimits_{j=1}^N \partial P^*_j\big) \le  C\mathcal{H}^1(J_u)$. 

\smallskip

\emph{Definition of the partition:} For each $P_j^* \in \mathcal{P}_{\rm good}$ choose a set $\mathcal{T}(P_j^*) \subset \mathcal{P}_{\rm small}$ with $P^*_{k,\mathcal{C}} \subset P^*_{j,\mathcal{C}}$ for $P^*_{k} \in \mathcal{T}(P^*_{j})$   (cf. \eqref{eq:lastlast}(iv)) such that $\mathcal{P}_{\rm small} = \dot{\bigcup}_{P_j^* \in \mathcal{P}_{\rm good}} \mathcal{T}(P_j^*)$  is a decomposition of $\mathcal{P}_{\rm small}$. Set $\mathcal{T}(P^*_j) = \emptyset$ for $P^*_j \in \mathcal{P}_R$.  Let $(P_j)^n_{j=0}$ be the partition of $Q_\mu$ with  $P_0 = \emptyset$ consisting of the connected components $(X^l_k)_{l,k}$ of $(Z^l_u)_{l \ge 8}$  and the sets
\begin{align}\label{eq:corresponding}
P_j := \big(P^*_j \cup \bigcup\nolimits_{P^*_k \in \mathcal{T}(P^*_j)} P^*_k\big) \setminus \bigcup\nolimits_{l \ge 8} Z^l_u, \ \ \text{ for }  P_j^* \in \mathcal{P}_{R} \cup \mathcal{P}_{\rm good}.
\end{align}

\emph{Proof of \eqref{eq: sharpi2}:}  First, the definition of the partition directly implies \eqref{eq: sharpi2}(iii). By \eqref{eq: main proXXX1}(ii) and $\mathcal{H}^1(J_u) \le 2\sqrt{2}\mu\theta^{-2}$   we see $\mathcal{H}^1(\bigcup_{j=1}^n\partial P_j) \le C\mathcal{H}^1(J_u) \le C\mu$. Then \eqref{eq: sharpi2}(i) follows from the assumption $\mathcal{H}^1(S_u) \ge \mu \theta^2$, where $S_u = J_{u}  \setminus \bigcup_{Q \in \mathcal{C}_u} Q'$. Likewise, the definition of $S_u$ and the fact that  $Q' \subset Q_\mu$ for all $Q \in \mathcal{C}_u$ (see \eqref{eq: main pro}(i)) imply $Q' \subset Q_\mu \setminus S_u$ and particularly \eqref{eq: sharpi2}(ii). 

\smallskip

\emph{Proof of \eqref{eq: sharpi3}:} The first part of \eqref{eq: sharpi3}(i) follows from \eqref{eq: except}. Applying \eqref{eq: partial W}(ii) with $\frac{\theta}{13}\eta$ in place of $\eta$, we get $|R_u| \le \eta$ by \eqref{eq.prep}.  Likewise, \eqref{eq: sharpi3}(iii),(iv) are consequences of Lemma \ref{lemma: bad sets}, \BBB see \eqref{eq: main proXXX1}(i) and \eqref{eq: main proXXX}(i), and $\mathcal{H}^1(J_u) \le 2\sqrt{2}\mu\theta^{-2}$. \EEE Recalling \BBB Remark \ref{rem: excep} \EEE we find  $E_u \setminus R \subset F$ with $F$ as defined in Theorem \ref{th: modifica}  and $R$ as in \eqref{eq: partial W}(ii). Thus, Remark \ref{rem:modi}(i) implies $|Q' \cap (E_u \setminus R)| \le \theta|Q|$ for all $Q\cap \bigcup_{l \ge 8} Z_u^l = \emptyset$. Then the fact that $|Q' \cap R| \le \theta|Q|$  for all $Q \in \mathcal{C}_u$  with $Q \not\subset R_u$ (see \eqref{eq.prep}) yields \eqref{eq: sharpi3}(ii).

\smallskip

\emph{Proof of \eqref{eq: sharpi2.b}:}  \BBB We distinguish the three cases of (a) isolated components, (b) components $\mathcal{P}_R$, and (c) components $\mathcal{P}_{\rm good}$. \EEE

 (a) First, each component $P_j \in (X^l_k)_{l,k}$ of $Z_u := \bigcup_{l \ge 8} Z^l_u$ consists of a union of squares in $\mathcal{Q}^l \cap \mathcal{C}_u$ and thus $P_{j,{\rm cov}} \subset P_j \subset F \subset E_u$ (cf. Remark \ref{rem:modi}(i) and \BBB Remark \ref{rem: excep}).  \EEE

(b) Now consider $P_j$ such that $P_j^* \in \mathcal{P}_{R}$ with $P_j^* $ as in \eqref{eq:corresponding}. As 
$P_j \subset P_j^*\subset R_u$ by \eqref{eq:lastlast}(ii), \eqref{eq: specialdef} implies   $P_{j,{\rm cov}} = \emptyset$. 

(c) For $P_j$ with $P^*_j \in \mathcal{P}_{\rm good}$ we recall that $P^*_j$  is a $c$-John domain ($c \ge 1$). Let $Q \in \mathcal{Q}(P_j;  \mathcal{C}_u, R_u)$ and note that by \eqref{eq:lastlast}(iv), \eqref{eq:corresponding}, and the choice of $\mathcal{T}(P^*_j)$ we have $Q \cap P_j^* \neq \emptyset$. Since $d(Q) \le   4\sqrt{2}d(P_j^*)$ by  \eqref{eq:lastlast}(iii), \UUU we find some $x \in P_j^* \cap Q$  such that $P_j^* \setminus B_{r}(x) \neq \emptyset$, where $r= d(Q)/(9\sqrt{2})$. Therefore, by Lemma \ref{lemma: plump} there exists $z \in \overline{B_{r}(x)}$ such that $B_{r/(2c)}(z) \subset P_j^*$. Also note that $B_{r/(2c)}(z) \subset Q'$. \EEE This  implies that  $|Q' \cap P^*_j| \UUU \ge |B_{r/(2c)}(z)| \EEE \ge c'd(Q)^2$ for some $c'$ small only depending on the John constant $c$. As for   $Q \in \mathcal{Q}(P_j;  \mathcal{C}_u, R_u)$ we have $Q \cap   Z_u = \emptyset$ (see \eqref{eq:corresponding}) and $Q \not\subset R_u$, this yields
\begin{align}\label{eq: also}
|(P^*_j \cap Q') \setminus E_u| \ge c|Q| 
\end{align}
for $\theta$ small by   \eqref{eq: sharpi3}(ii). \BBB In particular,  this also implies that each  $Q \in \mathcal{C}_u$ is contained only in a bounded number of different $\mathcal{Q}(P_j; \mathcal{C}_u,R_u)$. \EEE For   $Q \in \mathcal{Q}(P_j;  \mathcal{C}_u, R_u)$ we have $Q \cap   Z_u = \emptyset$   and thus \BBB Remark \ref{rem:modi}(ii)  holds for $Q$. \EEE This together with Lemma \ref{lemma: rigid}, $\BBB F \EEE \subset E_u$ (see Remark \ref{rem: excep}),  and \eqref{eq: also} yields $A_Q \in \R^{2 \times 2}_{\rm skew}$  such that 
\begin{align*}
 \Vert \nabla u - A_j \Vert^p_{L^p(Q \setminus E_u)}    &\le C\Vert \nabla u - A_Q \Vert^p_{L^p(Q \setminus \BBB E_u \EEE )} +  C\Vert A_Q - A_j \Vert^p_{L^p((P^*_j \cap Q') \setminus E_u)} \\ & \le C\Vert \nabla u - A_Q \Vert^p_{L^p(Q' \setminus \BBB F \EEE )} + C\Vert  \nabla u - A_j \Vert^p_{L^p((P^*_j \cap Q') \setminus E_u)} \\&\le Cd(Q)^{2-p}\Vert e(u) \Vert^p_{L^2(Q')}+ C\Vert \nabla u - A_j \Vert^p_{L^p((P^*_j \cap Q') \setminus E_u)}.
\end{align*}
Summing over all $P_j$ with $P_j^* \in \mathcal{P}_{\rm good}$ and $Q \in \mathcal{Q}(P_j; \mathcal{C}_u,R_u)$  we finally get \eqref{eq: sharpi2.b}(i). Indeed, for the term on the right we use \BBB \eqref{eq:lastlast}(i)    \EEE and the fact that each $x \in Q_\mu$ is contained in at most 12   different $Q'$, $Q \in \mathcal{C}_u$, \BBB see \eqref{eq: main pro}(iii). \EEE For the left term we \BBB use that each square is contained only in a bounded number of different $\mathcal{Q}(P_j; \mathcal{C}_u,R_u)$ (see below \eqref{eq: also}) and \EEE apply the discrete H\"older inequality to compute (cf. \eqref{eq: main holder} for a similar argument)
 \begin{align}\label{eq: main holderXXX}
\sum\nolimits_{Q \in \mathcal{C}_u} (d(Q))^{2-p}\Vert e(u) \Vert^p_{L^2(Q')} \le \big(\sum\nolimits_{Q \in \mathcal{C}_u} |Q|\big)^{1-\frac{p}{2}}\Vert e(u) \Vert^p_{L^2(Q_\mu)}.  
 \end{align}
\BBB This concludes the proof of \eqref{eq: sharpi2.b}(i). Note that (a),(b) and \eqref{eq: also} also show that \EEE each $Q \in \mathcal{C}_u$ intersects only a bounded number of different $(P_{j,{\rm cov}})_{j=1}^n$, which establishes \eqref{eq: sharpi2.b}(ii).

\smallskip

(2) (i) Now suppose that  \BBB $\mathcal{H}^1(J_{u}) \ge \mu \theta^2$ and \EEE  $\mathcal{H}^1(J_{u} \setminus \bigcup\nolimits_{Q \in \mathcal{C}} Q') < \mu \theta^2$. In particular, we have $\mathcal{H}^1(J_{\bar{u}}) \le \mu \theta^2$, \BBB where $\bar{u}$ is the  approximation  from   Theorem \ref{th: modifica}.  We now apply   Theorem \ref{th: kornSBDsmall} on $\bar{u}$. \EEE   Let $E_u = E' \cup F$ with $E'$ as in \eqref{eq: R2main} and $F$ as given by Theorem \ref{th: modifica}. Moreover, we let $P_0 = Q_\mu$, $a_0$ as in \eqref{eq: main estmain},   $R_u = \emptyset$,  $\mathcal{C}_u = \lbrace Q \in \mathcal{Q}^1: Q \subset Q_\mu \rbrace$ and  $Z_u^l = \emptyset$ for all $l \ge 8$.  \BBB This also implies $S_u = (J_{u} \cap Q_\mu)  \setminus \bigcup_{Q \in \mathcal{C}_u} Q'= \emptyset$. \EEE

First, \eqref{eq: sharpi2}(i),(iii) are trivially satisfied since the partition only consists of $P_0$. Moreover, \eqref{eq: sharpi2}(ii) follows from the fact that $\bigcup_{l \ge 8} Z^l_u = \emptyset$ and $S_u = \emptyset$. Properties \eqref{eq: sharpi3}(iii),(iv) are trivial. By $\mathcal{H}^1(J_u) \le 2\sqrt{2}\mu \theta^{-2}$, \BBB $\mathcal{H}^1(J_{\bar{u}}) \le \mu \theta^2$, \EEE \eqref{eq: R2main} and \eqref{eq: except***} we get
$$|E_u| \le |E'| + |F| \le c(\mathcal{H}^1(J_{\bar{u}}))^2 + c\mu\theta^5 \mathcal{H}^1(J_u) \le c\mu^2\theta^3.$$
Since $\mathcal{H}^1(J_{\bar{u}}) \le \mathcal{H}^1(J_{u})$ and \BBB $\mathcal{H}^1(J_{\bar{u}}) \le \mu \theta^2$, \EEE this implies \eqref{eq: sharpi3}(i). Also \eqref{eq: sharpi3}(ii) holds as $|Q| = 4\mu^2\theta^2$ for all $Q \in \mathcal{C}_u$.  As $P_{0,{\rm cov}} = P_0$,  \eqref{eq: sharpi2.b} follows from  \eqref{eq: main estmain} \BBB (applied on $\bar{u}$) and \eqref{eq: modifica prop**}(ii). \EEE Finally, \eqref{eq: main estmainXX} applied on $\bar{u}$ together with $\mathcal{H}^1(J_{\bar{u}}) \le \mu \theta^2$ and  $Tu = T\bar{u}$ on $\partial Q_\mu$ (see Theorem \ref{th: modifica}) yields \BBB a set $\Gamma_u \subset \partial Q_\mu$ with $\mathcal{H}^1(\Gamma_u) \le c\theta\mu$ such that \EEE \eqref{eq: traceestimate} holds. 

\BBB  (ii) Finally, if $\mathcal{H}^1(J_{u}) < \mu \theta^2$, we proceed as in (i), but apply Theorem \ref{th: kornSBDsmall} directly  on $u$. \EEE  \eop

\UUU Let us conclude with a brief remark about the choice of $\mathcal{C}_u$ in case (2). In case (i) it would also have been natural to choose the covering from Theorem \ref{th: bad part}, but since this is not possible in case (ii) (Theorem \ref{th: bad part} is not applicable), the present choice is the most convenient. \EEE

 \section{Derivation of the piecewise Korn inequality}\label{sec: main-proof}  

This section is devoted to the proof of our piecewise Korn inequality.   We first assume that $\Omega = (-\mu_0,\mu_0)^2$ is a square with $\mu_0 >0$. Similarly as in \eqref{eq: enlarged squaresXXX}-\eqref{eq: enlarged squares}, we define $\mathcal{Q}^i := \BBB \tilde{\mathcal{Q}}^{t_i} \EEE$ for $t_i = \mu_0 \theta^{i}$, $\theta \BBB \in 2^{-\N} \EEE $, and denote by  $Q'= \frac{3}{2}Q$  the enlarged square corresponding to $Q \in \mathcal{Q}^i$. (We will apply the results of Section \ref{sec: small set} on squares $(-\mu,\mu)^2$ of different sizes and thus in general $s_i \neq t_i$, where $s_i$ as defined below \eqref{eq: enlarged squaresXXX}.)    Recall the definition of the set $\mathcal{W}(\Omega)$ before Theorem \ref{th: korn-small set}. We first establish the main result for  functions in $\mathcal{W}(\Omega)$ which additionally satisfy 
\begin{align}\label{eq: general assumption}
\mathcal{H}^1(Q \cap J_u) \le \theta^{-1} d(Q) \ \ \ \ \text{ for all } \ Q \in \bigcup\nolimits_{i\ge 1} \mathcal{Q}^i.
\end{align}
Afterwards, we drop   assumption \eqref{eq: general assumption} and show the result for all functions in $\mathcal{W}(\Omega)$. Finally, we pass to general domains with Lipschitz boundary and Theorem \ref{th: main korn} will be derived using a density argument (see Section \ref{sec: main***}). 

\subsection{Proof for configurations with regular jump set}\label{sec: main*}

In this section we let $\Omega =  Q_{\mu_0} = (-\mu_0,\mu_0)^2$. We first prove Theorem \ref{th: main korn} for configurations in $\mathcal{W}(Q_{\mu_0})$ satisfying \eqref{eq: general assumption}, which already represents the core of the argument. If \eqref{eq: general assumption} is violated on a square $Q$, we add $\partial Q$ to the boundary of the partition and treat the problem on $Q$ as a separate problem independent from the analysis on the rest of the domain (see Theorem \ref{th: korn-small-nearly2} below for details).

\BBB

Recall that in Section \ref{sec: small set} we have already established a piecewise Korn inequality which holds outside of a small exceptional set. The derivation of the main result without exceptional set is based on an iterative application of  Lemma \ref{rem: large components}   on various mesoscopic scales. The strategy is to apply  first Lemma \ref{rem: large components} on $Q_{\mu_0}$ to obtain a partition, corresponding infinitesimal rigid motions, and a covering of  $Q_\mu$  with pairwise disjoint dyadic squares. Then Lemma \ref{rem: large components} is again applied on each square of this covering. This procedure is iterated to obtain a progressively smaller exceptional set. For the following we will assume that the reader is already familiar with the statement and the notions of Lemma  \ref{rem: large components}. Note, however, that in this section we will exclusively use Lemma \ref{rem: large components} and no other details from Section \ref{sec: small set} are needed.

Before we come to the rigorous proof, \UUU we first briefly give an intuition why an iteration scheme seems to be indispensable and we describe the main strategy of this procedure.

\smallskip

\textbf{Necessity of an iteration scheme:} First, we observe that it is basically enough to establish the result up to a small exceptional set whose area depends on the function $u$: Indeed, as $u \in \mathcal{W}(Q_{\mu_0})$ and thus $\nabla u \in L^p( Q_{\mu_0})$, we find $\eps=\eps(u)>0$  small such that  for all Borel sets $B \subset Q_{\mu_0}$ with $|B| \le \eps$ one has $\Vert \nabla u \Vert_{L^p(B)} \le \mu_0^{2/p-1}\Vert e(u) \Vert_{L^2( Q_{\mu_0})}$. Based on that, one can show that it is enough to control   $\nabla u$ on the complement of such a set of small measure.

The first naive idea could be to apply Theorem \ref{th: korn-small set} for $\theta$ sufficiently small such that the exceptional set $E$ in \eqref{eq: except} indeed satisfies $|E| \le \eps$ so that the above estimate applies. Observe, however, that then $\theta$ depends on $\eps$ (and thus on $u$), i.e., the constant $C=C(p,\theta)$ in \eqref{eq: small set main} depends on $u$ entailing nonuniform estimates for the partition and $\nabla u$. Therefore, it appears that $\theta$ has to be fixed and the arguments have to be reliant upon  an iterative application of the results presented in Section \ref{sec: small set}.

A first idea of an iterative scheme could be to apply Theorem \ref{th: korn-small set} on $Q_{\mu_0}$, then to cover suitably the resulting exceptional set with smaller squares and to again apply Theorem \ref{th: korn-small set} in these squares. Without going into detail, we mention that an iterative application of this strategy will indeed lead to a sequence of exceptional sets $(E_k)_k$ ($k$ denotes the iteration step) with $|E_k| \to 0 $ as $k \to \infty$. At the same time, however, the constant for the estimate on the partition and $\nabla u$ (cf.  \eqref{eq: small set main}) might be of the form $C^k$ for some $C>1$, eventually blowing up as the number of iteration steps $k$ tends to infinity.    

Consequently, this strategy is non expedient in general  and a more elaborated iteration scheme has to be applied, which is based on the refined version of Theorem \ref{th: korn-small set} and ensures that the relevant constants do not blow up as $k \to \infty$. We will now proceed by describing the cornerstones of this iteration scheme, concentrating on the underlying ideas and without giving exact statements and notations. For the actual proof we refer to Theorem \ref{th: korn-small-nearly} below.  \EEE

\smallskip

\textbf{Description of the iteration scheme:} To fix the main ideas and to give some intuition about the algorithmic procedure, we now illustrate one iteration step in detail.  In particular, we describe  how the partition, the infinitesimal rigid motions, the covering, and the exceptional sets are updated and we  explain how uniform estimates on the boundary of the partition \eqref{eq: main1} and the Korn-type inequality \eqref{eq: main korn2}(ii)  can be obtained. Recall that, due to some technical constructions in Section \ref{sec: small set}, in Lemma \ref{rem: large components} we have to take a set of so-called isolated components into account (we refer to Lemma \ref{lemma: bad sets}  and the comments thereafter for more details). Since in these sets a slightly different construction has to be performed, we  first neglect them. The necessary adaptions will be described afterwards.

 \smallskip

\textbf{The case without isolated components:} Let us suppose that in an iteration step we have given a partition $(P_j)_j$ of $Q_{\mu_0}$, corresponding infinitesimal rigid motions $(a_j)_j$, as well as a covering $\mathcal{C}$, an exceptional set $E$, and a rest set $R$. \UUU We  remark that in general squares of $\mathcal{C}$ might intersect various components of the partition   $(P_j)_j$. \EEE We assume that for each square of $\mathcal{C}$ which is not completely contained in $R$ only a small portion is covered by  $E$.   Moreover, suppose that there is a subset $S \subset J_u$ such that the length of the boundary of the partition can be controlled in terms of $S$ and $S$ does not intersect the squares of the covering. More precisely, we have
\begin{align}\label{eq: intuition1} 
\mathcal{H}^1\Big( \bigcup\nolimits_j\partial P_j  \setminus \partial Q_{\mu_0}\Big) \le C\mathcal{H}^1(S), \ \ \ \ \ \ S \cap Q = \emptyset \ \ \ \text{for all} \ \ \ Q \in \mathcal{C}. 
\end{align}
Finally, we assume that we have a Korn-type estimate of the form
\begin{align}\label{eq: intuition2} 
\sum\nolimits_j\Vert \nabla u - A_j\Vert^{p'}_{L^{p'}(P_{j,{\rm cov}} \setminus E)}  \le C\mu_0^{ 2-p'} \Vert e(u) \Vert^{p'}_{L^2(Q_{\mu_0})},
\end{align}
for $p'<2$, where the covering-adapted sets $P_{j,{\rm cov}}$ are defined below \eqref{eq: specialdef} with respect to $\mathcal{C}$ and $R$. In particular, we remark that in the first iteration step all properties are guaranteed by an application of Lemma \ref{rem: large components} on $Q_{\mu_0}$.

We now pass to the next iteration step. To this end, we apply Lemma \ref{rem: large components} on each $Q \in \mathcal{C}$. (Note that Lemma \ref{rem: large components}  is applicable due to   \eqref{eq: general assumption}.) In view of \eqref{eq:either}, we distinguish between two types of squares. For squares of \emph{type 1}, denoted by $\mathcal{C}_1$, one has a partition $\mathcal{P}^Q$ of $Q$ satisfying \eqref{eq: sharpi2}(i) and associated infinitesimal rigid motions (see   \eqref{eq: sharpi2.b}).  For squares of \emph{type 2}, denoted by $\mathcal{C}_2$,  the partition given by Lemma \ref{rem: large components} is trivial and consists only of $Q$. Correspondingly, there is only one infinitesimal rigid motion $a_Q$, see also Remark \ref{rem: case}.   Apart from the partition and the affine mappings, Lemma \ref{rem: large components} also yields a covering $\mathcal{C}^Q$ and sets $E^Q$, $R^Q$. The relevant notions are now updated as follows:

\begin{itemize}
\item[(a)] \emph{Covering:} We define $\mathcal{C}_{\rm new} = \bigcup_{Q\in \mathcal{C}} \mathcal{C}^Q$. In particular, $\mathcal{C}_{\rm new}$ is a refinement of $\mathcal{C}$.
\item[(b)] \emph{Exceptional set:} We define $E_{\rm new} = \bigcup_{Q\in \mathcal{C}} E^Q$. Note that in general  this set has no relation to the exceptional  set $E$ of the previous iteration step. 
\item[(c)] \emph{Rest set:} We define $R_{\rm new} = R \cup \bigcup_{Q\in \mathcal{C}} R^Q$.
\item[(d)] \emph{Partition:} The new partition contains all partitions $\mathcal{P}^Q$ given by squares $Q \in \mathcal{C}_1$ of type 1 and \UUU it contains  the components $(P_j^{\rm new})_j$ which \EEE are obtained from the partition of the previous iteration step \UUU by cutting away the squares of type 1, \EEE i.e., $P_j^{\rm new} := P_j \setminus \bigcup_{Q \in \mathcal{C}_1} Q$. \UUU We observe that each $P_j^{\rm new}$ is covered by squares in $\mathcal{C}_2$ but in general squares in $\mathcal{C}_2$ may intersect various components $P_j^{\rm new}$. \EEE
\item[(e)] \emph{Infinitesimal rigid motions:} On the components of the partitions given on squares  of type 1 we use the corresponding infinitesimal rigid motions given by Lemma \ref{rem: large components}, see \eqref{eq: sharpi2.b}. For $P_j^{\rm new}$ we choose the mapping $a_j$ from the previous iteration step. 
\end{itemize}

As $\mathcal{C}_{\rm new}$ is a refinement of $\mathcal{C}$, we particularly have that the squares of the covering are smaller, i.e., $\max_{Q \in \mathcal{C}_{\rm new}}  d(Q)\le \theta \max_{Q \in \mathcal{C}}  d(Q)$. From \eqref{eq: sharpi3}(i) we get that on each square $|E^Q|$ scales like $d(Q) \mathcal{H}^1(J_u \cap Q)$. Consequently, $|E_{\rm new}|$ scales like $\max_{Q \in \mathcal{C}_{\rm new}}  d(Q) \mathcal{H}^1(J_u)$ and, due to the refinement of the covering, we will see that the area of the exceptional set will become progressively smaller along the iteration steps. Moreover, although the rest set increases during the iteration scheme, its area can still be controlled uniformly by some (arbitrarily small) $\eps>0$ since each $|R^Q|$ can be chosen arbitrarily small, cf. \eqref{eq: sharpi3}(i). (Indeed,  we will see that Theorem \ref{th: main korn} already holds if it holds up to an exceptional set of arbitrarily small size.) We now come to the fundamental points (d), (e) and  explain how a uniform control on the boundary of the partition  and estimates of Korn-type   can be obtained.

\smallskip

\emph{Partition:} \UUU Let us first note that the procedure described in (d) leads to a partition of $Q_{\mu_0}$: indeed, the components of $\mathcal{P}^Q$, $Q \in \mathcal{C}_1$, form a partition of $\bigcup_{Q \in \mathcal{C}_1} Q$ and the sets $(P_j^{\rm new})_j$ form a partition of $\bigcup_{Q \in \mathcal{C}_2} Q$ (up to sets of negligible measure). Concerning \eqref{eq: intuition1}, \EEE the only delicate squares are those of type 1  since  for them new components are added to the partition. In view of  \eqref {eq: sharpi2}(i), the length of the boundary of each partition $\mathcal{P}^Q$, $Q \in \mathcal{C}_1$, is controlled in terms of $S^Q := (J_u \cap Q) \setminus \bigcup_{\hat{Q} \in \mathcal{C}^Q} \hat{Q}'$. Since by \eqref{eq: intuition1} the set $S^Q$ has empty intersection with $S$, we see that the parts of $J_u$ contained in $S$ are \emph{not `used'} to control the new parts of the partition. By \eqref {eq: sharpi2}(ii) we also get that the squares of $\mathcal{C}_{\rm new}$ do not intersect $S^Q$, $Q \in \mathcal{C}_1$. Thus, \eqref{eq: intuition1} also  holds in the new iteration step for the new partition and the new covering with  $S_{\rm new} = S \cup \bigcup_{Q \in \mathcal{C}_1} S^Q$. Consequently, we  conclude that  in the iteration scheme  each part of $J_u$ is only `used' in \emph{one single iteration step} to control the length of the boundary of the partition, which will eventually ensure a uniform bound independent of the number of iteration steps.

\smallskip

\emph{Korn inequalities:} Here, the squares of type 1 are less delicate since by the choice of the infinitesimal rigid motions in (e), we obtain an estimate of Korn-type directly from \eqref{eq: sharpi2.b}. We therefore concentrate on a square $Q \in \mathcal{C}_2$ of type 2 which is not completely contained in the rest set $R_{\rm new}$ and intersects a component $P_j^{\rm new}$. 

 \UUU In view of the definition of   the affine mappings in (e), our goal is to show that 
\begin{align}\label{eq: intuition11}
\Vert \nabla u - A_j\Vert^{p'}_{L^{p'}((Q \cap P^{\rm new}_{j,{\rm cov}})  \setminus E_{\rm new})} \le cd(Q)^{2-p'}\Vert  e(u)\Vert^{p'}_{L^{2}(Q)} + c\Vert \nabla u - A_j\Vert^{p'}_{L^{p'}((Q \cap P_{j,{\rm cov}}) \setminus E)}.
\end{align}
Then by summing over all squares of type 2 and all components $P_{j,{\rm cov}}^{\rm new}$ we indeed obtain an estimate of the form  \eqref{eq: intuition2} in the new iteration step, where the rightmost term in \eqref{eq: intuition11} is controlled  by estimate  \eqref{eq: intuition2} from the previous iteration step.  

We now explain how to derive \eqref{eq: intuition11}.  By the definition of the partition and the covering we have $P^{\rm new}_{j,{\rm cov}} \subset P_{j,{\rm cov}}$. Thus, the control on $\Vert \nabla u - A_j\Vert_{L^{p'}((Q \cap P^{\rm new}_{j,{\rm cov}})  \setminus E)}$  in terms of the right hand side of \eqref{eq: intuition11} is immediate and to establish \eqref{eq: intuition11} it suffices to control $\Vert \nabla u - A_j\Vert_{L^{p'}(Q \cap (E \setminus E_{\rm new}))}$.

To this end,  from Lemma \ref{rem: large components}(2) applied on $Q$ and the definition of $E_{\rm new}$ in (b)  we get the estimate 
\begin{align}\label{eq: intuition12}
\Vert \nabla u - A_Q\Vert^{p'}_{L^{p'}(Q \setminus E_{\rm new})}= \Vert \nabla u - A_Q\Vert^{p'}_{L^{p'}(Q \setminus E^Q)} \le cd(Q)^{2-p'}\Vert  e(u)\Vert^{p'}_{L^{2}(Q)}.
\end{align}  
(Since $Q$ is of type 2, the partition is trivial and only consists of $Q$ itself, see also Remark \ref{rem: case}). We now aim at estimating the difference of $A_j$ and $A_Q$. \EEE We recall that $Q \not\subset R_{\rm new}$. Thus, by the assumption  stated before \eqref{eq: intuition1}, $E$ covers only a small part of $Q$. Likewise, $E^Q$ covers only small portion of $Q$ by \eqref{eq: sharpi3}(ii). (This property was one of the motivations to establish Lemma \ref{rem: large components}, see refinement (A) in Section \ref{sec: refined}.) By the triangle inequality this allows to estimate the difference of $A_j$ and $A_Q$ and we are able to show (we refer to \eqref{eq: esti1} below for details)
\begin{align}\label{eq: intuition3}
\Vert \nabla u - A_j\Vert^{p'}_{L^{p'}( Q \cap (E \setminus E_{\rm new}))} & \le  c\Vert \nabla u - A_Q\Vert^{p'}_{L^{p'}((Q \cap E) \setminus E^Q)} + c|Q \setminus (E \cup E^Q)| |A_j -A_Q|^{p'}\notag\\
&\le c\Vert \nabla u - A_Q\Vert^{p'}_{L^{p'}(Q \setminus E^Q)} + c\Vert \nabla u - A_j\Vert^{p'}_{L^{p'}(Q \setminus E)}\notag \\
& \le \UUU cd(Q)^{2-p'} \Vert   e(u)\Vert^{p'}_{L^{2}(Q)} + c\Vert \nabla u - A_j\Vert^{p'}_{L^{p'}(Q \setminus E)}. \EEE
\end{align}
\UUU where the last step follows from \eqref{eq: intuition12}. This indeed yields the desired estimate \eqref{eq: intuition11} as by the definition of the  \emph{covering-adapted sets} we have $Q \subset P_{j, \rm cov}$. Note that at this point we fundamentally use the definition of $P_{j, \rm cov}$: it allows us to estimate the difference of $A_j$ and $A_Q$ in terms of $\nabla u - A_j$ and $\nabla u - A_Q$ on the set $Q \setminus (E \cup E^Q) = (Q \cap P_{j, \rm cov}) \setminus (E \cup E^Q) $ instead of only on $(Q \cap P_{j}) \setminus (E \cup E^Q)$. This is essential as  $|(Q \cap P_{j}) \setminus (E \cup E^Q)|$ might be arbitrarily small. \EEE (Recall that this  was one of the motivations to establish Lemma \ref{rem: large components}, see refinement (B) in Section \ref{sec: refined}.)

\GGG

Although this procedure leads to an  estimate outside $E_{\rm new}$, we have to face the additional difficulty that the constant $C$ in \eqref{eq: intuition2} is not uniform but blows up when \eqref{eq: intuition3} is repeatedly applied in different iteration steps. In the proof, this issue is reflected by a decomposition $Q_{\mu_0} = \bigcup_l D_l$,   where, roughly speaking, $D_l$ stands for the set where estimates of the  form \eqref{eq: intuition3} have already been applied (at most) $l$ times. To be more precise, on each of the sets $D_l$ the constant in \eqref{eq: intuition2} can be controlled in terms of  $(\theta^{-\lambda})^l$ for some small $\lambda>0$ to be specified in \eqref{lambdaXXX}, i.e., we have (cf. \eqref{eq: L4} below)
\begin{align}\label{LLL}
\sum\nolimits_j\Vert \nabla u - A_j\Vert^{p'}_{L^{p'}( (P^{\rm new}_{j,{\rm cov}} \cap D_l) \setminus E_{\rm new})}  \le C\theta^{-\lambda l}\mu_0^{ 2-p'} \Vert e(u) \Vert^{p'}_{L^2(Q_{\mu_0})}.
\end{align}
To overcome the difficulty that $\theta^{-\lambda l} \to \infty$ as $l \to \infty$, we use an idea from the proof of Theorem \ref{th: modifica}: the problem that the estimate on $D_l$ blows up as $l \to \infty$ can be compensated by the fact that $|D_l| \to 0$ as $l \to \infty$. Indeed, as discussed above, in each iteration step we need to establish an estimate of the  form \eqref{eq: intuition3} only on $E \setminus E_{\rm new}$. Consequently, as $|E|$ scales like $\max_{Q \in \mathcal{C}} d(Q) \mathcal{H}^1(J_u)$ and thus becomes smaller along the iteration steps, also $|D_l|$ becomes smaller as $l \to \infty$. More precisely, we will see that $|D_l|$ is controlled in terms of $\mu_0^2(\theta^{1-r})^l$ with $r$ from Lemma \ref{rem: large components}, see \eqref{eq: Ddec}.

As in the proof of Theorem \ref{th: modifica}, the strategy is then to reduce the exponent $p'$ to some smaller $p$ and to apply H\"older's inequality: using \eqref{LLL} and $|D_l| \le C\mu_0^2\theta^{l(1-r)}$ we calculate (see \eqref{eq: for later2} below for details)
 \begin{align*}
\sum\nolimits_{l} \sum\nolimits_j\Vert \nabla u - A_j\Vert^{p}_{L^{p}( (P^{\rm new}_{j,{\rm cov}} \cap D_l) \setminus E_{\rm new})}   & \le \sum\nolimits_{l}  |D_l|^{1- \frac{p}{p'}} \big( \sum\nolimits_j\Vert \nabla u - A_j\Vert^{p'}_{L^{p'}( (P^{\rm new}_{j,{\rm cov}} \cap D_l) \setminus E_{\rm new})}\big)^\frac{p}{p'} \\
& \le  C\mu_0^{2-p}  \Vert  e(u) \Vert^{p}_{L^2( Q_{\mu_0})}  \sum\nolimits_{l\ge 0} \theta^{l(1-r)(1-p/p') -l\lambda p /p'}  .
\end{align*}
 Thus, for $\lambda>0$ sufficiently small such that $(1-r)(1-p/p') -\lambda p /p'> 0$ (see \eqref{lambdaXXX} for its choice), we obtain a uniform control \emph{independently of $l$}, i.e., independently of the number of iteration steps in which estimates of the form  \eqref{eq: intuition3} are applied. Here, we emphasize that this H\"older-type argument and the corresponding reduction of the exponent is not performed iteratively in each iteration step but \emph{only once} after the final iteration step. 
 
\EEE

%
%

\smallskip

\textbf{Adaptions for isolated components:} As mentioned above, the construction is slightly different in the presence of isolated components $\bigcup_{l \ge 8} Z^l$ as given in \eqref{eq: sharpi3}. In this case, for the squares of the covering outside of these components we proceed exactly as before. Fix a component $X^l_k$ of $Z^l$ and recall that by \eqref{eq: sharpi2}(iii) the set $X^l_k$ is a component $P_j$ of the partition. We indicate the necessary adaptions inside $X^l_k$. The covering, exceptional set and rest set are updated as in (a)-(c). For the partition and infinitesimal rigid motions we proceed as follows:

\begin{itemize}
\item[(d')] \emph{Partition:} The  partition of $X^l_k$ contains all partitions given by squares   of type 1 contained in $X^l_k$ and the remaining components are given by the connected components of  $X^l_k \setminus \bigcup_{Q \in \mathcal{C}_1} \overline{Q}$, denoted by $(P^{l,k}_j)_j$,  which  consist exclusively of squares of type 2. 
\item[(e')] \emph{Infinitesimal rigid motions:} On the components of the partitions given on squares  of type 1 we use the corresponding infinitesimal rigid motions given by Lemma \ref{rem: large components}, see \eqref{eq: sharpi2.b}. For $(P^{l,k}_j)_j$ we select an arbitrary square $Q$  contained in $P^{l,k}_j$ and as infinitesimal rigid motion  we choose the  mapping $a_Q$, see Remark \ref{rem: case}. 
\end{itemize}

The control on the length of the boundary of the partition can be obtained as   above. Moreover, for squares of type 1 we may repeat the above arguments   for the Korn-type inequality. For the sets $P^{l,k}_j$, however, \eqref{eq: intuition3} cannot be repeated since  $E$ might cover  $X^l_k$ completely, cf. \eqref{eq: sharpi3}(ii). 

Recall that each square contained in $P^{l,k}_j$ is of type 2 and we have a corresponding infinitesimal rigid motion $a_Q$ such that  $\Vert \nabla u - A_Q \Vert_{L^{p'}(Q\setminus E^Q)}$ is controlled (cf.   Remark \ref{rem: case}). Proceeding as in the proof of Theorem \ref{th: modifica}, drawing some ideas from  \cite[Section 4]{FrieseckeJamesMueller:02}, we  estimate the difference of these affine mappings using the trace estimate \eqref{eq: traceestimate}. (In this context, it is convenient that all squares have the same sidelength $2t_l$.) Due to  $d(X^l_k) \le \theta^{-lr} t_l$ (see  \eqref{eq: L5}(iii) below), the involved constant may  blow up in $l$, but as in the proof of Theorem \ref{th: modifica} this can be compensated by using a H\"older-type estimate and the fact that the area of these components is small, cf. \eqref{eq: extendi2} and \eqref{eq: discrHol}(ii) below. 

In order to carry out all Korn-type estimates  simultaneously, in the proof below also the isolated components will be incorporated suitably in the sets $(D_l)_l$.

\EEE


\begin{theorem}\label{th: korn-small-nearly}
Let $p \in [1,2)$. Then Theorem \ref{th: main korn} holds for all $u \in \mathcal{W}(Q_{\mu_0})$ satisfying \eqref{eq: general assumption}.
\end{theorem}

\Proof  Let  $p \in [1,2)$  and $u \in \mathcal{W}(Q_{\mu_0})$   be given satisfying  \eqref{eq: general assumption}. Note that this particularly implies $\mathcal{H}^1(J_u) \le 2\sqrt{2}\theta^{-2}\mu_0$. A classical result states that $u$ is piecewise rigid if $\Vert e(u)\Vert_{L^2(Q_{\mu_0})} = 0$ (see also Remark \ref{rem: remark}(i)), so we can concentrate on the case $\Vert e(u) \Vert_{L^2(Q_{\mu_0})}>0$. Since $\nabla u \in L^p( Q_{\mu_0})$, we can select $\eps$ sufficiently small such that  for all Borel sets $B \subset Q_{\mu_0}$ with $|B| \le 4\eps$ one has 
\begin{align}\label{eq: epsdef}
\Vert \nabla u \Vert^p_{L^p(B)} \le \mu_0^{2-p}\Vert e(u) \Vert^p_{L^2( Q_{\mu_0})}.
\end{align}
We introduce $p' = 1 + \frac{p}{2} \in (p,2)$. Later we will pass from $p'$  to $p$ by  H\"older's inequality. We  let 
\begin{align}\label{lambdaXXX}
r = \frac{2-p}{24} \ \ \ \text{and} \ \ \ \lambda= \frac{(1-r)(2-p)}{3p+2}.
\end{align}
 In the following  $c>0$ stands for a universal constant and  $C=C(\theta,p)>0$  for a generic constant independent of $\eps$ and $\mu_0$.  We can suppose that $\theta$ is small with respect to  $\lambda$, particularly that $\BBB 8 \EEE \theta^{\lambda} \le \frac{1}{6}$ (cf. \eqref{eq: extra} below). At the end of the proof we will fix $\theta$ depending on $p$ such that eventually the constant $C$ depends only on $p$. 

\BBB 
We first formulate the induction hypothesis (Step I) and then discuss the induction step as described above (Step II - Step V). The iteration is repeated a finite number of times until the exceptional set is  smaller than $\eps$. We then conclude the proof by exploiting \eqref{eq: epsdef} (Step VI). 

\EEE

\smallskip
\textbf{Step I (Induction hypothesis):}  \BBB In this step we \EEE formulate the induction hypothesis for step $i$ \BBB and by applying Lemma \ref{rem: large components} we check that it holds for $i=0$. \EEE 

\smallskip

\emph{Hypothesis for step $i$:}  Assume that we have a partition  $\mathcal{P}^i = (P^i_j)_{j=1}^{\infty}$  of $Q_{\mu_0}$ (up to a set of negligible measure),   a covering $\mathcal{C}_i$ of $Q_{\mu_0}$ with $\mathcal{C}_i \subset \bigcup_{j \ge  i+1} \mathcal{Q}^j$ consisting of pairwise disjoint dyadic squares,  exceptional sets $E_i, R_i \subset Q_{\mu_0}$  with $|R_i|\le \frac{\eps}{2}  \sum_{j=0}^i 2^{-j}$,   a set $S_i \subset J_u$  and \BBB pairwise disjoint, closed \EEE sets $(Z_i^l)_{l \ge i+8}$ such that 
\begin{align}\label{eq: L3}
\begin{split}
(i)& \ \ \mathcal{H}^1\Big(\bigcup\nolimits_{j=1}^{\infty}\partial P^i_j \setminus \partial Q_{\mu_0}\Big) \le C_1 \mathcal{H}^1(S_i),\\
(ii) & \ \ Q \subset Q_{\mu_0} \setminus S_i   \ \text{for all}  \ Q \in \mathcal{C}_i, \ \ \ \ \  \overline{Q} \subset Q_{\mu_0} \setminus S_i  \  \text{for all} \  Q \subset \bigcup\nolimits_{l \ge i+8} Z_i^l,\\
(iii) & \ \ \bigcup\nolimits_{l \ge i+8} \bigcup\nolimits_k \partial X^{l,i}_k \subset \bigcup\nolimits_{j=1}^\infty \partial P^i_j
\end{split}
\end{align} 
for some $C_1 = C_1(\theta,p)>0$, where $(X^{l,i}_k)_k$ denote the connected components  of $Z^l_i$. Each set $Z^l_i$ is \BBB a \EEE union of squares in $\mathcal{Q}^l \cap \mathcal{C}_i$ up to a set of measure zero.  Moreover, we have for  $l\ge i+8$ 
\begin{align}\label{eq: L5}
\begin{split}
(i)& \ \  |E_i| \le C_1t_i \mu_0, \\
(ii)& \ \ |Z^l_i| \le C_1\theta^{-rl} t_l\mu_0,  \\
(iii)& \ \  d(X^{l,i}_k) \le \theta^{-rl} t_l \ \ \text{ for all $X^{l,i}_k$,}\\
(iv)  & \ \ |Q \cap E_i| \le c\theta |Q|  \ \text{ for } Q \in \mathcal{C}_i \text{ with } \  Q \not\subset R_i,\ \  Q \cap   \bigcup\nolimits_{l \ge i+8}{Z}^l_i = \emptyset.
\end{split}
\end{align}  
For each $P^i_j$ we set $P^i_{j,{\rm cov}} = \bigcup\nolimits_{Q   \in \mathcal{Q}(P_j; \mathcal{C}_i,R_i)} Q$ with  $\mathcal{Q}(P_j; \mathcal{C}_i, R_i)$ as defined in \eqref{eq: specialdef}. \BBB  We suppose that  there is a decomposition  $ Q_{\mu_0} = \bigcup^{i}_{l=0} D_l^{i}$ of $Q_{\mu_0}$ satisfying
\begin{align}\label{eq: Ddec}
|D_l^i| \le C_2\theta^{-rl} t_l\mu_0\  \ \ \text{ for all } 0 \le l \le i
\end{align}
for some $C_2 = C_2(\theta,p)>0$, and that \EEE for each $P^{i}_j$  there is $A_j^i \in \R^{2 \times 2}_{\rm skew}$ such that for all $0 \le l \le i$
\begin{align}\label{eq: L4}
\begin{split}
(i) & \ \ 
\sum\nolimits_{j=1}^{\infty} \Vert \nabla u - A^{i}_j \Vert^{p'}_{L^{p'}((P^{i}_{j,{\rm cov}} \cap D^i_l) \setminus  E_i )} \le C_3 \theta^{-\lambda l} \mu_0^{2-p'} \Vert  e(u) \Vert^{p'}_{L^2( Q_{\mu_0})},\\
(ii) & \ \ \#\lbrace P^i_{j,{\rm cov}}: x \in P^i_{j,{\rm cov}} \rbrace \le N_0 \ \ \text{ for all } \ \ x \in Q_{\mu_0},
\end{split}
\end{align}
for some $N_0 \in \N$ and a constant $C_3 = C_3(\theta,p)>0$ large enough, which will eventually be specified in \eqref{eq: lll} and \eqref{eq: extra}. 

\BBB Let us briefly comment on the conditions \eqref{eq: L3}-\eqref{eq: L4}. Properties of the form \eqref{eq: L3}-\eqref{eq: L5} appeared already in \eqref{eq: sharpi2}-\eqref{eq: sharpi3} and we refer to the discussion below Lemma \ref{rem: large components} for more details. Moreover, \eqref{eq: L4} is similar to \eqref{eq: sharpi2.b}   with an additional presence of the decomposition $ (D_l^{i})^{i}_{l=0}$. As motivated above, the sets $ (D_l^{i})^{i}_{l=0}$ represent regions where Korn-type estimates have been applied at most $l$ times.  Note that the constant $C_3 \theta^{-\lambda l}$  in \eqref{eq: L4}(i) blows up for $l \to \infty$. This, however, will eventually be compensated by using a H\"older-type estimate  and the fact that $|D^i_l| \to 0$ for $l \to \infty$.   \EEE

\smallskip

\emph{Step $i=0$:} To see that the hypothesis holds for $i=0$, we apply Lemma \ref{rem: large components} on $Q_{\mu_0}$ with $\eta =\frac{\eps}{2}$   to find a partition  $Q_{\mu_0} = \bigcup_{j =0}^{n_0}P_j^0$ and define $\mathcal{C}_0= \mathcal{C}_u$, $E_0 = E_u$,  $R_0 = R_u$,  $Z_0^l= Z_u^l$ for $l \ge 8$ as well as $S_0 := J_u \setminus \bigcup_{Q \in \mathcal{C}_0} Q'$. (Observe that Lemma \ref{rem: large components} is applicable as $\mathcal{H}^1(J_u) \le 2\sqrt{2}\theta^{-2} \mu_0$ by \eqref{eq: general assumption}.) 

Since $\eta = \frac{\eps}{2}$, we get $|R_0| \le \frac{\eps}{2} $ and \eqref{eq: L3}-\eqref{eq: L5} follow from  \eqref{eq: sharpi2}-\eqref{eq: sharpi3} \BBB (for $\mu = \mu_0$ and $s_l= t_l = \mu_0\theta^l$). \EEE Moreover, \eqref{eq: L4} is a  consequence of \eqref{eq: sharpi2.b} with $D^0_0 := Q_{\mu_0}$ and $p'$ in place of $p$  if we choose $C_3$ \BBB and $N_0$ \EEE large enough.

\smallskip
\textbf{Step II (Induction step: Application of Lemma \ref{rem: large components}):}  We now pass from step $i-1$ to $i$. \BBB In Step II we apply Lemma \ref{rem: large components} on each square of the covering. In Step III - Step V we then define all relevant notions and show \eqref{eq: L3}-\eqref{eq: L4} for iteration step $i$. The reader may also first proceed with the conclusion in Step VI in order to see how the result can be deduced from the induction hypothesis. \EEE

\smallskip

For each $Q \in \mathcal{C}_{i-1}$ we apply Lemma \ref{rem: large components} with $p'$ in place of $p$, $\eta = \eps   2^{-i-1} |Q_{\mu_0}|^{-1}|Q|$ and $\mu = t_k$ \BBB with  $k \ge i$ \EEE such that $Q \in \mathcal{Q}^k$. Note that Lemma \ref{rem: large components}  is applicable \BBB since by    \eqref{eq: general assumption} we have $\mathcal{H}^1(Q \cap J_u) \le \theta^{-1} d(Q) = 2\sqrt{2}\theta^{-1}t_k = 2\sqrt{2}\theta^{-1}\mu$. \EEE We obtain a partition    $Q =   \bigcup^{n_Q}_{j=0}  P^Q_j$ up to a set of negligible measure satisfying \eqref{eq:either}, a set $R^Q \subset Q$ with 
\begin{align}\label{eq: RQ}
|R^Q| \le  \eps   2^{-i-1} |Q_{\mu_0}|^{-1}|Q|, 
\end{align}
an exceptional set $E^Q$, a covering $\mathcal{C}^Q$ with $\mathcal{C}^Q \subset \bigcup_{j \ge i+1} \mathcal{Q}^j$, the set $S^Q = (J_u \cap Q) \setminus \bigcup_{\hat{Q} \in \mathcal{C}^Q} \hat{Q}'$,  sets $(Z^l_Q)_{l \ge 8}$,  and infinitesimal rigid motions  $(a^Q_j)_{j=0}^{n_Q}$ such that \eqref{eq: sharpi2}-\eqref{eq: sharpi2.b} hold (with $s_l = t_k \theta^l$ for $l\in \N$ and $\mu = t_k$). Recall particularly that, depending on the case in \eqref{eq:either}, \BBB a partition of $Q$ is formed  either by $\bigcup_{j=1}^{n_Q} P_j^Q$ (square of \emph{type 1}) or by  $P^Q_0$ (square of \emph{type 2}). In the following we will use the notions introduced here without further reference to Step II.   \EEE 

\smallskip

\textbf{Step III (Induction step: Partition and covering}):} \BBB In this step we construct the partition and the covering and we confirm \eqref{eq: L3}(i),(ii). For the underlying ideas we refer to (d),(d') in the description of the iteration scheme. \EEE
 
\smallskip

\emph{Decomposition of the covering:} We split up the covering $\mathcal{C}_{i-1}$ into different sets as follows. Let $\mathcal{C}_{i-1}'$ be the squares $Q \in \mathcal{C}_{i-1}$ with $P^Q_0 =\emptyset$, define
\begin{align}\label{eq:extradef}
\mathcal{C}_{i-1}''' := \big\{ Q \in \mathcal{C}_{i-1} \setminus \mathcal{C}_{i-1}': Q \subset \bigcup\nolimits_{l \ge i+7} Z_{i-1}^l, \ \  \mathcal{H}^1(J_u \cap \partial Q) \ge \theta d(Q)\big\} 
\end{align} 
and $\mathcal{C}_{i-1}'' = \mathcal{C}_{i-1} \setminus (\mathcal{C}_{i-1}' \cup \mathcal{C}_{i-1}''')$. \BBB The squares $\mathcal{C}_{i-1}'$ are of type 1. The squares $\mathcal{C}''_{i-1}$ are of type 2, i.e., the associated partition is trivial (see Lemma \ref{rem: large components}\,(2) and Remark \ref{rem: case}). We remark that the introduction of $\mathcal{C}_{i-1}'''$ is only a technical point \BBB needed for the application of a trace estimate in Step V below. \EEE The core of the argument lies in the separation of $\mathcal{C}_{i-1}'$ and $\mathcal{C}_{i-1}''$. \BBB For the construction of the partition as described in (d),(d') above, it is convenient to treat the squares $\mathcal{C}_{i-1}'''$ as squares of type 1 although the associated partition is trivial. \EEE

\smallskip

\emph{Definition of the partition:}  Recall that the connected components $(X^{l,i-1}_k)_{k \ge 1}$  of $Z^l_{i-1}$ consist of squares $\mathcal{Q}^l \cap \mathcal{C}_{i-1}$, i.e. $\lbrace Q \in \mathcal{Q}^l: Q \subset X_k^{l,i-1} \rbrace = \lbrace Q \in \mathcal{C}_{i-1}: Q \subset X_k^{l,i-1} \rbrace$.  Then by $\mathcal{P}_k^l$ we denote the connected components of  the set
 \begin{align}\label{eq: Y}
 Y^l_k:= {\rm int}\Big( \bigcup\nolimits_{Q \in \mathcal{C}''_{i-1}, Q \subset X^{l,i-1}_k} \overline{Q} \Big)  =  {\rm int} \Big(\bigcup\nolimits_{Q \in \mathcal{Q}^l \setminus ( \mathcal{C}'_{i-1} \cup \mathcal{C}'''_{i-1}), Q \subset X^{l,i-1}_k} \,  \overline{Q}  \Big).
 \end{align}
Note that by \eqref{eq: sharpi2}(i) we have $\mathcal{H}^1(\partial Q) \le \mathcal{H}^1(\bigcup_{j=1}^{n_Q} \partial  P_j^Q) \le C\mathcal{H}^1(S^Q)$ for all $Q \in \mathcal{C}_{i-1}'$. Then 
 \begin{align}\label{eq: clear2}
 \sum\nolimits_{P \in \mathcal{P}^l_k} \mathcal{H}^1(\partial P \setminus \partial X^{l,i-1}_k) &\le \sum\nolimits_{Q \in \mathcal{C}'_{i-1} \cup \mathcal{C}'''_{i-1}, Q \subset X^{l,i-1}_k} \mathcal{H}^1(\partial Q)  \\ 
 &\le C\sum_{Q \in \mathcal{C}'_{i-1}, Q \subset X^{l,i-1}_k}  \mathcal{H}^1(S^Q) + \sum_{Q \in  \mathcal{C}'''_{i-1} , Q \subset X^{l,i-1}_k} \mathcal{H}^1(\partial Q).\notag
  \end{align}
We now introduce the partition $\mathcal{P}^i$ for iteration step $i$. We set
\begin{align}\label{eq: def parti} 
\begin{split}
&\mathcal{P}_a = \lbrace P^Q_j:  Q \in \mathcal{C}'_{i-1}, \ j=1,\ldots,n_Q  \rbrace \cup \lbrace P^Q_0 =Q:  Q \in \mathcal{C}'''_{i-1} \rbrace , \\
& \mathcal{P}_b =  \lbrace   P \in \mathcal{P}^l_k: l \ge i+7, k \ge 1 \rbrace \\
& P'_j := P^{i-1}_j \setminus \ \bigcup\nolimits_{P \in {\mathcal{P}}_a \cup  {\mathcal{P}}_b} \overline{P} \ \  \ \text{for all} \ \ \ P_j^{i-1} \in \mathcal{P}^{i-1}, \ \ \ \mathcal{P}_c = (P'_j)_{j=1}^{\infty}.
\end{split}
\end{align}
Define  $\mathcal{P}^i=  \mathcal{P}_a \cup \mathcal{P}_b \cup \mathcal{P}_c  $ and denote the sets also by $\mathcal{P}^i=(P^i_j)_{j=1}^{\infty}$. We observe that $\mathcal{P}^i$ is a partition of $Q_{\mu_0}$ up to a set of negligible measure and that the sets in $\mathcal{P}_c$ do not intersect $\bigcup_{l\ge i+7} Z^l_{i-1}$. We define the covering $\mathcal{C}_i = \bigcup_{Q \in \mathcal{C}_{i-1}} \mathcal{C}^Q \subset \bigcup_{j \ge i+1} \mathcal{Q}^j$ and now show  \eqref{eq: L3}(i),(ii).

\smallskip

\emph{Proof of \eqref{eq: L3}(i):} \BBB By \eqref{eq: sharpi2}(i) we have $\mathcal{H}^1(\bigcup_{j=1}^{n_Q} \partial  P_j^Q) \le C\mathcal{H}^1(S^Q)$ for all $Q \in \mathcal{C}_{i-1}'$. Moreover, by \EEE \eqref{eq:extradef} with $\hat{S}^Q := J_u \cap \partial Q$  we find $\BBB\mathcal{H}^1(\partial P_0^Q) = \EEE \mathcal{H}^1(\partial Q) =2\sqrt{2}d(Q) \le \BBB 2\sqrt{2}\theta^{-1} \EEE \mathcal{H}^1(\hat{S}^Q)$ for $Q \in \mathcal{C}'''_{i-1}$. Consequently, applying \eqref{eq: L3}(iii) for step $i-1$ as well as \eqref{eq: clear2} we get 
\begin{align*}
\mathcal{H}^1\Big( \bigcup\nolimits^{\infty}_{j=1} \partial P^i_j  \setminus \big( \bigcup\nolimits^{\infty}_{j=1} \partial P^{i-1}_j   \cup \partial Q_{\mu_0} \big)\Big) &\le \BBB\mathcal{H}^1\Big(\bigcup_{ P \in \mathcal{P}_a} \partial P \Big) +  \mathcal{H}^1\Big(\bigcup_{ P \in \mathcal{P}_b} \partial P \setminus \bigcup_{l \ge i+7} \bigcup_k \partial X^{l,i-1}_k \Big) \EEE \\
& \le  C\sum\nolimits_{Q \in \mathcal{C}'_{i-1}} \mathcal{H}^1(S^Q) + C\sum\nolimits_{Q \in \mathcal{C}'''_{i-1}} \mathcal{H}^1(\hat{S}^Q).
\end{align*} 
Letting $S_i = S_{i-1} \cup S^*_i$ with $S^*_i :=\bigcup_{Q \in \mathcal{C}'_{i-1}} S^Q  \cup \bigcup_{Q \in \mathcal{C}'''_{i-1}} \hat{S}^Q$ we obtain $S_i \subset J_u$ and confirm \eqref{eq: L3}(i) as follows. Since   $S^Q \subset Q  \subset Q_{\mu_0} \setminus S_{i-1}$ for all $Q \in \mathcal{C}_{i-1}'$ and $\hat{S}^Q \subset \overline{Q} \subset Q_{\mu_0} \setminus S_{i-1}$ for all $Q \in \mathcal{C}_{i-1}'''$  by \eqref{eq: L3}(ii) for step $i-1$, we get $S_{i-1} \cap S^*_i = \emptyset$. Moreover, each  $x \in Q_{\mu_0}$ is contained in at most  four different $\hat{S}^Q$. Then  the claim follows from \eqref{eq: L3}(i) for step $i-1$ (if we choose $C_1 \ge C$).

\smallskip

\emph{Proof of \eqref{eq: L3}(ii):}   Consider $Q \in \mathcal{C}_i$. First, the fact that $\mathcal{C}_i$ is a refinement of $\mathcal{C}_{i-1}$ yields some $\hat{Q} \in \mathcal{C}_{i-1}$ such that $Q \subset \hat{Q}$ . Then  \eqref{eq: L3}(ii) for step $i-1$ implies  $Q \subset Q_{\mu_0} \setminus S_{i-1}$. Likewise, recalling the definition of the set $S_i^*$, we get $Q \cap S^*_{i} =\emptyset$ by \eqref{eq: sharpi2}(ii) and therefore $Q \cap S_i = \emptyset$.  \BBB This shows the left statement in \eqref{eq: L3}(ii). \EEE

\BBB Now consider the special case that $Q \subset \bigcup_{l \ge 8} Z_{\hat{Q}}^l$. (Recall that this particularly implies $\hat{Q} \in \mathcal{C}_{i-1}'$, see Lemma \ref{rem: large components}\,(2).) \EEE We get $\overline{Q} \subset \hat{Q}$ by \eqref{eq: sharpi2}(ii) and thus as before by  \eqref{eq: L3}(ii), \eqref{eq: sharpi2}(ii) we also derive  $\overline{Q} \subset  Q_{\mu_0} \setminus  S_{i-1}$ and  $\overline{Q} \cap S^*_{i} =\emptyset$. \BBB Anticipating the definition of  $(Z_i^l)_{l \ge i+ 8}$ from \eqref{eq: E def},  we see that the right statement in \eqref{eq: L3}(ii) holds. \EEE

\smallskip

\textbf{Step IV (Induction step: Exceptional sets and isolated components):} \BBB In this step we define the other relevant sets and show \eqref{eq: L3}(iii), \eqref{eq: L5}-\eqref{eq: Ddec}. \EEE Let $R_i = R_{i-1} \cup \bigcup_{Q \in \mathcal{C}_{i-1}} R^Q$ and observe that $|R_i| \le \frac{\eps}{2} \sum_{j=0}^i 2^{-j}$ by \eqref{eq: RQ}. For $l \ge i+8$  we set 
\begin{align}\label{eq: E def}
Z^l_i = \bigcup\nolimits^{l-8}_{k= i} \bigcup\nolimits_{Q \in \mathcal{C}'_{i-1} \cap \mathcal{Q}^k} Z^{l-k}_Q, \ \ \ \ \ \  E_i =  \bigcup\nolimits_{Q \in \mathcal{C}_{i-1}} E^Q.  
\end{align}
Each $Z^l_i$ is a union of squares in $\mathcal{Q}^l \cap \mathcal{C}_i$ up to a set of negligible measure since, as noted below \eqref{eq: sharpi2}, for each $Q \in \mathcal{C}_{i-1} \cap \mathcal{Q}^k$ the sets $Z_Q^j$ consist of squares with sidelength $2s_j$, where according to the notation in the previous section we have $s_j = t_k \theta^j$  for $j \in \N$. 

\smallskip

\emph{Proof of \eqref{eq: L3}(iii):}  \BBB From \eqref{eq: sharpi2}(iii) we get that $\bigcup_{l \ge 8} \partial Z^l_Q \subset \bigcup_{j=1}^{n_Q} \partial P_j^Q$ for each $Q \in \mathcal{C}_{i-1}'$. \EEE Thus, property  \eqref{eq: L3}(iii)  follows from  the definition of $\mathcal{P}_a$ in \eqref{eq: def parti}.

\smallskip

\emph{Proof of \eqref{eq: L5}:}  By \eqref{eq: E def} we get for each square $Q \in \mathcal{C}_{i-1} \cap \mathcal{Q}^k$, $k\ge i$, and $l \ge i+8$ \BBB that $Z_i^l \cap Q  = Z_Q^{l-k} $  and thus \EEE by \eqref{eq: sharpi3}(iv) 
\begin{align*}
|Z_i^l \cap Q|  &\le C  \theta^{-(l-k)r} s_{l-k}\mathcal{H}^1(J_u \cap Q)  = C t_k \theta^{(l-k)(1-r)} \mathcal{H}^1(J_u \cap Q) \le C t_{l} \theta^{-rl} \mathcal{H}^1(J_u \cap Q),
\end{align*}
where as before we have $s_j = t_k \theta^j$ for $j \in \N$.  Then taking the sum over all squares we derive \eqref{eq: L5}(ii), where we use the fact that $\mathcal{H}^1(J_u) \le 2\sqrt{2}\theta^{-2}\mu_0$.  Likewise, \eqref{eq: L5}(iii) follows from  \eqref{eq: sharpi3}(iii) and a similar calculation. Moreover, \eqref{eq: sharpi3}(ii) yields \eqref{eq: L5}(iv) by the definition in \eqref{eq: E def} and the fact that $R^Q \subset R_i$. Finally, to see \eqref{eq: L5}(i), we recall that for all $Q \in \mathcal{C}_{i-1}$ 
$$
|E^Q| \le c\theta^2 d(Q) \mathcal{H}^1(J_u \cap Q) \le Ct_i\mathcal{H}^1(J_u \cap Q)
$$
by \eqref{eq: sharpi3}(i). Consequently, in view of \eqref{eq: E def} and $\mathcal{H}^1(J_u) \le 2\sqrt{2}\theta^{-2}\mu_0$, the result  follows when we take the sum over all $Q \in \mathcal{C}_{i-1}$ (if we choose $C_1 \ge C$). 

\smallskip

\emph{Decomposition $(D_l^i)_l$:} Recalling the definition  of the partition in \eqref{eq: def parti},  we define 
\begin{align}\label{eq: lll2}
D^i_{i} = E_{i-1} \cup \bigcup\nolimits_{P \in \mathcal{P}_a \cup \mathcal{P}_b} P \ \ \text{ and } \ \ D^i_l = D^{i-1}_l \setminus D^i_i \ \ \text{ for } \ \ 0 \le l \le i-1,
\end{align}
and note that $(D^i_l)_{l=0}^i$ is a decomposition of $Q_{\mu_0}$. \BBB As motivated in the description of the iteration scheme, $D^i_i$ represents the region of the domain where we will apply estimates of Korn-type in the current iteration step. We conclude this step with the proof of  \eqref{eq: Ddec}. \EEE First,  \eqref{eq: Ddec} for $0 \le l \le i-1$ follows directly. \BBB Note that each $Q\in \mathcal{C}_{i-1}$ satisfies $|Q| = \frac{1}{16}\mathcal{H}^1(\partial Q)^2 \le \frac{1}{2}t_i\mathcal{H}^1(\partial Q) \le t_i\mathcal{H}^1(\partial Q \setminus \partial Q_{\mu_0})$. \EEE Recalling \eqref{eq: L3}(i), \eqref{eq: def parti} and $\mathcal{H}^1(J_u) \le 2\sqrt{2}\theta^{-2}\mu_0$ we compute
\begin{align*}
\begin{split}
\big|\bigcup\nolimits_{P \in \mathcal{P}_a} P \big| &= \big|\bigcup\nolimits_{Q \in \mathcal{C}'_{i-1} \cup \mathcal{C}'''_{i-1}} Q\big| \le t_i\sum\nolimits_{Q \in \mathcal{C}'_{i-1} \cup \mathcal{C}'''_{i-1}} \mathcal{H}^1(\partial Q \setminus \partial Q_{\mu_0}) \\&\le C_1 t_i\mathcal{H}^1(S_i) \le  C_1 t_i\mathcal{H}^1(J_u) \le C C_1 t_i\mu_0.
\end{split}
\end{align*}
Likewise, by \eqref{eq: def parti} and \eqref{eq: L5}(ii) we see  $|\bigcup\nolimits_{P \in \mathcal{P}_b} P| \le \sum_{l \ge i+7}|Z_{i-1}^l| \le C  C_1\theta^{-ri} t_i\mu_0$. This and the estimate for $|E_{i-1}|$ in \eqref{eq: L5}(i) yields \eqref{eq: Ddec} for $D^i_i$ if one chooses $C_2 =C_2(C_1)$ large enough.

\smallskip

\textbf{Step V (Induction step: Korn inequalities):} \BBB In the previous steps we have already defined the partition in \eqref{eq: def parti},  the exceptional set and the isolated components (see \eqref{eq: E def}). We also recall the definition of the squares $\mathcal{C}_{i-1}', \mathcal{C}_{i-1}'',\mathcal{C}_{i-1}'''$ in \eqref{eq:extradef} and  the decomposition in \eqref{eq: lll2}. It remains to define matrices  $(A^i_j)_{j=1}^{\infty} \subset \R_{\rm skew}^{2 \times 2}$ and to show \eqref{eq: L4}. For the underlying ideas we refer to (e),(e') in the description of the iteration scheme. \EEE 

We recall that  by \eqref{eq: sharpi2.b} we have for all $Q \in \mathcal{C}_{i-1}$
\begin{align}\label{eq: L6.2}
\begin{split}
\sum\nolimits_{j=0}^{n_Q}\Vert  \nabla u - A^Q_j\Vert^{p'}_{L^{p'}(P^Q_{j,{\rm cov}} \setminus E^Q)}\le C(d(Q))^{2-p'} \Vert e(u)\Vert^{p'}_{L^2(Q)} 
\end{split}
\end{align} 
for suitable $(a^Q_j)_j = (a_{A^Q_j,b^Q_j})_j$, where $P^Q_{j,{\rm cov}}  = \bigcup_{\hat{Q}   \in \mathcal{Q}(P^Q_j;  \mathcal{C}^Q,R^Q)} \hat{Q}$   is defined below \eqref{eq: specialdef}. \BBB (Note that for squares of type 2 the estimate simplifies as described in Remark \ref{rem: case}.) \EEE

\smallskip

\enlargethispage{\baselineskip}
\emph{Isolated components:} \BBB As mentioned above in (e'), we need a special construction inside of isolated components, which we discuss here as a preparatory step. \EEE  Fix some $X^{l,i-1}_k$ with $l \ge i+7$ and recalling \eqref{eq: Y} we consider a corresponding connected component  $P \in \mathcal{P}^l_k$.  By definition there are squares $\mathcal{Q}(P) \subset \lbrace Q \in \mathcal{Q}^l: Q \subset X^{l,i-1}_k \rbrace$  such that $P = {\rm int}(\bigcup_{Q \in \mathcal{Q}(P) } \overline{Q})$. 

As $Q \in \mathcal{C}_{i-1}''$ for $Q \in \mathcal{Q}(P)$,  a single infinitesimal rigid motion $a^Q_0$ is given such that \eqref{eq: L6.2} holds with $P_{0,{\rm cov}}^Q = Q$  (see Lemma \ref{rem: large components}\,(2) \BBB and  Remark \ref{rem: case}) \EEE and  by  \eqref{eq: traceestimate}
$$\int_{\partial Q \setminus \Gamma^Q} |Tu-a_0^Q|^2 \, d\mathcal{H}^1  \le ct_l\Vert e(u) \Vert^2_{L^2(Q)}$$
for a set $\Gamma^Q$ with $\mathcal{H}^1(\Gamma^Q) \le c\theta t_l$. We can estimate the difference of the infinitesimal rigid motions on adjacent squares as follows. Consider $Q_1,Q_2 \in \mathcal{Q}(P)$  such that $\partial Q_1$ and $\partial Q_2$ have a common edge. Then the previous estimate implies with $\Gamma := \Gamma^{Q_1} \cup \Gamma^{Q_2} \cup (J_u \cap \partial Q_1 \cap \partial Q_2)$
 $$ \int_{(\partial Q_1 \cap \partial Q_2) \setminus \Gamma }|a_0^{Q_1} - a_0^{Q_2}|^2\,d\mathcal{H}^1\le c\sum\nolimits_{k=1,2} \int_{\partial Q_k \setminus \Gamma}|Tu-a_0^{Q_k}|^2\,d\mathcal{H}^1  \le ct_l\Vert e(u) \Vert^2_{L^2(Q_1 \cup Q_2)},$$
 where we used that the traces $Tu$ on $\partial Q_1$ and $\partial Q_2$ coincide on $(\partial Q_1 \cap \partial Q_2) \setminus \Gamma $. By \eqref{eq:extradef} we get $\mathcal{H}^1(\Gamma) \le c\theta t_l$ as $Q_1,Q_2 \notin \mathcal{C}'''_{i-1}$. Therefore, by Remark \ref{rem: slice}  we get for $\theta$ small $t_l^2|A_0^{Q_1} - A_0^{Q_2}|^2 \le C \Vert e(u) \Vert^2_{L^2(Q_1 \cup Q_2)}$. Consequently, recalling that $P$ is connected, $\#  \mathcal{Q}(P) \le c\theta^{-2rl}$ by \eqref{eq: L5}(iii) and using the discrete H\"older inequality   we find 
$$
\max_{Q_1,Q_2 \in \mathcal{Q}(P)} \ t_l^2|A_0^{Q_1} - A_0^{Q_2}|^{p'}\le ct_l^{2-p'} \theta^{-2 r l(p'-1)} \sum\nolimits_{Q \in \mathcal{Q}(P)}\Vert e(u) \Vert^{p'}_{L^2(Q)}.
$$
Applying \eqref{eq: L6.2}, $\# \mathcal{Q}(P) \le c\theta^{-2rl}$, $p'\le 2$ and recalling \eqref{eq: E def} we derive for each $\hat{Q} \in   \mathcal{Q}(P)$  
\begin{align}\label{eq: extendi2}
\Vert \nabla u- A_0^{\hat{Q}}\Vert^{p'}_{L^{p'}(P \setminus E_i)} \le (C + c\theta^{-4 r l}) t_l^{2-p'}  \sum\nolimits_{Q \in {\mathcal{Q}}(P)}\Vert e(u) \Vert^{p'}_{L^2(Q)}.
\end{align}

\smallskip

\emph{Definition of infinitesimal rigid motions:}   We  now  define   matrices associated to $\mathcal{P}^i$. Recalling the definition of the partition $\mathcal{P}^i$ in \eqref{eq: def parti}, we distinguish   the following cases: (a) For $P^i_j = P^Q_k $ for some $Q \in \mathcal{C}'_{i-1}$ and $k =1,\ldots,n_Q$ we let $A^i_j = A^Q_k$ and for $P^i_j = P^Q_0 $ for some $Q \in \mathcal{C}'''_{i-1}$ we let $A^i_j = A^Q_0$  as given by  \eqref{eq: L6.2}. (b) If $P^i_j \in \mathcal{P}^l_k$, we set $A^i_j = A_0^{\hat{Q}}$ for an (arbitrary) $\hat{Q} \in \mathcal{Q}(P^i_j)$, \BBB see also \eqref{eq: extendi2}. \EEE (c) Finally, for all $P^i_j$ with $P^i_j = P'_k$ for some $P^{i-1}_k  \in \mathcal{P}^{i-1}$  we let $A^i_j = A^{i-1}_{k}$.

\smallskip

\emph{Proof of \eqref{eq: L4}(i):} \BBB We first consider the sets $(D^i_l)_{l=0}^{i-1}$. \EEE By the definition in \eqref{eq: lll2}, we see $D^i_l \subset \bigcup_{P \in \mathcal{P}_c}P$ for $0 \le l \le i-1$. Consider a component $P_j^i = P'_k \in \mathcal{P}_c$. \BBB By \eqref{eq: def parti} \EEE note $P_k' \subset P^{i-1}_k$ with $P^{i-1}_k \in \mathcal{P}^{i-1}$. Recall the definition of the sets $P^i_{j,{\rm cov}}, P^{i-1}_{k,{\rm cov}}$ before \BBB \eqref{eq: Ddec}. \EEE As  $\mathcal{C}_{i}$ is a refinement of $\mathcal{C}_{i-1}$ and $R_i \supset R_{i-1}$,  we  obtain 
\begin{align}\label{eq: *def}
P^i_{j,{\rm cov}} \subset P^*_j:= \bigcup\lbrace Q \in \mathcal{C}_{i-1}:  \ Q \cap P^i_j \neq \emptyset, \ Q \subset P^{i-1}_{k,{\rm cov}} \rbrace \subset P^{i-1}_{k,{\rm cov}}.
\end{align} 
Recall also from \eqref{eq: lll2} that  $D^i_l \setminus E_i \subset D^{i-1}_l \setminus E_{i-1}$ for $0 \le l \le i-1$.   Consequently, \eqref{eq: L4}(i) for $0 \le l \le i-1$ follows directly from the corresponding estimates in step $i-1$ and \eqref{eq: *def}. 

\BBB It remains to consider the set $D^i_i$  as defined in \eqref{eq: lll2}. \EEE It suffices to show that for $k=a,b,c$
\begin{align}\label{eq: lll}
\sum\nolimits_{P^i_j \in \mathcal{P}_k} \Vert \nabla u - A^{i}_j \Vert^{p'}_{L^{p'}((P^{i}_{j,{\rm cov}} \cap D^i_i) \setminus  E_i )} \le \tfrac{1}{3}C_3 \theta^{-\lambda i} \mu_0^{2-p'} \Vert  e(u) \Vert^{p'}_{L^2( Q_{\mu_0})}.
\end{align}
\BBB  As a preparation, \EEE we recall that the discrete H\"older inequality together with $\# \lbrace Q \in \mathcal{C}_i: Q \subset Z^l_i \rbrace \le t_l^{-2} |Z^l_i|  \le C_1 t_l^{-1}\theta^{-rl}\mu_0$   (see  \eqref{eq: L5}(ii)),  $r = \frac{1}{24}(2-p)$,  and  $t_l = \mu_0\theta^l$ implies   
\begin{align}\label{eq: discrHol}
(i)& \ \ \sum\nolimits_{Q \in \mathcal{C}_{i-1}}(d(Q))^{2-p'} \Vert e(u)\Vert^{p'}_{L^2(Q)} \le   c\mu_0^{2-p'}\Vert e(u) \Vert^{p'}_{L^2(Q_{\mu_0})}\notag,\\
(ii) & \ \ t_l^{2-p'} \theta^{-4rl} \sum\nolimits_{Q \in \mathcal{C}_i: Q\subset Z_i^l} \Vert e(u)\Vert^{p'}_{L^2(Q)}  \le  C_1\theta^{-4rl} (\theta^{-rl}t_l\mu_0)^{1-p'/2} \Vert e(u) \Vert^{p'}_{L^2(Q_{\mu_0})} \notag \\ & \quad\quad \quad\quad\quad\quad\quad \quad\quad\quad\quad\quad\quad \quad  \quad \quad \ \le  C_1\theta^{rl} \mu_0^{2-p'} \Vert e(u) \Vert^{p'}_{L^2(Q_{\mu_0})}. 
\end{align}
For (i) we also refer to \eqref{eq: main holderXXX} and for (ii) to \eqref{eq: main holder}-\eqref{eq: main holderXXXX} for  similar arguments,   where we particularly use $1-\frac{p'}{2} = \frac{1}{2}(1-\frac{p}{2}) = 6r$ \BBB and $t_l = \mu_0\theta^l$. \EEE 

\smallskip

\emph{Proof of \eqref{eq: lll} for $\mathcal{P}_a$:} For the components $\mathcal{P}_a$ the property in \eqref{eq: lll}  follows  directly from \eqref{eq: L6.2} and \eqref{eq: discrHol}(i) provided that $C_3$ is chosen large enough, where in this case the additional factor $\theta^{-\lambda  i}$  is not needed. In fact, for all $P^i_j \in \mathcal{P}_a$ we have $P^i_j = P^Q_{k}$ for some $Q\in \mathcal{C}'_{i-1} \cup \mathcal{C}'''_{i-1}$  and  $P^i_{j,{\rm cov}} \BBB \subset \EEE P^Q_{k,{\rm cov}}$ due to the definition of $\mathcal{C}_i$ \BBB and $R_i \supset R^Q$ \EEE with  $P^i_{j,{\rm cov}}$ as defined before \eqref{eq: Ddec}. 

\smallskip

\emph{Proof of \eqref{eq: lll} for $\mathcal{P}_b$:}  Likewise, for  sets $P^i_j \in \mathcal{P}_b$ we get \eqref{eq: lll}  by  \eqref{eq: extendi2} and \eqref{eq: discrHol}(ii) using $P^i_{j,{\rm cov}} \subset \bigcup_{Q \in \mathcal{Q}(P^i_j)}Q \subset  P^i_j$ as well as the fact that the sets $\mathcal{Q}(P)$ for $P \in \mathcal{P}_b$ are pairwise disjoint. Note that again the additional factor $\theta^{-\lambda  i}$ is not needed, but the choice of $C_3$  also depends  on $C_1$ and $\sum_l \theta^{rl} < + \infty$. (A very similar computation has been performed in \eqref{eq: main holderXXXX}.)

\smallskip

\emph{Proof of \eqref{eq: lll} for $\mathcal{P}_c$:} Consider a component $P_j^i = P'_k \in \mathcal{P}_c$ and note $P_k' \subset P^{i-1}_k$ with $P^{i-1}_k \in \mathcal{P}^{i-1}$. \BBB Recall \eqref{eq: *def}. \EEE Fix $Q \in \mathcal{C}_{i-1} \cap \mathcal{Q}^l$ with $Q \subset P^*_j \BBB \subset P^{i-1}_{k,{\rm cov}} \EEE $. As $P_j^i \cap Q \neq \emptyset$,  \eqref{eq: def parti}  implies $Q \cap  Z^l_{i-1} = \emptyset$ and  \BBB that $Q$ is a square of type 2, i.e., \EEE $P_0^Q =P_{0, {\rm cov}} ^Q  = Q$. We observe 
$$(i) \ \ |E_{i-1}  \cap Q| \le c\theta t_l^2, \ \ \ \ (ii) \ \ |E^Q | \le   c\theta t_l^2.$$
In fact, (i) follows from \eqref{eq: L5}(iv) for step $i-1$ and the fact that \BBB $Q \not\subset R_{i-1}$ since \EEE $Q \subset P^{i-1}_{k,{\rm cov}}$ (cf. \eqref{eq: specialdef}). By \eqref{eq: sharpi3}(i)  and \eqref{eq: general assumption} we obtain (ii). \BBB Thus, for $\theta$ sufficiently small with respect to $c$ we obtain $|Q \setminus (E_{i-1} \cup E^Q) | \ge \frac{1}{2}|Q|$ and thus  by the triangle inequality
\begin{align*}
|Q| |A_0^Q - A^{i-1}_k|^{p'} &\le 4\Vert \nabla u - A_0^Q \Vert^{p'}_{L^{p'}(Q \setminus (E^Q \cup E_{i-1}))} +  4\Vert \nabla u - A^{i-1}_k \Vert^{p'}_{L^{p'}(Q \setminus (E^Q \cup E_{i-1})},
  \end{align*}
  where $A^Q_0$ is the matrix given in \eqref{eq: L6.2} for the square $Q$, which is of type 2. By \eqref{eq: E def} and the triangle inequality we find
  \begin{align*}
\Vert \nabla u- A^{i-1}_k \Vert^{p'}_{L^{p'}(Q \setminus E_i)} & \le 2\Vert \nabla u- A_0^Q \Vert^{p'}_{L^{p'}(Q \setminus E^Q)} + 2|Q| |A_0^Q - A^{i-1}_k|^{p'} \\
& \le 10\Vert \nabla u - A_0^Q \Vert^{p'}_{L^{p'}(Q \setminus E^Q)} +  8\Vert \nabla u - A^{i-1}_k \Vert^{p'}_{L^{p'}(Q \setminus   E_{i-1})},
  \end{align*}
and  therefore, using  \eqref{eq: L6.2} and recalling $P_{0, {\rm cov}} ^Q  = Q$, \EEE we derive
\begin{align}\label{eq: esti1}
\Vert \nabla u- A^{i-1}_k \Vert^{p'}_{L^{p'}(Q \setminus E_i)} \le   C(d(Q))^{2-p'}\Vert e(u) \Vert^{p'}_{L^2(Q)} + \BBB 8 \EEE \Vert \nabla u - A^{i-1}_k \Vert^{p'}_{L^{p'}(Q \setminus E_{i-1})}.
  \end{align}
We now confirm \eqref{eq: lll}. Note that each $Q \in \mathcal{C}_{i-1}$ is contained  in at most $N_0$ different $P^*_{j}$ by \eqref{eq: *def} and \eqref{eq: L4}(ii) \BBB for step $i-1$. \EEE   Summing over all squares $Q \subset P^*_{j}  $ and all   $P^{i}_j = P_k' \in \mathcal{P}_c$, we deduce   using  \eqref{eq: L4}(i) for step $i-1$ as well as \eqref{eq: *def}, \eqref{eq: discrHol}(i),  and \eqref{eq: esti1} 
\begin{align}\label{eq: extra}
&\sum\nolimits_{P^i_j = P'_k \in \mathcal{P}_c}  \Vert \nabla u - A^{i-1}_k \Vert^{p'}_{L^{p'}((P^{i}_{j,{\rm cov}} \cap  D^i_i) \setminus E_i)} \le \sum\nolimits_{P^i_j  = P'_k \in \mathcal{P}_c}  \Vert \nabla u - A^{i-1}_k \Vert^{p'}_{L^{p'}((P^{*}_{j} \cap  D^i_i) \setminus E_i)}\notag \\
& \  \ \  \BBB \le   \sum\nolimits_{P^i_j  = P'_k \in \mathcal{P}_c} \sum\nolimits_{Q \in \mathcal{C}_{i-1}, Q \subset P^*_{j}} \big(C(d(Q))^{2-p'}\Vert e(u) \Vert^{p'}_{L^2(Q)} +  8  \Vert \nabla u - A^{i-1}_k \Vert^{p'}_{L^{p'}(Q \setminus E_{i-1})} \big) \EEE \notag \\
& \  \ \  \le N_0\cdot C\mu_0^{2-p'} \Vert e(u) \Vert^{p'}_{L^2(Q_{\mu_0})}   + 8\sum\nolimits_{P^{i-1}_k \in \mathcal{P}^{i-1}}\Vert \nabla u - A^{i-1}_k \Vert^{p'}_{L^{p'}(P^{i-1}_{k,{\rm cov}} \setminus E_{i-1})}\notag \\
&\ \  \ \le (CN_0+ 8C_3 \theta^{-(i-1) \lambda })  \mu_0^{2-p'}  \Vert e(u) \Vert^{p'}_{L^2( Q_{\mu_0})}  \le \tfrac{1}{3}C_3\theta^{-i \lambda } \mu_0^{2-p'}  \Vert e(u) \Vert^{p'}_{L^2( Q_{\mu_0})}, 
\end{align} 
where the last step holds for $C_3=C_3(\theta,p)$ large enough  and $\theta$ small (depending on $\lambda$) such that $CN_0 \le \frac{1}{6}C_3$ and $\BBB 8 \EEE \theta^\lambda \le \frac{1}{6}$. This concludes the proof of \eqref{eq: lll}.

\smallskip

\emph{Proof of \eqref{eq: L4}(ii):} It   suffices to treat the cases (a)-(c) separately. Indeed, by the definition of $\mathcal{C}_i$ we have $P^i_{j_1, \rm cov}$ and $P^i_{j_2, \rm cov}$ are disjoint if they lie in different sets $\mathcal{P}_a,\mathcal{P}_b,\mathcal{P}_c$. For $\mathcal{P}_a$ the desired  property follows directly from \eqref{eq: sharpi2.b}(ii). For $\mathcal{P}_b$ it is obvious as the sets $ P^i_{j,{\rm cov}} \subset P^i_{j} \in \mathcal{P}_b$ are pairwise disjoint. Finally, by \eqref{eq: *def} for each $P^i_j \in \mathcal{P}_c$ we have a (different) $P^{i-1}_{k}$ with  $P^i_{j,{\rm cov}} \subset P^{i-1}_{k,{\rm cov}}$ and the property follows from \eqref{eq: L4}(ii) for step $i-1$.

\smallskip
\textbf{Step VI (Conclusion):} We finally  define the partition and the infinitesimal rigid motions such that the assertion of Theorem \ref{th: main korn} holds.  \BBB For this step of the proof we need the properties stated in Step I. \EEE  As $\limsup_{i \to \infty}|E_i| =0$  by \eqref{eq: L5}(i), we can choose $I \in \N$ so large that  $|E_I| \le \eps$. By $\mathcal{P}^I = (P^I_j)_{j=1}^\infty $ we denote the partition  given in iteration step $I$ and get by \eqref{eq: L3}(i) that 
\begin{align}\label{eq: for later}
\sum\nolimits^\infty_{j=1}\mathcal{H}^1(\partial P^I_j \setminus \partial Q_{\mu_0}) \le C\mathcal{H}^1(S_I)  \le C\mathcal{H}^1(J_u),
\end{align}
where in the last step we used $S_I \subset J_u$.     In view of the definition \eqref{eq: specialdef}, we get
$$\sum\nolimits_{j=1}^\infty |P^I_j \setminus P^I_{j,{\rm cov}}| \le \sum\nolimits_{j=1}^\infty \sum\nolimits_{Q \in \mathcal{C}_I} |P^I_j \cap Q \cap R_I| \le  |R_I|   \le \eps.$$
Thus,  letting $E' = E_I \cup \bigcup_{j=1}^\infty (P^I_j \setminus  P^I_{j,{\rm cov}})$ we get $|E'| \le 2\eps$. We set $\mathcal{P}^* = \lbrace P^I_j \in \mathcal{P}^I: |P^I_j \setminus E'| \ge \frac{1}{2}|P^I_j| \rbrace$  and define 
\begin{align}\label{eq: vdef}
v':= u - \sum\nolimits_{P^I_j \in \mathcal{P}^*} (A^I_j \,x) \chi_{P^I_j}, \ \ \ \ E'' = E' \cup \bigcup\nolimits_{P^I_j \notin \mathcal{P}^*} P^I_j,
\end{align} 
where $(A^I_j)_j \subset \R^{2 \times 2}_{\rm skew}$ are the matrices corresponding to $(P^I_j)_j$ (cf. \eqref{eq: L4}). We now  show \eqref{eq: main korn2}(ii) for $v'$. We deduce from \eqref{eq: L4} \BBB and the definitions of $E', E''$ \EEE that
\begin{align}\label{eq: Dsets2}
 \Vert \nabla v' \Vert^{p'}_{L^{p'}(( Q_{\mu_0}\setminus E'')\cap D^I_l)} \le C\mu_0^{2-p'}\theta^{-l\lambda}  \Vert  e(u) \Vert^{p'}_{L^2( Q_{\mu_0})}
\end{align}
for all $0 \le l \le I$ for a constant $C=C(\theta,p)>0$. Moreover, by \eqref{eq: Ddec}
$$
|D_l^I| \le C \theta^{-rl}t_{l}\mu_0 \le C\theta^{l(1-r)}\mu_0^2
$$
for all $0 \le l \le I$.  Recall the definition $p' = 1+ \frac{p}{2}$, $r = \frac{2-p}{24}$ and $\lambda= \frac{(1-r)(2-p)}{3p+2}$. Passing from $p'$  to $p$ and  using \eqref{eq: vdef}, \eqref{eq: Dsets2} we obtain by H\"older's inequality  
\begin{align}\label{eq: for later2}
\Vert \nabla v'\Vert^p_{L^p( Q_{\mu_0} \setminus E'')} &= \sum\nolimits_{l=0}^{I}  \Vert \nabla v' \Vert^p_{L^p(( Q_{\mu_0}\setminus E'') \cap D^I_l)}  \le \sum\nolimits_{l=0}^I  |D^I_l|^{1- p / p'}\Vert \nabla v' \Vert^p_{L^{p'}(( Q_{\mu_0} \setminus E'') \cap D^I_l)}\notag \\
& \le  C\mu_0^{2-p}\sum\nolimits_{l\ge 0} \theta^{l(1-r)(1-p/p') -l\lambda p /p'}   \Vert  e(u) \Vert^{p}_{L^2( Q_{\mu_0})} \notag \\&\le C\mu_0^{2-p}\sum\nolimits_{l \ge 0}  \theta^{l\lambda}  \Vert  e(u) \Vert^{p}_{L^2( Q_{\mu_0})} \le C\mu_0^{2-p}\Vert  e(u) \Vert^p_{L^2( Q_{\mu_0})}  
\end{align}
with $C=C(\theta,p,\lambda) = C(\theta,p)$.   Now we  have to analyze the behavior in $E''$. To this end, we fix some $P^I_j \in \mathcal{P}^I$. If $P^I_j \in \mathcal{P}^*$, we can choose a measurable set $F_j$ with  $P^I_j \cap E' \subset F_j \subset P^I_j$ and  $|F_j| = 2 |P^I_j \cap E'|$. We    find by $|F_j| \BBB = \EEE  2|F_j \setminus E'| = 2|F_j \setminus E''|$ and the triangle inequality
\begin{align}\label{eq: in betw}
|F_j||A^I_j|^p \le c\Vert \nabla u - A^I_j \Vert^p_{L^p(F_j \setminus E'')} + c \Vert \nabla u \Vert^p_{L^p(F_j \setminus E'')}.
\end{align}
If $P^I_j \notin \mathcal{P}^*$, we set $F_j =P^I_j$. Then we see that $F := \bigcup_{j=1}^\infty F_j$ satisfies $|F| \le 2|E'| \le 4\eps$ and thus by \eqref{eq: epsdef} we derive $\Vert \nabla u\Vert^p_{L^p(F)} \le  \mu_0^{2 - p }\Vert e(u) \Vert^p_{L^2( Q_{\mu_0})}$. Consequently, as $E'' \subset F$ we obtain by \BBB \eqref{eq: vdef} and \EEE \eqref{eq: for later2}-\eqref{eq: in betw}
\begin{align}\label{eq:latref}
\Vert \nabla v' \Vert^p_{L^p( Q_{\mu_0})} &\le \Vert \nabla v' \Vert^p_{L^p( Q_{\mu_0} \setminus E'')} +\Vert \nabla v' \Vert^p_{L^p(F)} \\&\le  C\mu_0^{2-p}\Vert e(u) \Vert^p_{L^2( Q_{\mu_0})} +c\Vert \nabla u \Vert^p_{L^p(F)} + c\sum\nolimits_{P^I_j \in \mathcal{P}^*} |F_j|\vert A^I_j \vert^p \notag \\& \le C\mu_0^{2-p}\Vert e(u) \Vert^p_{L^2( Q_{\mu_0})}  +  c \Vert \nabla v' \Vert^p_{L^p(Q_{\mu_0}\setminus E'')} \le  C\mu_0^{2-p}\Vert e(u) \Vert^p_{L^2( Q_{\mu_0})}. \notag
\end{align}
All arguments hold for a fixed $\theta$ sufficiently small and thus the constant $C$ eventually only depends on $p$. To obtain \eqref{eq: main korn2}(i), we apply Theorem \ref{th: main poinc} on $v'$ and $\rho = \mathcal{H}^1(J_u) + \mathcal{H}^1(\partial Q_{\mu_0})$ to get a Caccioppoli partition $(P'_j)_j$ and corresponding translations $(b_j)_j$ such that with $v := v' -\sum_j b_j \chi_{P_j'}$ by \eqref{eq: for later}, \eqref{eq: vdef}, and H\"older's inequality
\begin{align*}
&\sum\nolimits_j \mathcal{H}^1(\partial^* P'_j) \le \mathcal{H}^1(J_{v'}) + \mathcal{H}^1(\partial Q_{\mu_0}) +c\rho \le C(\mathcal{H}^1(J_u) + \mathcal{H}^1(\partial Q_{\mu_0})),\\
& \Vert v \Vert_{L^\infty(Q_{\mu_0})} \le c\rho^{-1} \Vert \nabla v' \Vert_{L^1(Q_{\mu_0})} \le c \mu_0 \rho^{-1}\Vert e(u) \Vert_{L^2( Q_{\mu_0})}.
\end{align*} 
(For the proof of Theorem \ref{th: main poinc}, which is completely independent from the arguments in Section \ref{sec: small set} and   Section \ref{sec: main-proof}, we refer to Section \ref{subsec: main**1} below.)  Let $(P_j)_j$ be the Caccioppoli partition consisting of the sets $(P^I_i \cap P'_j)_{i,j}$ and observe that  \eqref{eq: main1} follows. Choosing corresponding $(a_j)_j$ such that $v = u - \sum_j a_j \chi_{P_j}$, we obtain \eqref{eq: main korn2} as $\nabla v' = \nabla v$.   Finally, since $u \in \mathcal{W}( Q_{\mu_0})$, we also have that $v \in SBV^2( Q_{\mu_0}) \cap L^\infty(Q_{\mu_0})$. \eop

We now drop assumption \eqref{eq: general assumption}. The idea is to add the boundary of the squares violating \eqref{eq: general assumption} to the boundary of the partition. As in each of these squares the jump set is large, \BBB a careful bookkeeping of their number eventually shows that \EEE the length of the added boundary can be controlled, cf. \eqref{eq:below} below.  On each of these squares we then \BBB apply separately  \EEE Theorem \ref{th: korn-small-nearly}.  For the following proof we introduce the notation $\mathcal{F}^k := \lbrace F = {\rm int}(\bigcup_{Q \in \mathcal{Q}_F} \overline{Q}): \mathcal{Q}_F \subset \mathcal{Q}^k \rbrace$ for $k \in \N$.

\begin{theorem}\label{th: korn-small-nearly2}
Let $p \in [1,2)$. Then Theorem \ref{th: main korn} holds for all $u \in \mathcal{W}(Q_{\mu_0})$.
\end{theorem}

\Proof We reduce the problem to Theorem \ref{th: korn-small-nearly} by passing to suitable modifications of $u$.  Since $u \in \mathcal{W}(Q_{\mu_0})$, we have that $J_u = \bigcup^n_{j=1}  \Gamma^u_j$ consists of closed segments and $\nabla u \in L^p(Q_{\mu_0})$. As in Theorem \ref{th: korn-small-nearly}, it suffices to treat the case $\Vert e(u) \Vert_{L^2(Q_{\mu_0})}>0$.   We choose $I \in \N$ such that
for $\mathcal{R}_0 = \lbrace Q \in \mathcal{Q}^I: Q \cap J_{u} \neq \emptyset \rbrace $
\begin{align}\label{eq: R1}
(i)& \ \ \Vert \nabla u \Vert^p_{L^p(A)} \le \mu_0^{2-p}\Vert e(u) \Vert^p_{L^2(Q_{\mu_0})} \  \text{ for } A \subset \R^2 \text{ with } |A| \le \theta^{-1}\mathcal{H}^1(J_{u}) t_I, \notag\\
(ii) & \ \ \sum\nolimits_{Q \in \mathcal{R}_0} \mathcal{H}^1(\partial Q) \le c\mathcal{H}^1(J_u).
\end{align}
Let  $R_0 = \bigcup_{Q \in \mathcal{R}_0} Q$. Clearly, we have $|R_0| \le c\mathcal{H}^1(J_{u}) t_I\le \theta^{-1}\mathcal{H}^1(J_{u}) t_I$ for $\theta$ small (with respect to $c$). Moreover,  $J_u \cap Q = \emptyset$  for each $Q \in \bigcup_{j \ge I}\mathcal{Q}^{j}$ with $Q \cap R_0 = \emptyset$. For later purpose, we also define $J_u' = J_u \cup \partial R_0$. 

\smallskip

\emph{Identification of bad squares:} Assume that $\mathcal{R}_j \subset \mathcal{Q}^{I-j}$, $0 \le j \le i$, have already been constructed and let $R_j := \bigcup_{Q \in \mathcal{R}_j} Q$ and \BBB $u_j := u \chi_{Q_{\mu_0} \setminus \bigcup_{k=0}^{j} \overline{R_{k}}} \in \mathcal{W}(Q_{\mu_0})$. \EEE For all $j=1,\ldots,i$ we define the sets $S^j_\tau$, $\tau = (\tau_0,\ldots,\tau_{j-1}) \in J_j := \lbrace 0,1\rbrace^{j}$, by 
$$S^j_\tau =  \big( \overline{R_j} \cap    \bigcup\nolimits_{k: \tau_k =1} \overline{R_k}   \big) \setminus   \bigcup\nolimits_{k: \tau_k =0} \overline{R_k} \text{ for } \tau \neq 0, \  \ \ \ \ S^j_0 = \overline{R_j} \setminus \bigcup\nolimits_{k=0}^{j-1} \overline{R_k}.$$
For $j=0$ we let $J_0 = \lbrace 0 \rbrace$ and $S^0_0 = \overline{R_0}$. Denoting $\# \tau = \sum_{k=0}^{j-1} \tau_k$ for $\tau \in J_j$ we assume that for all $j=0,\ldots,i$
\begin{align}\label{eq: R2}
(i)& \ \ J_{u_j}  \subset \Gamma_j := \bigcup\nolimits_{k=0}^j \big(\partial R_k \setminus \bigcup\nolimits_{l = k+1}^j \overline{R_{l}}\big) \cup \big(J_u\setminus \bigcup\nolimits_{k=0}^j \overline{R_k}\big),\notag \\ 
(ii) & \ \ \mathcal{H}^1(J_{u_j} \cap Q) \le \theta^{-1} d(Q) \ \ \text{ for all } Q \in \bigcup\nolimits^\infty_{k = I-j}\mathcal{Q}^k,\\ \notag 
(iii)& \ \ \mathcal{H}^1(\partial R_j \cap F) \le \sum\nolimits^{j}_{k=0} \sum\nolimits_{\tau \in J_j: \#\tau=k} 2^{-k} \  \mathcal{H}^1\big(J'_{u} \cap S^j_\tau \cap F\big) \text{ for all $F \in \mathcal{F}^{I-(j+1)}$}. \notag
\end{align}
From the above discussion we get that the properties hold for $i=0$ where  particularly  (iii) follows from the fact that $J_u' \supset \partial R_0$. We pass to step $i+1$ as follows. Let 
\begin{align}\label{eq: cubbond}
\mathcal{R}_{i+1} = \{ Q\in \mathcal{Q}^{I-i-1}: \mathcal{H}^1(\Gamma_i \cap Q) > 16 t_{I-i-1} \}
\end{align}
and $R_{i+1}  = \bigcup_{Q \in \mathcal{R}_{i+1}} Q$. Set $u_{i+1} = u_i \chi_{Q_{\mu_0} \setminus \overline{R_{i+1}}}$. Then (for $\theta$ small) \eqref{eq: R2}(i),(ii) hold by construction and \eqref{eq: R2}(i),(ii) for step $i$. We now show (iii). Fix  $ F \in \mathcal{F}^{I-(i+2)}$ and recall that $F$ is open. Then \eqref{eq: cubbond} and  the fact that $\mathcal{H}^1(\partial Q) = 8t_{I-i-1}$ yield by \eqref{eq: R2}(i)
\begin{align}\label{eq: first addendXXXXX}
\begin{split}
\mathcal{H}^1(\partial R_{i+1} \cap F)  \le \frac{1}{2} \mathcal{H}^1(\Gamma_i \cap F \cap \overline{R_{i+1}})   
&\le \frac{1}{2}\mathcal{H}^1\Big(J_{u} \cap F \cap \Big(\overline{R_{i+1}} \setminus \bigcup\nolimits_{k=0}^{i} \overline{R_k}\Big)\Big) \\& \ \ \ + \frac{1}{2}\sum\nolimits_{j=0}^{i} \mathcal{H}^1\Big( (\partial  R_j \cap F) \setminus \bigcup\nolimits_{l = j+1}^i \overline{R_l}\Big).
\end{split}
\end{align}
We then note that by definition of $S^{i+1}_0$
\begin{align}\label{eq: first addend}
\mathcal{H}^1\big(J_{u} \cap F \cap \big(\overline{R_{i+1}} \setminus \bigcup\nolimits_{k=0}^{i} \overline{R_k}\big)\big) = \mathcal{H}^1(J_{u}  \cap F \cap S^{i+1}_{0}  ).
\end{align}
Moreover,  for each $\tau \in J_{j}$ we have $S^j_\tau \setminus \bigcup\nolimits_{l = j+1}^i \overline{R_l} = S^{i+1}_{\tau'} = S^{i+1}_{\tau'} \cap (\overline{R_j} \setminus \bigcup\nolimits_{l = j+1}^i \overline{R_l})$ with $\tau' = (\tau,1,0\ldots,0) \in J_{i+1}$. Consequently, since $\#\tau' = \#\tau + 1$ we find by \eqref{eq: R2}(iii) for step $i-1$ and the fact that $F  \setminus \bigcup\nolimits_{l = j+1}^i \overline{R_{l}} \in \mathcal{F}^{I-(j+1)}$  
\begin{align*}
\frac{1}{2}&\sum\nolimits\nolimits_{j=0}^{i} \mathcal{H}^1\big( (\partial  R_j \cap F) \setminus \bigcup\nolimits_{l = j+1}^i \overline{R_{l}}\big)\\
&\le    \sum\nolimits\nolimits_{j=0}^{i}\sum\nolimits^{j}_{k=0} \sum\nolimits_{\tau \in J_j: \#\tau=k} 2^{-k-1} \  \mathcal{H}^1\big(J'_{u} \cap S^j_\tau \cap  (F \setminus \bigcup\nolimits_{l = j+1}^i \overline{R_l})\big)\\
&\le    \sum\nolimits\nolimits_{j=0}^{i}\sum\nolimits^{j+1}_{k=1} \sum\nolimits_{\tau' \in J_{i+1}: \#\tau'=k} 2^{-k} \  \mathcal{H}^1\big(J'_{u} \cap S^{i+1}_{\tau'} \cap F \cap (\overline{R_j} \setminus \bigcup\nolimits_{l = j+1}^i \overline{R_l})\big)\\
& \le  \sum\nolimits^{i+1}_{k=1} \sum\nolimits_{\tau': \#\tau'=k} 2^{-k} \  \mathcal{H}^1\big(J'_{u} \cap S^{i+1}_{\tau'} \cap F \big).
\end{align*}
This together with \eqref{eq: first addendXXXXX}-\eqref{eq: first addend} yields \eqref{eq: R2}(iii). 

For later purpose,  note that  $\big|Q \cap \bigcup\nolimits_{k \le i} R_k \big| \le 2t_{I-i} \mathcal{H}^1(\Gamma_i \cap \overline{Q})$ for all $Q \in \mathcal{Q}^{I-i-1}$,  where we used that $\bigcup_{k \le i} R_k$ consists of squares with sidelength smaller than $2t_{I-i}$ whose boundaries are contained in $\Gamma_i$. \BBB Thus, for all $Q \in \mathcal{Q}^{I-i-1} \setminus \mathcal{R}_{i+1}$ we get for $\theta$ small by \eqref{eq: cubbond} \EEE
\begin{align}\label{eq: thetasmall}
\big|Q \cap \bigcup\nolimits_{k \le i} R_k \big| \le 2t_{I-i} \mathcal{H}^1(\Gamma_i \cap \overline{Q}) \le 2t_{I-i} (16 t_{I-i-1} + 8t_{I-i-1}) \le 2t^2_{I-i-1} = \tfrac{1}{2}|Q|.
\end{align}

We proceed in this way until step $i = I-1$ and finally let $R_I = Q_{\mu_0}$, $\mathcal{R}_I = \lbrace Q_{\mu_0} \rbrace$. By  \eqref{eq: R2}(iii) and the fact that $(S^I_{\tau'})_{\tau' \in J_I}$ is a partition of $Q_{\mu_0}$ we derive
\begin{align*}
\sum\nolimits_{j=0}^{I-1}\mathcal{H}^1(\partial  R_j \cap Q_{\mu_0}) &\le \sum\nolimits_{j=0}^{I-1}\sum\nolimits^{j}_{k=0} \sum\nolimits_{\tau \in J_j: \#\tau=k} 2^{-k} \  \mathcal{H}^1\big(J'_{u} \cap S^j_\tau\big)\notag\\
&= \sum\nolimits_{\tau' \in J_{I}}\sum\nolimits_{j=0}^{I-1}\sum\nolimits^{j}_{k=0} \sum\nolimits_{\tau \in J_j: \#\tau=k} 2^{-k} \  \mathcal{H}^1\big(J'_{u} \cap S^j_\tau \cap S^{I}_{\tau'}\big).
\end{align*}
Note that $S^j_\tau \cap S^{I}_{\tau'} \neq \emptyset$ if and only if $\tau$ coincides with the first $j$ entries of $\tau'$ and $\tau'_j=1$. Consequently, with $\mathcal{H}^1(J_u') \le c\mathcal{H}^1(J_u)$ we find
\begin{align}\label{eq:below}
\sum\nolimits_{j=0}^{I-1}\mathcal{H}^1(\partial  R_j \cap Q_{\mu_0}) & \le  \sum\nolimits_{\tau' \in J_{I}}\sum\nolimits_{j: \tau'_j = 1}  \hspace{-0.1cm} 2^{-(\sum^{j-1}_{l=0} \tau'_l)} \  \mathcal{H}^1\big(J'_{u}   \cap S^{I}_{\tau'}\big) +  \mathcal{H}^1\big(J'_{u}   \cap S^{I}_{0}\big)\notag\\ & \le c\sum\nolimits_{\tau' \in J_{I}}  \mathcal{H}^1\big(J'_{u}   \cap S^{I}_{\tau'}\big)\le c\mathcal{H}^1(J_u).
\end{align}

\smallskip

\emph{Application of Theorem \ref{th: korn-small-nearly}:} We now apply Theorem \ref{th: korn-small-nearly} separately on each square $Q \in \mathcal{R}_j$, $j=1,\ldots,I$, for the function $u_{j-1}$. In fact, in view of   \eqref{eq: R2}(ii), condition \eqref{eq: general assumption} is satisfied (replacing $Q_{\mu_0}$ by $Q$ and $u$ by $u_{j-1}$). Hereby, by \eqref{eq: for later} we obtain a partition of the square $Q \in \mathcal{R}_j$, which we can restrict to $Q \setminus \bigcup_{k=0}^{j-1}\overline{R_k}$. Consequently, taking the union over all partitions defined  on each $Q \in \bigcup_{j=1}^I\mathcal{R}_j$  we obtain a partition $(P'_j)_{j=0}^\infty$ of $Q_{\mu_0}$ with $\bigcup_{j=0}^{I} \partial R_j \subset \bigcup_{j=0}^\infty \partial P'_j$, where we set $P'_0 = R_0$. By \eqref{eq: R2}(i), \eqref{eq:below} and the fact that $\partial R_{k} \cap  \bigcup_{l=0}^{k} R_l = 0$ we get  
\begin{align*}
\sum\nolimits_{j=1}^{I-1}&\sum\nolimits_{Q \in  \mathcal{R}_j}  \big( \mathcal{H}^1(\partial Q \cap Q_{\mu_0})  + \mathcal{H}^1(J_{u_{j-1}} \cap Q) \big) + \mathcal{H}^1(J_{u_{I-1}} \cap Q_{\mu_0}) \\
& \le \sum\nolimits_{j=1}^{I-1} \mathcal{H}^1(\partial R_j \cap Q_{\mu_0}) + \sum\nolimits_{j=1}^{I} \mathcal{H}^1(\Gamma_{j-1} \cap   R_j) \\
&  \le  c\mathcal{H}^1(J_u) + \sum_{j=1}^{I} \mathcal{H}^1\Big(  R_j \cap \bigcup\nolimits^{j-1}_{k=0} \big(  \partial R_k     \setminus \bigcup\nolimits_{l=k+1}^{j-1} \overline{R_l}\big) \Big) + \sum_{j=1}^{I} \mathcal{H}^1\Big(  J_u  \cap \big( R_j \setminus \bigcup\nolimits_{k=0}^{j-1} \overline{R_k} \big)\Big)\\
&  \le  c\mathcal{H}^1(J_u) + \sum_{j=1}^{I} \mathcal{H}^1\Big(  \bigcup\nolimits^{j-1}_{k=0} \partial R_k    \cap \big( R_j \setminus \bigcup\nolimits_{l=0}^{j-1} R_l \big)\Big) + \sum_{j=1}^{I} \mathcal{H}^1\Big(  J_u  \cap \big( R_j \setminus \bigcup\nolimits_{k=0}^{j-1} \overline{R_k} \big)\Big)\\
&\le \sum\nolimits_{j=0}^{I-1}  \mathcal{H}^1(\partial R_j \cap Q_{\mu_0}) + c\mathcal{H}^1(J_u) \le  c\mathcal{H}^1(J_u).
\end{align*}
Consequently, applying \eqref{eq: for later} for  each $Q \in \bigcup_{j=1}^I \mathcal{R}_j$,  we find that $\sum\nolimits^\infty_{j=0}\mathcal{H}^1(\partial P'_j \setminus \partial Q_{\mu_0}) \le C\mathcal{H}^1(J_u)$. Moreover,  Theorem \ref{th: korn-small-nearly} yields piecewise constant $\R^{2 \times 2}_{\rm skew}$-valued functions $\bar{A}_Q$ on each $Q \in \bigcup_{j=1}^I \mathcal{R}_j$ (being constant on each component of the \BBB corresponding \EEE partition) such that by \eqref{eq:latref} and the \BBB fact that $u_{j-1}= u \chi_{Q_{\mu_0} \setminus \bigcup_{k=0}^{j-1} \overline{R_{k}}}$ \EEE we have for $Q \in \mathcal{R}_j$
\begin{align*}
\Vert \nabla u - \bar{A}_Q \Vert^p_{L^p(Q \setminus \bigcup_{k <j} \overline{R_k})}
&\le \Vert \nabla u_{j-1} - \bar{A}_Q \Vert^p_{L^p(Q)} \le C t_{I-j}^{2-p} \Vert e(u_{j-1})\Vert^p_{L^2(Q)} \\ & = C t_{I-j}^{2-p} \Vert e(u)\Vert^p_{L^2(Q \setminus \bigcup_{k <j} \overline{R_k})}.
\end{align*}
If $|Q \setminus \bigcup_{k <j} \overline{R_k}|>0$, we find some $\hat{Q} \in \mathcal{Q}^{I-j+1} \setminus \mathcal{R}_{j-1}$ with $\hat{Q} \subset Q$. Then \eqref{eq: thetasmall} (for $i=j-2$) implies $|Q \setminus \bigcup_{k <j} \overline{R_k}| \ge |\hat{Q} \setminus \bigcup_{k <j-1} \overline{R_k}| = |\hat{Q} \setminus \bigcup_{k \le j-2} {R_k}| \ge 2t^2_{I-j+1} = 2\theta^2 t^2_{I-j}$. Thus, for each $Q \in \mathcal{R}_j$, $j=1,\ldots,I$,
\begin{align*}
\Vert \nabla u - \bar{A}_Q \Vert^p_{L^p(Q \setminus \bigcup_{k <j} \overline{R_k})}
\le  C |Q \setminus \bigcup\nolimits_{k <j} \overline{R_k}|^{1-\frac{p}{2}} \Vert e(u)\Vert^p_{L^2(Q \setminus \bigcup_{k <j} \overline{R_k})},
\end{align*}
where $C=C(\theta,p)>0$. Therefore, for each $(P'_j)_{j \ge 1}$ we find a corresponding $A_j' \in \R^{2 \times 2}_{\rm skew}$ such that for $v' := u -\sum_{j=1}^\infty (A'_j\,x)\chi_{P'_j}$ we have by summing over all squares $Q \in \bigcup_{j=1}^I \mathcal{R}_j$
$$\Vert v' \Vert^p_{L^p(Q_{\mu_0} \setminus R_0)} \le C\mu_0^{2-p} \Vert e(u) \Vert^p_{L^2(Q_{\mu_0})},$$
where similarly as before  we have used the discrete H\"older inequality. (See, e.g., \eqref{eq: main holderXXX}, \eqref{eq: discrHol}. We omit the details.) It remains to analyze the behavior on $R_0$. Observe that Theorem \ref{th: korn-small-nearly} is not applicable in the squares contained in $R_0$ since property \eqref{eq: general assumption} might not hold. However, we may argue similarly as in Step VI of the previous proof. By \eqref{eq: R1}(i) and the fact that $|R_0| \le \theta^{-1}\mathcal{H}^1(J_{u})t_I$ we find $\Vert \nabla u \Vert^p_{L^p(R_0)} \le \mu_0^{2-p}\Vert e(u) \Vert^p_{L^2(Q_{\mu_0})}$. Since $u = v'$ on $R_0=P_0'$, we get $\Vert \nabla v' \Vert^p_{L^p( Q_{\mu_0})} \le C\mu_0^{2-p}\Vert e(u) \Vert^p_{L^2(Q_{\mu_0})}$ and therefore have re-derived \eqref{eq:latref}.   To conclude the proof of \eqref{eq: main korn2}, it now remains to apply Theorem \ref{th: main poinc} on $v'$ as in the previous proof.  \eop

\begin{rem}\label{rem:square}
{\normalfont

The proof of Theorem \ref{th: korn-small-nearly} and Theorem \ref{th: korn-small-nearly2}, in particular \eqref{eq: for later}, show that applying Theorem \ref{th: main poinc} with $\rho = \mathcal{H}^1(J_u)$ instead of  $\rho = \mathcal{H}^1(J_u) + \mathcal{H}^1(\partial Q_{\mu_0})$ yields a partition $(P_j)_j$ and a function $v$ such that \eqref{eq: main korn2}(ii) still holds and \eqref{eq: main1}, \eqref{eq: main korn2}(i) are replaced by 
$$\sum\nolimits_{j=1}^\infty \mathcal{H}^1(\partial^* P_j \cap  Q_{\mu_0}) \le C\mathcal{H}^1(J_u), \ \ \ \ \Vert v\Vert_{L^\infty(Q_{\mu_0})} \le C(\mathcal{H}^1(J_u))^{-1} \Vert e(u) \Vert_{L^2(\Omega)}.$$ 

}
\end{rem}

\subsection{Density arguments}\label{sec: main***}
 
This section is devoted to the proof of the general version of Theorem \ref{th: main korn} and to the proof of Theorem \ref{th: fried-iurlano}. The strategy is to approximate a given function $u \in GSBD^2(\Omega)$ by a sequence $(u_n)_n \subset \mathcal{W}(\Omega)$. (Recall the definition before Theorem \ref{th: korn-small set}.) We then apply Theorem \ref{th: korn-small-nearly2} on $(u_n)_n$ and show that the desired properties can be recovered for the limiting function $u$. We first present a variant of Theorem \ref{th: iurlano} which yields an approximation result for every $GSBD^2$ function at the expense of a nonoptimal estimate for the surface energy.  The reader only interested in the main result for functions in $(G)SBD^2(\Omega) \cap L^2(\Omega)$ may skip the following lemma and corollary and may replace Corollary \ref{cor: iurlano*} by Theorem \ref{th: iurlano} in the proof of Theorem \ref{th: main korn} below. \BBB We also refer  to the recent generalization \cite{Crismale} where density results in $GSBD^p$ are obtained for all $1 < p < \infty$. \EEE

\begin{lemma}\label{lemma: iurlano*}
Let $\Omega \subset \R^d$ open, bounded. Let $u \in GSBD^2(\Omega)$ and  $\Omega' \subset \subset \Omega$ with Lipschitz boundary. Then there exists a sequence $(u_k)_k \subset  \mathcal{W}(\Omega')$ such that  for a universal constant $c>1$ one has
\begin{align}\label{eq: iurl}
(i) & \ \  u_k \to  u \text{ a.e. in $\Omega'$ as } k \to \infty,\\
(ii) & \ \ \Vert e(u_k) \Vert_{L^2(\Omega')} \le c\Vert e(u) \Vert_{L^2(\Omega)}, \  \mathcal{H}^{d-1}(J_u)\le \mathcal{H}^{d-1}(J_{u_k})  \le c\mathcal{H}^{d-1}(J_u), \ k  \in \N.\notag
\end{align}
\end{lemma}
 
\Proof We follow closely the proofs in \cite[Theorem 1]{Chambolle:2004}, \cite[Theorem 3.5]{Iurlano:13} and only indicate the necessary changes. By \cite[Lemma 2.2]{Iurlano:13} one can find a basis $\e_1,\ldots,\e_d$ of $\R^d$ such that 
$$\mathcal{H}^{d-1}(\lbrace x \in J_u: [u](x) \cdot e = 0\rbrace) = 0$$
for all $e \in D: = \lbrace \e_i, 1\le i \le d, \e_i + \e_j, 1 \le i < j \le d \rbrace$. For $h>0$ small and each $y \in [0,1)^d$ we introduce the discretized function of $u$ defined by 
$$u^y_h(\xi) = u(hy+\xi),  \ \ \ \xi \in h\Z^d \cap (\Omega - hy). $$
Moreover, letting $\triangle(x) = \prod^d_{i=1} \max \lbrace (1-|x_i|),0\rbrace$ we define the continuous interpolation 
$$w^y_h(x) = \sum\nolimits_{\xi \in h\Z^d \cap \Omega} u^y_h(\xi) \triangle\Big(\frac{x - (\xi + hy)}{h} \Big).$$
Note that $w^y_h \in W^{1,\infty}(\Omega')$ and for $h$ small enough $w^y_h$ is indeed well defined for all $x \in \Omega'$ since $\Omega' \subset\subset \Omega$. We now show that $w^y_h \to u$ in measure on $\Omega'$ as $h \to 0$ for a.e. $y\in [0,1)^d$ which is equivalent to 
\begin{align}\label{eq: densi1}
\int_{y \in [0,1)^d}\Vert \varphi_1 (|w^y_h -u|) \Vert_{L^1(\Omega')} \to 0,
\end{align}
where $\varphi_t(s) = \min\lbrace s,\frac{1}{t}\rbrace$ for $s \ge 0$ and $t>0$. Using the subadditivity of $\varphi_1$,  a change of variables and the fact that $\sum\nolimits_{\xi \in h\Z^d \cap \Omega}  \triangle\Big(\frac{x - \xi}{h} -  y \Big) = 1$ for all $x \in \Omega'$, $y \in [0,1)^d$, we deduce
\begin{align*}
&\int_{y \in [0,1)^d} \,dy \int_{\Omega'} \varphi_1(|w^y_h(x)-u(x)|)\,dx \\
& \ \ = \int_{y \in [0,1)^d} \,dy \int_{\Omega'} \varphi_1\Big(\Big| \sum\nolimits_{\xi \in h\Z^d \cap \Omega}  \triangle\Big(\frac{x - \xi}{h} -  y \Big) (u(x)-u(hy + \xi))\Big|\Big)\,dx \\ 
& \ \ \le \int_{y \in [0,1)^d} \,dy \int_{\Omega'} \sum\nolimits_{\xi \in h\Z^d \cap \Omega}  \varphi_1\Big(\Big|  \triangle\Big(\frac{x - \xi}{h} -  y \Big) (u(x)-u(hy + \xi))\Big|\Big)\,dx \\ 
& \ \ \le \sum\nolimits_{\xi \in h\Z^d \cap \Omega}  \int_{\frac{x-\xi}{h}-[0,1)^d} \int_{\Omega'} \varphi_1(|  \triangle(z) (u(x)-u(x -hz))|)\,dx \,dz\\
& \ \ \le \int_{(-1,1)^d} \int_{\Omega'} \varphi_1(|  \triangle(z) (u(x)-u(x -hz))|)\,dx \,dz.
\end{align*} 
In the last step we used that the sets $\frac{x-\xi}{h}-[0,1)^d$ are pairwise disjoint for $\xi \in h\Z^d \cap \Omega$.
Since 
$$\int_{\Omega'} \varphi_1(|  \triangle(z) (u(x)-u(x -hz))|)\,dx = \int_{\Omega'} \triangle(z) \varphi_{\triangle(z)}(|   u(x)-u(x -hz)|)\,dx \to 0$$
for $h \to 0$ and is uniformly bounded by $|\Omega'|$ for all $z \in (-1,1)^d$, we obtain \eqref{eq: densi1} by dominated convergence. We now follow closely the proof in \cite{Iurlano:13}. For a specific choice of $y \in [0,1)^d$ and a set of `bad cubes' (see (3.13)-(3.15) in \cite{Iurlano:13} for details)
\begin{align}\label{eq:badsquaresiurl}
\mathcal{Q}_h\subset \lbrace Q_{\xi}=\xi + hy + [0,h)^d: \xi \in h\Z^d \rbrace
\end{align}
with $\# \mathcal{Q}_h \le Ch^{-(d-1)}$ one defines the function $u_h = w_h^y \chi_{\Omega' \setminus \bigcup_{Q \in \mathcal{Q}_h}Q}$. Clearly, $u_h \in \mathcal{W}(\Omega')$ as the jump set is given by the boundary of the cubes $\mathcal{Q}_h$. (Strictly speaking $J_{u_h}$ is only contained in the boundary of the cubes. However, the desired property may always be achieved by an infinitesimal perturbation of $u_h$, see \cite[Remark 5.3]{Chambolle:2004}). In view of \eqref{eq: densi1} and $|\bigcup_{Q \in \mathcal{Q}_h}Q| \le Ch$, we get, possibly passing to a not relabeled subsequence,  $u_h \to u$ a.e. on $\Omega'$ which gives \eqref{eq: iurl}(i). 

Following the lines of \cite{Iurlano:13},  for a suitable choice of $y \in [0,1)^d$ we get $\Vert e(u_h)\Vert_{L^2(\Omega')} \le c\Vert e(u)\Vert_{L^2(\Omega)}$ and $\mathcal{H}^{d-1}(J_{u_h}) \le c\mathcal{H}^{d-1}(J_u)$ for a universal constant $c>1$ after passage to a not relabeled subsequence. In fact, in the estimates for the elastic energy, which are based on a slicing technique, the assumption $u \in L^2(\Omega)$ is not needed. (This assumption was only needed to define an extension of $u$, which is not necessary in our setting.)   To conclude the proof of \eqref{eq: iurl}(ii), we remark that the inequality $\mathcal{H}^{d-1}(J_u) \le \mathcal{H}^{d-1}(J_{u_h})$ has not been stated explicitly in \cite{Iurlano:13}, but can always be achieved by introducing arbitrary additional `bad squares' in \eqref{eq:badsquaresiurl}.   \eop

The proof of the following corollary is now straightforward.  

\begin{corollary}\label{cor: iurlano*}
Let $\Omega \subset \R^d$ open, bounded with Lipschitz boundary and let $Q \subset \R^d$ be a cube with $\Omega \subset \subset Q$. Let $u \in GSBD^2(\Omega)$. Then there exists a sequence $(u_k)_k \subset \mathcal{W}(Q)$ such that  for a universal constant $c>1$ one has
\begin{align}\label{eq: iurlX}
\begin{split}
(i) & \ \  u_k \to  u \text{ a.e. in $\Omega$ as } k \to \infty,\\
(ii) & \ \ \Vert e(u_k) \Vert_{L^2(Q)} \le c\Vert e(u) \Vert_{L^2(\Omega)},\\
(iii) & \ \  \mathcal{H}^{d-1}(J_u) + \mathcal{H}^{d-1}(\partial \Omega) \le \mathcal{H}^{d-1}(J_{u_k})  \le c\mathcal{H}^{d-1}(J_u) +  c\mathcal{H}^{d-1}(\partial \Omega).
\end{split}
\end{align}
\end{corollary}

\Proof Choose a sequence of Lipschitz sets $\Omega_n \subset \subset \Omega$ with $\mathcal{H}^{d-1}(\partial \Omega) \le \mathcal{H}^{d-1}(\partial \Omega_n) \le c\mathcal{H}^{d-1}(\partial \Omega)$  and $|\Omega \setminus \Omega_n| \le \frac{1}{n}$ for all $n \in \N$ such that $\partial \Omega_n$ consists of a finite number of closed $(d-1)$-simplices. For each $n \in \N$ we apply Lemma \ref{lemma: iurlano*} on $\Omega_n \subset \subset \Omega$ and obtain sequences $(v^n_l)_l \subset \mathcal{W}(\Omega_n)$ converging in measure to $u$ on $\Omega_n$ such that $\Vert e(v^n_l) \Vert_{L^2(\Omega_n)}$ and $\mathcal{H}^{d-1}(J_{v^n_l})$ are uniformly controlled by $c\Vert e(u) \Vert_{L^2(\Omega)}$ and $c\mathcal{H}^{d-1}(J_u)$, respectively. Define the extensions $\hat{v}^n_l = v^n_l \chi_{\Omega_n} \in \mathcal{W}(Q)$. Possibly replacing $\hat{v}^n_l$ by $\hat{v}^n_l + t^n_l$ for a suitable $t^n_l \in \R^d$, $|t^n_l| \le \frac{1}{l}$, we obtain (not relabeled) sequences still converging to $u$ on $\Omega_n$ such that $\mathcal{H}^{d-1}(\partial \Omega_n \setminus J_{\hat{v}^n_l})=0$ and  thus
$$\mathcal{H}^{d-1}(J_u) + \mathcal{H}^{d-1}(\partial \Omega) \le \mathcal{H}^{d-1}(J_{\hat{v}^n_l})  \le c\mathcal{H}^{d-1}(J_u) +  c\mathcal{H}^{d-1}(\partial \Omega). $$
Consequently, by a standard diagonal sequence argument taking into account that convergence in measure is metrizable (take $(f,g) \mapsto \int_\Omega \min\lbrace |f-g|,1\rbrace$) we get a sequence $(w_k)_k \subset \mathcal{W}(Q)$ satisfying \eqref{eq: iurlX}.   \eop

We are now in the position to give the proof of Theorem \ref{th: main korn}. 

\noindent {\em Proof of Theorem \ref{th: main korn}.}   Let $p\in (1,2)$  and let $u \in GSBD^2(\Omega)$ be given. Without restriction we assume that $\Omega$ is connected as otherwise the following arguments are applied on each connected component of $\Omega$. Choose a square $Q_{\mu_0} \supset \supset \Omega$ with $\mathcal{H}^1(\partial Q_{\mu_0}) \le c\mathcal{H}^1(\partial \Omega)$ for a universal constant $c>0$. 

By Corollary \ref{cor: iurlano*} we find a sequence $(u_k)_k \subset \mathcal{W}(Q_{\mu_0})$  satisfying \eqref{eq: iurlX}. (The reader only interested in the case $(G)SBD^2(\Omega) \cap L^2(\Omega)$ can instead apply Theorem \ref{th: iurlano} and Remark \ref{rem: iur}.)  We apply Theorem \ref{th: korn-small-nearly2} on each $u_k$ and obtain a sequence of (ordered) Caccioppoli partitions $(P^k_j)_j$ of $Q_{\mu_0}$ and corresponding infinitesimal rigid motions $(a^k_j)_j = (a_{A^k_j,b^k_j})_j$ such that 
$v_k := u_k - \sum\nolimits^\infty_{j=1} a^k_j \chi_{P^k_j} \in SBV^p(Q_{\mu_0}) \cap L^\infty(Q_{\mu_0})$
satisfies by \eqref{eq: main1} and \eqref{eq: main korn2} 
\begin{align}\label{eq: R3}
 (i)& \ \ \sum\nolimits^\infty_{j=1}\mathcal{H}^1(\partial^* P^k_j) \le c(\mathcal{H}^1(J_{u_k})+ \mathcal{H}^1(\partial Q_{\mu_0}))\le c(\mathcal{H}^1(J_{u})+ \mathcal{H}^1(\partial \Omega)),\notag \\  
(ii) & \ \  \Vert \nabla v_k \Vert_{L^p(Q_{\mu_0})} \le C \Vert e(u_k) \Vert_{L^2(Q_{\mu_0})} \le C \Vert e(u) \Vert_{L^2(\Omega)},\\
(iii)  & \ \ \Vert v_k \Vert_{L^\infty(Q_{\mu_0})} \le C(\mathcal{H}^1(J_{u_k}))^{-1} \Vert e(u_k) \Vert_{L^2(Q_{\mu_0})}  \le C(\mathcal{H}^1(J_{u})+ \mathcal{H}^1(\partial \Omega))^{-1} \Vert e(u) \Vert_{L^2(\Omega)}.\notag
 \end{align}
Possibly passing to a (not relabeled) refinement of the partition (consisting of the sets $(P_j^k \cap \Omega)_j \cup (P_j^k \setminus \Omega)_j$) we may assume that each component $P^k_j$ satisfies $P^k_j \subset \Omega$ or $P^k_j \cap \Omega = \emptyset$ and \eqref{eq: R3}(i) still holds. Therefore, Theorem \ref{th: comp cacciop} implies the existence of a Caccioppoli partition $(P_j)_{j=1}^\infty$ of $\Omega$  and a (not relabeled) subsequence such that $\chi_{P^k_j} \to \chi_{P_j}$ in $L^1(\Omega)$, when $k \to \infty$, for all $j \in \N$ and such that,  passing to the limit in \eqref{eq: R3}(i) via the lower semicontinuity of the perimeter  
$$
\sum\nolimits_{j=1}^\infty \mathcal{H}^1( \partial^* P_j  ) \le c (\mathcal{H}^1(J_u) + \mathcal{H}^1(\partial \Omega)).
$$
This gives  \eqref{eq: main1}. Applying Ambrosio's compactness theorem (Theorem \ref{clea-th: compact})  on the sequence $(v_k)_k$ we find  $v \in SBV^p(Q_{\mu_0}) \cap L^\infty(Q_{\mu_0})$  such that $v_k \to v$ a.e. and $\nabla v_k \rightharpoonup \nabla v$ weakly in $L^p(Q_{\mu_0})$ up to a not relabeled subsequence. In particular, we have considering the restriction of $v$ to $\Omega$
\begin{align}\label{eq: R3***}
\Vert v \Vert_{ L^\infty(\Omega)} \le C(\mathcal{H}^1(J_{u})+ \mathcal{H}^1(\partial \Omega))^{-1} \Vert e(u) \Vert_{L^2(\Omega)},\ \ \ \ \  \Vert \nabla v\Vert_{ L^p(\Omega)} \le C \Vert  e(u) \Vert_{L^2(\Omega)}.
\end{align}
Since $(P_j)_j$ is a partition of $\Omega$, it now suffices to show the existence of infinitesimal rigid motions $(a_j)_j = (a_{A_j,b_j})_j$ such that 
\begin{equation}\label{claim}
(u-v)\chi_{P_j}=a_j \chi_{P_j}
\end{equation}
a.e. in $P_j$ for all $j\in \N$. Indeed, \eqref{eq: main korn2} then is immediate by \eqref{eq: R3***}.  Clearly, if $|P_j|=0$ it suffices to set $a_j=0$. If instead $|P_j|>0$, then it exists $\delta>0$ independently of $k$ such that $|P^k_j|\ge \delta$. As $u_k \chi_{P^k_j} \to u \chi_{P_j}$ a.e. and $v_k \chi_{P^k_j}  \to v\chi_{P_j} $ a.e., the sequence 
$$
a^k_j \chi_{P^k_j}=(u_k-v_k)\chi_{P^k_j}$$
converges in measure to $(u-v)\chi_{P_j}$. Therefore, it exists a  positive nondecreasing continuous function $\psi$ with $\lim_{s \to \infty} \psi(s) = +\infty$    such that $\sup_{k \ge 1}\int_{P^k_j}\psi(|a^k_j|)\,\mathrm{d}x \le 1$ (see e.g. \cite[Remark 2.2]{FriedrichSolombrino}). By Lemma \ref{lemma: rigid motionXXX} we infer that $(a^k_j)_k$ are bounded in $W^{1,\infty}(\Omega)$ for a constant independent of $k$.  Consequently, we find an infinitesimal rigid motion $a_j$ such that
$
a^k_j \chi_{P^k_j} \to a_j \chi_{P_j}
$ in $L^{1}(\Omega)$. By the convergence of $a^k_j \chi_{P^k_j}$ to $(u-v)\chi_{P_j}$, this implies \eqref{claim} and concludes the proof. \eop

\UUU

\begin{rem}[Piecewise Korn inequality for $q \neq 2$]\label{rem: general version}
{\normalfont

Having concluded the proof of Theorem \ref{th: main korn}, we remark that the result can be extended to a version in the space $GSBD^q(\Omega)$ for general $1 < q <\infty$ and exponents $1 \le p < q$. Specifically, for  each $u \in GSBD^q(\Omega)$ there is a Caccioppoli partition $(P_j)^\infty_{j=1}$ of  $\Omega$ with $\sum\nolimits^\infty_{j=1} \mathcal{H}^1( \partial^* P_j) \le c(\mathcal{H}^1(J_u) + \mathcal{H}^1(\partial \Omega))$ and   infinitesimal rigid motions $(a_j)_j = (a_{A_j,b_j})_j$ such that $v:= u - \sum\nolimits_{j=1}^\infty a_j \chi_{P_j} \in SBV^p(\Omega;\R^2) \cap L^\infty(\Omega;\R^2)$ and
\begin{align*}
\begin{split}
\Vert v \Vert_{L^\infty(\Omega)}\le C(\mathcal{H}^1(J_u) + \mathcal{H}^1(\partial \Omega))^{-1} \Vert e(u) \Vert_{L^q(\Omega)}, \ \ \  \Vert \nabla v \Vert_{L^p(\Omega)} \le C \Vert e(u) \Vert_{L^q(\Omega)}.
\end{split}
\end{align*}
We emphasize that we do not formulate this as a main result of the paper in Section \ref{sec: main}, but only on a level of a remark since we do not give a complete proof. We only provide some hints on how the arguments in the proof of  Theorem \ref{th: main korn} have to be adapted to the general case. Note that the main adaptions concern the replacement of $L^2$ by $L^q$ norms and adjusting corresponding H\"older-type estimates. (Setting $q=2$ below, we see that the parameters and exponents are consistent to the ones given along the proof of  Theorem \ref{th: main korn}.) Moreover, we need an $L^q$-version of Theorem  \ref{th: kornSBDsmall} (see  \cite{Chambolle-Conti-Francfort:2014, Conti-Iurlano:15.2}) and use the density result \cite{Crismale} instead of Corollary \ref{cor: iurlano*}. More specifically, we have the following adaptions, ordered by sections:

\smallskip

\emph{Section \ref{sec: 3sub1}:} The results in Section \ref{sec: 3sub1} (Theorem \ref{th: bad part} and Lemma \ref{lemma: bad sets}) remain completely unchanged as they rely exclusively on the geometry of $J_u$. 

\smallskip

\emph{Section \ref{sec: 3sub2}:} Theorem \ref{th: modifica} holds with $r = (q-p)/(6q^2)$ in place of $(2-p)/24$ and  in \eqref{eq: modifica prop**} we need to replace $\mu^{2/p}$ by $\mu^{2/p - 2/q + 1}$ and   $\Vert e(u) \Vert_{L^2(Q_\mu)}$ by $\Vert e(u) \Vert_{L^q(Q_\mu)}$. In Step I of the proof, the application of the Korn inequality for $SBD^2$ functions with small jump set (see \eqref{eq: poincare estim}) has to be replaced by the corresponding version for $SBD^q$ functions \cite{Conti-Iurlano:15.2}. (Again replace the exponent $2/p$ by $2/p -2/q + 1$.)  In the rest of the proof we only need to  adjust some exponents as we indicate now: In Step I-II, in the estimates on the differences of rigid motions \eqref{eq: diff1*}-\eqref{eq:33} we replace  $\Vert e(u) \Vert_2$ by $\Vert e(u) \Vert_q$ and $d(Q)^2,s_l^2$ by $d(Q)^{2- 2p/q + p}, s_l^{2-2p/q + p}$. The argument leading to \eqref{eq: main holderXXXX} is essentially the same, where we have  two  slightly different H\"older-type estimates: First, we need 
\begin{align}\label{eq: newHolder}
\sum\nolimits_{Q \in \mathcal{C}}  (d(Q))^{2-\frac{2p}{q}} \Vert e(u) \Vert_{L^q(N_Q)}^{p} & \le \Big(\sum\nolimits_{Q \in \mathcal{C}}  (d(Q))^{ (2-\frac{2p}{q})\frac{q}{q-p}} \Big)^{1-\frac{p}{q}} 
\Big(\sum\nolimits_{Q \in \mathcal{C}} \Vert e(u) \Vert_{L^q(N_Q)}^q \Big)^{\frac{p}{q}}\notag\\
&\le C \mu^{2 - \frac{2p}{q}} \Vert e(u) \Vert_{L^q(Q_{\mu})}^{p}.
\end{align}
Moreover, using $M_l= \sum_k \# \mathcal{X}^l_k \le C\mu s_l^{-1}$ and  $\# {\mathcal{X}}^l_k \le c \theta^{-2lr}$, we have
\begin{align}\label{eq: newHolder2}
\sum_{l \ge 8} \sum_k \# {\mathcal{X}}^l_k \, s_l^{2-\frac{2p}{q}} \theta^{-2(p+1)lr} \Vert e(u)\Vert^p_{L^q(N^l_k)} & \le C \sum_{l \ge 8} \sum_k  (\# {\mathcal{X}}^l_k)^{1-\frac{p}{q}} \ s_l^{2-\frac{2p}{q}} \theta^{-2(p+1+\frac{p}{q})lr} \Vert e(u)\Vert^p_{L^q(N^l_k)} \notag \\
& \le C\sum_{ l \ge 8} M_l^{1-\frac{p}{q}} s_l^{2-\frac{2p}{q}} \theta^{-2(p+1+\frac{p}{q})lr} \big(\sum_{k}\Vert e(u)\Vert^q_{L^q(N^l_k)}\big)^{\frac{p}{q}}\notag\\
& \le C\mu^{2-\frac{2p}{q}} \Vert e(u) \Vert^p_{L^q(Q_\mu)}  \sum_{l \ge 8} \theta^{l(1 - \frac{p}{q}) -2(p+1+\frac{p}{q})lr }   .
\end{align}
Recalling $r = (q-p)/(6q^2)$, one can check that $1 - p/q -2(p+1+p/q)r \ge 2(q-p)r$, i.e., the sum indeed converges. Finally, we note that the construction of $\bar{u}$ in Step III and the derivation of \eqref{eq: except***}-\eqref{eq: modifica prop**} remain unchanged.

\smallskip

\emph{Section \ref{sec: 3sub5}:} In the statement of Theorem \ref{th: partiparti} we again  replace $\mu^{2/p-1}$ by $\mu^{2/p - 2/q}$ and   $\Vert e(u) \Vert_{L^2(Q_\mu)}$ by $\Vert e(u) \Vert_{L^q(Q_\mu)}$. The proof remains unchanged as it only relies on geometric arguments and an application of the Korn inequality in John domains (Theorem \ref{th: kornsobo}) for $1 < p < q$.

\smallskip

\emph{Section \ref{sec: refined}:} In the statement of Lemma \ref{cor: partiparti} we only need to adapt the exponent in \eqref{eq:lastlast} replacing $\mu^{2-p}$ by $\mu^{2 - 2p/q}$. The proof relying entirely on geometric arguments remains the same. In Lemma  \ref{rem: large components} we define $r = (q-p)/(6q^2)$, the exponents in \eqref{eq: sharpi2.b} are adjusted as before and instead of \eqref{eq: traceestimate} we have 
\begin{align}\label{eq: newtrace}
\int_{\partial Q_\mu \setminus \Gamma_u} |Tu - a_0|^q\, d\mathcal{H}^1 \le c\mu^{q-1}\Vert e(u) \Vert^q_{L^q(Q_\mu)}.
\end{align}
In part (1) of the proof we only need to change the H\"older-type estimate \eqref{eq: main holderXXX} in the sense of  \eqref{eq: newHolder}. In part (2) of the proof, in the application of  Theorem \ref{th: kornSBDsmall} we observe that there is also a general version of \eqref{eq: main estmainXX} in terms of $L^q$-norms, see \cite{Chambolle-Conti-Francfort:2014}. 

\smallskip

\emph{Section \ref{sec: main*}:} At the beginning of the  proof of Theorem \ref{th: korn-small-nearly} we define the parameters \begin{align*}
r = (q-p)/(6q^2), \ \ \ \ \ p' = (q+p)/2, \ \ \ \  \   \lambda =  (1-r)(q-p)/(q+3p)
\end{align*}
and make the usual adaptions of the exponents replacing $2-p'$ by $2-2p'/q$ (see e.g. \eqref{eq: epsdef}, \eqref{eq: L4} and  \eqref{eq: L6.2} among others). Step II - Step IV of the proof do not address Korn-type estimates and remain unchanged. In Step V we need to adapt  the estimate for isolated components \eqref{eq: extendi2} and the H\"older-type estimates \eqref{eq: discrHol}. The derivation of \eqref{eq: extendi2} relies on the scaling in the trace estimate \eqref{eq: newtrace} and an application of Remark \ref{rem: slice} for exponent $q$. To see \eqref{eq: discrHol}, we follow the lines of the adaptions given in  \eqref{eq: newHolder}-\eqref{eq: newHolder2}. Finally, the essential change in Step VI concerns the H\"older-type estimate \eqref{eq: for later2}. Using $|D^I_l| \le C_2 \theta^{-rl}t_l \mu_0$ and $\Vert \nabla v' \Vert^{p'}_{L^{p'}(( Q_{\mu_0} \setminus E'') \cap D^I_l)} \le C  \mu_0^{2-2p'/q} \theta^{-\lambda l} \Vert e(u) \Vert_{L^q(Q_{\mu_0})}^{p'}$ (cf. \eqref{eq: Ddec} and \eqref{eq: Dsets2})  we calculate
\begin{align*}
\Vert \nabla v'\Vert^p_{L^p( Q_{\mu_0} \setminus E'')} &= \sum\nolimits_{l=0}^{I}  \Vert \nabla v' \Vert^p_{L^p(( Q_{\mu_0}\setminus E'') \cap D^I_l)}  \le \sum\nolimits_{l=0}^I  |D^I_l|^{1- \frac{p}{p'}}\Vert \nabla v' \Vert^p_{L^{p'}(( Q_{\mu_0} \setminus E'') \cap D^I_l)}\notag \\
& \le  C\mu_0^{2-\frac{2p}{q}}\sum\nolimits_{l} \theta^{l(1-r)(1-\frac{p}{p'}) -l\lambda \frac{p}{p'}}   \Vert  e(u) \Vert^{p}_{L^q( Q_{\mu_0})} = C\mu_0^{2-\frac{2p}{q}}\sum\nolimits_{l}  \theta^{l\lambda}  \Vert  e(u) \Vert^{p}_{L^q( Q_{\mu_0})}   
\end{align*}
where  in the last step we used the definition of $p'$ and $\lambda$. Theorem \ref{th: korn-small-nearly2} addresses the modification of the jump set and remains unchanged.

\smallskip

\emph{Section \ref{sec: main***}:} Finally, to extend the result from $\mathcal{W}(\Omega)$ to general functions $GSBD^q(\Omega)$, we proceed as in the proof of Theorem \ref{th: main korn} by employing the recently proved density result \cite{Crismale} instead of Corollary \ref{cor: iurlano*}. \eop

}

\end{rem}
\EEE

Now with Theorem \ref{th: main korn} at hand, the proof of our density result Theorem \ref{th: fried-iurlano} is straightforward.

\noindent {\em Proof of Theorem \ref{th: fried-iurlano}.} Let $u \in GSBD^2(\Omega)$. By Theorem \ref{th: main korn} applied for some $p \in [1,2)$ we obtain a Caccioppoli partition $(P_j)_j$ of $\Omega$ and corresponding infinitesimal rigid motions $(a_j)_j$ such that $v:= u - \sum_j a_j\chi_{P_j} \in SBV^p(\Omega) \cap L^\infty(\Omega)$. As motivated in \eqref{eq: sequence}, we consider the sequence 
$$
v_k = u - \sum\nolimits_{j\ge k} a_j \chi_{P_j} \in SBV^p(\Omega) \cap L^\infty(\Omega)
$$
and observe that $v_k \to u$ in measure on $\Omega$, $e(v_k) = e(u_k)$ for all $k \in \N$ and $\mathcal{H}^1(J_{v_k} \triangle J_u) \to 0$ when $k \to \infty$. Using Theorem \ref{th: iurlano}, each function $v_k$ can be approximated in $L^2(\Omega)$ by a sequence with the properties stated in Theorem \ref{th: iurlano} such that \eqref{eq: iuiu} holds. Now the assertion follows from a diagonal sequence argument. \eop

 \section{Proof of further results}\label{sec: main**}
 

 \subsection{Piecewise Poincar\'e inequality}\label{subsec: main**1}
 

  We start with the proof of the piecewise Poincar\'e inequality which is essentially based on the coarea formula for $BV$ functions and can be derived completely independently from the results discussed in the previous sections.

  \noindent {\em Proof of Theorem \ref{th: main poinc}.} Without restriction we assume $\Vert \nabla u \Vert_{L^1(\Omega)}>0$ as otherwise $u$ is piecewise constant (see \cite[Theorem 4.23]{Ambrosio-Fusco-Pallara:2000}, \cite{Chambolle-Giacomini-Ponsiglione:2007}) and there is nothing to show.  We start with the case $m=1$ and $u \in SBV(\Omega;\R)$. Following ideas in \cite{Ambrosio:90, Braides-Defranceschi, Francfort-Larsen:2003} we can use the coarea formula in $BV$ (see \cite[Theorem 3.40]{Ambrosio-Fusco-Pallara:2000}) to write
  $$\Vert \nabla u \Vert_{L^1(\Omega)} = |Du|(\Omega \setminus J_u) = \int_{-\infty}^\infty \mathcal{H}^{d-1}\big((\Omega \setminus J_u)\cap \partial^*\lbrace u>t \rbrace \big)\, dt.$$
Thus, \BBB for $\rho >0$ and \EEE $M:=\rho^{-1}\Vert \nabla u \Vert_{L^1(\Omega)}$ we find $t_i \in (iM, (i+1)M]$ for all $i \in \Z$ such that 
\begin{align}\label{eq: coarea}
\mathcal{H}^{d-1}\big((\Omega \setminus J_u) \cap \partial^*\lbrace u>t_i \rbrace \big) \le \frac{1}{M} \int_{iM}^{{(i+1)M}} \mathcal{H}^{d-1}\big((\Omega \setminus J_u) \cap \partial^*\lbrace u>t \rbrace \big)\, dt.
\end{align}
Let $E_i = \lbrace u > t_{i} \rbrace \setminus \lbrace u > t_{i+1} \rbrace$ and note that each set has finite perimeter in $\Omega$ since it is the difference of two sets of finite perimeter.  Now \eqref{eq: coarea} implies
\begin{align}\label{eq: coarea2}
\sum\nolimits_{i \in \Z} \mathcal{H}^{d-1}\big((\Omega \cap \partial^* E_i) \setminus J_u \big) \le \frac{2}{M}\Vert \nabla u \Vert_{L^1(\Omega)} = 2\rho.
\end{align} 
Since $|\Omega \setminus \bigcup_{i \in \Z} E_i|=0$, $(E_i)_i$ is a Caccioppoli partition of $\Omega$. Moreover, we note that the function $v:= u - \sum_i t_i\chi_{E_i}$ lies in $L^\infty(\Omega)$ and satisfies $\Vert v \Vert_{L^\infty(E_i)} \le 2M$. This implies \eqref{eq: kornpoinsharp0}.  Consider the sequence $v_n = \sum_{|i| \le n} (u- t_i)\chi_{E_i}$ for $n \in \N$ with $v_n \in SBV(\Omega;\R)$ by \cite[Theorem 3.84]{Ambrosio-Fusco-Pallara:2000}. Since $v_n \to v$ in $BV$ norm and $SBV$ is a closed subspace of $BV$, we conclude $v \in SBV(\Omega;\R) \cap L^\infty(\Omega)$.  This   concludes the proof in the case  $m=1$ and $u \in SBV(\Omega;\R)$.

If $u \in GSBV(\Omega;\R)$, we apply the analog of the coarea for $GSBV$ functions  (see \cite[Theorem 4.34]{Ambrosio-Fusco-Pallara:2000}) and again obtain $(t_i)_{i \in \Z}$, $(E_i)_{i \in \Z}$ such that $v: = u- \sum_i t_i\chi_{E_i} \in L^\infty(\Omega)$ and \eqref{eq: kornpoinsharp0} holds. (Note that $\nabla u$, $J_u$ have to be understood in a weaker sense, cf. \cite[Section 4.5]{Ambrosio-Fusco-Pallara:2000}.) The characterization of scalar $GSBV$ functions (see \cite[Remark 4.27]{Ambrosio-Fusco-Pallara:2000}) together with $\nabla u \in L^1(\Omega)$, $\mathcal{H}^{d-1}(J_u) < \infty$, yields  $\min\lbrace  \max \lbrace -n, u \rbrace, n \rbrace \in SBV(\Omega;\R)$ for all $n \in \N$. Consequently, similarly as before, the sequence $v_n =  \sum_{|i| \le n} (u -t_i)\chi_{E_i}$ lies in $SBV(\Omega;\R)$, converges to $v$ in $BV$ norm and thus $v \in SBV(\Omega;\R) \cap L^\infty(\Omega)$. 

Now let $u \in (GSBV(\Omega;\R))^m$ for $m \ge 2$. We repeat the above argumentation for each component $u^j$, $j=1,\ldots,m$, and obtain corresponding $(E^j_i)_{i \in \Z}$ and $(t^j_i)_{i \in \Z}$ such that  \eqref{eq: coarea2} holds with $E^j_i$ in place of $E_i$ and $v^j := u^j - \sum_i t_i^j \chi_{E^j_i}$ satisfies $\Vert v^j \Vert_\infty \le 2\rho^{-1} \Vert \nabla u \Vert_{L^1(\Omega)}$. We introduce the Caccioppoli partition  of $\Omega$, denoted by $(P_j)_{j \ge 1}$,  consisting of the (nonempty) sets in
$$ \lbrace E^1_{i_1} \cap \ldots \cap E^m_{i_m}: i_1,\ldots,i_m \in \Z\rbrace.$$
In view of  \eqref{eq: coarea2}, we get that \eqref{eq: kornpoinsharp0}(i) holds for a constant also depending on the dimension $m$. Fix $P_j = E^1_{i_1} \cap \ldots \cap E^m_{i_m}$ and define the corresponding translation $b_j \in \R^m$ by $b_j = (t_{i_1},\ldots,t_{i_m})$. Then we conclude that $v:= u - \sum_j b_j \chi_{P_j}$ satisfies  \eqref{eq: kornpoinsharp0}(ii) for a constant depending on $m$.  Arguing as before for each component, we find $v\in SBV(\Omega;\R^m) \cap L^\infty(\Omega; \R^m)$.   \eop

 \subsection{Embedding results}\label{subsec: main**2}
 

 We now prove the results stated in Section \ref{sec: main2}.

\noindent {\em Proof of Theorem \ref{th: embed}.} We first prove \eqref{eq: embed} for $u \in \mathcal{W}(Q)$ for a square $Q \subset \R^2$ in dimension two. Afterwards, we use a slicing and density argument to derive the result for domains in $\R^d$. By Theorem \ref{th: korn-small-nearly2} for $p=1$   we find a Caccioppoli partition $(P_j)_j$ of $Q \subset \R^2$ and corresponding infinitesimal rigid motions $(a_j)_j$ such that $v := u - \sum_j a_j \chi_{P_j} \in SBV(Q) \cap L^\infty(Q)$ and
\begin{align}\label{eq: all together}
\begin{split}
&\Vert v \Vert_\infty \le C( \mathcal{H}^1(J_u) + \mathcal{H}^1(\partial Q))^{-1}\Vert e(u) \Vert_{L^2(Q)}, \  \  \Vert \nabla v \Vert_{L^1(Q)} \le C\Vert e(u) \Vert_{L^2(Q)}, \\
& \sum\nolimits_j\mathcal{H}^1(\partial^* P_j) \le C\mathcal{H}^1(J_u)+C\mathcal{H}^1(\partial Q),
\end{split}
\end{align}
where $C$ only depends on the diameter of $Q$. To conclude the proof of \eqref{eq: embed} for $d=2$, it now suffices to show that \begin{align}\label{eq: all together2}
\sum\nolimits_j  |A_j||P_j| \le C(\mathcal{H}^1(J_u)+1)\Vert u \Vert_\infty  + C\Vert e(u) \Vert_{L^2(Q)}
\end{align}
for $C=C(Q)$. To this end, we use Lemma \ref{lemma: rigid motion} and the isoperimetric inequality   to obtain
$$
\sum\nolimits_j  |A_j||P_j| \le c\sum\nolimits_j|P_j|^{\frac{1}{2}} \Vert a_j \Vert_{L^\infty(P_j)} \le c(\Vert u \Vert_{\infty} + \Vert v \Vert_{\infty})\sum\nolimits_j  \mathcal{H}^1(\partial^* P_j).
$$
Then  \eqref{eq: all together2} follows from  \eqref{eq: all together}.

We now treat the case $Q=(-\mu,\mu)^d$ and $u \in \mathcal{W}(Q) \subset SBV(Q) \cap L^\infty(Q)$.  To prove the assertion, we need to control $\Vert \partial_j u_i \Vert_{L^1(\Omega)}$ for each $1 \le i,j\le d$. For notational convenience we only treat the case  $i,j \in \lbrace 1,2 \rbrace$. The other terms follow analogously due to the symmetry of the problem. For $x \in Q$ we write $x = (x_1,x_2,y)$ with $y \in (-\mu,\mu)^{d-2}$ and introduce the functions $w^{y} : (-\mu,\mu)^2 \to \R^2$, $w^{y}(x_1,x_2) = (u_1(x), u_2(x))$ for $y \in (-\mu,\mu)^{d-2}$. Applying the result in $d=2$ we obtain for a.e. $y \in (-\mu,\mu)^{d-2}$
$$\sum\nolimits_{1\le i,j \le 2} \Vert \partial_j w^{y}_i \Vert_{L^1((-\mu,\mu)^2;\R)}   \le C\Vert e(w^{y})\Vert_{L^2((-\mu,\mu)^2;\R^{2 \times 2}_{\rm sym})} + C\Vert u \Vert_\infty (\mathcal{H}^1(J_{w^{y}}) + 1),$$
where $C=C(\mu)$. Once we have proved
\begin{align}\label{eq: BV}
 \int_{(-\mu,\mu)^{d-2}} \mathcal{H}^1(J_{w^{y}})\, d\mathcal{H}^{d-2}(y) \le C\mathcal{H}^{d-1}(J_u),
\end{align}
we take the integral over $(-\mu,\mu)^{d-2}$, use Fubini's theorem and H\"older's  inequality to conclude
 $$\sum\nolimits_{1\le i,j \le 2} \Vert \partial_j u_i \Vert_{L^1(Q;\R)} \le C\Vert e(u)\Vert_{L^2(Q;\R^{d \times d}_{\rm sym})} + C\Vert u \Vert_\infty (\mathcal{H}^{d-1}(J_{u}) + 1).$$
To see \eqref{eq: BV}, we apply  slicing techniques for $BV$ functions  (see \cite[Section 3.11]{Ambrosio-Fusco-Pallara:2000}): for $f \in SBV((-\mu,\mu)^n;\R^m)$ and $j=1,\ldots,n$ we define $f_{j,s}: (-\mu,\mu) \to \R^m$ by $f_{j,s}(t) = f(s+t \e_j)$ for $s \in \Pi^n_j := \lbrace s \in (-\mu,\mu)^{n}: s \cdot \e_j = 0\rbrace$. We obtain
$$\int_{\Pi^n_j} \# J_{f_{j,s}} \,d\mathcal{H}^{n-1}(s) = \int_{J_f} |\nu_f \cdot \e_j|\, d\mathcal{H}^{n-1},$$
where $\nu_f$ denotes a normal of the jump set $J_f$. First, applying this estimate on $w^y \in SBV((-\mu,\mu)^2;\R^2)$ for $j=1,2$ and $y \in (-\mu,\mu)^{n-2}$ a.e. we obtain
 $$\int_{\Pi^2_1} \# J_{w^{y}_{1,s}} \,d\mathcal{H}^{1}(s) + \int_{\Pi^2_2} \# J_{w^{y}_{2,s}} \,d\mathcal{H}^{1}(s) = \int_{J_{w^y}} (|\nu_{w^y} \cdot \e_1| + |\nu_{w^y} \cdot \e_2|)\, d\mathcal{H}^1 \ge \mathcal{H}^1(J_{w^y}).$$
Repeating the argument for the function $w=(u_1,u_2) \in SBV((-\mu,\mu)^d;\R^2)$ for $j=1,2$, we also get  
 \begin{align*}
 \int_{\Pi^d_1} \# J_{w_{1,s}} \,d\mathcal{H}^{d-1}(s) + \int_{\Pi^d_2} \# J_{w_{2,s}} \,d\mathcal{H}^{d-1}(s)  &= \int_{J_{w}} (|\nu_{w} \cdot \e_1| + |\nu_{w} \cdot \e_2|)\, d\mathcal{H}^{d-1} \\&\le 2\mathcal{H}^{d-1}(J_{w})\le 2\mathcal{H}^{d-1}(J_{u}).
 \end{align*}
Taking the integral over $(-\mu,\mu)^{d-2}$ we derive \eqref{eq: BV} from the last two estimates.

It remains to consider general Lipschitz domains $\Omega \subset \R^d$ and  $u \in SBD^2(\Omega) \cap L^\infty(\Omega)$. First, we choose a cube $Q$ containing $\Omega$ and define the extension $\bar{u} = u \chi_\Omega \in SBD^2(Q) \cap L^\infty(Q)$. Clearly, the choice of $Q$ depends only on $\Omega$. Note that $\mathcal{H}^{d-1}(J_{\bar{u}}) \le \mathcal{H}^{d-1}(J_u) + \mathcal{H}^{d-1}(\partial \Omega)$. By Theorem \ref{th: iurlano}, Remark \ref{rem: iur}  we find a sequence $(u_k)_k \subset \mathcal{W}(Q)$ with $u_k \to \bar{u}$ in $L^2(Q)$  and $\Vert u_k \Vert_\infty \le \Vert u \Vert_\infty$ such that by the above arguments
 \begin{align}\label{eq: embed2}
 \begin{split}
  |Du_k|(Q) &= \Vert \nabla u_k \Vert_{L^1(Q)} +  \int_{J_{u_k}}|[u_k]|\,d\mathcal{H}^{d-1} \le C\Vert e(u)\Vert_{L^2(\Omega)} + C\Vert u \Vert_\infty (\mathcal{H}^{d-1}(J_{u}) + 1)
  \end{split}
 \end{align}
 for $C=C(\Omega)$, where we used $|[u_k]| \le 2 \Vert u \Vert_\infty$ almost everywhere. As $(u_k)_k$ is uniformly bounded in  $BV$ norm, we deduce   $u \in BV(\Omega)$ and that \eqref{eq: embed}  holds by lower semicontinuity. Finally, by  Alberti's rank one property $|D^c u| \le \sqrt{2}|E^c u|$ and the fact that $u \in SBD(\Omega)$ we conclude $u \in SBV(\Omega)$.\eop

 \noindent {\em Proof of Theorem \ref{th: embed2}.} Let $u \in GSBD^2(\Omega)$. We show that each component $u_i$, $i=1,\ldots,d$,  satisfies \eqref{eq: embed0} for the truncation $u_i^M = \min \lbrace\max\lbrace u_i,-M\rbrace,M\rbrace$. Herefrom we particularly deduce $u_i \in GBV(\Omega; \R)$ since in the scalar case this property is equivalent to $u_i^M \in BV_{\rm loc}(\Omega)$ for all $M > 0$ (see \cite[Remark 4.27]{Ambrosio-Fusco-Pallara:2000}). As in the previous proof it essentially suffices to show 
 \begin{align}\label{eq: embed1}
 \begin{split}
|Du^M_i|(Q) &\le  \Vert \nabla u_i^M \Vert_{L^1(Q)} + 2M\mathcal{H}^{d-1}(J_{u^M_i}) \le CM(\mathcal{H}^{d-1}(J_u)+1) + C\Vert e(u) \Vert_{L^2(Q)},
\end{split} 
 \end{align}
 where $\Omega =Q$ is a cube, $C=C(Q)$ and $u \in \mathcal{W}(Q) \subset SBV(Q)$.  In fact, we then establish the general case  using the approximation of $u$ given by Corollary \ref{cor: iurlano*} and repeating the argument in  \eqref{eq: embed2}. (Note, however, that in contrast to \eqref{eq: embed2}, we cannot apply Alberti's theorem and therefore only obtain $u^M_i \in BV(\Omega)$.) Finally, in view of the slicing argument \eqref{eq: BV}, it is enough to treat the planar case $u \in \mathcal{W}(Q)$ for a square $Q \subset \R^2$.

  By Theorem \ref{th: korn-small-nearly2} for $p=1$ we get that $v := u - \sum_j a_j\chi_{P_j} \in SBV(Q) \cap L^\infty(Q)$  satisfying $\Vert v \Vert_{L^\infty(Q)}\le  M'$ with $M' := C(\mathcal{H}^1(J_u) + \mathcal{H}^1(\partial Q))^{-1} \Vert e(u) \Vert_{L^2(Q)}$ by \eqref{eq: main korn2}. Let $a = \sum_j a_j\chi_{P_j}$ and $a_i$ be the $i$-th component for $i = 1,2$. From \cite[Theorem 2.2]{Conti-Iurlano:15} we obtain $a \in (GSBV(Q))^2$, in particular it is shown that  $|D a_i^M|(Q) \le  cM \sum\nolimits_j\mathcal{H}^1(\partial^* P_j)$ for a universal $c>0$. Thus, 
 \begin{align}\label{eq: embed3}
 |D a_i^M|(Q) \le  cM \sum\nolimits_j\mathcal{H}^1(\partial^* P_j) \le cM(\mathcal{H}^1(J_u) + \mathcal{H}^1(\partial Q))
 \end{align}
 by \eqref{eq: main1}.  Up to sets of negligible $\mathcal{L}^2$-measure we have
 \begin{align*}
T:= \lbrace \nabla u_i^M \neq 0 \rbrace \subset \lbrace |u_i| \le M \rbrace \subset \lbrace |a_i| \le M+M' \rbrace
 \end{align*}
since $\Vert v \Vert_{L^\infty(Q)}\le  M'$. Therefore, we compute using \eqref{eq: embed3} and  \eqref{eq: main korn2}
  \begin{align}\label{eq: viewargu}
  \begin{split}
\Vert \nabla u_i \Vert_{L^1(T)} &\le \Vert \nabla v_i \Vert_{L^1(T)} + \Vert \nabla a_i \Vert_{L^1(T)}   \le \Vert \nabla v_i \Vert_{L^1(Q)} +  \Vert \nabla a^{M+M'}_i \Vert_{L^1(Q)} \\&\le c(M+M')(\mathcal{H}^1(J_u)+\mathcal{H}^1(\partial Q)) + C\Vert e(u) \Vert_{L^2(Q)}\\ &\le CM(\mathcal{H}^1(J_u)+1) + C\Vert e(u) \Vert_{L^2(Q)}.
\end{split}
 \end{align}
 As $\Vert \nabla u_i^M \Vert_{L^1(Q)} = \Vert \nabla u_i \Vert_{L^1(T)}$ and $\mathcal{H}^1(J_{u_i^M} \setminus  J_{u}) = 0$, we obtain the second inequality in \eqref{eq: embed1}. The first inequality follows from the decomposition of the distributional derivative and the fact that $\Vert u_i^M \Vert_\infty \le M$. This concludes the proof.   \eop

We also obtain the following variant of \eqref{eq: embed0} needed in Section \ref{subsec: main**3}.

\begin{corollary}\label{cor: vari}
Let $Q \subset \R^d$ be a cube. Then there is a constant $C=C(Q)>0$ such that for all $u \in \mathcal{W}(Q)$ and Borel sets $F \subset Q$ one has
$$|Du_i^M|(Q) \le CM\big(\mathcal{H}^{d-1}(J_u) + |F|\big) + C\big(\Vert e(u) \Vert_{L^2(Q)}  + \Vert u \Vert_{L^1(Q \setminus F)}\big).$$
for all $i=1,\ldots,d$ and $M>0$, where $u_i^M=: \min \lbrace\max\lbrace u_i,-M\rbrace,M\rbrace$.   
\end{corollary}

\Proof Similarly as in the proof of Theorem \ref{th: embed2}, it suffices to show this estimate in the planar setting $d=2$ for a square $Q \subset \R^2$ as the general case then follows by Fubini and the slicing argument \eqref{eq: BV} for a constant depending on the dimension. Let $\bar{c} = \bar{c}(Q)>0$ to be specified below. If $\mathcal{H}^1(J_u) + |F| \ge \bar{c}$, the result follows directly from \eqref{eq: embed1} for a constant depending  on $\bar{c}$.

Now let $\mathcal{H}^1(J_u) + |F| \le \bar{c}$. By Theorem \ref{th: korn-small-nearly2} for $p=1$ in the version of Remark \ref{rem:square} we obtain an ordered Caccioppoli partition $(P_j)_j$ and $v := u - \sum_j a_j\chi_{P_j}$ such that for $C=C(Q)>0$
\begin{align}\label{eq:likewise-new}
\Vert \nabla v \Vert_{L^1(Q)}  + \mathcal{H}^1(J_u)\Vert v \Vert_{L^\infty(Q)}\le C\Vert e(u) \Vert_{L^2(Q)},   \ \ \ \  \sum\nolimits_j\mathcal{H}^1(\partial^* P_j \cap Q) \le C\mathcal{H}^1(J_u).
\end{align}
The essential step is now to show 
\begin{align}\label{eq:likewise}
\Vert\nabla a_i^{M+M'}\Vert_{L^1(Q)} \le C\big((M+M')\mathcal{H}^{1}(J_u) + \Vert u \Vert_{L^1(Q \setminus F)} + \Vert e(u) \Vert_{L^2(Q)}\big)
\end{align}
for $a : =\sum_j a_j\chi_{P_j}$ and $M' := C(\mathcal{H}^1(J_u))^{-1} \Vert e(u) \Vert_{L^2(Q)}$. Indeed, the claim then follows by repeating the argument in \eqref{eq: viewargu} noting that $\lbrace \nabla u_i^M \neq 0 \rbrace \subset \lbrace |a_i| \le M+M'\rbrace$ since $\Vert v \Vert_\infty \le M'$.

\BBB We now confirm \eqref{eq:likewise}. \EEE  For notational convenience we set  ${\bar{M}} = M+M'$. Since $\sum_j\mathcal{H}^1(\partial^* P_j \cap Q) \le C\bar{c}$, for $\bar{c}$ small enough the relative isoperimetric inequality implies $|P_1| > \frac{1}{2}|Q|$ and $|P_j| \le C(\mathcal{H}^1(\partial^* P_j \cap Q))^2$ for $j \ge 2$. Without restriction we assume that the sets $(P_j)_{j \ge 2}$ are connected (more precisely indecomposable, see \cite[Example 4.18]{Ambrosio-Fusco-Pallara:2000}) as otherwise we consider the indecomposable components. By \cite[Lemma 4.6]{FriedrichSolombrino} we get that the diameter of each $P_j$, $j \ge 2$,  is controlled in terms of $C\mathcal{H}^1(\partial^* P_j \cap Q)$ for $C=C(Q)$, which by \BBB \cite[Lemma 4.3]{FriedrichSolombrino} \EEE then also yields $\mathcal{H}^1(\partial^* P_j) \le C\mathcal{H}^1(\partial^* P_j \cap Q)$ for all $j \ge 2$.  Then again by \cite[Theorem 2.2]{Conti-Iurlano:15} (cf. also \eqref{eq: embed3}) \BBB and \eqref{eq:likewise-new} \EEE
\begin{align}\label{eq: diffi}
\begin{split}
\sum_{j \ge 2} \Vert \nabla a_i^{\bar{M}} \Vert_{L^1(P_j)} &\le c{\bar{M}}\sum_{j \ge 2}\mathcal{H}^1(\partial^* P_j) \le C{\bar{M}}\sum_{j \ge 2}\mathcal{H}^1(\partial^* P_j \cap Q) \le C{\bar{M}}\mathcal{H}^1(J_u).
\end{split}
\end{align}
By Theorem \ref{th: kornSBDsmall} we find an infinitesimal rigid motion $a' = a_{A',b'}$ and an exceptional set $F' \subset Q$ with 
\begin{align*}
\Vert \nabla u - A' \Vert_{L^1(Q \setminus F')} +  \Vert  u - a' \Vert_{L^1(Q \setminus F')}\le C\Vert e(u) \Vert_{L^2(Q)}, \  \ \ \ \ \  |F'| \le C(\mathcal{H}^1(J_u))^2 \le C\bar{c}^2.
\end{align*}
Consequently, for $\bar{c}$  small we have $|F\cup F'| \le \BBB \bar{c} + C\bar{c}^2 \EEE \le\frac{1}{4}|Q|$  and derive   using Lemma \ref{lemma: rigid motion}, \BBB \eqref{eq:likewise-new} \EEE  and the triangle inequality with $G := F \cup F'$
\begin{align*}
\Vert \nabla a_i^{\bar{M}} \Vert_{L^1(P_1)}& \le |P_1||A_1|\le 2|P_1 \setminus G||A_1| \\&\le 2|Q \setminus G||A'| + 2\Vert \nabla u - A_1 \Vert_{L^1(P_1)} +  2\Vert \nabla u - A' \Vert_{L^1(Q \setminus F')} \\
&\le C|Q\setminus G|^{-\frac{1}{2}} \Vert a' \Vert_{L^1(Q \setminus  G)} + C\Vert e(u)\Vert_{L^2(Q)} \\
& \le C\Vert u - a' \Vert_{L^1(Q \setminus F')} + C\Vert u \Vert_{L^1(Q \setminus F)} + C\Vert e(u)\Vert_{L^2(Q)} \\&\le C\Vert u \Vert_{L^1(Q \setminus F)} + C\Vert e(u)\Vert_{L^2(Q)}.
\end{align*}
This together with \eqref{eq: diffi} concludes the proof of \eqref{eq:likewise}. \eop

 We close this section with the proof of   Lemma \ref{lemma: lq}.

 \noindent {\em Proof of Lemma \ref{lemma: lq}.}  Let $q \in [1,\infty)$. Consider the function defined in \eqref{eq: cif} with $r_k = k^{-p}$ and $d_k = k^{-1 + dp}$, where $p = \frac{1}{d-1} + \frac{1}{q(d-1)^2}$. (Note that the existence of pairwise disjoint balls with this property is guaranteed by \cite[Lemma 12.2]{DalMaso:13} and the approximate differential of $u$ in the sense \cite[Definition 3.70]{Ambrosio-Fusco-Pallara:2000} exists a.e.)  We first see that $\mathcal{H}^{d-1}(J_u) < + \infty$ as $\sum_k r_k^{d-1} < \infty$ due to the fact that $p > \frac{1}{d-1}$. To see that $u \in L^q(\Omega)$ and $\nabla u \notin L^1(\Omega)$, it suffices to show
$$\sum\nolimits_k r_k^d (r_k d_k)^q < \infty, \ \ \ \ \ \sum\nolimits_k r_k^d d_k = \infty.$$
The latter is immediate since $-pd - 1 + pd  = -1$. To see the first property we calculate  
$-p(d+q) -q + dpq = \tfrac{1}{d-1} - pd< -1,$ where in the last step we used $p > \frac{1}{d-1}$.

It remains to show that $u \in GSBD^2(\Omega)$. To this end, we consider the sequence of functions 
$$u_j = \sum\nolimits^j_{k=1} (A_k \,(x-x_k)) \chi_{B_k}(x) \in GSBD^2(\Omega) \cap L^q(\Omega)$$
converging to $u$ and by the compactness theorem for $GSBD$ (see \cite[Theorem 11.3]{DalMaso:13}) we see that $u \in GSBD^2(\Omega)$ since $e(u_j)=0$ a.e., $\sup_j \mathcal{H}^{d-1}(J_{u_j}) < \infty$ and $\sup_j \Vert u\Vert_{L^q(\Omega)} < \infty$. \eop

 \subsection{A Korn-Poincar\'e inequality for functions with small jump set}\label{subsec: main**3}
 

We finally give the proof of the Korn-Poincar\'e inequality for functions with small jump set. 

 \noindent {\em Proof of Theorem \ref{th: kornpoin-d}.} As in the previous sections we first treat the case $u \in \mathcal{W}(Q) \subset SBV(Q)$ with $\mathcal{H}^{d-1}(J_u) >0$ and at  the end of the proof we  indicate the adaptions if $u \in GSBD^2(Q)$. We may assume that $\mathcal{H}^{d-1}(J_u) \le \mathcal{H}^{d-1}(\partial Q)$ since otherwise the theorem trivially holds with $E = Q$ for $C=C(Q)>0$ sufficiently large. We apply the Korn-Poincar\'e inequality due to Chambolle, Conti, and Francfort (see \cite[Theorem 1]{Chambolle-Conti-Francfort:2014}) which is as the assertion of Theorem \ref{th: kornpoin-d} with the difference that only the volume of the exceptional set can be controlled. We find $F \subset Q$ with $|F| \le C(\mathcal{H}^{d-1}(J_u))^{\frac{d}{d-1}}$ for $C=C(Q)>0$ and an infinitesimal rigid motion $a$ such that with $q = \frac{2d}{d-1}$ and $v := u-a$ 
\begin{align}\label{eq: ccf1}
\Vert v \Vert_{L^q(Q \setminus F)} \le C \Vert e(u) \Vert_{L^2(Q)}.
\end{align}
We now define 
$$M := (\mathcal{H}^{d-1}(J_u))^{-\frac{d}{q(d-1)}} \Vert e(u) \Vert_{L^2(Q)} = (\mathcal{H}^{d-1}(J_u))^{-\frac{1}{2}} \Vert e(u) \Vert_{L^2(Q)}$$
and observe that by \eqref{eq: ccf1} 
\begin{align}\label{eq: eta2}
|\lbrace |v| > M \rbrace| \le  |F| + \tfrac{1}{M^{q}}\Vert v \Vert^{q}_{L^{q}(Q \setminus F)} \le C(\mathcal{H}^{d-1}(J_u))^{\frac{d}{d-1}}.
\end{align}
 We now consider the truncation $v'_i := \min \lbrace\max\lbrace v_i,-2M\rbrace,2M\rbrace$ for $i=1,\ldots,d$. By Corollary \ref{cor: vari} we get   for $i=1,\ldots,d$
\begin{align}\label{eq: calcul}
\begin{split}
|D v'_i|(Q) &\le CM(\mathcal{H}^{d-1}(J_u)+|F|) + C\Vert e(v) \Vert_{L^2(Q)} + C\Vert v \Vert_{L^1(Q\setminus F)} \\
& \le CM\mathcal{H}^{d-1}(J_u) + C\Vert e(u) \Vert_{L^2(Q)},
\end{split}
\end{align}
where in the last step we used \eqref{eq: ccf1}  and  $|F| \le C(\mathcal{H}^{d-1}(J_u))^{\frac{d}{d-1}} \le C\mathcal{H}^{d-1}(J_u)$. (Recall $\mathcal{H}^{d-1}(J_u) \le \mathcal{H}^{d-1}(\partial Q)$.)  Now arguing similarly as in the proof of Theorem \ref{th: main poinc} (see Section \ref{subsec: main**1}), using the coarea formula \cite[Theorem 3.40]{Ambrosio-Fusco-Pallara:2000}, we can find $t_i \in (M,2M)$ such that $G_i := \lbrace -t_i < v_i' < t_i \rbrace$ is a set of finite perimeter and
\begin{align}\label{eq: R2mainXXX}
\mathcal{H}^{d-1}(Q \cap \partial^*G_i)  &\le \frac{1}{M} \int_{-\infty}^\infty \mathcal{H}^{d-1}\big(Q \cap \partial^*\lbrace v'_i>t \rbrace \big)\, dt \le \frac{1}{M}|Dv'_i|(Q) \\
 &\le  \frac{1}{M} \big( CM\mathcal{H}^{d-1}(J_{u}) + C\Vert e(u) \Vert_{L^2(Q)}\big) \le C(\mathcal{H}^{d-1}(J_{u}))^{\frac{1}{2}},\notag
 \end{align}
 where in the last step we used $\mathcal{H}^{d-1}(J_{u})\le C(\mathcal{H}^{d-1}(J_{u}))^{\frac{1}{2}} $ for a constant $C=C(Q)$. 
 
  We define $E = Q \setminus \bigcap_{i=1}^d G_i$. As  $Q \setminus G_i \subset \lbrace |v'_i| > M \rbrace$ and $\lbrace |v'_i| > M \rbrace \subset \lbrace |v_i| > M \rbrace$, \eqref{eq: eta2} yields the second part of \eqref{eq: R2mainx^2}.  By  \eqref{eq: R2mainXXX} we get $\mathcal{H}^{d-1}(\partial^* E \cap Q) \le C(\mathcal{H}^{d-1}(J_{u}))^{\frac{1}{2}}$, which shows the first part of \eqref{eq: R2mainx^2}. Since $v'_i = v_i$ on $\lbrace |v'_i| < 2M \rbrace$ for $i=1,\ldots,d$ and 
  $$Q \setminus E = \bigcap\nolimits_{i=1}^d G_i \subset \bigcap\nolimits_{i=1}^d \lbrace |v'_i| < 2M  \rbrace,$$
  we observe $v'  = v$ on $Q \setminus E$, where $v' := (v_1',\ldots,v_d')$. Consequently, we have $\Vert v \Vert_{L^\infty(Q \setminus E)} \le CM = C(\mathcal{H}^{d-1}(J_u))^{-\frac{1}{2}}\Vert e(u) \Vert_{L^2(Q)}$ by  the definition of $M$. This  gives \eqref{eq: main estmainx^2}(ii) and  together with \eqref{eq: ccf1}, $|F| \le C(\mathcal{H}^{d-1}(J_u))^{\frac{d}{d-1}}$ we obtain
\begin{align*}
\Vert v \Vert_{L^q(Q \setminus E)} & \le \Vert v \Vert_{L^q(Q \setminus F)} + \Vert v \Vert_{L^q(F \setminus E)} \le C \Vert e(u) \Vert_{L^2(Q)} + CM|F|^{\frac{1}{q}} \le C\Vert e(u) \Vert_{L^2(Q)},
\end{align*}
which establishes  \eqref{eq: main estmainx^2}(i). Finally, $\bar{u} := (u - a) \chi_{Q \setminus E} =  v \chi_{Q \setminus E}=  v' \chi_{Q \setminus E}$ clearly lies in $SBV^2(Q) \cap L^\infty(Q)$ since $u \in \mathcal{W}(Q)$. We use \eqref{eq: calcul} and \cite[Theorem 3.84]{Ambrosio-Fusco-Pallara:2000}  to calculate
\begin{align*}
|D\bar{u}|(Q) &\le |Dv'|(Q) + CM\mathcal{H}^{d-1}(\partial^* E \cap Q) \le CM(\mathcal{H}^{d-1}(J_u))^{\frac{1}{2}} + C\Vert e(u) \Vert_{L^2(Q)},
\end{align*}
which by the definition of $M$ gives \eqref{eq: main estmain2X}. The variant announced below Theorem \ref{th: kornpoin-d} with $p \in [1,2]$ follows by replacing $q= \frac{2d}{d-1}$ by $q'= \frac{2d}{p(d-1)}$ in the above proof. 

Finally, we treat the case $u \in GSBD^2(Q)$ with $\mathcal{H}^{d-1}(J_u)>0$. We proceed similarly as in the density argument in the proof of Theorem \ref{th: main korn} (see Section \ref{sec: main***}) and refer therein for details. Let $(u^j)_j$ be a sequence of rescaled versions of $u$ defined on squares $Q_j \supset \supset Q$ with $|Q_j \setminus Q| \to 0$ for $j \to \infty$.  Using Lemma \ref{lemma: iurlano*} we approximate $(u^j)_j$ with a (diagonal) sequence $(u_k)_k \subset \mathcal{W}(Q)$ such that \eqref{eq: iurl} holds with $u_k \to u$ a.e. on $Q$.  

We then obtain infinitesimal rigid motions $(a_k)_k$ and exceptional sets $(E_k)_k$ such that \eqref{eq: R2mainx^2}-\eqref{eq: main estmain2X} hold for each $k \in \N$. By compactness we find a set of finite perimeter $E$ satisfying \eqref{eq: R2mainx^2} such that $\chi_{E_k} \to \chi_E$ in measure after extracting a not relabeled subsequence. Moreover, applying Lemma \ref{lemma: rigid motionXXX} we find that also the infinitesimal rigid motions converge to some $a$ and \eqref{eq: main estmainx^2} follows by lower semicontinuity. 

The sequence $\bar{u}_k := (u_k -a_k) \chi_{Q \setminus E_k}$ is uniformly bounded in $BV$ norm and we therefore deduce $\bar{u} := (u - a) \chi_{Q \setminus E} \in BV(Q)$ satisfies \eqref{eq: main estmain2X} by lower semicontinuity. Likewise, a compactness result in $SBD^2$ together with the uniform bound on $\Vert \bar{u}_k \Vert_\infty$ gives $\bar{u} \in SBD^2(Q)$ (see \cite{Bellettini-Coscia-DalMaso:98}). Finally, Alberti's rank one property $|D^c \bar{u}| \le \sqrt{2}|E^c \bar{u}|$ also yields $\bar{u} \in SBV(Q)$.  \eop

 \smallskip

\noindent \textbf{Acknowledgements} This work has been funded by the Vienna Science and Technology Fund (WWTF)
through Project MA14-009. The support by the Alexander von Humboldt Stiftung is gratefully acknowledged. I would like to thank {\sc Gianni Dal Maso} for turning my attention to this problem by conjecturing Theorem \ref{th: embed}. Moreover, I am very grateful to {\sc Francesco Solombrino} for  stimulating discussions on the content of this paper. \BBB Finally, I am gratefully indebted to the referee for her/his careful reading of the manuscript and many helpful suggestions. \EEE


 \typeout{References}

\end{document}